\documentclass[10pt,oneside,reqno]{amsart}

\pdfoutput=1

\usepackage{accents}
\usepackage{afterpage}
\usepackage{amscd}
\usepackage{amsfonts}
\usepackage{amsmath}
\usepackage{amssymb}
\usepackage{amsthm}
\usepackage{appendix}

\usepackage{bbm}
\usepackage{bm}
\usepackage{braket}

\usepackage{calc}
\usepackage{caption}
\usepackage{changepage}
\usepackage{color}
\usepackage{comment}

\usepackage{dsfont}

\usepackage{enumitem}
\usepackage{etoolbox}
\usepackage{extpfeil}

\usepackage{graphicx}

\usepackage{import}

\usepackage{mathbbol}
\usepackage{mathrsfs}
\usepackage{mathtools}

\usepackage{pifont}

\usepackage{physics}

\usepackage{relsize}

\usepackage{scalerel}
\usepackage{stackengine}
\usepackage{stackrel}
\usepackage{stmaryrd}
\usepackage{subcaption}

\usepackage{tikz}
    \usetikzlibrary{arrows.meta}
    \usetikzlibrary{shapes}
    \usetikzlibrary{decorations.shapes}
    \usetikzlibrary{decorations.pathmorphing}

\usepackage{tikz-cd}
\usepackage{transparent}

\usepackage{wasysym}
\usepackage{wrapfig}

\usepackage{xcolor}
\usepackage{xfrac}
\usepackage{xparse}

\makeatletter
\let\old@tocline\@tocline
\let\section@tocline\@tocline
\newcommand{\subsection@dotsep}{4.5}
\newcommand{\subsubsection@dotsep}{4.5}
\patchcmd{\@tocline}
  {\hfil}
  {\nobreak
     \leaders\hbox{$\m@th
        \mkern \subsection@dotsep mu\hbox{.}\mkern \subsection@dotsep mu$}\hfill
     \nobreak}{}{}
\let\subsection@tocline\@tocline
\let\@tocline\old@tocline

\patchcmd{\@tocline}
  {\hfil}
  {\nobreak
     \leaders\hbox{$\m@th
        \mkern \subsubsection@dotsep mu\hbox{.}\mkern \subsubsection@dotsep mu$}\hfill
     \nobreak}{}{}
\let\subsubsection@tocline\@tocline
\let\@tocline\old@tocline

\let\old@l@subsection\l@subsection
\let\old@l@subsubsection\l@subsubsection

\def\@tocwriteb#1#2#3{%
  \begingroup
    \@xp\def\csname #2@tocline\endcsname##1##2##3##4##5##6{%
      \ifnum##1>\c@tocdepth
      \else \sbox\z@{##5\let\indentlabel\@tochangmeasure##6}\fi}%
    \csname l@#2\endcsname{#1{\csname#2name\endcsname}{\@secnumber}{}}%
  \endgroup
  \addcontentsline{toc}{#2}%
    {\protect#1{\csname#2name\endcsname}{\@secnumber}{#3}}}%

\newlength{\@tocsectionindent}
\newlength{\@tocsubsectionindent}
\newlength{\@tocsubsubsectionindent}
\newlength{\@tocsectionnumwidth}
\newlength{\@tocsubsectionnumwidth}
\newlength{\@tocsubsubsectionnumwidth}
\newcommand{\settocsectionnumwidth}[1]{\setlength{\@tocsectionnumwidth}{#1}}
\newcommand{\settocsubsectionnumwidth}[1]{\setlength{\@tocsubsectionnumwidth}{#1}}
\newcommand{\settocsubsubsectionnumwidth}[1]{\setlength{\@tocsubsubsectionnumwidth}{#1}}
\newcommand{\settocsectionindent}[1]{\setlength{\@tocsectionindent}{#1}}
\newcommand{\settocsubsectionindent}[1]{\setlength{\@tocsubsectionindent}{#1}}
\newcommand{\settocsubsubsectionindent}[1]{\setlength{\@tocsubsubsectionindent}{#1}}

\renewcommand{\l@section}{\section@tocline{1}{\@tocsectionvskip}{\@tocsectionindent}{}{\@tocsectionformat}}%
\renewcommand{\l@subsection}{\subsection@tocline{2}{\@tocsubsectionvskip}{\@tocsubsectionindent}{}{\@tocsubsectionformat}}%
\renewcommand{\l@subsubsection}{\subsubsection@tocline{3}{\@tocsubsubsectionvskip}{\@tocsubsubsectionindent}{}{\@tocsubsubsectionformat}}%
\newcommand{\@tocsectionformat}{}
\newcommand{\@tocsubsectionformat}{}
\newcommand{\@tocsubsubsectionformat}{}
\expandafter\def\csname toc@1format\endcsname{\@tocsectionformat}
\expandafter\def\csname toc@2format\endcsname{\@tocsubsectionformat}
\expandafter\def\csname toc@3format\endcsname{\@tocsubsubsectionformat}
\newcommand{\settocsectionformat}[1]{\renewcommand{\@tocsectionformat}{#1}}
\newcommand{\settocsubsectionformat}[1]{\renewcommand{\@tocsubsectionformat}{#1}}
\newcommand{\settocsubsubsectionformat}[1]{\renewcommand{\@tocsubsubsectionformat}{#1}}
\newlength{\@tocsectionvskip}
\newcommand{\settocsectionvskip}[1]{\setlength{\@tocsectionvskip}{#1}}
\newlength{\@tocsubsectionvskip}
\newcommand{\settocsubsectionvskip}[1]{\setlength{\@tocsubsectionvskip}{#1}}
\newlength{\@tocsubsubsectionvskip}
\newcommand{\settocsubsubsectionvskip}[1]{\setlength{\@tocsubsubsectionvskip}{#1}}

\patchcmd{\tocsection}{\indentlabel}{\makebox[\@tocsectionnumwidth][l]}{}{}
\patchcmd{\tocsubsection}{\indentlabel}{\makebox[\@tocsubsectionnumwidth][l]}{}{}
\patchcmd{\tocsubsubsection}{\indentlabel}{\makebox[\@tocsubsubsectionnumwidth][l]}{}{}

\newcommand{\@sectypepnumformat}{}
\renewcommand{\contentsline}[1]{%
  \expandafter\let\expandafter\@sectypepnumformat\csname @toc#1pnumformat\endcsname%
  \csname l@#1\endcsname}
\newcommand{\@tocsectionpnumformat}{}
\newcommand{\@tocsubsectionpnumformat}{}
\newcommand{\@tocsubsubsectionpnumformat}{}
\newcommand{\setsectionpnumformat}[1]{\renewcommand{\@tocsectionpnumformat}{#1}}
\newcommand{\setsubsectionpnumformat}[1]{\renewcommand{\@tocsubsectionpnumformat}{#1}}
\newcommand{\setsubsubsectionpnumformat}[1]{\renewcommand{\@tocsubsubsectionpnumformat}{#1}}
\renewcommand{\@tocpagenum}[1]{%
  \hfill {\mdseries\@sectypepnumformat #1}}

\let\oldappendix\appendix
\renewcommand{\appendix}{%
  \leavevmode\oldappendix%
  \addtocontents{toc}{%
    \protect\settowidth{\protect\@tocsectionnumwidth}{\protect\@tocsectionformat\sectionname\space}%
    \protect\addtolength{\protect\@tocsectionnumwidth}{2em}}%
}
\makeatother

\makeatletter
\settocsectionnumwidth{2em}
\settocsubsectionnumwidth{2.5em}
\settocsubsubsectionnumwidth{3em}
\settocsectionindent{1pc}%
\settocsubsectionindent{\dimexpr\@tocsectionindent+\@tocsectionnumwidth}%
\settocsubsubsectionindent{\dimexpr\@tocsubsectionindent+\@tocsubsectionnumwidth}%
\makeatother

\settocsectionvskip{10pt}
\settocsubsectionvskip{3pt}
\settocsubsubsectionvskip{0pt}
    
\settocsectionformat{\bfseries}
\settocsubsectionformat{\mdseries}
\settocsubsubsectionformat{\mdseries}
\setsectionpnumformat{\bfseries}
\setsubsectionpnumformat{\mdseries}
\setsubsubsectionpnumformat{\mdseries}

\let\oldtableofcontents\tableofcontents
\renewcommand{\tableofcontents}{%
  \vspace*{-\linespacing}% Default gap to top of CONTENTS is \linespacing.
  \oldtableofcontents}

\DeclareSymbolFontAlphabet{\amsmathbb}{AMSb}% to have different mathbb -- we need them for greek letters if I'm not mistaken

\DeclareFontFamily{U}{stix2bb}{\skewchar\font127 }
\DeclareFontShape{U}{stix2bb}{m}{n} {<-> stix2-mathbb}{}
\DeclareMathAlphabet{\mathbbnew}{U}{stix2bb}{m}{n}

\usepackage[
backref=true,
backrefstyle=none,
backend=biber,
style=alphabetic,
citestyle=alphabetic,
maxnames=50,
maxalphanames=7
]{biblatex}

\DefineBibliographyStrings{english}{
  backrefpage = {cited on page},% originally "cit. on p."
  backrefpages = {cited on pages},% originally "cit. on pp."
}

\usepackage[linktocpage=true]{hyperref}

\mathtoolsset{showonlyrefs,showmanualtags}

\allowdisplaybreaks

\newtheorem{theorem}{Theorem}
\newtheorem{proposition}[theorem]{Proposition}
\newtheorem{corollary}[theorem]{Corollary}
\newtheorem{lemma}[theorem]{Lemma}

\newtheorem*{notation}{Notation}

\theoremstyle{definition}
\newtheorem{definition}[theorem]{Definition}
\newtheorem{remark}[theorem]{Remark}
\newtheorem{example}[theorem]{Example}

\numberwithin{theorem}{section}
\numberwithin{equation}{section}
\DeclareDocumentCommand{\faktor}{s m O{0.5} m O{-0.5}}{% \newfaktor[*]{#2}[#3]{#4}[#5] -> #2/#4
  \setbox0=\hbox{\ensuremath{#2}}% Store numerator
  \setbox1=\hbox{\ensuremath{\diagup}}% Store slash /
  \setbox2=\hbox{\ensuremath{#4}}% Store denominator
  \raisebox{#3\ht1}{\usebox0}% Numerator
  \mkern-5mu% Slash /
    {\rotatebox{-44}{\rule[#5\ht2]{0.4pt}{-#5\ht2+#3\ht0+\ht0}}}% tilted rule as a slash
  \mkern-4mu%
  \raisebox{#5\ht2}{\usebox2}% Denominator
}
\DeclareDocumentCommand{\doublefaktor}{s m O{0.5} m O{-0.5}}{% \newdoublefaktor[*]{#2}[#3]{#4}[#5] -> #2//=#4
  \setbox0=\hbox{\ensuremath{#2}}% Store numerator
  \setbox1=\hbox{\ensuremath{\diagup}}% Store slash /
  \setbox2=\hbox{\ensuremath{#4}}% Store denominator
  \raisebox{#3\ht1}{\usebox0}% Numerator
  \mkern-5mu% Space between numerator and double slashes
  {\rotatebox{-44}{\rule[#5\ht2]{0.4pt}{-#5\ht2+#3\ht0+\ht0}}}% First tilted rule as part of double slash
  \mkern-14mu% Space between slashes
  {\rotatebox{-44}{\rule[#5\ht2]{0.4pt}{-#5\ht2+#3\ht0+\ht0}}}% Second tilted rule as part of double slash
  \mkern-4mu% Space between slash and denominator
  \raisebox{#5\ht2}{\usebox2}% Denominator
}
\makeatletter
\newcommand{\pushright}[1]{\ifmeasuring@#1\else\omit\hfill$\displaystyle#1$\fi\ignorespaces}
\newcommand{\pushleft}[1]{\ifmeasuring@#1\else\omit$\displaystyle#1$\hfill\fi\ignorespaces}
\makeatother
\input pdf-trans
\newbox\qbox
\def\usecolor#1{\csname\string\color@#1\endcsname\space}
\newcommand\bordercolor[1]{\colsplit{1}{#1}}
\newcommand\fillcolor[1]{\colsplit{0}{#1}}
\newcommand\colsplit[2]{\colorlet{tmpcolor}{#2}\edef\tmp{\usecolor{tmpcolor}}%
  \def\tmpB{}\expandafter\colsplithelp\tmp\relax%
  \ifnum0=#1\relax\edef\fillcol{\tmpB}\else\edef\bordercol{\tmpC}\fi}
\def\colsplithelp#1#2 #3\relax{%
  \edef\tmpB{\tmpB#1#2 }%
  \ifnum `#1>`9\relax\def\tmpC{#3}\else\colsplithelp#3\relax\fi
}
\newcommand\outline[1]{\leavevmode%
  \def\maltext{#1}%
  \setbox\qbox=\hbox{\maltext}%
  \boxgs{Q q 2 Tr \thickness\space w \fillcol\space \bordercol\space}{}%
  \copy\qbox%
}
\newcommand\mathcalbb[2][1]{%
   \stackengine{0pt}{\outline{$\mathcal{#2}$}}{\kern.3pt\outline{$\mathcal{#2}$}}{O}{l}{F}{F}{L}}
\bordercolor{black}
\fillcolor{white}
\def\thickness{0.35}% TO CHANGE THICKNESS OF SHADOW
\newcommand{\dstar}{\mathrel{\star\star}}
\newcommand{\orep}{\mathcal{O}(\operatorname{Rep}_N(A))}
\newcommand{\rep}{\operatorname{Rep}_N(A)}
\newcommand{\GL}{\operatorname{GL}_N}
\newcommand{\Mat}{\operatorname{Mat}_N^*(\Bbbk)}
\newcommand{\Matt}{\operatorname{Mat}_N(\Bbbk)}
\newcommand{\TMat}{\operatorname{T}\left(\operatorname{Mat}_N^*(\Bbbk)\right)}
\newcommand{\s}{\operatorname{S}(A_{\natural})}
\renewcommand{\S}{\amsmathbb{S}}
\newcommand{\SSym}{\operatorname{Sym_{\S}}}
\renewcommand{\O}{\mathcalbb{O}}
\newcommand{\U}{\mathcalbb{U}}
\newcommand{\E}{\mathcalbb{E}}
\newcommand{\ored}{\O^{\operatorname{red}}}

\renewcommand{\tr}{\operatorname{tr}_N}
\renewcommand{\a}{\alpha}
\renewcommand{\b}{\beta}
\newcommand{\cc}{\gamma}
\newcommand{\R}{\mathfrak{R}}
\newcommand{\1}{\widehat{1}}
\newcommand{\vect}{\operatorname{vect}}
\newcommand{\covect}{\operatorname{covect}}
\newcommand{\longtwoheadrightarrow}{\xtwoheadrightarrow[\phantom{a}]{}}
\newcommand{\longhookrightarrow}{\xhookrightarrow[\phantom{aaa}]{}}
\newcommand{\End}{\E nd}
\newcommand{\Endred}{\E nd^{\operatorname{red}}}
\newcommand{\Com}{\operatorname{Com}}

\newcommand{\stimes}{\mathbin{\otimes_{\scalebox{0.6}{$\S$}}}}
\newcommand{\Sh}{\operatorname{Sh}}

\newcommand{\id}{\operatorname{id}}
\newcommand{\tildeD}{\widetilde{\mathcal{D}}}
\newcommand{\D}{\mathcal{D}}

\newcommand{\A}{\mathcal{A}}
\newcommand{\mathbbp}{\mathbbnew{p}}
\DeclareMathSymbol{\mh}{\mathord}{operators}{`\-}
\newcommand{\Ldi}{L_{di\mh As}}
\newcommand{\Ldistar}{L_{di\mh star}}
\newcommand{\Ldinorm}{L_{1,di\mh As}}
\newcommand{\Ldistarnorm}{L_{1,di\mh star}}
\newcommand{\Lpoisdi}{L_{di\mh Pois}}
\newcommand{\kh}{\Bbbk[[\hbar]]}
\newcommand{\g}{\mathfrak{g}}
\newcommand{\Hom}{\operatorname{Hom}}
\newcommand{\Homadm}{\operatorname{Hom}_{\S}^{\operatorname{adm}}}
\newcommand{\Hompadm}[1]{\operatorname{Hom}_{\S}^{\operatorname{#1-adm}}}

\newcommand{\bulletrn}{\bullet_{\mathsmaller{RN}}}
\newcommand{\ord}{\operatorname{ord}}

\newcommand{\ii}{\mathbbnew{i}}

\definecolor{blue-violet}{rgb}{0.54, 0.17, 0.89}
\definecolor{darkpastelgreen}{rgb}{0.01, 0.75, 0.24}
\definecolor{darkpowderblue}{rgb}{0.0, 0.2, 0.6}
\definecolor{dodgerblue}{rgb}{0.12, 0.56, 1.0}
\definecolor{denim}{rgb}{0.08, 0.38, 0.74}
\definecolor{amber}{rgb}{1.0, 0.75, 0.0}
\definecolor{applegreen}{rgb}{0.55, 0.71, 0.0}
\title{What is a double star-product?}

\date{}

\author{Nikita Safonkin}

\address{Institute of Mathematics, Leipzig University, Augustusplatz 10, 04109 Leipzig, Germany.}

\email{safonkin.nik@gmail.com}

\begin{document}

\begin{abstract}
   Double Poisson brackets, introduced by M. Van den Bergh in 2004, are noncommutative analogs of the usual Poisson brackets in the sense of the Kontsevich-Rosenberg principle: for any $N$, they induce Poisson structures on the space of $N$-dimensional representations $\rep$ of an associative algebra $A$. The problem of deformation quantization of double Poisson brackets was raised by D. Calaque in 2010 and has remained open since then.
   
   In this paper, we solve this problem by presenting a structure on $A$ that induces a star-product under the representation functor and therefore, according to the Kontsevich-Rosenberg principle, can be viewed as an analog of star-products in noncommutative geometry. We also provide an explicit example for the noncommutative affine space $A=\Bbbk\langle x_1,\ldots,x_d\rangle$ and prove an analog of the famous formality theorem of M.Kontsevich. Along the way, we invert the Kontsevich-Rosenberg principle by introducing the notion of a double algebra over an arbitrary operad.

   \hfill\textit{Preliminary draft.}
\end{abstract}

\maketitle

    \setcounter{tocdepth}{2}
    \tableofcontents
    
    \section{Introduction}

In this paper, we work within a realm of noncommutative geometry guided by the Kontsevich-Rosenberg principle (see \cite{kontsevich2000noncommutative}; see also \cite{kontsevich1993formal}), a principle that helps identify ``geometric'' structures on associative unital algebras. For such an algebra $A$, it reads as follows:
\begin{equation}
    \hspace{-0ex}\text{\parbox{0.85\textwidth}{\textit{A noncommutative structure of some kind on $A$ should give an analogous ``commutative'' structure on all schemes $\rep$ for $N \ge 1$.}}}
    %\label{KR}
\end{equation}

An example of such a noncommutative structure is the double Poisson structure introduced by M. Van den Bergh in 2004 (see \cite{van2008double}). Roughly speaking, a double Poisson bracket is a linear map $A \otimes A \longrightarrow A \otimes A$ that is skew-symmetric in a certain sense and satisfies analogs of the Leibniz rule and the Jacobi identity. Any double Poisson bracket on $A$ induces a (conventional) Poisson bracket on the representation spaces $\rep$ for any $N \ge 1$. Another noncommutative structure closely related to double Poisson brackets, and also satisfying the Kontsevich-Rosenberg principle, is the bi-symplectic structure introduced by W. Crawley-Boevey, P. Etingof, and V. Ginzburg (see \cite{crawley2007noncommutative}).

The problem of quantizing conventional Poisson brackets is well-studied; see, for instance, \cite{bordemann2008deformation}, \cite{cattaneo2005deformation}, along with references therein. We only recall that star-products were introduced in \cite{bayen1978deformation1}, and in \cite{kontsevich2003deformation} M. Kontsevich proved that any finite-dimensional Poisson manifold can be canonically quantized.

Once a definition of what one should consider a noncommutative Poisson bracket is given, it is natural to ask how to quantize this structure. This question was explicitly formulated by D. Calaque and posted on MathOverflow in 2010 (see \cite{mathoverflowCalaque}). This question was also mentioned in Section~4 of \cite{odesskii2013double}, see also a comment after Definition~4.1 in \cite{schedler2009poisson}.

In the present paper, we work with the category of diagonal $\S$-bimodules, which are graded vector spaces equipped with a bimodule action of the symmetric group $S(n)$ on the $n$-th graded component. All the structures we consider are ``di-twisted'' in the sense that they are twisted with respect to the adjoint $\S$-module structure; see \cite{barratt1978twisted, stover1993equivalence}. The idea of working with $\S$-modules in a similar context goes back to T. Schedler. It appeared in the first version of the paper \cite{schedler2009poisson}, which was published on arXiv in 2006 -- only two years after double Poisson brackets were introduced by M. Van den Bergh. T. Schedler found a connection between double Poisson brackets on $A$ and twisted Poisson brackets on the tensor algebra $\operatorname{T}_{\Bbbk}(A)$, viewed as a free commutative algebra in the category of $\S$-modules. This connection was pushed further in \cite{ginzburg2010differentialoperatorsbvstructures}, where V. Ginzburg and T. Schedler worked with diagonal $\S$-bimodules equipped with an additional structure; they called these \textit{wheelspaces}. Among other results, they introduced a commutative algebra in the category of wheelspaces $\mathcal{F}(A)$, which is commutative viewed as an algebra in the category of diagonal $\S$-bimodules, and mentioned, without proof, that any double Poisson bracket on $A$ can be canonically extended to a di-twisted Poisson bracket on $\mathcal{F}(A)$. This was proved in a recent paper \cite{fernández2025symplecticwheelgebrasnoncommutativegeometry} by D. Fernandez and E. Herscovich; see Proposition~6.8 therein. The authors of \cite{ginzburg2010differentialoperatorsbvstructures} also defined differential operators on $\mathcal{F}(A)$; we will adopt their definition with minor adjustments.

In the present paper, we denote the algebra $\mathcal{F}(A)$ by $\O(A)$ and refer to it as \textit{the double coordinate ring of $A$}, see Section~\ref{section_noncommutative_functions} for details. We first came across $\O(A)$ through lengthy computations with representation spaces presented in Section \ref{section_example}, and only later realized that it had already appeared in the literature. This double coordinate ring satisfies the Kontsevich-Rosenberg principle, i.e. there is a way to produce elements of $\orep$ for any $N$ out of elements of $\O(A)$. Then we can consider the ideal of elements vanishing in any finite dimensional representation, denoted by $\R(A)$. We mainly work with the quotient of the double coordinate ring by the ideal $\R(A)$, which we treat as the reduced version of the double coordinate ring, and denote it by $\ored(A)$. Note that for the Weyl algebra $A=\Bbbk\langle x,\partial_x\rangle$ we have $\R(A)=\O(A)$ as expected, and for the free polynomial algebra $A=\Bbbk\langle x_1,\ldots,x_d\rangle$ the ideal $\R(A)$ is zero, hence the double coordinate ring $\O(A)$ is reduced. The main advantage of the reduce double coordinate ring $\ored(A)$ over the non-reduced $\O(A)$ is that the former satisfies the following property: if an element of $\ored(A)$ is zero in any finite dimensional representation, then it's zero as an element of $\ored(A)$, cf. Proposition 11.1.1 in \cite{crawley2007noncommutative}, see also \cite{ginzburg2006moyal}. 

Another reason to treat $\O(A)$ and $\ored(A)$ as noncommutative analogs of coordinate rings $\orep$ is that they naturally replace $\orep$ in noncommutative analogs of many $\GL$-equivariant structures on $\rep$; this is precisely true at least when the $\GL$-equivariant structure on $\rep$ is given by a $\GL$-equivariant $\mathcal{P}$-algebra structure on $\orep$ for an operad $\mathcal{P}$. This is due to the fact that any $\GL$-equivariant linear map $\orep^{\otimes n}\longrightarrow \orep$ arises under the action of the representation functor from a linear map $\ored(A)^{\otimes n}\longrightarrow\ored(A)$ satisfying some additional properties, see Theorem \ref{th3} below. We will call these linear maps \textit{admissible}. Moreover, this construction satisfies the following property: if an admissible map $\ored(A)^{\otimes n}\longrightarrow\ored(A)$ is zero in any finite dimensional representation, i.e. the induced maps $\orep^{\otimes n}\longrightarrow\orep$ are zero for any $N$, then the initial admissible map is zero too. We will refer to this phenomenon as \textit{asymptotical injectivity}. We borrowed this terminology from \cite{ginzburg2006moyal}. 

We will use asymptotical injectivity for $\ored(A)$ and for admissible linear maps $\ored(A)^{\otimes n}\longrightarrow\ored(A)$ many times in the paper and not only for inspiration -- almost all our proofs are based on it. For example, we do not use the Kontsevich-Rosenberg principle only to define noncommutative structures, but also to prove properties about them (after the structures in question have been defined, of course). On the level of computations it looks like we go along the Kontsevich-Rosenberg principle in the opposite direction -- not passing from the noncommutative side to the commutative one via the representation functor, but lifting our commutative structures to the noncommutative world, and once we've lifted them, all the necessary properties of the lift are proven automatically due to the asymptotical injectivity. For instance, if we know that the multiplication in $\orep$ lifts to the multiplication in $\ored(A)$, we don't need to check that the latter is associative, because the associator of the multiplication in $\ored(A)$ maps under the representation functor to the associator of the multiplication in $\orep$, which is zero for any $N$. Hence by the asymptotical injectivity, the associator of the multiplication in $\ored(A)$ must be zero. Of course, it may not be possible to lift each and every structure on $\rep$ to $\ored(A)$, because the initial structure must be "universal" in a sense, by which we mean "stable in $N$". That's why we usually lift the structures we're interested in by hand or explain why such a lift exists and how to obtain it. Then, as mentioned above, all the expected properties of the lift follow automatically and effortlessly.

Note that the point of view described above is in agreement with a remark that led M.Van den Bergh to the notion of double Poisson brackets. Namely, he noticed that when computations with certain quiver varieties are suitably organized, explicit matrices appear only at the very last step, as an afterthought, see the third paragraph below Theorem 1.1 in \cite{van2008double} and compare this to Theorem \ref{th3} below. In this sense, asymptotical injectivity provides a mathematical tool for systematically studying more general effects of this kind.

The reduced double coordinate ring $\ored(A)$ seems to share some properties of the coordinate ring of the reduction of an affine scheme, which makes us think that the following analogy is appropriate. It, particularly, suggests that the noncommutative geometry in the sense of the Kontsevich-Rosenberg principle is actually commutative but in a different symmetric monoidal category -- not in the category of vector spaces but in the category of diagonal $\S$-bimodules.

\renewcommand{\arraystretch}{1.7}
\setlength{\arrayrulewidth}{0.25mm}
\setlength{\tabcolsep}{10pt}
\begin{table}[h]
    \centering
    \begin{tabular}{||c||c||}
    \hline
        Affine scheme of finite type $X$ over $\Bbbk=\overline{\Bbbk}$ & $\scalebox{1.3}{?}$ Noncommutative scheme associated to $A$ $\scalebox{1.3}{?}$ \\ [5pt]
       \hline
       \hline
        Coordinate ring $\mathcal{O}(X)$ & Double coordinate ring $\O(A)$  \\
       \hline
        Closed points of $X$ & Representation schemes $\big\{\rep\big\}_{N\geq 1}$ \\
       \hline
        Nilradical of $\mathcal{O}(X)$ & The ideal $\R(A)$\\
       \hline
       Coordinate ring of the reduction $\mathcal{O}(X_{\operatorname{red}})$ &Reduced double coordinate ring $\ored(A)$ \\
       \hline
    \end{tabular}
    %\caption{An analogy between commutative and noncommutative pictures}
    \label{tab1}
\end{table}

Recall that closed points of $X$ are in bijection with $\Hom_{\operatorname{Alg}}(R,\Bbbk)$, where $X=\operatorname{Spec(R)}$ for a finitely generated $\Bbbk$-algebra $R$, and $f\in \mathcal{O}(X_{\operatorname{red}})$ vanishes at any closed point if and only if $f=0\in \mathcal{O}(X_{\operatorname{red}})$, which is similar to the picture with $\O(A)$, $\ored(A)$, and representation spaces $\rep$.

We did not emphasize this earlier, but there are two different conventions about representation spaces: some authors consider only unital representations, i.e. unital homomorphisms from a unital associative algebra $A$ to the matrix algebra, while others consider nonunital representations of nonunital algebras. Usually, the difference is minor: the unit of $A$ does not explicitly enter into any formulas or structures of interest, nor does it play a role on the side of representation spaces. So, it's usually easy to pass from the unital setting to the nonunital one. However, this is not the case in the present paper. We do need the algebra $A$ to be unital and we do need to work with unital representations. The reason is Theorem \ref{th3} -- its proof explicitly uses the unit of $A$. The claim boils down to describing $\GL$-equivariant linear maps $\Mat^{\otimes n}\longrightarrow \Mat^{\otimes m}$, which, of course, follows from the invariant theory for $\GL$. This description is unavoidably tied to the counit $\varepsilon:\Mat\longrightarrow\Bbbk$, which gives one of the simplest examples of a such $\GL$-equivariant linear map for $n=1$, $m=0$; see formula \eqref{f98} below. The key point is that this counit is responsible for the relation in $\orep$ involving the unit of $A$, see formula \eqref{f36} below. So, we cannot give up the unit.

\subsection{Main results}
There are three main results of the present paper -- two definitions and one theorem:

\begin{itemize}
    \item the definition of a double algebra over an arbitrary operad, Definition \ref{def8};
    
    \item the definition of a double star-product, which we call a \textit{di-twisted} star-product, Definition~\ref{def11};

    \item the double formality theorem for the noncommutative affine space $A=\Bbbk\langle x_1,\ldots,x_d\rangle$, Theorem \ref{th4}.
\end{itemize}
All of them satisfy the Kontsevich-Rosenberg principle. Let us briefly explain them. Set
\begin{align}
        \O(A):=\SSym(A_{\natural}\oplus A[-1])=\bigoplus\limits_{n\geq 0}A^{\otimes n}\otimes \s\otimes \Bbbk[S(n)],
    \end{align}
    where $A_{\natural}:=\faktor{A}{[A,A]}$ and $[A,A]$ is the vector space spanned by $ab-ba$ for $a,b\in A$. The subscript $\S$ stands for the category of (diagonal) $\S$-bimodules, and $[-1]$ is the shift to the right by $1$. 
    
    This $\O(A)$ comes equipped with two homogeneous linear maps $\pi,\widehat{1}:\O(A)\longrightarrow\O(A)$ of degrees $-1$ and $+1$ induced by the multiplication and unity in $A$. In a sense, these maps are responsible for the relations in $\orep$.
    
    Next, one can consider a two-sided ideal of $\O(A)$ that consists of all elements vanishing in any representation of $A$. We denote this ideal by $\R(A)$. There are two extreme case: when $A=\Bbbk\langle x_1,\ldots,x_d\rangle$, this ideal is zero, see Proposition \ref{prop5} in the text below, and when $A$ is the Weyl algebra in any number of variables, this ideal is the whole algebra $\O(A)$. The quotient of $\O(A)$ by this ideal $\R(A)$ is the central object of this paper, we denote it by $\ored(A)$.

    \subsubsection{Double \texorpdfstring{$\mathcal{P}$}{P}-algebras}
    We denote by $\Homadm\big(\ored(A)^{\stimes n},\ored(A)\big)$ the subset of $\S$-bimodule homomorphisms $\ored(A)^{\stimes n}\longrightarrow\ored(A)$ satisfying certain compatibility conditions with $\pi$ and $\1$. We call these $\S$-bimodule homomorphisms \textit{admissible}. The compatibility conditions are algebraic relations of some kind which literally mean that the $\S$-bimodule homomorphism induces a well-defined linear map $\orep^{\otimes n}\longrightarrow\orep$, which is (automatically) $\GL$-equivariant. Theorem \ref{th3} in the text below shows also that (almost) any $\GL$-equivariant linear map $\orep^{\otimes n}\longrightarrow\orep$ is of such form; see also Theorem \ref{th2}, which generalizes Theorem \ref{th3} to $\GL$-equivariant linear maps $\orep^{\otimes n}\longrightarrow\orep\otimes \Matt^{\otimes p}$. These two theorems are our milestones, because they show that if anyone wants to work with $\GL$-equivariant structures on $\rep$, then they don't have much choice but to work with $\ored(A)$ (or $\O(A)$).
    
    Admissible $\S$-bimodule homomorphisms form an operad 
    \begin{align}
    \End_A^{\operatorname{red}}:=\Bigg(\Homadm\left(\ored(A)^{\stimes n},\ored(A)\right)\Bigg)_{n\geq 0}
    \end{align}
    which we call \textit{the double endomorpism operad}, because the representation functor induces a homomorphism of operads 
    \begin{equation}
        \End_A^{\operatorname{red}}\longrightarrow \mathcal{E}nd_{\orep}^{\operatorname{GL}_N}
    \end{equation}
    where on the right-hand side one has the $\GL$-equivariant version of the usual endomorphism operad. Moreover, by the very definition of the ideal $\R(A)$, these operad homomorphisms considered together for all $N$ give an injective operad homomorphism
    \begin{align}\label{f92}
        \End_A^{\operatorname{red}}\longhookrightarrow\prod\limits_{N\geq 1}\mathcal{E}nd_{\orep}^{\operatorname{GL}_N}.
    \end{align}

    Thus, we say that an associative algebra $A$ is a double algebra over an operad $\mathcal{P}$ if, roughly speaking, there is an operad homomorphism $\mathcal{P}\longrightarrow\End_A^{\operatorname{red}}$, see Definition \ref{def8}. Corollary \ref{cor3} of Theorem \ref{th3} explains why this definition is in a sense maximal, see Remark \ref{rem7}.

    \subsubsection{Double star-products and the double formality theorem}

    Now we will apply the discussion above to deformations and star-products. Let us list the main objects we work with
    \begin{itemize}[leftmargin=3ex]
        \item by a \textit{di-twisted Poisson bracket} we mean
        an admissible $\S$-bimodule homomorphism 
    \begin{align}
        \ored(A)\stimes \ored(A)\longrightarrow\ored(A)
    \end{align}
    satisfying certain diagonally-twisted versions of diagonally cyclic anti-symmetry, Leibniz rule, and Jacobi identity, see Definition \ref{def12};
        
        \item by a \textit{di-twisted deformation} of the associative multiplication in $\ored(A)$ we mean an associative admissible $\S$-bimodule homomorphism 
    \begin{align}
        \dstar:\ored(A)[[\hbar]]\stimes \ored(A)[[\hbar]]\longrightarrow\ored(A)[[\hbar]],
    \end{align}
    which coincides with the standard multiplication in $\ored(A)$ modulo $\hbar$, see Definition \ref{def7}. We consider these di-twisted deformations modulo the gauge equivalence realized by admissible $\S$-bimodule homomorphism. One can then easily check that the di-twisted commutator 
    \begin{align}
        [\a,\b]_{\dstar}:=\a\dstar\b-\operatorname{Ad}((12)^{|\b|,|\a|})\b\dstar\a
    \end{align}vanishes in the zero-th order and defines a di-twisted Poisson bracket in the first order, making it possible to speak about deformations of a given di-twisted Poisson bracket;
    
    \item by differential operators on $\ored(A)$ we, roughly speaking, mean those admissible $\S$-bimodule homomorphism $\ored(A)\longrightarrow\ored(A)$ that induce differential operators in any $\rep$. Such a definition is rather implicit, see Definition \ref{def9} for an explicit one, and Proposition \ref{prop21} for equivalence of the definitions. Note that Definition \ref{def9} almost repeats the definition of differential operators from \cite{ginzburg2010differentialoperatorsbvstructures} with only minor changes related to $\O(A)\rightsquigarrow\ored(A)$ and the fact that, in contrast to \cite{ginzburg2010differentialoperatorsbvstructures}, we work with unital representations.

    \item by $n$-poly-differential operators on $\ored(A)$ we mean those admissible $\S$-bimodule homomorphisms $\ored(A)^{\stimes n}\longrightarrow\ored(A)$ that induce $n$-poly-differential operators in any $\rep$. An explicit version of this definition is Definition \ref{def10}, see also Proposition \ref{prop20} for the equivalence. Note that V. Ginzburg and T. Schedler didn't work with poly-differential operators in \cite{ginzburg2010differentialoperatorsbvstructures}.

    \item by a double star-product, which we call \textit{di-twisted star-product}, we mean a di-twisted deformation of the standard associative multiplication in $\ored(A)$ realized by bi-differential operators, see Definition \ref{def11}. We consider the di-twisted star-products modulo the gauge equivalence realized by differential operators.
    \end{itemize}

    Due to the discussion in the previous subsection, all these structures satisfy the Kontsevich-Rosenberg principle.

    Next, we introduce two dg Lie algebras $\Ldi(\ored(A))$ and $\Ldistar(\ored(A))$ which control di-twisted deformations and di-twisted star-products, Definition \ref{def13}. Roughly speaking, the first one consists of admissible $\S$-bimodule homomorphism $\ored(A)^{\stimes n}\longrightarrow\ored(A)$ for all $n$, and the second one consists of poly-differential operators. They both satisfy the Kontsevich-Rosenberg principle, since there are dg Lie algebra homomorphisms 
    \begin{align}
        \Ldi(\ored(A))&\longrightarrow L_{As}\Big(\orep\Big)^{\GL}\\
        \Ldistar(\ored(A))&\longrightarrow L_{star}\Big(\orep\Big)^{\GL}
    \end{align}
    induced by the representation functor.

    It is well-known that any deformation of an associative multiplication is gauge equivalent to a unital deformation. The same is true for di-twisted deformations and di-twisted star-products, we prove this in Section \ref{section_normalized_cochains}.

    Next, for any $\{-,-\}$ di-twisted Poisson bracket on $\ored(A)$, we introduce a dg Lie algebra $\Lpoisdi(\ored(A))$ consisting of di-twisted poly-vector fields, which are skew-symmetric di-twisted derivations, see Definition \ref{def14}. This dg Lie algebra is responsible for the deformations of $\{-,-\}$ modulo gauge equivalence infinitesimally realized by di-twisted derivations. Recall that by results of V. Ginzburg, T. Schedler and D. Fernandez, E. Herscovich any double bracket on $A$ can be canonically extended to a di-twisted Poisson bracket on $\O(A)$ and, obviously, on $\ored(A)$. So, $\Lpoisdi(\ored(A))$ is well-defined for any double Poisson bracket on $A$ as well. Moreover, these poly-vector fields satisfy the Kontsevich-Rosenberg principle, as there is a homomorphism of dg Lie algebras 
    \begin{align}
    L_{di\mh Pois}(\ored(A))\longrightarrow L_{Pois}\Big(\orep\Big)^{\GL}
    \end{align}
    induced by the representation functor. On the right-hand side the differential is given by the induced Poisson bracket $\{-,-\}_N$ corresponding to the di-twisted Poisson bracket $\{-,-\}$ on $\ored(A)$.

    Finally, we consider a special case $A=\amsmathbb{C}\langle x_1,\ldots,x_d\rangle$, for which we don't have to distinguish $\ored(A)$ and $\O(A)$ due to Proposition \ref{prop5}. We prove the double formality theorem, Theorem \ref{th4}, which states that $\Ldistar(\O(A))$ is $L_{\infty}$-quasi-isomorphic to its cohomology, which is equal to $\Lpoisdi(\O(A))$ with the zero differential. We prove it by presenting a noncommutative lift of Kontsevich's universal formula, i.e. an $L_{\infty}$-quasi-isomorphism $\U$ making the following diagram commute
    \begin{equation*}
        \begin{tikzcd}
            \Lpoisdi(\O(A))\arrow[r,rightsquigarrow,"\U"]\arrow[d]&[20pt] \Ldistar(\O(A))\arrow[d]\\[10pt]
            L_{Pois}\Big(\orep\Big)^{\GL} \arrow[r,rightsquigarrow,"\mathcal{U}"]& L_{star}\Big(\orep\Big)^{\GL},
        \end{tikzcd}
    \end{equation*}
    where vertical arrows are the dg Lie algebra homomorphisms induced by the representation functor, and $\mathcal{U}$ is the $L_{\infty}$-quasi-isomorphism constructed by M.Kontsevich in 1997, see \cite{kontsevich2003deformation}. More precisely, the $L_{\infty}$-morphism $\mathcal{U}$ is given by its Taylor components
    \begin{align}
        \mathcal{U}_n:=\sum\limits_{m\geq 0}\sum\limits_{\Gamma\in G_{n,m}}W_{\Gamma}\,\mathcal{U}_{\Gamma},
\end{align}
where $W_{\Gamma}\in\amsmathbb{R}$ are the universal coefficients given by integrals over configuration spaces of points on the upper half-plane, and $\mathcal{U}_{\Gamma}$ are certain explicitly defined differential-algebraic expressions. Our formula for $\U$ reads then as follows
\begin{align}
    \U_n:=\sum\limits_{m\geq 0}\sum\limits_{\Gamma\in G_{n,m}}W_{\Gamma}\,\U_{\Gamma},
\end{align}
so we didn't change the coefficients $W_{\Gamma}$ at all, but constructed $\U_{\Gamma}$ as a noncommutative lift of $\mathcal{U}_{\Gamma}$.
    
This is essentially the whole proof of the $L_{\infty}$-part -- we don't have to perform any computations with Taylor coefficients of the lift to prove that they satisfy the $L_{\infty}$-quadratic relation due to the asymptotical injectivity phenomenon. The quasi-isomorphism part follows from a noncommutative analog of an explicit homotopy constructed in \cite{de1995homotopy}. As a corollary, we have that any di-twisted Poisson bracket on $\O(A)$ and, particularly, any double Poisson bracket on $A$ can be canonically quantized.
    
\subsection{Outline of the paper.}
In Section \ref{section_preliminaries}, we fix the notation. In Section \ref{section_s_bimodules}, we recall the main definitions concerning diagonal $\S$-bimodules. In Section \ref{section_noncommutative_functions}, we discuss the algebra $\O(A)$ and its connection to representation spaces $\rep$. In Section \ref{section_inverting_KR}, we discuss the Kontsevich-Rosenberg principle in more detail and propose a way of inverting it by giving a definition of double algebras over an arbitrary operad. In Section \ref{section_equivariant_linear_maps}, we define and examine the main structures on $\O(A)$ we are working with: di-twisted Poisson brackets, di-twisted deformations, di-twisted star-products, and dg Lie algebras of di-twisted poly-vector fields and di-twisted poly-differential operators controlling the deformations, along with their normalized (vanishing at unity) analogs controlling unital deformations. In Section \ref{section_example}, we provide an example of a di-twisted star-product on $\O(\Bbbk\langle x_1,\ldots,x_d\rangle)$ given by a very explicit but long combinatorial formula. Finally, in Section \ref{section_double_formality}, we prove a double formality theorem for $A=\Bbbk\langle x_1,\ldots,x_d\rangle$ by presenting an $L_{\infty}$-quasi-isomorphism between dg Lie algebras of di-twisted poly-vector fields and di-twisted poly-differential operators on the noncommutative affine space $N\amsmathbb{A}^d$.

\textbf{Acknowledgments.}
I would like to thank Anton Khoroshkin for a helpful discussion on Hochschild cohomology and Maxime Fairon for many discussion on double Poisson brackets. These conversations were very helpful to me. I am also grateful to Victor Ginzburg for pointing out the importance of the unit in $A$ for Theorem \ref{th3}, and to Will Sawin for Lemma \ref{lemma20} and his clarifications. I would like to thank Michael Pevzner as well for providing useful references, and Ivan Kaygorodov for Remark \ref{rem6}. I am especially grateful to Grigori Olshanski for his constant support and many valuable discussions, and to Alexey Bufetov for his support and for providing excellent working conditions at Leipzig University. The author was partially supported by the European Research
Council (ERC), Grant Agreement No. 101041499.

\section{Preliminaries and notation}\label{section_preliminaries}

Throughout the paper we assume that $\Bbbk$ is an algebraically closed field of characteristic zero. We denote the tensor product of vector spaces over $\Bbbk$ by $\otimes$. We will also need the following Hadamard tensor product of graded algebras.

\begin{definition}
    Let $X=\bigoplus\limits_{n\geq 0} X_n$ and $Y=\bigoplus\limits_{n\geq 0} Y_n$ be graded vector spaces. By their \textit{Hadamard tensor product} $X\odot Y$ we mean the following graded subspace $X\odot Y=\bigoplus\limits_{n\geq 0} X_n\otimes Y_n$ of the usual tensor product $X\otimes Y$. If $X$ and $Y$ are graded algebras, so is $X\odot Y$.
\end{definition}

The main property of this Hadamard tensor product for us is
\begin{align}
    \operatorname{T}(V\otimes W)=\operatorname{T}(V)\odot\operatorname{T}(W),
\end{align}
where $V$ and $W$ are vector spaces and $\operatorname{T}(-)$ stands for the tensor algebra. 

\subsection{Permutations and symmetric groups}
Recall that the symmetric group $S(n)$ embeds in $S(n+1)$ as the subgroup fixing the last node. More generally, for any $n,m$ there is a homomorphism $S(n)\times S(m)\hookrightarrow S(n+m)$ that writes the second permutation from $S(m)$ next to the first one from $S(n)$. Namely, the resulting permutation from $S(n+m)$ permutes the first $n$ nodes as the first permutation would do, and the rest $m$ nodes as the second permutation. We will denote this homomorphism by $(w,\tau)\mapsto w\times \tau$. 

For any permutation $\tau\in S(n)$ and any integers $k_1,\ldots,k_n\in\mathbb{Z}_{\geq 0}$, we denote by $\tau^{k_1,\ldots,k_n}\in S(k_1+\ldots+k_n)$ the permutation that splits all $k_1+\ldots+k_n$ nodes in $n$ blocks of lengths $k_1,\ldots,k_n$, and permutes these blocks according to $\tau$. In other words, $\tau^{k_1,\ldots,k_n}$ is the permutation obtained from $\tau$ by blowing up its edges: the $i$-th edge of $\tau$, connecting $i$ and $\tau(i)$, splits in $k_i$ parallel lines.  

For instance, for $n=3$, $\tau=(12)\in S(3)$, and $k_1=5$, $k_2=2$, $k_3=3$ one has 
$$(12)^{5,2,3}=
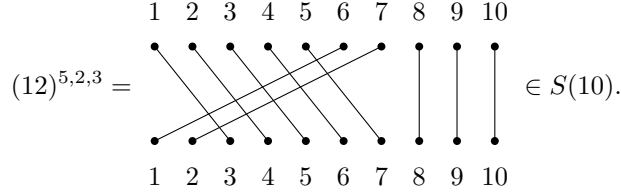
\begin{tikzpicture}[baseline=(A.base), scale=1, every node/.style={transform shape},
dot/.style = {circle, fill, minimum size=#1, inner sep=0pt, outer sep=5pt},
emptydot/.style = {circle, fill=white, minimum size=#1, inner sep=0pt, outer sep=5pt}]

%preliminaries
\def\radd{3pt}% radius of the dots
\def\height{1.25}% the distance between two rows of dots
\def\width{0.5}% the distance between two consecutive dots in the same row

%initializing the top row vertices
\coordinate (1up) at (0,0);
\path (1up) ++(\width,0) coordinate (2up);
\path (2up) ++(\width,0) coordinate (3up);
\path (3up) ++(\width,0) coordinate (4up);
\path (4up) ++(\width,0) coordinate (5up);
\path (5up) ++(\width,0) coordinate (6up);
\path (6up) ++(\width,0) coordinate (7up);
\path (7up) ++(\width,0) coordinate (8up);
\path (8up) ++(\width,0) coordinate (9up);
\path (9up) ++(\width,0) coordinate (10up);

%initializing the bottom row vertices
\path (1up) ++(0,-\height) coordinate (1down);
\path (2up) ++(0,-\height) coordinate (2down);
\path (3up) ++(0,-\height) coordinate (3down);
\path (4up) ++(0,-\height) coordinate (4down);
\path (5up) ++(0,-\height) coordinate (5down);
\path (6up) ++(0,-\height) coordinate (6down);
\path (7up) ++(0,-\height) coordinate (7down);
\path (8up) ++(0,-\height) coordinate (8down);
\path (9up) ++(0,-\height) coordinate (9down);
\path (10up) ++(0,-\height) coordinate (10down);

%the top row dots
\node[dot=\radd,label=above:{ $1$}] at (1up) {};
\node[dot=\radd,label=above:{ $2$}] at (2up) {};
\node[dot=\radd,label=above:{ $3$}] at (3up) {};
\node[dot=\radd,label=above:{ $4$}] at (4up) {};
\node[dot=\radd,label=above:{ $5$}] at (5up) {};
\node[dot=\radd,label=above:{ $6$}] at (6up) {};
\node[dot=\radd,label=above:{ $7$}] at (7up) {};
\node[dot=\radd,label=above:{ $8$}] at (8up) {};
\node[dot=\radd,label=above:{ $9$}] at (9up) {};
\node[dot=\radd,label=above:{ $10$}] at (10up) {};

%the bottom row dots
\node[dot=\radd,label=below:{ $1$}] at (1down) {};
\node[dot=\radd,label=below:{ $2$}] at (2down) {};
\node[dot=\radd,label=below:{ $3$}] at (3down) {};
\node[dot=\radd,label=below:{ $4$}] at (4down) {};
\node[dot=\radd,label=below:{ $5$}] at (5down) {};
\node[dot=\radd,label=below:{ $6$}] at (6down) {};
\node[dot=\radd,label=below:{ $7$}] at (7down) {};
\node[dot=\radd,label=below:{ $8$}] at (8down) {};
\node[dot=\radd,label=below:{ $9$}] at (9down) {};
\node[dot=\radd,label=below:{ $10$}] at (10down) {};

%the black lines
\draw (1up) -- (3down);
\draw (2up) -- (4down);
\draw (3up) -- (5down);
\draw (4up) -- (6down);
\draw (5up) -- (7down);
\draw (6up) -- (1down);
\draw (7up) -- (2down);
\draw (8up) -- (8down);
\draw (9up) -- (9down);
\draw (10up) -- (10down);

\node[](A) at (0,-0.5*\height) {};

\end{tikzpicture}\in S(10).$$

\begin{lemma}\label{lemma8}
    For any $\tau,w\in S(n)$, $k_1,\ldots,k_n\in\mathbb{Z}_{\geq 0}$, and $\sigma_1\in S(k_1),\ldots, \sigma_n\in S(k_n)$ one has
    \begin{align}
        (\tau w)^{k_1,\ldots,k_n}&=\tau^{k_{w^{-1}(1)},\ldots,k_{w^{-1}(n)}}w^{k_1,\ldots,k_n},\\
        w^{k_1,\ldots,k_n}(\sigma_1\times\ldots\times\sigma_n)&=(\sigma_{w^{-1}(1)}\times\ldots\times \sigma_{w^{-1}(n)})w^{k_1,\ldots,k_n}.
    \end{align}
\end{lemma}
\begin{proof}
    The first identity is clear. The second states that there is no difference between permuting blocks of nodes and then permuting nodes within each block, or first permuting nodes within each block and then permuting the blocks, which is also clear.
\end{proof}

\begin{definition}
    For any two permutations $w\in S(n)$, $\tau\in S(m)$, and integers $l,r\in[0,n]$ such that $l+r=n$ we denote by $w'\underset{l}{*}\tau\underset{r}{*}w''$ the permutation obtained by splitting $w$ into two parts of sizes $l$ and $r$, pulling them apart preserving all the edges, and inserting the permutation $\tau$ in the middle.
\end{definition}

For example, if

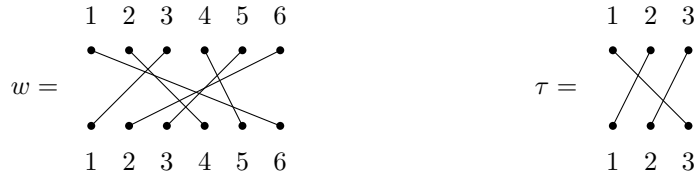
\begin{figure}[ht]
    \centering 
    \begin{tikzpicture}[scale=1, every node/.style={transform shape},
dot/.style = {circle, fill, minimum size=#1, inner sep=0pt, outer sep=5pt},
emptydot/.style = {circle, fill=white, minimum size=#1, inner sep=0pt, outer sep=5pt}]
%%%%%%%%%%%%%%%%%%%%%%%%%%%%%%%% preliminaries %%%%%%%%%%%%%%%%%%%%%%%%%%%%%%%%
\def\radd{3pt}% radius of the dots
\def\height{1}% the distance between two rows of dots
\def\width{0.5}% the distance between two consecutive dots in the same row

%initializing the top row vertices
\coordinate (1up) at (0,0);
\path (1up) ++(\width,0) coordinate (2up);
\path (2up) ++(\width,0) coordinate (3up);
\path (3up) ++(\width,0) coordinate (4up);
\path (4up) ++(\width,0) coordinate (5up);
\path (5up) ++(\width,0) coordinate (6up);

%initializing the bottom row vertices
\path (1up) ++(0,-\height) coordinate (1down);
\path (1down) ++(\width,0) coordinate (2down);
\path (2down) ++(\width,0) coordinate (3down);
\path (3down) ++(\width,0) coordinate (4down);
\path (4down) ++(\width,0) coordinate (5down);
\path (5down) ++(\width,0) coordinate (6down);

%%%%%%%%%%%%%%%%%%%%%%%%%%% the top row dots %%%%%%%%%%%%%%%%%%%%%%%%%%%%%%%%%%

\node[dot=\radd,label={ $1$}] at (1up) {};
\node[dot=\radd,label={ $2$}] at (2up) {};
\node[dot=\radd,label={ $3$}] at (3up) {};
\node[dot=\radd,label={ $4$}] at (4up) {};
\node[dot=\radd,label={ $5$}] at (5up) {};
\node[dot=\radd,label={ $6$}] at (6up) {};

%%%%%%%%%%%%%%%%%%%%%%%%% the bottom row dots %%%%%%%%%%%%%%%%%%%%%%%%%%%%%%%%%

\node[dot=\radd,label=below:{ $1$}] at (1down) {};
\node[dot=\radd,label=below:{ $2$}] at (2down) {};
\node[dot=\radd,label=below:{ $3$}] at (3down) {};
\node[dot=\radd,label=below:{ $4$}] at (4down) {};
\node[dot=\radd,label=below:{ $5$}] at (5down) {};
\node[dot=\radd,label=below:{ $6$}] at (6down) {};

%%%%%%%%%%%%%%%%%%%%%%%%%%%%% the black lines %%%%%%%%%%%%%%%%%%%%%%%%%%%%%%%%%

\draw (1up) -- (6down);
\draw (2up) -- (4down);
\draw (3up) -- (1down);
\draw (4up) -- (5down);
\draw (5up) -- (3down);
\draw (6up) -- (2down);

\node[] at (-0.75,-0.5) {$w=$};

\end{tikzpicture}\hspace{80pt}
    \begin{tikzpicture}[scale=1, every node/.style={transform shape},
dot/.style = {circle, fill, minimum size=#1, inner sep=0pt, outer sep=5pt},
emptydot/.style = {circle, fill=white, minimum size=#1, inner sep=0pt, outer sep=5pt}]
%%%%%%%%%%%%%%%%%%%%%%%%%%%%%%%% preliminaries %%%%%%%%%%%%%%%%%%%%%%%%%%%%%%%%
\def\radd{3pt}% radius of the dots
\def\height{1}% the distance between two rows of dots
\def\width{0.5}% the distance between two consecutive dots in the same row

%initializing the top row vertices
\coordinate (1up) at (0,0);
\path (1up) ++(\width,0) coordinate (2up);
\path (2up) ++(\width,0) coordinate (3up);

%initializing the bottom row vertices
\path (1up) ++(0,-\height) coordinate (1down);
\path (1down) ++(\width,0) coordinate (2down);
\path (2down) ++(\width,0) coordinate (3down);

%%%%%%%%%%%%%%%%%%%%%%%%%%% the top row dots %%%%%%%%%%%%%%%%%%%%%%%%%%%%%%%%%%

\node[dot=\radd,label={ $1$}] at (1up) {};
\node[dot=\radd,label={ $2$}] at (2up) {};
\node[dot=\radd,label={ $3$}] at (3up) {};

%%%%%%%%%%%%%%%%%%%%%%%%% the bottom row dots %%%%%%%%%%%%%%%%%%%%%%%%%%%%%%%%%

\node[dot=\radd,label=below:{ $1$}] at (1down) {};
\node[dot=\radd,label=below:{ $2$}] at (2down) {};
\node[dot=\radd,label=below:{ $3$}] at (3down) {};

%%%%%%%%%%%%%%%%%%%%%%%%%%%%% the black lines %%%%%%%%%%%%%%%%%%%%%%%%%%%%%%%%%

\draw (1up) -- (3down);
\draw (2up) -- (1down);
\draw (3up) -- (2down);

\node[] at (-0.75,-0.5) {$\tau=$};

\end{tikzpicture}
\end{figure}
then
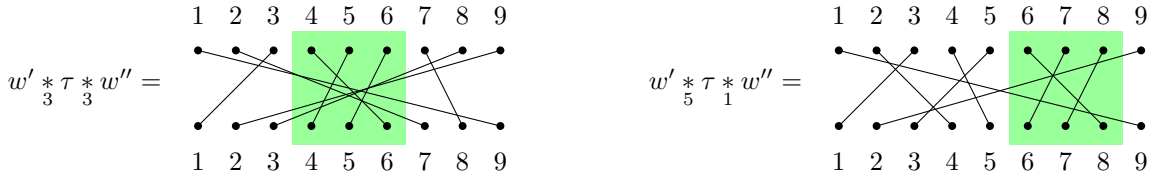
\begin{figure}[ht]
    \centering 
    \begin{tikzpicture}[scale=1, every node/.style={transform shape},
dot/.style = {circle, fill, minimum size=#1, inner sep=0pt, outer sep=5pt},
emptydot/.style = {circle, fill=white, minimum size=#1, inner sep=0pt, outer sep=5pt}]
%preliminaries
\def\radd{3pt}% radius of the dots
\def\height{1}% the distance between two rows of dots
\def\width{0.5}% the distance between two consecutive dots in the same row

%initializing the top row vertices
\coordinate (1up) at (0,0);
\path (1up) ++(\width,0) coordinate (2up);
\path (2up) ++(\width,0) coordinate (3up);
\path (3up) ++(\width,0) coordinate (4up);
\path (4up) ++(\width,0) coordinate (5up);
\path (5up) ++(\width,0) coordinate (6up);
\path (6up) ++(\width,0) coordinate (7up);
\path (7up) ++(\width,0) coordinate (8up);
\path (8up) ++(\width,0) coordinate (9up);

%initializing the bottom row vertices
\path (1up) ++(0,-\height) coordinate (1down);
\path (1down) ++(\width,0) coordinate (2down);
\path (2down) ++(\width,0) coordinate (3down);
\path (3down) ++(\width,0) coordinate (4down);
\path (4down) ++(\width,0) coordinate (5down);
\path (5down) ++(\width,0) coordinate (6down);
\path (6down) ++(\width,0) coordinate (7down);
\path (7down) ++(\width,0) coordinate (8down);
\path (8down) ++(\width,0) coordinate (9down);

%initializaing the colored area
\path (4up) ++(-0.5*\width,0.5*\width) coordinate (4upcolor);
\path (6up)++(0.5*\width,0.5*\width) coordinate (6upcolor);
\path (4down) ++(-0.5*\width,-0.5*\width) coordinate (4downcolor);
\path (6down) ++(0.5*\width,-0.5*\width) coordinate (6downcolor);

%filling the area with color
\fill[green!40] (4upcolor) -- (6upcolor) -- (6downcolor) -- (4downcolor) -- cycle;

% drawing the top row dots
\node[dot=\radd,label={ $1$}] at (1up) {};
\node[dot=\radd,label={ $2$}] at (2up) {};
\node[dot=\radd,label={ $3$}] at (3up) {};
\node[dot=\radd,label={ $4$}] at (4up) {};
\node[dot=\radd,label={ $5$}] at (5up) {};
\node[dot=\radd,label={ $6$}] at (6up) {};
\node[dot=\radd,label={ $7$}] at (7up) {};
\node[dot=\radd,label={ $8$}] at (8up) {};
\node[dot=\radd,label={ $9$}] at (9up) {};

%drawing the bottom row dots
\node[dot=\radd,label=below:{ $1$}] at (1down) {};
\node[dot=\radd,label=below:{ $2$}] at (2down) {};
\node[dot=\radd,label=below:{ $3$}] at (3down) {};
\node[dot=\radd,label=below:{ $4$}] at (4down) {};
\node[dot=\radd,label=below:{ $5$}] at (5down) {};
\node[dot=\radd,label=below:{ $6$}] at (6down) {};
\node[dot=\radd,label=below:{ $7$}] at (7down) {};
\node[dot=\radd,label=below:{ $8$}] at (8down) {};
\node[dot=\radd,label=below:{ $9$}] at (9down) {};

%drawing the black lines
\draw (1up) -- (9down);
\draw (2up) -- (7down);
\draw (3up) -- (1down);
\draw (4up) -- (6down);
\draw (5up) -- (4down);
\draw (6up) -- (5down);
\draw (7up) -- (8down);
\draw (8up) -- (3down);
\draw (9up) -- (2down);

\node[] at (-1.5,-0.5) {$w'\underset{3}{*}\tau\underset{3}{*}w''=$};

\end{tikzpicture}\hspace{40pt}
    \begin{tikzpicture}[scale=1, every node/.style={transform shape},
dot/.style = {circle, fill, minimum size=#1, inner sep=0pt, outer sep=5pt},
emptydot/.style = {circle, fill=white, minimum size=#1, inner sep=0pt, outer sep=5pt}]
%preliminaries
\def\radd{3pt}% radius of the dots
\def\height{1}% the distance between two rows of dots
\def\width{0.5}% the distance between two consecutive dots in the same row

%initializing the top row vertices
\coordinate (1up) at (0,0);
\path (1up) ++(\width,0) coordinate (2up);
\path (2up) ++(\width,0) coordinate (3up);
\path (3up) ++(\width,0) coordinate (4up);
\path (4up) ++(\width,0) coordinate (5up);
\path (5up) ++(\width,0) coordinate (6up);
\path (6up) ++(\width,0) coordinate (7up);
\path (7up) ++(\width,0) coordinate (8up);
\path (8up) ++(\width,0) coordinate (9up);

%initializing the bottom row vertices
\path (1up) ++(0,-\height) coordinate (1down);
\path (1down) ++(\width,0) coordinate (2down);
\path (2down) ++(\width,0) coordinate (3down);
\path (3down) ++(\width,0) coordinate (4down);
\path (4down) ++(\width,0) coordinate (5down);
\path (5down) ++(\width,0) coordinate (6down);
\path (6down) ++(\width,0) coordinate (7down);
\path (7down) ++(\width,0) coordinate (8down);
\path (8down) ++(\width,0) coordinate (9down);

%initializaing the colored area
\path (6up) ++(-0.5*\width,0.5*\width) coordinate (6upcolor);
\path (8up)++(0.5*\width,0.5*\width) coordinate (8upcolor);
\path (6down) ++(-0.5*\width,-0.5*\width) coordinate (6downcolor);
\path (8down) ++(0.5*\width,-0.5*\width) coordinate (8downcolor);

%filling the area with color
\fill[green!40] (6upcolor) -- (8upcolor) -- (8downcolor) -- (6downcolor) -- cycle;

%drawing the top row dots
\node[dot=\radd,label={ $1$}] at (1up) {};
\node[dot=\radd,label={ $2$}] at (2up) {};
\node[dot=\radd,label={ $3$}] at (3up) {};
\node[dot=\radd,label={ $4$}] at (4up) {};
\node[dot=\radd,label={ $5$}] at (5up) {};
\node[dot=\radd,label={ $6$}] at (6up) {};
\node[dot=\radd,label={ $7$}] at (7up) {};
\node[dot=\radd,label={ $8$}] at (8up) {};
\node[dot=\radd,label={ $9$}] at (9up) {};

%drawing the bottom row dots
\node[dot=\radd,label=below:{ $1$}] at (1down) {};
\node[dot=\radd,label=below:{ $2$}] at (2down) {};
\node[dot=\radd,label=below:{ $3$}] at (3down) {};
\node[dot=\radd,label=below:{ $4$}] at (4down) {};
\node[dot=\radd,label=below:{ $5$}] at (5down) {};
\node[dot=\radd,label=below:{ $6$}] at (6down) {};
\node[dot=\radd,label=below:{ $7$}] at (7down) {};
\node[dot=\radd,label=below:{ $8$}] at (8down) {};
\node[dot=\radd,label=below:{ $9$}] at (9down) {};

%drawing the black lines

\draw (1up) -- (9down);
\draw (2up) -- (4down);
\draw (3up) -- (1down);
\draw (4up) -- (5down);
\draw (5up) -- (3down);
\draw (6up) -- (8down);
\draw (7up) -- (6down);
\draw (8up) -- (7down);
\draw (9up) -- (2down);

\node[] at (-1.5,-0.5) {$w'\underset{5}{*}\tau\underset{1}{*}w''=$};

\end{tikzpicture}
    
    The part related to $\tau$ is highlighted in green.
\end{figure}
\begin{lemma}\label{lemma10}
For any two permutations $w\in S(n)$, $\tau\in S(m)$, and integers $l,r\in[0,n]$ such that $l+r=n$ one has
        \begin{align}
        w'\underset{l}{*}\tau\underset{r}{*}w''=\operatorname{Ad}\big((12)^{m,l,r}\big)\big[\tau\times w\big],
    \end{align}
    where $\operatorname{Ad}(\sigma)$ stand for the conjugate action $\operatorname{Ad}(\sigma)(u):=\sigma u\sigma^{-1}$; and
    \begin{align}
        \Big(w'\underset{l}{*}\tau\underset{r}{*}w''\Big)^{-1}=(w^{-1})'\underset{l}{*}\tau^{-1}\underset{r}{*}(w^{-1})''.
    \end{align}
\end{lemma}

\subsubsection{Canonical projections}

We will also need an operation on permutations called the canonical projection in \cite{kerov2004harmonic}. Let $u\in S(n+1)$ be any permutation. Consider $u$ as a collection of $n+1$ nodes in the first row, joined by $n+1$ edges to $n+1$ nodes in the second row. In the case $u=\operatorname{id}_1\times w\in S(n+1)$ we set by definition $u_*=w\in S(n)$. If $u(1)>1$, we add one extra edge connecting the first vertices in the upper and lower rows. Then $u_*\in S(n)$ is the permutation defined by the paths in this new graph. Another useful description of $u_*$ is as follows. First, we remove $1$ from the cycle that contains it, obtaining a bijection between $2,\ldots, n+1$ and $1,\ldots,u(1)-1,u(1)+1,\ldots,n+1$. Then we relabel the nodes and obtain a permutation from $S(n)$ denoted by $u_*$. More precisely, we decompose $u$ into a product of commuting cycles, and replace the cycle $(\ldots \mapsto i\mapsto 1\mapsto j\mapsto\ldots)$ containing $1$ with the cycle $(\ldots \mapsto i\mapsto j\mapsto\ldots)$. Then we relabel the nodes and denote the result by $u_*$. Note that our definition of $u_*$ slightly differs from the one given in \cite{kerov2004harmonic}, because we remove the number $1$ from a cycle but not $n+1$.

\begin{proposition}[Proposition 1.1.1 in \cite{kerov2004harmonic}]\label{prop9}
    The map $S(n+1)\longrightarrow S(n)$, $u\mapsto u_*$ is equivariant with respect to the two-sided action of the group $S(n)$ which acts on $S(n+1)$ via $S(n)\simeq\operatorname{id}_1\times S(n)\subset S(n+1)$. If $n\geq 4$, this is the only map $S(n+1)\longrightarrow S(n)$ with this property.
\end{proposition}

\subsubsection{Actions on matrices}
Below we work with linear functionals on matrices of a fixed size $\Mat$ and their tensor powers $\Mat^{\otimes n}$. We will not only permute the tensor factors but also permute ''vectors'' and ''covectors'' independently. Namely, we identify matrices $\operatorname{Mat}_N(\Bbbk)$ with $\big(\Bbbk^N\big)^*\otimes \Bbbk^N$, and the functionals on matrices with $\Bbbk^N\otimes\big(\Bbbk^N\big)^*$. We will refer to the first factor of the latter tensor product as ''vectors'' and to the second one as ''covectors''. Then we can identify $\Mat^{\otimes n}$ with $\big(\Bbbk^N\big)^{\otimes n}\otimes\Big(\big(\Bbbk^N\big)^*\Big)^{\otimes n}$ and it's now clear what it means to permute vectors or covectors. Note that these actions of $S(n)$ commute with the natural action of $\GL$ on $\Mat^{\otimes n}$.

Explicitly, for any $\sigma\in S(n)$ we denote by $\operatorname{vect}(\sigma)$ and $\operatorname{covect}(\sigma)$ the left actions of the symmetric group $S(n)$ on $\Mat^{\otimes n}$ by permutations of vectors and covectors, i.e.
\begin{align}
    \operatorname{vect}(\sigma)\left(E^*_{\mathbf{i}\mathbf{j}}\right)&:=E^*_{\sigma(\mathbf{i})\mathbf{j}}=E^*_{i_{\sigma^{-1}(1)}j_1}\otimes \ldots\otimes E^*_{i_{\sigma^{-1}(n)}j_n},\\
    \operatorname{covect}(\sigma)\left(E^*_{\mathbf{i}\mathbf{j}}\right)&:=E^*_{\mathbf{i}\sigma(\mathbf{j})}=E^*_{i_1j_{\sigma^{-1}(1)}}\otimes \ldots\otimes E^*_{i_nj_{\sigma^{-1}(n)}},
\end{align}
where $E^*_{\mathbf{i}\mathbf{j}}:=E^*_{i_1j_1}\otimes \ldots\otimes E^*_{i_nj_n}$ for any $n$-tuples of integers $\mathbf{i}=(i_1,\ldots,i_n)\in[1,N]^n$, and $\mathbf{j}=(j_1,\ldots,j_n)\in[1,N]^n$.

In this notation the usual action of $S(n)$ on $\Mat^{\otimes n}$ by permutations of tensor components becomes $\operatorname{vect}(\sigma)\operatorname{covect}(\sigma)$, but we will denote it simply by $\sigma$.

\subsection{Representation spaces \texorpdfstring{$\rep$}{RepNA} and double Poisson brackets}
 Let $A$ be a finitely generated algebra over $\Bbbk$. Below we will work with tensor powers of the algebra $A$ and use the following notation.

\begin{itemize}
\item An element $a=\sum\limits_i a_i'\otimes a_i''\in A^{\otimes 2}$ is usually written in the shorthand notation $a=a'\otimes a''$. Sometimes we will call these $a'$ and $a''$ \textit{the Sweedler's components} of $a$. Similarly for the tensor cube or any other tensor power of $A$.

\item For any $n$ we consider the natural action of the symmetric group on $n$ letters $S(n)$ on $A^{\otimes n}$ given by 
\begin{equation}
    \sigma\cdot (a_1\otimes\ldots\otimes a_n):=a_{\sigma^{-1}(1)}\otimes\ldots\otimes a_{\sigma^{-1}(n)},\ \ \ \sigma\in S(n).
\end{equation}

\item For $n=2$ we denote the action of the unique nontrivial element $(12)\in S(2)$ by $\circ$
\begin{equation}
    (a_1\otimes a_2)^{\circ}=a_2\otimes a_1.
\end{equation}
\end{itemize}

We will use two commuting $A$-bimodule structures on $A^{\otimes 2}$ --- \textit{outer and inner} which are respecively given by
\begin{gather}
    a(x'\otimes x'')b:=(ax')\otimes (x''b),\\
    a*(x'\otimes x'')*b:=(x'b)\otimes (ax''),
\end{gather}
where $a,b,x',x''\in A$.

Note that the permutation $(12)$ interchanges these two bimodule structures
\begin{equation}
    (axb)^{\circ}=a*x^{\circ}*b,\ \ \ \ \ \ \ a,b\in A,\ \ \ \ \ x\in A^{\otimes 2}.
\end{equation}

\begin{definition}[Van den Bergh,\cite{van2008double}]
A \emph{double Poisson bracket} on $A$ is a linear map 
$$
\{\!\!\{-,-\}\!\!\}: A\otimes A \longrightarrow  A\otimes A
$$
satisfying the following three conditions. 

1. \emph{Skew symmetry}: for $a,b\in A$,
\begin{equation}
\{\!\!\{a,b\}\!\!\}=-\{\!\!\{b,a\}\!\!\}^\circ.
\end{equation}

2. \emph{Leibniz rule}: for $a,b,c\in A$, 
\begin{equation}
\{\!\!\{a, bc\}\!\!\}=\{\!\!\{a,b\}\!\!\} c+ b\{\!\!\{a,c\}\!\!\}.
\end{equation}

3. \emph{Double Jacobi identity}: for $a,b,c\in A$, 
\begin{equation}
\Bigl\{\!\!\!\Bigl\{ a,\{\!\!\{b,c\}\!\!\}\Bigr\}\!\!\!\Bigr\}_{\mathbf L} +(123)\cdot \Bigl\{\!\!\!\Bigl\{ b,\{\!\!\{ c,a\}\!\!\}\Bigr\}\!\!\!\Bigr\}_{\mathbf L} +(123)^2\cdot \Bigl\{\!\!\!\Bigl\{ c,\{\!\!\{ a,b\}\!\!\}\Bigr\}\!\!\!\Bigr\}_{\mathbf L}=0,
\end{equation}
where $(123)$ is the cyclic permutation $1\to 2\to 3\to 1$, which pushes everything to the right
$$
(123)\cdot (x'\otimes x''\otimes x''')=x'''\otimes x'\otimes x'',
$$
and for $x=x'\otimes x''\in A^{\otimes 2}$,
$$
\{\!\!\{ a,x\}\!\!\}_{\mathbf L}:= \{\!\!\{ a,x'\}\!\!\}\otimes x''\in A^{\otimes 3}.
$$
\end{definition}

To any double Poisson bracket we can associate a map $\{-,-\}:A\otimes A\rightarrow A$ which is the composition of $\{\!\!\{-,-\}\!\!\}$ with the multiplication $m$ in $A$, i.e. $\{-,-\}=m\{\!\!\{-,-\}\!\!\}$. 

One can easily check that $\{ab,c\}=\{ba,c\}$ for any $a,b,c\in A$. Then the bracket $\{-,-\}$ degenerates to a map $A_{\natural}\otimes A\rightarrow A$, which we denote by $\{-,-\}$ too, where $A_{\natural}$ stands for the quotient vector space $\faktor{A}{[A,A]}$. Let us denote the natural map $A\twoheadrightarrow A_{\natural}$ by $a\mapsto\bar{a}$. 

\begin{proposition}[\cite{van2008double}, Corollaries 2.4.4 and 2.4.6]\label{prop4}
\phantom{a}
\begin{enumerate}[label=\arabic*)]
    \item The map $\{-,-\}:A_{\natural}\otimes A\rightarrow A$ gives rise to a well-defined skew-symmetric bilinear map $A_{\natural}\times A_{\natural}\rightarrow A_{\natural}$, which satisfies Jacobi identity, turning $A_{\natural}$ into a Lie algebra. We will denote this Lie bracket by the same symbol $\{-,-\}$.

    \item The map $\{-,-\}:A_{\natural}\otimes A\rightarrow A$ defines an action of the Lie algebra $A_{\natural}$ on $A$, turning $A$ into a representation of the Lie algebra $A_{\natural}$. The action of $\bar{a}\in A_{\natural}$ on $A$ is given by $\{\bar{a},-\}$.
\end{enumerate}
\end{proposition}

Let $N$ be a positive integer. 

\begin{definition}
    Representation space $\rep$ is by definition the affine $\Bbbk$-scheme representing the functor 
    \begin{align}
        &CAlg\longrightarrow Set,\\
        &R\mapsto \operatorname{Hom}_{\operatorname{Alg}}\left(A,\operatorname{Mat}_N(R)\right),
    \end{align}
    where $CAlg$ is the category of commutative unital algebras.
\end{definition}

The coordinate ring $\orep$ of $\rep$ is the commutative algebra generated by the symbols $a_{ij}$ for $a\in A$ and $i,j\in\{1,\dots,N\}$, which are $\Bbbk$-linear in $a$  and  are subject to the relations
\begin{equation}\label{f17}
(ab)_{ij}=\sum_{k=1}^N a_{ik}b_{kj}, \qquad a,b\in A, \quad i,j=1,\dots,N
\end{equation}
together with
\begin{equation}\label{f36}
    (1)_{ij}=\delta_{ij} \qquad i,j=1,\dots,N,
\end{equation}
where $1$ is the identity in $A$. 

\begin{proposition}[\cite{van2008double}, Proposition 1.2]
Let $A$ be an associative $\Bbbk$-algebra and $\{\!\!\{-,-\}\!\!\}$ be a double Poisson bracket on $A$. For each $N=1,2,\dots$, the formula
\begin{equation}
\{a_{ij},b_{kl}\}:=\{\!\!\{ a,b\}\!\!\}'_{kj}\{\!\!\{ a,b\}\!\!\}''_{il}, \qquad a,b\in A, \quad i,j,k,l\in\{1,\dots,N\},
\end{equation}
gives rise to a Poisson bracket on the commutative algebra $\orep$. 
\end{proposition}

\section{Diagonal \texorpdfstring{$\S$}{S}-bimodules}\label{section_s_bimodules}
%\addtocontents{toc}{\protect\setcounter{tocdepth}{1}}

Here we recall basic definitions related to diagonal $\S$-bimodules. This name first appeared in \cite{fernández2025symplecticwheelgebrasnoncommutativegeometry} and we refer the reader to this paper for a detailed treatment of the symmetric monoidal category of diagonal $\S$-bimodules.

Roughly speaking, one can recast definitions of $\S$-modules and of standard structures related to them replacing symmetric groups $S(n)$ with their squares $S(n)\times S(n)$, and obtain the definition of diagonal $\S$-bimodules and related structures. The word diagonal refers to the fact that we work with $S(n)\times S(n)$ but not $S(n)\times S(m)$. In our presentation we will mainly follow \cite{stover1993equivalence} and Section 5.1 of \cite{loday2012algebraic} adapting constructions to our bimodule context in an obvious way.

\begin{definition}
    A \textit{diagonal $\S$-bimodule} $M$ is a graded vector space $M=\bigoplus\limits_{n\geq 0}M^{(n)}$, where each $M^{(n)}$ is an $S(n)$-bimodule. The degree of a homogeneous element $m\in M^{(n)}$ is denoted by $|m|=n$. A \textit{morphism of diagonal $\S$-bimodules} $f:M\longrightarrow N$ is a collection of $S(n)$-bimodule homomorphisms $f_n:M^{(n)}\longrightarrow N^{(n)}$.
\end{definition}

Note that to any diagonal $\S$-bimodule $M$ one can canonically associate a (left) $\S$-module which is the same $M$ as graded vector spaces, but the $\S$-module structure comes from the adjoint action, i.e. $\sigma\in S(n)$ acts on $m\in M^{(n)}$ as $\operatorname{Ad}(\sigma)(m):=\sigma\cdot m\cdot\sigma^{-1}$. We will refer to this $\S$-module structure as \textit{the adjoint $\S$-module structure on $M$}. Roughly speaking, all ''twisted'' structures on diagonal $\S$-bimodules we are going to consider, are twisted with respect to this adjoint $\S$-module structure.

To any positive integer $k$ and any diagonal $\S$-bimodule $M$ one can associate the shifted (to the left) diagonal $\S$-bimodule $M[k]$ defined by
\begin{align}
    \big(M[k]\big)^{(n)}:=\operatorname{Res}^{S(n+k)}_{S(n)} M^{(n+k)},
\end{align}
where $\operatorname{Res}$ stands for the restriction of the $S(n+k)$-bimodule structure to $S(n)$ realized as the subgroup of $S(n+k)$ fixing the first $k$ nodes. In other words, $\big(M[k]\big)^{(n)}=M^{(n+k)}$ and the $S(n)$-bimodule structure is given by the usual bimodule action of $S(n+k)$ precomposed with the map $\operatorname{id}_k\times S(n)\hookrightarrow S(n+k)$.

For any vector space $V$ we denote by the same symbol $V$ the $\S$-bimodule concentrated in degree zero, and by $V[-1]$ the $\S$-bimodule concentrated in degree one.

\subsection{Tensor product \texorpdfstring{$\stimes$}{S}}
For any finite group $G$ and its subgroup $H\subset G$, and any $H$-bimodule $M$ we denote by $\operatorname{bi-Ind}_H^G(M)$ the induced $G$-bimodule defined by 
$$
\operatorname{bi-Ind}_H^G(M):=\Bbbk[G]\otimes_{H}M\otimes_{H}\Bbbk[G].
$$

For diagonal $\S$-bimodules $M$ and $N$ we define their tensor product $\stimes$ as follows 
\begin{align}
    (M\stimes N)^{(n)}:=\bigoplus\limits_{i+j=n}\operatorname{bi-Ind}_{S(i)\times S(j)}^{S(n)} \big(M^{(i)}\otimes N^{(j)}\big).
\end{align}

We will usually denote the element $\sigma\otimes \big(m\otimes n\big)\otimes \tau$ of $(M\stimes N)^{(n)}$ by $\sigma\circ \big(m\otimes n\big)\circ \tau$, and element $\operatorname{id}\otimes \big(m\otimes n\big)\otimes \operatorname{id}$ simply by $m\otimes n$.

This tensor product is commutative. The swap map $(12):M\stimes N\longrightarrow N\stimes M$ is given on homogeneous elements by 
\begin{align}
    \sigma\circ \big(m\otimes n\big)\circ \tau\mapsto \sigma(12)^{|n|,|m|}\circ \big(n\otimes m\big) \circ (12)^{|m|,|n|}\tau
\end{align}

\begin{proposition}
The category of diagonal $\S$-bimodules is symmetric monoidal with respect to the tensor product $\stimes$, swap map $(12)$ introduced above, and the unit $\big(\Bbbk,0,0,\ldots,\big)$, which is the diagonal $\S$-bimodule concentrated in degree zero and whose only non-trivial component is $\Bbbk$.  
\end{proposition}

The swap map $(12)$ extends to the natural left action of the symmetric group $S(k)$ on the $k$-fold tensor product. 

Say, we had diagonal $\S$-bimodules $M_1,\ldots,M_k$. One easily sees that their tensor product $M_1\stimes\ldots\stimes M_k$ is given by 
$$
\big(M_1\stimes\ldots\stimes M_k\big)^{(n)}=\bigoplus\limits_{i_1+\ldots+i_k=n}\operatorname{bi-Ind}_{S(i_1)\times\ldots\times S(i_k)}^{S(n)} \Big(M_1^{(i_1)}\otimes\ldots\otimes  M_k^{(i_k)}\Big).
$$

Then the permutation $\sigma\in S(k)$ acts by
\begin{align}\label{f70}
    \sigma:M_1\stimes\ldots\stimes M_k&\longrightarrow M_{\sigma^{-1}(1)}\stimes \ldots \stimes M_{\sigma^{-1}(k)},\\
    \tau_1\circ (m_1\otimes\ldots\otimes m_k)\circ \tau_2&\mapsto \tau_1(\sigma^{|m_1|,\ldots,|m_k|})^{-1}\circ(m_{\sigma^{-1}(1)}\otimes\ldots\otimes m_{\sigma^{-1}(k)})\circ\sigma^{|m_1|,\ldots,|m_k|}\tau_2,
\end{align}
where $m_1,\ldots,m_k$ are homogeneous and $\tau_1,\tau_2\in S(|m_1|+\ldots+|m_k|)$.

\begin{proposition}\label{prop10}
\phantom{a}
\begin{enumerate}[label=\theenumi), leftmargin=5ex]
    \item For any permutation $\sigma\in S(k)$ the map \eqref{f70} is a well-defined homomorphism of diagonal $\S$-bimodules.

    \item The map \eqref{f70} corresponding to a product of permutations equals the composite of similar maps corresponding to the factors, e.g. the map \eqref{f53} corresponding to $\sigma\tau$ equals the composite $\sigma\big(\tau(-)\big)$. 
\end{enumerate}
In the case when $M_1=\ldots=M_k=M$ this defines a left action of $S(k)$ on $M^{\stimes k}$ by diagonal $\S$-bimodule automorphisms.
\end{proposition}

\subsection{Di-twisted algebras}
\begin{definition}
    By a \textit{di-twisted algebra}\footnote{di is from ''diagonally''.} we mean an algebra in the category of diagonal $\S$-bimodules, i.e. an diagonal $\S$-bimodule $A$ together with diagonal $\S$-bimodule morphisms $\mu:A\stimes A\longrightarrow A$ and $\eta:\Bbbk\longrightarrow A$ such that the following well-known diagrams commute
    \begin{equation}
        \begin{tikzcd}
            A\stimes A\stimes A \arrow{r}[sloped, above]{\mu\otimes 1} \arrow{dd}[xshift=-4.5ex]{1\otimes \mu} & A\stimes A \arrow{dd}{\mu}& \\
            & \\
            A\stimes A \arrow{r}[sloped, below]{\mu}& A 
        \end{tikzcd}\hspace{30pt}
        \begin{tikzcd}
            & \Bbbk\stimes A\arrow{r}[above]{\eta\otimes 1} & A\stimes A \arrow{rd}[sloped, above]{\mu}&\\
            A \arrow[equal]{rrr} \arrow{ru}[sloped, above]{\sim}\arrow{rd}[sloped, above]{\sim} & & & A\\
            & A\stimes \Bbbk\arrow{r}[below]{1\otimes \eta} & A\stimes A\arrow{ru}[sloped, below]{\mu}& 
        \end{tikzcd}
    \end{equation}
\end{definition}

Note that any di-twisted algebra is automatically a graded algebra, and the graded components of the multiplication $\mu_{n,m}:A^{(n)}\otimes A^{(m)}\longrightarrow A^{(n+m)}$ are left and right $S(n)\times S(m)$-equivariant. Thus, any di-twisted algebra $A$ is also a twisted algebra (an algebra in the category of $\S$-modules) with respect to the adjoint $\S$-module structure.

We say that a di-twisted algebra $A$ is \textit{commutative}, if $\mu(12)=\mu$. This explicitly reads as 
$$
    \mu(a,b)=\operatorname{Ad}(12)^{|b|,|a|}\mu\big(b,a\big),
$$
for homogeneous $a,b\in A$, which means that the multiplication $\mu$ is twisted commutative with respect to the adjoint $\S$-module structure, see 3.5 in \cite{stover1993equivalence}, but note that the author works with right $\S$-modules.

\begin{definition}\label{def1}
    Let $A$ and $B$ be di-twisted algebras. Their tensor product $A\stimes B$ is naturally a di-twisted algebra with the multiplication
    \begin{equation}
        \begin{tikzcd}
            (A\stimes B)\stimes (A\stimes B)\arrow{r}{(23)} & (A\stimes A)\stimes (B\stimes B)\arrow{r}{\mu\otimes\mu}& A\stimes B
        \end{tikzcd}.
    \end{equation}
    Explicitly,
    \begin{align}
        &\big[\sigma_1\circ(a_1\otimes b_1)\circ\tau_1\big]\big[\sigma_2\circ(a_2\otimes b_2)\circ\tau_2\big]\\
        &\hspace{80pt}:=\Big((\sigma_1\times\sigma_2)((23)^{|a_1|,|b_1|,|a_2|,|b_2|})^{-1}\Big)\circ\big(a_1a_2\otimes b_1b_2\big)\circ\Big((23)^{|a_1|,|b_1|,|a_2|,|b_2|}(\tau_1\times\tau_2)\Big)
    \end{align}
    for homogeneous $a_1,a_2\in A$, $b_1,b_2\in B$, and $\sigma_i,\tau_i\in S(|a_i|+|b_i|)$, $i=1,2$.
\end{definition}

\begin{proposition}\label{prop12}
    \begin{enumerate}[label=\theenumi),leftmargin=3ex]
        \item If $A$ and $B$ are commutative di-twisted algebras, then so is their tensor product.

        \item If $A$ is a di-twisted algebra, then the natural action of the symmetric group $S(k)$ on $A^{\stimes k}$ is by di-twisted algebra automorphisms.
    \end{enumerate}
\end{proposition}

\subsubsection{Tensor, symmetric, and exterior algebras}
    For any diagonal $\S$-bimodule $M$ one can consider its tensor algebra $$\operatorname{T}_{\S}(M)=\bigoplus\limits_{k\geq 0} M^{\stimes k}.$$
    For instance, $\operatorname{T}_{\S}(M)^{(0)}=\operatorname{T}\big(M^{(0)}\big)$ and $\operatorname{T}_{\S}(M)^{(1)}=\operatorname{T}\big(M^{(0)}\big)\otimes M^{(1)}\otimes \operatorname{T}\big(M^{(0)}\big)$, where $T$ in the right-hand sides stands for the usual tensor algebra of a vector space. 
    
    Multiplication for $\operatorname{T}_{\S}(M)$ comes from the associativity isomorphisms $M^{\stimes k}\stimes M^{\stimes l}\simeq M^{\stimes k+l}$. Explicitly,
    $$
    \big[\sigma_1\circ(m_1\otimes\ldots\otimes m_k)\circ \tau_1\big]\big[\sigma_2\circ(m_{k+1}\otimes \ldots\otimes m_{k+l})\circ\tau_2\big]=(\sigma_1\times\sigma_2)\circ(m_1\otimes\ldots\otimes m_{k+l})\circ(\tau_1\times\tau_2)
    $$
    for homogeneous $m_1\ldots,m_{k+l}\in M$ and $\sigma_1,\tau_1\in S(|m_1|+\ldots+|m_k|)$, $\sigma_2,\tau_2\in S(|m_{k+1}|+\ldots+|m_{k+l}|)$.

One can also consider the symmetric algebra of $M$
\begin{align}\label{f60}
    \SSym\big(M\big):=\bigoplus\limits_{k\geq 0}\big(M^{\stimes k}\big)_{S(k)},
\end{align}
where $\big(-\big)_{S(k)}$ stands for the coinvariants of the $S(k)$-action on $M^{\stimes k}$ discussed in Proposition \ref{prop10}.

Warning: the sum \eqref{f60} is not the decomposition of a diagonal $\S$-bimodule into its graded components. The diagonal $\S$-bimodule structure on $\SSym(M)$ comes from the same structure on $M^{\stimes k}$'s. Due to the fact that the action of $S(k)$ on $M^{\stimes k}$ commutes with the diagonal $\S$-bimodule structure, the latter descends to the coinvariants. 

Next, one readily sees that the multiplication on the tensor algebra $\operatorname{T}_{\S}(M)$ descends to $\SSym(M)$ turning it into a di-twisted commutative algebra.

\begin{proposition}\label{prop11}
    Suppose that a diagonal $\S$-bimodule $M$ is concentrated in degrees zero and one, i.e. $M^{(n)}=0$ for $n\geq 2$. Then one has an isomorphism of di-twisted algebras
    \begin{align}\label{f61}
        \SSym(M)\simeq \bigoplus\limits_{n\geq 0}\big(M^{(1)}\big)^{\otimes n}\otimes \operatorname{S}(M^{(0)})\otimes \Bbbk[S(n)],
    \end{align}
    where the sum in the right-hand side is the decomposition of a diagonal $\S$-bimodule in its graded components, and on each of them 
    the left and right $S(n)$-actions are given by %the right action is the right multiplication of the group algebra, and the left action is the simultaneous permutation of $M^{(1)}$-components and the left multiplication of the group algebra.%, i.e. 
    \begin{align}
    \sigma\cdot(m\otimes f\otimes \tau)&:=\big(m_{\sigma^{-1}(1)}\otimes\ldots\otimes m_{\sigma^{-1}(n)}\big)\otimes f\otimes (\sigma\tau),\\
    (m\otimes f\otimes \tau)\cdot \sigma&:=m\otimes f\otimes (\tau\sigma),
    \end{align}
    for any $m=m_1\otimes\ldots\otimes m_n\in\big(M^{(1)}\big)^{\otimes n}$, $f\in \operatorname{S}(M^{(0)})$, $\tau\in S(n)$, and $\sigma\in S(n)$.
    
    The multiplication on the right-hand side of \eqref{f61} is as follows: 
    \begin{align}
        (m_1\otimes f_1\otimes \tau_1)(m_2\otimes f_2\otimes \tau_2):=(m_1\otimes m_2)\otimes (f_1f_2)\otimes (\tau_1\times\tau_2),
    \end{align}
    for any $m_1\in\big(M^{(1)}\big)^{\otimes n_1}$, $f_1\in \operatorname{S}(M^{(0)})$, $\tau_1\in S(n_1)$, and $m_2\in\big(M^{(1)}\big)^{\otimes n_2}$, $f_2\in \operatorname{S}(M^{(0)})$, $\tau_2\in S(n_2)$.
\end{proposition}
\begin{proof}
    By the very definition of the tensor product $\stimes$ we have
    \begin{align}\label{f62}
        &M^{\stimes k}=\bigoplus\limits_{i_1,\ldots,i_k=0,1}\Bbbk[S(i_1+\ldots+i_k)]\otimes\big(M^{(i_1)}\otimes\ldots\otimes M^{(i_k)}\big)\otimes \Bbbk[S(i_1+\ldots+i_k)].
    \end{align}
    Note that all tensor product are taken over $\Bbbk$ because $i_1,\ldots,i_k=0,1$, hence $S(i_1)\times \ldots\times S(i_k)$ is trivial. 
    
    Then for coinvariants one has
    \begin{align}
        \big(M^{\stimes k}\big)_{S(k)}=\bigoplus\limits_{n=0}^k \big(M^{(1)}\big)^{\otimes n}\otimes \operatorname{S}^{k-n}(M^{(0)})\otimes \Bbbk[S(n)].
    \end{align}
    To see this, one should act by $S(k)$ on \eqref{f62} and perform the following manipulations:
    \begin{enumerate}[label=\theenumi)]
        \item Move all $M^{(1)}$-components in \eqref{f62} to the left and all $M^{(0)}$-components to the right.

        \item Permute all the $M^{(0)}$-components -- this affects neither $M^{(1)}$-components nor the two permutations, hence we get an element of the symmetric algebra $\operatorname{S}(M^{(0)})$.

        \item Kill the left permutation by permuting the $M^{(1)}$-components.
    \end{enumerate}
    One can readily see that this gives an isomorphism of di-twisted algebras \eqref{f61} when the right-hand side is equipped with the aforementioned diagonal $\S$-bimodule structure and multiplication.
\end{proof}

One can similarly consider the exterior algebra of a diagonal $\S$-bimodule $M$
\begin{align}
    \bigwedge\nolimits_{\S}M=\bigoplus\limits_{k\geq 0}\big(M^{\stimes k}\big)_{S(k)}^-,
\end{align}
where $(-)_{S(k)}^-$ stands for anti-coinvariants, i.e. $(N)_{S(k)}^-=\faktor{N}{\big(\sigma.n-(-1)^{\sigma}n\ |\ n\in N,\ g\in S(k)\big)}$. 

One readily checks that the multiplication in $\operatorname{T}_{\S}(M)$ descends to $\bigwedge\nolimits_{\S}M$.

%\addtocontents{toc}{\protect\setcounter{tocdepth}{2}}

\section{Double coordinate ring \texorpdfstring{$\O(A)$}{OA}}\label{section_noncommutative_functions}

Here, we give a definition of the commutative di-twisted algebra $\O(A)$ and discuss its properties and connections to representation spaces $\rep$. Algebra $\O(A)$ first appeared in \cite{ginzburg2010differentialoperatorsbvstructures} as the right object to study differential operators in noncommutative geometry. This algebra was denoted by $\mathcal{F}(A)$; here $\mathcal{F}$ refers to ''Fock space''. Actually, $\O(A)=\mathcal{F}(A)$ is not merely a di-twisted algebra, but also a wheelgebra. This notion was introduced in Section 3.1 in \cite{ginzburg2010differentialoperatorsbvstructures} and studied in Section 4 in \cite{fernández2025symplecticwheelgebrasnoncommutativegeometry}. Roughly speaking, it reflects the fact that $A$ is equipped with a multiplication. This wheelgebra structure will not be sufficient for our needs as we will need an additional structure reflecting the fact that $A$ is also unital. So, the wheel-structure on $\O(A)$ will not be our primary interest. Instead of that we will focus on links between $\O(A)$ and representation spaces $\rep$. From this point of view, it would be easier for us to treat $\O(A)$ simply as a di-twisted algebra with a certain additional structure, without emphasizing exactly what this additional structure is. The notation $\O(A)$ reflects our point of view that $\O(A)$ should be regarded as a proper noncommutative analog of coordinate rings, at least for the purposes of this paper. We will explain this shortly.

Let us speculate a bit on the representation theoretic nature of the algebra $\O(A)$ before recalling its easy and concise definition. Suppose $A$ is an associative algebra and we are looking for a noncommutative analog of the coordinate rings $\orep$, i.e. we are looking for an associative algebra whose elements produce elements of $\orep$. The simplest solution would be to take $A$ itself. Indeed, having a pair of integers $i,j$ ranging from $1$ to $N$ and an element $a\in A$, one can produce a function $a_{ij}$ on the representation space $\rep$, which takes an $N$-dimensional representation $\rho$ of $A$ and returns the value of the $(ij)$-th matrix element of $\rho(a)$, so according to the Kontsevich-Rosenberg principle one can treat $A$ as an algebra of ''noncommutative functions'' on $\rep$. We argue that this is a noncommutative analog only of a generating set of $\orep$. In many cases this doesn't lead to any difficulties, because for many structures it's enough to define them on a system of generators, and then they can be uniquely extended to the whole algebra, e.g. Poisson brackets. However, differential operators and star-products are not uniquely determined by their values on generators. To treat those, a noncommutative analog of a set of generators is no longer enough: one needs an analog of a spanning set of $\orep$. Then the obvious solution is to consider the tensor algebra $\operatorname{T}(A)$, which obviously satisfies the Kontsevich-Rosenberg principle, as we can apply a pair of $n$-tuples of indices $\mathbf{i}=(i_1,\ldots,i_n)$ and $\mathbf{j}=(j_1,\ldots,j_n)$ to any element of $A^{\otimes n}$ in order to obtain an element of $\orep$:
\begin{align}
    (a_1\otimes\ldots\otimes a_n)_{\mathbf{i}\mathbf{j}}:=(a_1)_{i_1j_1}\ldots(a_n)_{i_nj_n}\in\orep.
\end{align}

The tensor algebra $\operatorname{T}(A)$ does not contain any information about the algebra structure in $A$, but one can implement it manually by considering an additional structure on $\operatorname{T}(A)$: two graded operators of degrees $+1$ and $-1$, which are responsible for the unit and multiplication in $A$. The former is the tensor multiplication by the unit, say, on the left, and the latter is the multiplication in algebra $A$ applied to the first two tensor factors in $A^{\otimes n}$ for $n\geq 2$ and, say, zero for $n=0,1$. The tensor algebra $\operatorname{T}(A)$ even contains an artifact of commutativity in $\orep$. Namely, there is a natural left action of the symmetric group $S(n)$ on the $n$-th graded component $A^{\otimes n}$. This action is responsible for commutativity in $\orep$ due to the following evident identity
\begin{align}\label{f63}
    \hspace{40pt} \Bigl(\sigma\cdot(a_1\otimes\ldots\otimes a_n)\Bigr)_{\mathbf{i}\mathbf{j}}=\left(a_1\otimes\ldots\otimes a_n\right)_{\sigma^{-1}(\mathbf{i})\sigma^{-1}(\mathbf{j})},\hspace{40pt} \sigma\in S(n).
\end{align}

This approach is basically correct. The only problem is that, if one wants to study $\GL$-invariant structures on $\rep$, the tensor algebra $\operatorname{T}(A)$ is hopelessly small. Suppose we wish to study $\GL$-equivariant linear maps $\orep\longrightarrow\orep$ without any additional constraints. This problem boils down to the study of $\GL$-equivariant linear maps between various tensor powers of $\Mat$, which can be done with the help of invariant theory. It's not very difficult to show that the vector space of $\GL$-equivariant maps from $\Mat^{\otimes n}$ to $\Mat^{\otimes m}$ is spanned by maps of the following form 
\begin{align}
        X_1\otimes\ldots\otimes X_n\mapsto \sigma_2\, D\, \operatorname{vect}(\sigma_1) (X_1\otimes\ldots\otimes X_n\otimes \tr^{\otimes s}),
    \end{align}
    where 
    \begin{itemize}
        \item $X_1,\ldots,X_n\in\Mat$;
        
        \item $\sigma_2$ is a permutation which permutes the tensor factors, i.e. $\vect(\sigma_2)\covect(\sigma_2)$;

        \item $D$ is a composition of several copies of the counit $\varepsilon$ and coproduct $\Delta$ taken in a certain order\footnote{Recall that $\varepsilon:\Mat\longrightarrow\Bbbk, \ \varepsilon(E_{ij}^*)=\delta_{ij}$ and $\Delta:\Mat\longrightarrow\Mat^{\otimes 2},\ \Delta(E_{ij}^*)=\sum\limits_{k=1}^NE_{ik}^*\otimes E_{kj}^*$.};

        \item $\sigma_1\in S(n)$, $\tr\in\Mat$ is the usual trace, and $s$ is an integer.
    \end{itemize}

Tracing this back to $\orep$, one sees that $\sigma_2$ can be eliminated by the commutativity of $\orep$ and $D$ by relations \eqref{f17} and \eqref{f36}. Thus, only $\operatorname{covect}(\sigma_1)$ and the product of traces remain. This suggests that the $n$-th graded component of the (future) algebra $\O(A)$ should be $A^{\otimes n}\otimes \s\otimes \Bbbk[S(n)]$, where $A^{\otimes n}$ corresponds to $X_1\otimes\ldots\otimes X_n$, $\s$ corresponds to traces $\tr$, and $\Bbbk[S(n)]$ corresponds to $\sigma_1$. Proposition \ref{prop11} then yields the following definition.

\begin{definition}\label{def4}
    Let $A$ be an associative algebra over $\Bbbk$. We set
    \begin{align}
        \O(A):=\SSym(A_{\natural}\oplus A[-1])=\bigoplus\limits_{n\geq 0}A^{\otimes n}\otimes \s\otimes \Bbbk[S(n)].
    \end{align}
\end{definition}

This definition appeared in exactly this form in Section 5.2 in \cite{fernández2025symplecticwheelgebrasnoncommutativegeometry} and in a bit different form in Section 3.2 in \cite{ginzburg2010differentialoperatorsbvstructures}, where it was first introduced.

From the discussion above it follows that $\O(A)$ is the smallest algebra one has to consider if they want to study $\GL$-equivariant linear maps $\orep\longrightarrow\orep$. All such maps arise then as images under the representation functor of linear maps $\O(A)\longrightarrow\O(A)$ satisfying certain conditions. The same is true for $\GL$-equivariant linear maps $\orep^{\otimes n}\longrightarrow\orep$, see Theorem \ref{th3} below. This is the reason why we treat this algebra as an analog of $\orep$. 

\begin{remark}
    Note that it's quite natural to take $A_{\natural}\oplus A[-1]$ in Definition \ref{def4}. The algebra $\orep$ is generated by $a_{ij}$ and seemingly, we do not need anything else other than $A$ placed in the first degree component. This is because we require not only an element of $A$, but also one additional piece of data (a pair of indices ranging from $1$ to $N$) in order to produce an element of $\orep$. However, there are important distinguished elements of $\orep$, such as $\tr(f)$ for $f\in A_{\natural}$, whose definition requires nothing more than an element of $A_{\natural}$, making it natural to place them in the zero degree component.
\end{remark}

Below we introduce two graded linear maps $\pi,\widehat{1}:\O(A)\longrightarrow\O(A)$ of degrees $-1$ and $+1$ respectively, which are responsible for the multiplication in $A$ and the unity. The operator $\pi$ is related to the wheelgebra structure on $\O(A)$, see \cite{ginzburg2010differentialoperatorsbvstructures}.

\begin{definition}
    The graded linear map $\pi:\O(A)\longrightarrow \O(A)$ is of degree $-1$, i.e. $\pi:\O(A)^{(n)}\longrightarrow \O(A)^{(n-1)}$, and it is given by
\begin{equation}
    \pi(a\otimes f\otimes u)=\begin{cases}
        (a_2\otimes\ldots\otimes a_n)\otimes \overline{a_1}\cdot f\otimes u_*, &\text{if}\ u(1)=1,\\[5pt]
        m_{1,u(1)}(a) \otimes f\otimes u_*, &\text{if}\ u(1)>1,
    \end{cases}
\end{equation}

where 
\begin{itemize}
    \item $a=a_1\otimes\ldots\otimes a_n\in A^{\otimes n}$, $f\in\s$, and $u\in S(n)$;
    
    \item the map $m_{1,k}:A^{\otimes n}\rightarrow A^{\otimes n-1}$ multiplies the $k$-th component by the first one, i.e. $$m_{1,k}(a_1\otimes\ldots\otimes a_n)=a_2\otimes\ldots\otimes a_{k-1}\otimes a_1a_k\otimes a_{k+1}\otimes\ldots\otimes a_n.$$
\end{itemize}

For instance, $\pi(a\otimes f\otimes \operatorname{id}_1)=\mathds{1}\otimes\overline{a}\cdot f\otimes\operatorname{id}_0$ and $\pi(\mathds{1}\otimes f\otimes\operatorname{id}_0)=0$ for $a\in A$, $f\in\s$.
\end{definition}

Note also that 
\begin{align}
    \pi(\a\b)=\pi(\a)\b,
\end{align}
if $\a$ is homogeneous of positive degree.

\begin{definition}
    Let the graded linear map $\widehat{1}:\O(A)\longrightarrow\O(A)$ of degree $-1$, i.e. $\widehat{1}:\O(A)^{(n)}\longrightarrow \O(A)^{(n+1)}$, be just multiplication on the left by $1\otimes\mathds{1}\otimes\operatorname{id}_1\in\O(A)^{(1)}$, where $1$ in the first tensor factor stands for the unit in $A$. Explicitly, we have
\begin{equation}
    \widehat{1}(\a)=(1\otimes a)\otimes f\otimes (\operatorname{id}_1\times u)\in\O(A)^{(n+1)},    
\end{equation}
where $\a=a\otimes f\otimes u\in\O(A)^{(n)}$. 
\end{definition}

Obviously, $\pi\circ \widehat{1}:\O(A)\longrightarrow\O(A)$ is the linear map of degree zero that multiplies the second tensor component of $\O(A)^{(n)}$ by $\overline{1}\in A_{\natural}$.

The half of the next proposition related to $\pi$ is known, see Section 3.2 in \cite{ginzburg2010differentialoperatorsbvstructures}. However, we provide the proof for completeness.

\begin{proposition}\label{prop14}
    Operators $\pi:\O(A)^{(n)}\longrightarrow \O(A)^{(n-1)}$ and $\1:\O(A)^{(n-1)}\longrightarrow\O(A)^{(n)}$ are $S(n-1)$ left and right invariant, where the action of $S(n-1)$ on $\O(A)^{(n)}$ is given as the composition of the $S(n)$-action with $S(n-1)\simeq \id_1\times S(n-1)\hookrightarrow S(n)$. Explicitly, this means
    \begin{align}
        \widehat{1} \left(w\cdot \a\right)&=(\operatorname{id}_{1}\times w)\cdot\widehat{1}(\a),\\
        \widehat{1} \left(\a\cdot w\right)&=\widehat{1}(\a)\cdot(\operatorname{id}_{1}\times w),\\
        \pi(\a)\cdot 
        w&=\pi\big(\a\cdot(\operatorname{id}_1\times w)\big),\\
        w\cdot\pi(\a)&=\pi\big((\operatorname{id}_1\times w)\cdot\a\big).
    \end{align}
\end{proposition}
\begin{proof}
The first two identities are clear. The third identity immediately follows from the definition of $\pi$ with the help of Proposition \ref{prop9}. Let us prove the fourth identity. By the very definition we have

\begin{align}
    \hspace{40pt}&\hspace{-40pt}w\cdot\pi(\a)=\begin{cases}
        (a_{w^{-1}(1)+1}\otimes\ldots\otimes a_{w^{-1}(n-1)+1})\otimes \overline{a_1}\cdot f\otimes wu_*, &\text{if}\ u(1)=1,\\[15pt]
        w(m_{1,u(1)}(a)) \otimes f\otimes wu_*, &\text{if}\ u(1)>1,
    \end{cases}\\[10pt]
    &=\begin{cases}
        \pi\bigl((\operatorname{id}_1\times w)\cdot\a\bigr), &\text{if}\ u(1)=1,\\[15pt]
        a_{w^{-1}(1)+1}\otimes\ldots\otimes a_{w^{-1}(i-1)+1}\otimes (a_1a_{u(1)})\otimes a_{w^{-1}(i+1)+1}\otimes \ldots\otimes a_{w^{-1}(n-1)+1}\\[5pt]
        \hspace{300pt}\otimes f\otimes wu_*, &\text{if}\ u(1)>1,
    \end{cases}
\end{align}
where in the second case $i=w(u(1)-1)$. In the first case the claim clearly holds. For the second case we note that $i=v(1)-1$, where $v=(\operatorname{id}_1\times w)u$, and
\begin{align}
    &a_{w^{-1}(1)+1}\otimes\ldots\otimes a_{w^{-1}(i-1)+1}\otimes (a_1a_{u(1)})\otimes a_{w^{-1}(i+1)+1}\otimes \ldots\otimes a_{w^{-1}(n-1)+1}\\
    &\hspace{320pt}=m_{1,i+1}\big((\operatorname{id}_1\times w)(a_1\otimes\ldots\otimes a_n)\big).
\end{align}
Then we apply Proposition \ref{prop9} to complete the proof.
\end{proof}

Below, we discuss a relation between operators $\pi$ and $\widehat{1}$. 

\begin{proposition}\label{prop8}
    For any $r\geq 0$ one has
        \begin{align}
            \widehat{1}^{\, r}\pi=\pi \operatorname{Ad}((12)^{r,1})\,\widehat{1}^{\, r}.
        \end{align}
\end{proposition}
\begin{proof}
Let us apply the left-hand side of the identity above to an element of $\O(A)$ of the form $\a=a\otimes f\otimes u$, where $a\in A^{\otimes n}$, $f\in\s$, and $u\in S(n)$. By the very definition we have
    \begin{align}
    \widehat{1}^{\, r}(\pi(a\otimes f\otimes u))&=\begin{cases}
        (\underbrace{1\otimes \ldots\otimes 1}_{r}\otimes a_2\otimes\ldots\otimes a_n)\otimes \overline{a_1}\cdot f\otimes \operatorname{id}_r\times u_*, &\text{if}\ u(1)=1,\\[15pt]
        (\underbrace{1\otimes \ldots\otimes 1}_{r}\otimes m_{1,u(1)}(a)) \otimes f\otimes \operatorname{id}_r\times u_*, &\text{if}\ u(1)>1,
    \end{cases}\\[10pt]
    &=\begin{cases}
        \pi\Bigl((a_1\otimes \underbrace{1\otimes \ldots\otimes 1}_{r}\otimes a_2\otimes\ldots\otimes a_n)\otimes f\otimes \operatorname{id}_{r+1}\times u_*\Bigr), &\text{if}\ u(1)=1,\\[15pt]
        m_{1,u(1)+r}(a_1\otimes \underbrace{1\otimes \ldots\otimes 1}_{r}\otimes a_2\otimes\ldots\otimes a_n) \otimes f\otimes \operatorname{id}_r\times u_*, &\text{if}\ u(1)>1,
    \end{cases}
    \end{align}

We will assume that $u(1)>1$, because the other case is trivial. Let us set $w=u'\underset{1}{*}\operatorname{id}_r\underset{n-1}{*} u''$. Then one has $w_*=\operatorname{id}_r\times u_*$.

Finally, we apply Lemma \ref{lemma10}. By this lemma in both cases, $u(1)=1$ and $u(1)>1$, we have
\begin{align}
    \widehat{1}^{\, r}(\pi(a\otimes f\otimes u))&=\pi\Bigl((a_1\otimes \underbrace{1\otimes \ldots\otimes 1}_{r}\otimes a_2\otimes\ldots\otimes a_n)\otimes f\otimes u'\underset{1}{*}\operatorname{id}_r\underset{n-1}{*} u''\Bigr)\\
    &=\pi \operatorname{Ad}((12\ldots r+1))\,\widehat{1}^{\, r} (\a)
\end{align}
\end{proof}

\subsection{\texorpdfstring{$\O(A)$}{OA} and the Kontsevich-Rosenberg principle}
The content of this subsection is a reformulation of certain results from Section 6 in \cite{ginzburg2010differentialoperatorsbvstructures}. See also Section 7 in \cite{fernández2025symplecticwheelgebrasnoncommutativegeometry}.

We would like to explain how to produce elements of $\orep$ out of elements of $\O(A)$. Assume that two tuples of indices ranging from $1$ to $N$ are given, say $\mathbf{i}=(i_1,\ldots,i_n)$ and $\mathbf{j}=(j_1,\ldots,j_n)$. Take any $a_1,\ldots, a_n\in A$, a permutation $u\in S(n)$, and any $f_1,\ldots, f_m\in A_{\natural}$. Then the element $\a=(a_1\otimes\ldots\otimes a_n) \otimes f_1\cdot\ldots\cdot f_m\otimes u$ belongs to $\O(A)^{(n)}$ and the element of the coordinate ring $\orep$ corresponding to $\mathbf{i}$, $\mathbf{j}$, and $\a$, denoted by $\a_{\mathbf{i}\mathbf{j}}$, is by definition
\begin{equation}\label{f28}
    \a_{\mathbf{i}\mathbf{j}}=(a_1)_{i_{u^{-1}(1)}j_1}\ldots (a_n)_{i_{u^{-1}(n)}j_n}\tr(f_1)\ldots\tr(f_m)\in\orep.
\end{equation}
In the case $m=0$ the last factor consisting of $m$ traces is missing. In the case when $\a\in\O(A)_0$ the first factor with $a$'s is missing and the tuples $\mathbf{i}$, $\mathbf{j}$ are of length zero, so we may simply write $\tr(\a)$ for $\tr(f_1)\ldots\tr(f_m)$. And also the unit of $\O(A)$ goes to the unit of $\orep$.

We will also use the following index-free variation of this notation. We define a linear map $(-,-):\O(A)\odot \TMat\longrightarrow \orep$ by setting
        \begin{equation}
            (\a|E_{i_1j_1}^*\otimes\ldots\otimes E^*_{i_nj_n}):=\a_{\mathbf{ij}}   
        \end{equation}
for any $\a\in\O(A)^{(n)}$, $\mathbf{i}=(i_1,\ldots,i_n)\in[1,N]$, $\mathbf{j}=(j_1,\ldots,j_n)\in[1,N]^n$, and extending this definition to arbitrary elements of $\Mat^{\otimes n}$ by linearity. That is, if 
\begin{align}
    X=\sum\limits_{\mathbf{i},\mathbf{j}\in[1,N]^n} X_{\mathbf{i}\mathbf{j}}E^*_{\mathbf{i}\mathbf{j}}\in\Mat^{\otimes n},\hspace{25pt} X_{\mathbf{i}\mathbf{j}}\in\Bbbk,
\end{align}
where $E^*_{\mathbf{i}\mathbf{j}}:=E^*_{i_1j_1}\otimes\ldots\otimes E^*_{i_nj_n}$, then for any $\a\in\O(A)^{(n)}$ we set
\begin{align}
    \left(\a\middle|X\right)=\sum\limits_{\mathbf{i},\mathbf{j}\in[1,N]^n} X_{\mathbf{i}\mathbf{j}}\cdot\a_{\mathbf{i}\mathbf{j}}.
\end{align}

\begin{remark}
It may seem that this definition of $\alpha_{\mathbf{i}\mathbf{j}}$ and $(\alpha|X)$ is not canonical, as we heavily relied on the isomorphism from Proposition \ref{prop11}. However, we can produce elements of $\orep$ even out of elements of the tensor algebra $\operatorname{T}_{\S}\big(A_{\natural}\oplus A[-1]\big)$ without passing to $\O(A)=\SSym\big(A_{\natural}\oplus A[-1]\big)$. Indeed, take $n$-tuples of indices $\mathbf{i},\mathbf{j}\in[1,N]^n$, and an element $\alpha=\sigma\circ a_1\otimes a_2\otimes \ldots\otimes a_{n+m}\circ\tau\in \operatorname{T}_{\S}\big(A_{\natural}\oplus A[-1]\big)^{(n)}$, where the number of $a_i$'s belonging to $A$ is $n$, all other $a_i$'s belong to $A_{\natural}$, and $\sigma,\tau\in S(n)$. Then we define $\alpha_{\mathbf{i}\mathbf{j}}$ as the product 
\begin{align}\label{f64}
    \prod\limits_{p=1}^{n+m}b_p,
\end{align}
where $b_p=\tr(a_p)$, if $a_p\in A_{\natural}$, and $b_p=(a_p)_{i_{\tau^{-1}(p)}j_{\sigma(p)}}$ otherwise. Since $\orep$ is commutative, one can reorder factors in \eqref{f64} in such a way that $j$-indices go in the increasing order from $1$ to $n$ and all the traces are on the right, which gives exactly \eqref{f28}. On the level of the tensor algebra $\operatorname{T}_{\S}\big(A_{\natural}\oplus A[-1]\big)$ this reordering corresponds to modding out by the action of the symmetric group that permutes the tensor factors.  
\end{remark}

Now we want to explain what corresponds to the multiplication and the diagonal $\S$-bimodule structure of $\O(A)$ on the side of $\orep$.

\subsubsection{Multiplication in \texorpdfstring{$\O(A)$}{OA} and \texorpdfstring{$\orep$}{ORepNA}}
The multiplication in $\O(A)$ is consistent with the multiplication in $\orep$ in the following sense
\begin{align}\label{f40}
    (\a|X)(\b|Y)=\left(\a\b\middle|X\otimes Y\right),
\end{align}
where $\a,\b\in\O(A)$ are homogeneous and $X\in\Mat^{\otimes |\a|}$, $Y\in\Mat^{\otimes |\b|}$.

\subsubsection{\texorpdfstring{$\S$}{S}-bimodule structure on \texorpdfstring{$\O(A)$}{OA} and commutativity of \texorpdfstring{$\orep$}{ORepNA}}

The diagonal $\S$-bimodule structure on $\O(A)$ is related to commutativity in $\orep$ through the following identities
\begin{align}
    &(\sigma\cdot\a)_{\mathbf{i}\mathbf{j}}=\a_{\mathbf{i}\,\sigma^{-1}(\mathbf{j})},\\[5pt]
    &(\a\cdot\sigma)_{\mathbf{i}\mathbf{j}}=\a_{\sigma(\mathbf{i})\mathbf{j}},
\end{align}
where $\a\in \O(A)^{(n)}$, $\sigma\in S(n)$, tuples of indices $\mathbf{i}$, $\mathbf{j}$ are of length $n$, and $\sigma(\mathbf{i})=(i_{\sigma^{-1}(1)},\ldots,i_{\sigma^{-1}(n)})$, $\sigma^{-1}(\mathbf{j})=(j_{\sigma(1)},\ldots,j_{\sigma(n)})$.

In the index-free notation this can be rewritten as
\begin{align}\label{f65}
    \left(\sigma\cdot\a\cdot\tau\middle|X\right)=\left(\a\middle|\operatorname{vect}(\tau)\operatorname{covect}(\sigma^{-1})(X)\right).
\end{align}
Particularly, this means that the pairing $(-,-):\O(A)\odot \TMat\longrightarrow \orep$ is invariant with respect to the adjoint $\S$-module structure on $\O(A)$, i.e.
\begin{align}
    \left(\operatorname{Ad}(\sigma)(\a)\middle| \sigma X\right)=\left(\a\middle|X\right)
\end{align}
for $\a\in\O(A)^{(n)}$, $X\in\Mat^{\otimes n}$, and $\sigma\in S(n)$.

\subsubsection{\texorpdfstring{$\O(A)^{\stimes k}$}{OA} and \texorpdfstring{$\orep^{\otimes k}$}{ORepNA}}

Below we deal not only with $\orep$ but with its tensor powers as well, so we will need to construct elements of $\orep^{\otimes k}$ out of $\O(A)^{\stimes k}$. Here we explain how to do it. 

Take $\a=\sigma\circ \alpha_1\otimes\ldots\otimes \a_k\circ\tau\in\O(A)^{\stimes k}$, where $\a_i\in \O(A)$ are homogeneous, and $\sigma,\tau\in S(|\a_1|+\ldots+|\a_k|)$. Suppose that $X\in \Mat^{|\a_1|+\ldots+|\a_k|}$ is a product of matrix units. Let $X_1\in\Mat^{\otimes |\a_1|}$ be the first $|\a_1|$ tensor factors from $\operatorname{vect}(\tau)\operatorname{covect}(\sigma^{-1})X$, $X_2\in\Mat^{\otimes |\a_2|}$ be the next $|\a_2|$ tensor factors and so on. Then obviously $\operatorname{vect}(\tau)\operatorname{covect}(\sigma^{-1})X=X_1\otimes \ldots\otimes X_k$ and we can define
\begin{align}
    (\a| X):=(\a_1| X_1)\otimes\ldots\otimes (\a_k| X_k)\in\orep^{\otimes k}.
\end{align}

Compare this definition to \eqref{f64}. One can easily check that \eqref{f40} and \eqref{f65} hold in this case as well. Recall that the product in $\O(A)^{\stimes k}$ is not just concatenation of tensors -- there is also a twist by a permutation, see Definition \ref{def1}.%If you need a proof just look at Image4.

There is a natural action of $S(k)$ on $\orep^{\otimes k}$ by permutations of tensor components. It is related to a similar action on $\O(A)^{\stimes k}$ via
\begin{align}
    (w(\a)|X)=w(\a|X),\hspace{50pt} w\in S(k).
\end{align}%proof is in Image5

\subsubsection{Linear maps \texorpdfstring{$\pi$}{pi}, \texorpdfstring{$\widehat{1}$}{1} and relations in \texorpdfstring{$\orep$}{ORepNA}}

The map $\pi$ is responsible for the relation in $\orep$ involving multiplication in $A$, relation \eqref{f17}, in the following sense.

\begin{proposition}\label{prop7}
For all $\a\in\O(A)^{(n)}$  with $n\geq 1$ and all $X\in\Mat^{\otimes n-1}$ one has
    \begin{align}\label{f2}
    \left(\pi(\a)\middle|X\right)=\left(\a\middle|\tr\otimes X\right).
\end{align}
\end{proposition}
\begin{remark}
    Note that if we take $\a=(a\otimes b)\otimes \mathds{1}\otimes (12)\in\O(A)_2$ in \eqref{f2}, where $\mathds{1}$ is the identity in $\s$, we get exactly \eqref{f17}.
\end{remark}
\begin{proof}[Proof of Proposition \ref{prop7}]
    Let $\a=a\otimes f\otimes u$. In the case $u(1)=1$ the verification is trivial. Let's assume that $u(1)>1$. Take $X=E^*_{\mathbf{i}\mathbf{j}}$ for simplicity. Then by the very definition we have
    \begin{equation}
        \sum\limits_{k=1}^N \a_{k\sqcup\mathbf{i},k\sqcup\mathbf{j}}=\sum\limits_{k=1}^N(a\otimes f\otimes\operatorname{id}_1)_{u(k\sqcup \mathbf{i}),k\sqcup\mathbf{j}}=\sum\limits_{k=1}^N(a\otimes f\otimes\operatorname{id}_1)_{(i_{u^{-1}(1)-1},\ldots,k,\ldots,i_{u^{-1}(n+1)-1}),k\sqcup\mathbf{j}},
    \end{equation}
    where $k$ is placed on the $u(1)$-st position and the indices of all $i$'s are shifted by minus one. Let us expand the right-hand side. The box below doesn't have a mathematical meaning, we use it to highlight the factors related to $k$.

\begin{align}\label{f43}
    \sum\limits_{k=1}^N \a_{k\sqcup\mathbf{i},k\sqcup\mathbf{j}}&=\sum\limits_{k=1}^N\boxed{(a_1)_{i_{u^{-1}(1)-1},k}}(a_2)_{i_{u^{-1}(2)-1},j_1}\ldots (a_{u(1)-1})_{i_{u^{-1}(u(1)-1)-1},j_{u(1)-2}}\notag\\[7pt]
    &\pushright{\boxed{(a_{u(1)})_{k,j_{u(1)-1}}}(a_{u(1)+1})_{i_{u^{-1}(u(1)+1)-1},j_{u(1)}}\ldots (a_{n+1})_{i_{u^{-1}(n+1)-1}j_n}\tr(f)}\notag\\[13pt]
    &=(a_2)_{i_{u^{-1}(2)-1},j_1}\ldots (a_{u(1)-1})_{i_{u^{-1}(u(1)-1)-1},j_{u(1)-2}}\\[10pt]
    &\phantom{aaaa}\boxed{(a_1a_{u(1)})_{i_{u^{-1}(1)-1},j_{u(1)-1}}}(a_{u(1)+1})_{i_{u^{-1}(u(1)+1)-1},j_{u(1)}}\ldots (a_{n+1})_{i_{u^{-1}(n+1)-1}j_n}\tr(f).
\end{align}

Note that
    \begin{equation}
        u_*^{-1}=(u^{-1}(2)-1,u^{-1}(3)-1,\ldots,u^{-1}(u(1)-1)-1,u^{-1}(1)-1,u^{-1}(u(1)+1)-1, \ldots,u^{-1}(n+1)-1),
    \end{equation}
    so the right-hand side of \eqref{f43} equals $\pi(\a)_{\mathbf{i}\mathbf{j}}$.
\end{proof}

The linear map $\widehat{1}$ is responsible for relation \eqref{f36} in $\orep$ in the following sense
\begin{align}
    \left(\widehat{1}(\a)\middle| X\right)=\left(\a\middle|\varepsilon_1(X)\right)    
\end{align}
for $\a\in\O(A)^{(n)}$, $X\in\Mat^{\otimes n+1}$, where $\varepsilon$ is the counit $\varepsilon:\Mat\longrightarrow \Bbbk$, and $\varepsilon_1$ means that we apply $\varepsilon$ to the first tensor factor, so $\varepsilon_1:\Mat^{\otimes n+1}\longrightarrow \Mat^{\otimes n}$.

\subsection{How to obtain \texorpdfstring{$\orep$}{ORepNA} from \texorpdfstring{$\O(A)$}{TA}?}

The answer is given in Proposition \ref{prop6} below. 

\begin{lemma}\label{lemma5}
    The subspace of $\O(A)\odot\TMat$ spanned by the elements of the form
        \begin{enumerate}[label=\alph*)]
            \item\label{cond14} $\pi(\a)\otimes X-\a\otimes(\tr\otimes X)$ for all $\a\in\O(A)^{(n)}$ with $n\geq 1$ and all $X\in\Mat^{\otimes n}$;
            \vspace{7pt}

            \item\label{cond15}  $\widehat{1}(\a)\otimes X-\a\otimes\varepsilon_1(X)$ for all $\a\in\O(A)^{(n)}$ and all $X\in\Mat^{\otimes n+1}$;
            \vspace{7pt}

            \item\label{cond10} $(\sigma\cdot \a)\otimes E^*_{\mathbf{i}\mathbf{j}}-\a\otimes E^*_{\mathbf{i}\sigma^{-1}(\mathbf{j})}$ and $(\a\cdot\sigma)\otimes E^*_{\mathbf{i}\mathbf{j}}-\a\otimes E^*_{\sigma(\mathbf{i})\mathbf{j}}$ for all $\a\in\O(A)^{(n)}$, $\sigma\in S(n)$, and all tuples of indices $\mathbf{i}$,$\mathbf{j}$
        \end{enumerate}
        is a two-sided ideal.
\end{lemma}
\begin{proof}
The subspace spanned by \ref{cond14} and \ref{cond15} is a right ideal due to $\pi(\a\b)=\pi(\a)\b$ and $\widehat{1}(\a\b)=\widehat{1}(\a)\b$ for all $\a\in\O(A)$ of positive degree. In fact the multiplication in $\O(A)\odot\TMat$ is commutative modulo the subspace spanned by the elements of the form \ref{cond10} and there is nothing to prove for \ref{cond14} and \ref{cond15}. Indeed, let us denote by $V$ the subspace spanned by \ref{cond10} and let $\a\in\O(A)^{(n)}$, $\b\in\O(A)^{(m)}$. Then one has
    \begin{multline}
        (\a\otimes E^*_{\mathbf{i_1}\mathbf{j_1}})(\b\otimes E^*_{\mathbf{i_2}\mathbf{j_2}})=(\a\b)\otimes E^*_{\mathbf{i_1}\sqcup\mathbf{i_2},\,\mathbf{j_1}\sqcup\mathbf{j_2}}=(\a\b)\otimes E^*_{(12)^{m,n}(\mathbf{i_2}\sqcup\mathbf{i_1}),\ (12)^{m,n}(\mathbf{j_2}\sqcup\mathbf{j_1})}\\
        =\left((12)^{n,m}\cdot\a\b\cdot (12)^{n,m}\right)\otimes E^*_{\mathbf{i_2}\sqcup\mathbf{i_1},\ \mathbf{j_2}\sqcup\mathbf{j_1}}\operatorname{mod}\ V=(\b\a)\otimes E^*_{\mathbf{i_2}\sqcup\mathbf{i_1},\,\mathbf{j_2}\sqcup\mathbf{j_1}}=(\b\otimes E^*_{\mathbf{i_2}\mathbf{j_2}})(\a\otimes E^*_{\mathbf{i_1}\mathbf{j_1}}).
    \end{multline}

Let us prove now that the vector space $V$ is a two-sided ideal. Let $\a=a\otimes f\otimes u\in\O(A)^{(n)}$ and $\b=b\otimes g\otimes v\in\O(A)_m$, and take any $\sigma\in S(n)$. Then 
    \begin{multline}
        (\b\otimes E^*_{\mathbf{k},\mathbf{l}})\left((\sigma\cdot\a)\otimes E^*_{\mathbf{i}\mathbf{j}}-\a\otimes E^*_{\mathbf{i}\,\sigma^{-1}(\mathbf{j})}\right)\\
        =(b\otimes \sigma(a))\otimes fg\otimes (v\times\sigma u))\otimes E^*_{\mathbf{k}\sqcup\mathbf{i},\ \mathbf{l}\sqcup\mathbf{j}}-\left((b\otimes a)\otimes fg\otimes (v\times u)\right)\otimes E^*_{\mathbf{k}\sqcup\mathbf{i},\ \mathbf{l}\sqcup\sigma^{-1}(\mathbf{j})}.
    \end{multline}
    Let us denote $\operatorname{id}_m\times \sigma\in S(n+m)$ by $\sigma_+$. Then the right-hand side equals
    \begin{align}
        \left(\sigma_+\cdot (\b\a)\right)\otimes E^*_{\mathbf{k}\sqcup\mathbf{i},\ \mathbf{l}\sqcup\mathbf{j}}-\left(\b\a\right)\otimes E^*_{\mathbf{k}\sqcup\mathbf{i},\sigma_+^{-1}(\mathbf{l}\sqcup\mathbf{j})}\in V
    \end{align}
    as desired. Similarly, for the remaining part of \ref{cond10}.
    %\begin{align}
        %\b_{\mathbf{k},\mathbf{l}}\left((\a\cdot\sigma)_{\mathbf{i}\mathbf{j}}-\a_{\sigma(\mathbf{i})\mathbf{j}}\right)&=\left(b\otimes a\otimes \imath(v,u\sigma)\otimes f\cdot g\right)_{\mathbf{k}\sqcup\mathbf{i},\ \mathbf{l}\sqcup\mathbf{j}}-\left(b\otimes a\otimes \imath(v,u)\otimes f\cdot g\right)_{\mathbf{k}\sqcup\sigma(\mathbf{i}),\ \mathbf{l}\sqcup\mathbf{j}}\\[5pt]
        %&=\left(b\otimes a\otimes \imath(v,u)\sigma_+\otimes f\cdot g\right)_{\mathbf{k}\sqcup\mathbf{i},\ \mathbf{l}\sqcup\mathbf{j}}-\left(b\otimes a\otimes \imath(v,u)\otimes f\cdot g\right)_{\sigma_+(\mathbf{k}\sqcup\mathbf{i}),\ \mathbf{l}\sqcup\mathbf{j}}\\[5pt]
        %&=\left((\b\a)\cdot\sigma_+\right)_{\mathbf{k}\sqcup\mathbf{i},\ \mathbf{l}\sqcup\mathbf{j}}-\left(\b\a\right)_{\sigma_+(\mathbf{k}\sqcup\mathbf{i}),\ \mathbf{l}\sqcup\mathbf{j}}\in V.
    %\end{align}   
\end{proof}

\begin{remark}\label{rem1}
    It's easy to see that the ideal in Lemma \ref{lemma5} coincides with the whole algebra $\O(A)\odot\TMat$ in the case of $A=\Bbbk[x,\partial]$ the (first) Weyl algebra. Indeed, one has
    \begin{align}
        \a\otimes X=\frac{1}{N}\cdot \a\otimes \varepsilon_1(\tr\otimes X)=\frac{1}{N}\cdot \big(\pi\,\widehat{1}(\a)\big)\otimes X\ \operatorname{mod}\ \text{the ideal}
    \end{align}
    But $\pi\,\widehat{1}$ is the multiplication of the second tensor component by $\overline{1}\in A_{\natural}$, which is zero as $1=[\partial,x]\subset [A,A]$. Note also that in this case $A_{\natural}=0$ because any element of $A$ belongs to the subspace of commutators $[A,A]$ due to the identity $x^n\partial^m=\left[\partial,\frac{x^{n+1}}{n+1}\partial^m\right]$. 
\end{remark}

\begin{proposition}\label{prop6}
    The algebra $\orep$ is isomorphic to the quotient of $\O(A)\odot\TMat$ by the two-sided ideal from Lemma \ref{lemma5}.
\end{proposition}
\begin{proof}
There is a surjective homomorphism from $\O(A)\odot\TMat$ to $\orep$ defined by $\a\otimes X\mapsto \left(\a\middle| X\right)$. It's clear that \ref{cond14}, \ref{cond15}, and \ref{cond10} from Lemma \ref{lemma5} belong to its kernel.

Let us denote by $R$ the quotient of $\O(A)\odot\TMat$ by the ideal from Lemma \ref{lemma5}. Then there is a surjective homomorphism, say $\phi$, from $R$ to $\orep$. On the other hand, $R$ is generated by images of the following elements of $\O(A)\odot\TMat$
\begin{align}
    \left(a\otimes\mathds{1}\otimes\operatorname{id}_1\right)\otimes E^*_{ij},  
\end{align}
where $a\otimes\mathds{1}\otimes\operatorname{id}_1\in \O(A)$ and $i,j$ are integers from $1$ to $N$. Let us (temporarily) denote these elements of $R$ simply by $a_{ij}$. One can easily see from \ref{cond14}, \ref{cond15}, and \ref{cond10} in Lemma \ref{lemma5} that they satisfy the following relations
\begin{align}
(ab)_{ij}=\sum_{r=1}^N a_{ir}b_{rj},&&(1)_{ij}=\delta_{ij}.
\end{align}
Hence there is a surjective homomorphism, say $\psi$, from $\orep$ to $R$. It's clear that $\phi$ and $\psi$ are mutually inverse, thus, $R$ and $\orep$ are isomorphic.
\end{proof}

\subsection{Ideals \texorpdfstring{$\R_N(A)$}{RNA} and \texorpdfstring{$\R(A)$}{RA} of \texorpdfstring{$\O(A)$}{TA}}

The main object of study in the present paper is not the associative algebra $A$ itself, but rather the collection of representation spaces $\rep$ for $N \geq 1$. Consequently, it is no surprise that elements of the algebra $A$ or $\O(A)$ that produce zero functions on all $\rep$ do not make any difference. In the present section we mention several properties of such elements.

By Proposition \ref{prop6} $\orep$ is the quotient of $\O(A)\odot\TMat$ by the ideal from Lemma \ref{lemma5}. Let us denote the quotient map by
\begin{align}
    \theta_N:\O(A)\odot\TMat\longrightarrow \orep.
\end{align}

\begin{definition}
    Denote by $\R_N(A)$ the subset of $\O(A)$ spanned by the graded elements $\a\in\O(A)^{(n)}$ such that for any $X\in\Mat^{\otimes n}$ one has $\a\otimes X\in\operatorname{ker}(\theta_N)$ or, equivalently, $(a|X)=0\in\orep$. Obviously, $\R_N(A)$ is a graded two-sided ideal.
\end{definition}

Note that $\R_N(A)$ is never zero because $N\a-\pi\widehat{1}(\a)\in\R_N(A)$ for any $\a\in\O(A)$ due to the fact that $\pi\widehat{1}$ is the multiplication of the second tensor component of $\O(A)$ by $\overline{1}\in A_{\natural}$, which is the multiplication by $N$ modulo $\operatorname{ker}(\theta_N)$:
    \begin{align}
        N\a\otimes X=\a\otimes \varepsilon_1(\tr\otimes X)=\big(\pi\,\widehat{1}(\a)\big)\otimes X\ \operatorname{mod}\ \operatorname{ker}(\theta_N).
    \end{align}

    \begin{example}
        From Remark \ref{rem1} it follows that, if $A$ is the Weyl algebra $A=\Bbbk[x,\partial]$, then $\R_N(A)=\O(A)$ for any $N$.  
    \end{example}

    \begin{proposition}\label{prop13}
        The ideal $\R_N(A)$ is invariant under the action of $\pi$ and $\widehat{1}$, and each its graded component $\R_N(A)_n$ is invariant under the natural $S(n)$-bimodule action.
    \end{proposition}
    \begin{proof}
        \textit{$\pi$-invariance:} Let $\a\in\R_N(A)_n$. Then by the very definition of $\R_N(A)$ we have $(a|X)=0\in\orep$ for any $X\in\Mat^{\otimes n}$. Particularly, $(a|\tr\otimes X)=0$ for any $X\in\Mat^{\otimes n-1}$, and by relation \ref{cond14} from Lemma \ref{lemma5} we have $(\pi(a)|X)=0$, so $\R_N(A)$ is invariant under $\pi$.

        \textit{$\1$-invariance:} Take the same $\a$ and $X$ as above and consider the expression $(a|\varepsilon(E^*_{ij})\otimes X)$. It equals zero, if $i\neq j$, and if $i=j$ then it equals zero by the very definition of $\R_N(A)$. So we conclude that $(a|\varepsilon_1(X))=0$ for any $X\in\Mat^{\otimes n+1}$. Thus, by relation \ref{cond15} from Lemma \ref{lemma5} $(\widehat{1}(\a)|X)=0$ and $\R_N(A)$ is $\widehat{1}$-invariant.

        \textit{$S(n)$-bimodule action:} Let $\a\in\R_N(A)_n$ be a graded element of degree $n$. Take any $\sigma,\tau\in S(n)$. Then by the very definition of $\R_N(A)$ we have $(a|\vect(\tau)\covect(\sigma^{-1})X)=0$ for any $X\in\Mat^{\otimes n}$. Then by relations \ref{cond10} from Lemma \ref{lemma5} we have $(\sigma\cdot a\cdot\tau|X)=0$ and $\R_N(A)_n$ is closed under the $S(n)$-bimodule action.
    \end{proof}

    So, if we were interested in only one particular representation space $\rep$ of a fixed dimension $N$, then we should actually work with the quotient $\faktor{\O(A)}{\R_N(A)}$, rather than the algebra $\O(A)$ itself. Note this this quotient is a di-twisted algebra as well. Since we wish to study all the representation spaces at once, we should consider the ideal 
    \begin{align}
        \R(A):=\bigcap\limits_{N\geq 1}\R_N(A)\subset \O(A).    
    \end{align}
    Obviously, for the Weyl algebras this ideal coincides with $\O(A)$, but in contrast to ideals $\R_N(A)$ the ideal $\R(A)$ can be zero. 

    \begin{proposition}\label{prop5}
    Let $A$ be the algebra of noncommutative polynomials in $d$ variables, i.e. $A=\Bbbk\langle x_1,\ldots,x_d\rangle$. Then $\R(A)$ is the zero ideal.
    \end{proposition}
    \begin{proof}
Below we will prove that if $\a\in\O(A)^{(n)}$ and $\a_{\mathbf{i}\mathbf{j}}=0\in\orep$ for any $N$ and any $n$-tuples $\mathbf{i}$, $\mathbf{j}$ of integers ranging from $1$ to $N$, then $\a=0\in \O(A)$. The proof is based on the answer of Will Sawin to a question on MathOverflow, see \cite{mathoverflowSawin}. Namely, we first reduce the general claim to the case $n=0$, and then apply the main argument of \cite{mathoverflowSawin}.

    Let us choose a basis in $\Bbbk\langle x_1,\ldots,x_d\rangle$ consisting of the words in $d$ letters. Denote the set of such words by $W_d$ and the subset of words of length not greater than $k$ by $W_d^{\leq k}$. Consider the cyclic words in $d$ letters, which are equivalence classes of the ordinary words under the cyclic permutations. For any cyclic word in $d$ letters we can pick a representative from $W_d$. Let us denote by $\widetilde{W_d}\subset W_d$ such set of representatives of all cyclic words that none of the elements of $\widetilde{W_d}$ is the cyclic permutation of the other. Such set of representatives is not unique, but we fix it now and will not change later. By $\widetilde{W_d}^{\leq k}$ we denote elements of $\widetilde{W_d}$ of length not greater that $k$. Cyclic words provide a basis of the vector space $A_{\natural}$ and we may treat each element of $\s$ as a polynomial in $\widetilde{W_d}$.

    Any $\a\in \O(A)^{(n)}$ can be uniquely written in the form
    \begin{equation}
        \a=\sum\limits_{\sigma\in S(n)}\sum\limits_{w^{\sigma}_1,\ldots, w^{\sigma}_n\in W_d^{\leq k}} \big(w_1^{\sigma}\otimes \ldots \otimes w_n^{\sigma}\big)\otimes f_{w^{\sigma}_1,\ldots, w^{\sigma}_n}^{\sigma}\otimes \sigma,
    \end{equation}
    where each $f_{w^{\sigma}_1,\ldots, w^{\sigma}_n}^{\sigma}$ is a polynomial in $\widetilde{W_d}^{\leq k}$, for some fixed $k$ which depends on $\a$. Our main goal is to prove that if $\a_{\mathbf{i}\mathbf{j}}=0$ for any $\mathbf{i}$, $\mathbf{j}$, then all $f_{w^{\sigma}_1,\ldots, w^{\sigma}_n}^{\sigma}\in\s$ are zero.

    Below we treat $\a_{\mathbf{i}\mathbf{j}}$ as a function on $\operatorname{Hom}_{\mathbf{Alg}}(A,\operatorname{Mat}_N(\Bbbk))$, which is just the affine space $\operatorname{Mat}_N(\Bbbk)^d\simeq\Bbbk^{d^2N}$.

    One must be careful with the empty cyclic word, because it behaves completely differently compared to all non-empty words. Namely, its trace differs under the usual inclusion of matrices of size $N$ into the matrices of size $N+1$ which adds $2N+1$ zeros to the matrix, while the traces of non-empty words do not change. Below we will assume that $f_{w^{\sigma}_1,\ldots, w^{\sigma}_n}^{\sigma}$ is a polynomial only in non-empty words, because we can expand $\a$ as a polynomial in $\overline{1}\in\s$ and consider the aforementioned inclusions to see that each coefficient of this polynomial must belong to $\R(A)$.

    \textit{Reduction to the case $n=0$.}

    Let $N$ be a positive integer. Consider the left $A$-module $V_N=\faktor{A}{I_N}$, where $I_N$ is the two-sided ideal of all polynomials of total degree at least $N+1$. The vector space $V_N$ is finite dimensional, so the images of $x_1,\ldots,x_d$ in this representation of $A$ are given by lower triangular matrices in the natural basis $W_d^{\leq N}$. Let us denote them by $L_1^{(N)},\ldots,L_d^{(N)}$. 

    Let $X_1,\ldots,X_d$ be arbitrary $\Bbbk$-valued matrices of a common size. For any $p=1,\ldots,d$ we set 
    \begin{equation}\label{f37}
        \widehat{X_p}=X_p\oplus \left(L_p^{(N)}\right)^{\oplus n}.
    \end{equation}
    For any $v\in W_d$ one clearly has $\operatorname{tr}(v(\widehat{X_1},\ldots,\widehat{X_d}))=\operatorname{tr}(v(X_1,\ldots,X_d))$.

    Take any permutation $\tau\in S(n)$. Suppose that $\mathbf{i}=(i_1,\ldots,i_n)$ and $\mathbf{j}=(j_1,\ldots,j_n)$ are such that each $j_l$ ranges over the $l$-th block in the second summand of \eqref{f37} and each $i_l$ ranges over the $\tau(l)$-th block of the same summand. Let us denote the $l$-th block of the second summand of \eqref{f37} by $L_{p,l}^{(N)}$. Then one has
    \begin{multline}
        \a_{\mathbf{i}\mathbf{j}}(\widehat{X_1},\ldots,\widehat{X_d})\\
        =\sum\limits_{w^{\tau}_1,\ldots, w^{\tau}_n\in W_d^{\leq k}} \left(w_1^{\tau}(L_{1,1}^{(N)},\ldots,L_{d,1}^{(N)})\right)_{i_{\tau^{-1}(1)} j_1} \ldots \left(w_n^{\tau}(L_{1,n}^{(N)},\ldots,L_{d,n}^{(N)})\right)_{i_{\tau^{-1}(n)} j_n}\\
        \times\operatorname{tr}(f_{w^{\tau}_1,\ldots, w^{\tau}_n}^{\tau})(X_1,\ldots,X_d).
    \end{multline}
    The summands with $\sigma\neq\tau$ vanish because $i_{\sigma^{-1}(l)}$ and $j_l$ belong to a common block of the second summand in \eqref{f37} if and only if $\tau(\sigma^{-1}(l))=l$. Now we can relabel the $\mathbf{i}$-indices and remove them at all by rewriting the identity $\a_{\mathbf{i}\mathbf{j}}=0$ as an identity for operators acting on $V_N^{\otimes n}$:
    \begin{equation}
        \sum\limits_{w^{\tau}_1,\ldots, w^{\tau}_n\in W_d^{\leq k}} \operatorname{tr}(f_{w^{\tau}_1,\ldots, w^{\tau}_n}^{\tau})(X_1,\ldots,X_d)\cdot w_1^{\tau}(L_{1,1}^{(N)},\ldots,L_{d,1}^{(N)})\otimes \ldots \otimes w_n^{\tau}(L_{1,n}^{(N)},\ldots,L_{d,n}^{(N)})=0.
    \end{equation}
    Then we can apply this operator to $1\otimes\ldots\otimes 1\in (V_N)^{\otimes n}$. Then
    \begin{equation}\label{f33}
        \sum\limits_{w^{\tau}_1,\ldots, w^{\tau}_n\in W_d^{\leq k}} \operatorname{tr}(f_{w^{\tau}_1,\ldots, w^{\tau}_n}^{\tau})(X_1,\ldots,X_d)\cdot w_1^{\tau}\otimes \ldots \otimes w_n^{\tau}=0\in (V_N)^{\otimes n}.
    \end{equation}
    Since $(V_N)^{\otimes n}\simeq \faktor{A^{\otimes n}}{I_N^n}$, where $I_N^n=\sum\limits_{l=1}^{n}A^{\otimes l-1}\otimes I_N\otimes A^{n-l}$, the element of $A^{\otimes n}$ defined by the left-hand side of \eqref{f33} must belong to the ideal $I_N^n$ for any $N$. Taking $N=k+1$ we see that this element is zero in $A^{\otimes n}$, hence $\operatorname{tr}(f_{w^{\tau}_1,\ldots, w^{\tau}_n}^{\tau})(X_1,\ldots,X_d)=0$ for any $\tau$, any $w^{\tau}_1,\ldots, w^{\tau}_n\in W_d^{\leq k}$ and any matrices $X_1,\ldots,X_d$ of a common size.

    \textit{Case \texorpdfstring{$n=0$}{n}}
    The desired claim immediately follows from the next lemma. 
    \end{proof}
    Here we reproduce the proof given by Will Sawin \cite{mathoverflowSawin} without any substantial changes. We are grateful to him for this nice proof and publish it here with his kind permission.
    
    \begin{lemma}[W. Sawin \cite{mathoverflowSawin}]\label{lemma20}
        For any collection of constants from $\Bbbk$ indexed by the words from $\widetilde{W_d}^{\leq k}$, say $\{a_{\widetilde{w}}\}$, there exist $\Bbbk$-valued matrices $X_1,\ldots,X_d$ of a certain size such that $\operatorname{tr}(\widetilde{w})(X_1,\ldots,X_d)=a_{\widetilde{w}}$ for any $\widetilde{w}\in\widetilde{W_d}^{\leq k}$.
    \end{lemma}
    \begin{proof}
        Let $p\in\{1,\ldots,d\}$ and $w=w_1\ldots w_l$ be the decomposition of a word $w\in W_d$ into the product of letters. So, $w_1,\ldots, w_l\in\{1,\ldots,d\}$. We define $Y_{p,w}$ as the matrix of the following operator on $\Bbbk^l$, which is almost the cyclical shift,
        \begin{equation}
            Y_{p,w}(e_i)=\delta_{p,w_i}\, e_{i-1\operatorname{mod} l},\ \ \ \ \ \ \ \ \ \ \ \ i=1,\ldots,l,
        \end{equation}
        where $e_1,\ldots,e_l$ is the standard basis of $\Bbbk^l$ and $\delta$ is the Kronecker symbol.
    
        Let $v,w$ be arbitrary words in $d$ letters with the following decompositions into product of letters $v=v_1\ldots v_r$ and $w=w_1\ldots w_l$. Then one has
        \begin{equation}
            v(Y_{1,w},\ldots, Y_{d,w})(e_i)=\delta_{v_1,w_{i-r+1 \operatorname{mod} l}}\delta_{v_2,w_{i-r+2 \operatorname{mod} l}}\ldots \delta_{v_r,w_{i\operatorname{mod} l}}\, e_{i-r\operatorname{mod} l},\ \ \ \ \ \ \ i=1,\ldots,l.
        \end{equation}
        So, $\operatorname{tr}(v(Y_{1,w},\ldots, Y_{d,w}))=0$, unless $r$ is divisible by $l$. For, instance, $\operatorname{tr}(v(Y_{1,w},\ldots, Y_{d,w}))=0$, if $l>r$. If $r=l$, then the trace is non-zero, if and only if $v$ is a cyclic permutation of $w$.

        For any $p=1,\ldots, d$ we set 
    \begin{equation}
        X_p=\bigoplus\limits_{\widetilde{w}\in \widetilde{W_d}^{\leq k}}c_{\widetilde{w}}^p Y_{p,\widetilde{w}},
    \end{equation} 
    where $c_{\widetilde{w}}^p$ are certain constants from $\overline{\Bbbk}$ which are defined below. In fact, we will not determine each of the coefficients $c_{\widetilde{w}}^p$ uniquely, but we will show that the coefficients with the desired property exist. 
    
    Let $\widetilde{w}\in\widetilde{W_d}^{\leq k}$. From the discussion above we have
    \begin{multline}
        \operatorname{tr}(\widetilde{w})(X_1,\ldots,X_d)
        =(c_{\widetilde{w}}^1)^{\operatorname{deg_1}(\widetilde{w})}\ldots (c_{\widetilde{w}}^d)^{\operatorname{deg_d}(\widetilde{w})}\operatorname{tr}(\widetilde{w})(Y_{1,\widetilde{w}},\ldots, Y_{d,\widetilde{w}})\\
        \phantom{aaaaaaaaaaaaaaaaaaaaaa}+\sum\limits_{\substack{\widetilde{u}\in\widetilde{W_d}^{\leq k}\\ l(\widetilde{u})<l(\widetilde{w})}}(c_{\widetilde{u}}^1)^{\operatorname{deg_1}(\widetilde{w})}\ldots (c_{\widetilde{u}}^d)^{\operatorname{deg_d}(\widetilde{w})}\operatorname{tr}(\widetilde{w})(Y_{1,\widetilde{u}},\ldots, Y_{d,\widetilde{u}}),
    \end{multline}
    where the sum runs over all $\widetilde{u}\in\widetilde{W_d}^{\leq k}$ of length strictly smaller than the length of $\widetilde{w}$, and $\operatorname{deg_i}(\widetilde{w})$ is the number of times $i$ occurs in $\widetilde{w}$. There is only one non-zero term corresponding to a word of the same length as $\widetilde{w}$, due to the fact that $\widetilde{W_d}^{\leq k}$ does not contain distinct words that differ by a cyclic permutation. So, we can express the product $(c_{\widetilde{w}}^1)^{\operatorname{deg_1}(\widetilde{w})}\ldots (c_{\widetilde{w}}^d)^{\operatorname{deg_d}(\widetilde{w})}$ in terms of $a_{\widetilde{w}}$ and $c_{\widetilde{u}}^p$ with $l(\widetilde{u})<l(\widetilde{w})$. We do not impose any other restrictions on the coefficients $c_{\widetilde{w}}^p$'s, so the system of equations we obtain after repeating the same procedure for any $\widetilde{w}\in\widetilde{W_d}^{\leq k}$ certainly has a $\overline{\Bbbk}$-valued solution for $c_{\widetilde{w}}^p$'s.
    \end{proof}

\section{Inverting the Kontsevich-Rosenberg principle}\label{section_inverting_KR}

An equivariant version of the Kontsevich-Rosenberg principle could be formulated as follows
\begin{equation}
    \hspace{-0ex}\text{\parbox{0.85\textwidth}{\textit{A noncommutative structure of some kind on $A$ should give an analogous ``commutative'' $\GL$-equivariant structure on all schemes $\rep$ for $N \ge 1$.}}}
\end{equation}

If one wishes to invert this principle, and asks \textit{"How can we define a non-commutative analog of structure $\mathcal{S}$?"}, then one would need to consider an equivariant version of the structure $\mathcal{S}$ for all representation spaces $\rep$ and try to lift it to $A$, or to a related object. In this section we focus on the case when the structure $\mathcal{S}$ is the algebra structure over an operad $\mathcal{P}$, i.e. $\mathcal{S}$ is an operad homomorphism $\mathcal{P}\longrightarrow \mathcal{E}nd_{\orep}^{\operatorname{GL}_N}$, where $\mathcal{E}nd_{\orep}^{\operatorname{GL}_N}$ is the suboperad of $\GL$-invariants of the usual endomorphism operad:
\begin{align}
    \mathcal{E}nd_{\orep}^{\operatorname{GL}_N}=\Bigg(\operatorname{Hom}_{\GL}\left(\orep^{\otimes n},\orep\right)\Bigg)_{n\geq 0}
\end{align}

We claim that in this case the object related to $A$ on which we should define the noncommutative analog of the structure $\mathcal{S}$ is $\O(A)$, or, strictly speaking, $\ored(A):=\faktor{\O(A)}{\R(A)}$. More precisely, a non-commutative analog of the structure $\mathcal{S}$ in this case is an operad homomorphism $\mathcal{P}\longrightarrow \E nd_{A}$, where $\E nd_{A}$ is a certain suboperad of the usual endomorphism operad of $\O(A)$ (or $\ored(A)$), see Definition \ref{def8} below.

\subsection{\texorpdfstring{$\GL$}{GLNk}-equivariant linear maps}
Note that in Theorem \ref{th3} below we use the usual tensor product of vector spaces $\otimes$ instead of the tensor product of diagonal $\S$-bimodules $\stimes$.

\begin{theorem}\label{th3}
A homogeneous linear map $\Phi:\O(A)^{\otimes n}\longrightarrow \O(A)$ of degree zero gives rise to a well-defined $\GL$-equivariant linear map $\varphi_N:\orep^{\otimes n}\longrightarrow \orep$ by the rule
    \begin{align}\label{f57}
        \varphi_N\big((\a_1|X_1),\ldots,(\a_n|X_n)\big)=\Big(\Phi(\a_1,\ldots,\a_n)\Big|X_1\otimes\ldots\otimes X_n\Big),
    \end{align}
    where $\a_1,\ldots,\a_n\in\O(A)$ and $X_1,\ldots,X_n\in\TMat$ are homogeneous elements such that $|\a_i|=|X_i|$, if and only if for any homogeneous elements $\a_1,\ldots,\a_n\in\O(A)$, any permutations $w_1\in S(|\a_1|),\ldots,w_n\in S(|\a_n|)$, and any $i$ the map $\Phi$ satisfies 
    \begin{align}
        &\hspace{-20pt}\Phi\big(\a_1,\ldots,\a_{i-1},\pi(\a_i),\a_{i+1},\ldots,\a_n\big)\equiv\pi\operatorname{Ad}((12)^{|\a_1|+\ldots+|\a_{i-1}|,1})\Phi(\a_1,\ldots,\a_n)\ 
        \operatorname{mod}\ \R_N(A),\label{f25}\\[10pt]
        &\hspace{-20pt}\Phi(\a_1,\ldots,\a_{i-1},\widehat{1}(\a_i),\a_{i+1},\ldots,\a_n)\equiv\operatorname{Ad}((12)^{1,|\a_1|+\ldots+|\a_{i-1}|})\widehat{1}\Phi(\a_1,\ldots,\a_n)\ \operatorname{mod}\ \R_N(A),\label{f39}\\[10pt]
        &\hspace{-20pt}\Phi\big(w_1\cdot \a_1, \ldots,w_n\cdot \a_n\big)\equiv(w_1\times\ldots\times w_n)\cdot\Phi(\a_1,\ldots,\a_n)\ 
        \operatorname{mod}\ \R_N(A),\label{f42}\\[5pt]
        &\hspace{-20pt}\Phi(\a_1\cdot w_1,\ldots,\a_n\cdot w_n)\equiv\Phi(\a_1\ldots,\a_n)\cdot \big(w_1\times\ldots\times w_n\big)\ 
        \operatorname{mod}\ \R_N(A).\label{f58}
    \end{align}

Moreover, any $\GL$-equivariant linear map $\varphi_N:\orep^{\otimes n}\longrightarrow \orep$ is of this form.

\end{theorem}
\begin{corollary}\label{cor3}
     A homogeneous linear map $\Phi$ as above gives rise to a well-defined linear map $\orep^{\otimes n}\longrightarrow \orep$ by rule \eqref{f57} for any $N\geq 1$ if and only if \eqref{f25}-\eqref{f58} are satisfied modulo $\R(A)$. 
\end{corollary}
\begin{proof}[Proof of Theorem \ref{th3}]
First of all note that if map $\varphi_N$ given by \eqref{f57} is well-defined, then it's always $\GL$-equivariant.

Let's assume that \eqref{f57} defines a linear map $\orep^{\otimes n}\longrightarrow\orep$ and prove \eqref{f25}-\eqref{f58}. 

Let us prove \eqref{f25}. In $\orep$ we have $(\pi(\a)|X)=(\a|\tr\otimes X)$ for homogeneous $\a\in\O(A)$ and $X\in\TMat$ such that $|\a|=|X|+1$. Hence in $\orep^{\otimes n}$ one has
\begin{align}
    &(\a_1|X_1)\otimes \ldots \otimes (\a_{i-1}|X_{i-1})\otimes (\pi(\a_i)|X_i)\otimes (\a_{i+1}|X_{i+1})\otimes\ldots\otimes (\a_n|X_n)\\[5pt]
    &\hspace{100pt}=(\a_1|X_1)\otimes \ldots \otimes (\a_{i-1}|X_{i-1})\otimes (\a_i|\tr\otimes X_i)\otimes (\a_{i+1}|X_{i+1})\otimes\ldots\otimes (\a_n|X_n)
\end{align}
for homogeneous elements $\a_1,\ldots,\a_n\in\O(A)$ and $X_1,\ldots,X_n\in\TMat$ such that for any $j=1,\ldots,n$, $j\neq i$ one has $|\a_j|=|X_j|$, and $|\a_i|=|X_i|+1$. 

Apply $\varphi_N$ to the both sides of this identity and transform the right-hand side
\begin{align}
    \hspace{40pt}&\hspace{-40pt}\Big(\Phi(\a_1,\ldots,\a_{i-1},\pi(\a_i),\a_{i+1},\ldots,\a_n)\Big|X_1\otimes \ldots\otimes X_n\Big)\\[5pt]
    &=\Big(\Phi(\a_1,\ldots,\a_n)\Big|X_1\otimes\ldots\otimes X_{i-1}\otimes \tr\otimes X_i\otimes X_{i+1}\otimes\ldots\otimes X_n\Big)\\[5pt]
    &=\Big(\Phi(\a_1,\ldots,\a_n)\Big|(12)^{1,|\a_1\ldots\a_{i-1}|}\tr\otimes X_1\otimes\ldots\otimes X_n\Big)\\[5pt]
    &=\Big(\operatorname{Ad}((12)^{|\a_1\ldots\a_{i-1}|,1})\Phi(\a_1,\ldots,\a_n)\Big|\tr\otimes X_1\otimes\ldots\otimes X_n\Big)\\[5pt]
    &=\Big(\pi\operatorname{Ad}((12)^{|\a_1\ldots\a_{i-1}|,1})\Phi(\a_1,\ldots,\a_n)\Big|X_1\otimes\ldots\otimes X_n\Big).
\end{align}
Hence 
\begin{align}
    \Phi(\a_1,\ldots,\a_{i-1},\pi(\a_i),\a_{i+1},\ldots,\a_n)-\pi\operatorname{Ad}((12)^{|\a_1\ldots\a_{i-1}|,1})\Phi(\a_1,\ldots,\a_n)\in\R_N(A).
\end{align}

Let us prove \eqref{f39}. In $\orep$ we have $(\widehat{1}(\a)|X)=(\a|\varepsilon_1(X))$ for homogeneous $\a\in\O(A)$ and $X\in\TMat$ such that $|\a|=|X|-1$. Hence in $\orep^{\otimes n}$ one has
\begin{align}
    &(\a_1|X_1)\otimes \ldots \otimes (\a_{i-1}|X_{i-1})\otimes (\widehat{1}(\a_i)|X_i)\otimes (\a_{i+1}|X_{i+1})\otimes\ldots\otimes (\a_n|X_n)\\[5pt]
    &\hspace{100pt}=(\a_1|X_1)\otimes \ldots \otimes (\a_{i-1}|X_{i-1})\otimes (\a_i|\varepsilon_1(X_i))\otimes (\a_{i+1}|X_{i+1})\otimes\ldots\otimes (\a_n|X_n)
\end{align}
for homogeneous elements $\a_1,\ldots,\a_n\in\O(A)$ and $X_1,\ldots,X_n\in\TMat$ such that for any $j=1,\ldots,n$, $j\neq i$ one has $|\a_j|=|X_j|$, and $|\a_i|=|X_i|-1$. Let us use the following Sweedler's notation $X_i=X_i'\otimes X_i''$, where $X_i'\in\Mat$. Then $\varepsilon_1(X_i)=\varepsilon(X_i')X_i''$.

Apply $\varphi_N$ to the both sides of the identity above and transform the right-hand side
\begin{align}
    \hspace{40pt}&\hspace{-40pt}\Big(\Phi(\a_1,\ldots,\a_{i-1},\widehat{1}(\a_i),\a_{i+1},\ldots,\a_n)\Big|X_1\otimes\ldots\otimes X_n\Big)\\[5pt]
    &=\Big(\Phi(\a_1,\ldots,\a_n)\Big|\varepsilon(X_i')X_1\otimes\ldots\otimes X_{i-1}\otimes X_i''\otimes X_{i+1}\otimes \ldots\otimes X_n\Big)\\[5pt]
    &=\Big(\Phi(\a_1,\ldots,\a_n)\Big|\varepsilon_1(X_i'\otimes X_1)\otimes\ldots\otimes X_{i-1}\otimes X_i''\otimes X_{i+1}\otimes \ldots\otimes X_n\Big)\\[5pt]
    &=\Big(\1\Phi(\a_1,\ldots,\a_n)\Big|X_i'\otimes X_1\otimes\ldots\otimes X_{i-1}\otimes X_i''\otimes X_{i+1}\otimes \ldots\otimes X_n\Big)\\[5pt]
    &=\Big(\1\Phi(\a_1,\ldots,\a_n)\Big|(12)^{|\a_1|+\ldots+|\a_{i-1}|,1} X_1\otimes\ldots\otimes X_{i-1}\otimes X_i'\otimes X_i''\otimes X_{i+1}\otimes \ldots\otimes X_n\Big)\\[5pt]
    &=\Big(\operatorname{Ad}((12)^{1,|\a_1|+\ldots+|\a_{i-1}|})\1\Phi(\a_1,\ldots,\a_n)\Big|X_i'\otimes X_1\otimes\ldots\otimes X_{i-1}\otimes X_i''\otimes X_{i+1}\otimes \ldots\otimes X_n\Big).
\end{align}
Hence 
\begin{align}
    \Phi(\a_1,\ldots,\a_{i-1},\widehat{1}(\a_i),\a_{i+1},\ldots,\a_n)-\operatorname{Ad}((12)^{1,|\a_1|+\ldots+|\a_{i-1}|})\1\Phi(\a_1,\ldots,\a_n)\in \R_N(A).
\end{align}

Let us prove \eqref{f42} and \eqref{f58}. In $\orep$ we have $(\sigma\cdot\a|X)=(\a|\covect(\sigma^{-1})X)$ and $(\a\cdot\sigma|X)=(\a|\vect(\sigma)X)$ for any homogeneous $\a\in\O(A)$ and $X\in\TMat$ such that $|\a|=|X|$ and any permutation $\sigma$ of rank $|\a|$. Hence in $\orep^{\otimes n}$ one has
\begin{align}
    &(\a_1|X_1)\otimes \ldots \otimes (\a_{i-1}|X_{i-1})\otimes (\sigma\cdot\a_i|X_i)\otimes (\a_{i+1}|X_{i+1})\otimes\ldots\otimes (\a_n|X_n)\\[5pt]
    &\hspace{60pt}=(\a_1|X_1)\otimes \ldots \otimes (\a_{i-1}|X_{i-1})\otimes (\a_i|\covect(\sigma^{-1})X_i)\otimes (\a_{i+1}|X_{i+1})\otimes\ldots\otimes (\a_n|X_n)\\[10pt]
    &(\a_1|X_1)\otimes \ldots \otimes (\a_{i-1}|X_{i-1})\otimes (\a_i\cdot\sigma|X_i)\otimes (\a_{i+1}|X_{i+1})\otimes\ldots\otimes (\a_n|X_n)\\[5pt]
    &\hspace{60pt}=(\a_1|X_1)\otimes \ldots \otimes (\a_{i-1}|X_{i-1})\otimes (\a_i|\vect(\sigma)X_i)\otimes (\a_{i+1}|X_{i+1})\otimes\ldots\otimes (\a_n|X_n)
\end{align}
for homogeneous elements $\a_1,\ldots,\a_n\in\O(A)$ and $X_1,\ldots,X_n\in\TMat$ such that for any $j=1,\ldots,n$ one has $|\a_j|=|X_j|$, and any permutation $\sigma$ of rank $|\a_i|$. 

Apply $\varphi_N$ to the both sides of these identities and transform the right-hand side. For the first identity we have
\begin{align}
    \hspace{40pt}&\hspace{-40pt}\Big(\Phi(\a_1,\ldots,\a_{i-1},\sigma\cdot\a_i,\a_{i+1},\ldots,\a_n)\Big|X_1\otimes\ldots\otimes X_n\Big)\\[5pt]
    &=\Big(\Phi(\a_1,\ldots,\a_n)\Big|X_1\otimes\ldots\otimes X_{i-1}\otimes \covect(\sigma^{-1})(X_i)\otimes X_{i+1}\otimes \ldots\otimes X_n\Big)\\[5pt]
    &=\Big(\Phi(\a_1,\ldots,\a_n)\Big|\covect(\operatorname{id}_{|\a_1\ldots\a_{i-1}|}\times \sigma^{-1}\times \operatorname{id}_{\a_{i+1}\ldots\a_n})\big(X_1\otimes\ldots\otimes X_n\big)\Big)\\[5pt]
    &=\Big(\big(\operatorname{id}_{|\a_1\ldots\a_{i-1}|}\times \sigma\times \operatorname{id}_{\a_{i+1}\ldots\a_n}\big)\cdot\Phi(\a_1,\ldots,\a_n)\Big|X_1\otimes\ldots\otimes X_n\Big).
\end{align}
Hence 
\begin{align}
    \Phi(\a_1,\ldots,\a_{i-1},\sigma\cdot\a_i,\a_{i+1},\ldots,\a_n)-\big(\operatorname{id}_{|\a_1\ldots\a_{i-1}|}\times \sigma\times \operatorname{id}_{\a_{i+1}\ldots\a_n}\big)\cdot\Phi(\a_1,\ldots,\a_n)\in\R_N(A),
\end{align}
which immediately implies \eqref{f42}. For the second identity we argue exactly the same and that leads us to \eqref{f58}.

Let's assume now that $\Phi$ satisfies \eqref{f25}-\eqref{f58} and let's prove that \eqref{f57} defines a linear map $\orep^{\otimes n}\longrightarrow\orep$. Let us denote by $\widetilde{\varphi_N}$ the map induced by $\Phi$
\begin{align}
    \widetilde{\varphi_N}:\Big[\O(A)\odot\TMat\Big]^{\otimes n}&\longrightarrow\O(A)\odot\TMat\\
    \big(\a_1\otimes X_1\big)\otimes\ldots\otimes \big(\a_n\otimes X_n\big)&\mapsto \Phi(\a_1,\ldots,\a_n)\otimes \big(X_1\otimes\ldots\otimes X_n\big)
\end{align}

Recall that by Proposition \ref{prop6} $\orep$ is the quotient of $\O(A)\odot\TMat$ by the ideal spanned by
\begin{enumerate}[label=\alph*)]
            \item\label{cond111} $\pi(\a)\otimes X-\a\otimes(\tr\otimes X)$ for all $\a\in\O(A)_k$ with $k\geq 1$ and all $X\in\Mat^{\otimes k}$;
            \vspace{7pt}

            \item\label{cond222}  $\widehat{1}(\a)\otimes X-\a\otimes\varepsilon_1(X)$ for all $\a\in\O(A)_k$ and all $X\in\Mat^{\otimes k+1}$;
            \vspace{7pt}

            \item\label{cond333} $(\sigma\cdot \a)\otimes E^*_{\mathbf{i}\mathbf{j}}-\a\otimes E^*_{\mathbf{i}\sigma^{-1}(\mathbf{j})}$ and $(\a\cdot\sigma)\otimes E^*_{\mathbf{i}\mathbf{j}}-\a\otimes E^*_{\sigma(\mathbf{i})\mathbf{j}}$ for all $\a\in\O(A)_k$, $\sigma\in S(k)$, and all tuples of indices $\mathbf{i}$,$\mathbf{j}$.
        \end{enumerate}

Let us denote the quotient map by
\begin{align}
    \theta_N:\O(A)\odot\TMat\longtwoheadrightarrow \orep.
\end{align}

Then, obviously, $\orep^{\otimes n}$ is the quotient of $\Big[\O(A)\odot\TMat\Big]^{\otimes n}$ by the ideal 
\begin{align}\label{f115}
    \operatorname{ker}(\theta_N^{\otimes n})=\sum\limits_{i=1}^n\Big[\O(A)\odot\TMat\Big]^{\otimes i-1}\otimes \operatorname{ker}(\theta_N)\otimes \Big[\O(A)\odot\TMat\Big]^{\otimes n-1-i}.
\end{align}

Then it's enough to check that $\widetilde{\varphi_N}\big(\operatorname{ker}(\theta_N^{\otimes n})\big)\subset \operatorname{ker}(\theta_N)$. Let's show that the image of the $i$-th summand in \eqref{f115} under $\widetilde{\varphi_N}$ belongs to $\operatorname{ker}(\theta_N)$. We will do this only in case when the $i$-th component is of the form \ref{cond111}; cases \ref{cond222} and \ref{cond333} can be done similarly. In fact, we have already done all the essential computations and now we just need to repeat the computation above but in the reverse order. 

So, we compute
\begin{align}
    \hspace{20pt}&\hspace{-20pt}\widetilde{\varphi_N}\Big((\a_1\otimes X_1)\otimes\ldots\otimes (\a_{i-1}\otimes X_{i-1})\otimes (\pi(\a_i)\otimes X_i)\otimes (\a_{i+1}\otimes X_{i+1})\otimes\ldots\otimes (\a_n\otimes X_n)\Big)\\[5pt]
    &=\Phi\big(\a_1,\ldots,\a_{i-1},\pi(\a_i),\a_{i+1},\ldots,\a_n\big)\otimes \big(X_1\otimes \ldots\otimes X_n\big)\\[5pt]
    &\equiv \Big(\pi\operatorname{Ad}((12)^{|\a_1\ldots\a_{i-1}|,1})\Phi(\a_1,\ldots,\a_n)\Big)\otimes \big(X_1\otimes\ldots\otimes X_n\big)\ \operatorname{mod}\ \ker(\theta_N)\tag{\text{by}\ \eqref{f25}}\\[5pt]
    &\equiv\Big(\operatorname{Ad}((12)^{|\a_1\ldots\a_{i-1}|,1})\Phi(\a_1,\ldots,\a_n)\Big)\otimes \big(\tr\otimes X_1\otimes\ldots\otimes X_n\big) \ \operatorname{mod}\ \ker(\theta_N)\tag{\text{by}\ \ref{cond111}\ \text{above}}\\[5pt]
    &\equiv\Big(\Phi(\a_1,\ldots,\a_n)\Big)\otimes \big((12)^{1,|\a_1\ldots\a_{i-1}|}\tr\otimes X_1\otimes\ldots\otimes X_n\big) \ \operatorname{mod}\ \ker(\theta_N)\tag{\text{by}\ \ref{cond333}\ \text{above}}\\[5pt]
    &\equiv\Big(\Phi(\a_1,\ldots,\a_n)\Big)\otimes \big(X_1\otimes\ldots\otimes X_{i-1}\otimes \tr\otimes X_i\otimes\ldots\otimes X_n\big) \ \operatorname{mod}\ \ker(\theta_N)
\end{align}
and see that the right-hand side equals 
\begin{align}
    \widetilde{\varphi_N}\Big((\a_1\otimes X_1)\otimes\ldots\otimes (\a_{i-1}\otimes X_{i-1})\otimes (\a_i\otimes (\tr\otimes X_i))\otimes (\a_{i+1}\otimes X_{i+1})\otimes\ldots\otimes (\a_n\otimes X_n)\Big) 
\end{align}
modulo $\ker(\theta_N)$ as desired.

Let us now prove that any $\GL$-equivariant linear map $\orep^{\otimes n}\longrightarrow\orep$ is of the form \eqref{f57} for some linear map $\Phi$. Then $\Phi$ would automatically satisfy \eqref{f25}-\eqref{f58} by the argument above. Note that the map
\begin{align}
    \theta^{\otimes n}_N:\Big[\O(A)\odot\TMat\Big]^{\otimes n}\longtwoheadrightarrow \orep^{\otimes n}    
\end{align}
is $\GL$-equivariant, so $\operatorname{ker}(\theta^{\otimes n}_N)$ is a $\GL$-submodule of $\Big[\O(A)\odot\TMat\Big]^{\otimes n}$. 

Next, $\Big[\O(A)\odot\TMat\Big]^{\otimes n}$ is a semisimple $\GL$-module due to the following $\GL$-module isomorphism
\begin{align}
    \Big[\O(A)\odot\TMat\Big]^{\otimes n}&=\Big[\bigoplus\limits_{p\geq 0}\O(A)^{(p)}\otimes \Mat^{\otimes p}\Big]^{\otimes n}\\[5pt]
    &\simeq \bigoplus\limits_{p_1\ldots,p_n\geq 0}\Big[\O(A)^{(p_1)}\otimes\ldots\otimes\O(A)^{(p_n)}\Big]\otimes \Mat^{\otimes p_1+\ldots+p_n}\\[5pt]
    &\simeq \O(A)^{\otimes n}\odot\TMat\label{f123}
\end{align}
and the fact that the latter is a direct sum, infinite of course, of finite dimensional $\GL$-modules, hence a sum of simple $\GL$-modules. So, $\operatorname{ker}(\theta^{\otimes n}_N)$ admits a complement $\GL$-module, say $P_n$, i.e. 
\begin{equation}\label{f109}
    \Big[\O(A)\odot\TMat\Big]^{\otimes n}\simeq\operatorname{ker}(\theta^{\otimes n}_N)\oplus P_n,
\end{equation}
see for instance Theorem in 6.2.1 in \cite{procesi2007lie}. Obviously, $P_n$ is isomorphic to $\orep^{\otimes n}$ as $\GL$-modules. Below we do not distinguish them. Let us denote by $\imath_n$ the inclusion of $\GL$-modules $\orep^{\otimes n}\simeq P_n\xhookrightarrow{\hspace{15pt}} \Big[\O(A)\odot\TMat\Big]^{\otimes n}$ provided by the right-hand side of \eqref{f109}. Recall that $\theta^{\otimes n}_N\circ\imath_n=\operatorname{id}$, which implies 
\begin{align}\label{f112}
    \imath_n\Big((\a_1|X_1)\otimes\ldots\otimes (\a_n|X_n)\Big)=(\a_1\otimes X_1)\otimes\ldots\otimes(\a_n\otimes X_n)\ \ \operatorname{mod}\ \operatorname{ker}(\theta^{\otimes n}_N)
\end{align}
for any homogeneous elements $\a_1,\ldots,\a_n\in\O(A)$, $X_1,\ldots,X_n\in\TMat$ such that $|\a_i|=|X_i|$.

We will also need similar objects corresponding to $n=1$, i.e. $P_1$ and $\imath_1$.

Next, for any $\GL$-equivariant linear map $\varphi_N:\orep^{\otimes n}\longrightarrow\orep$ we can define a linear map $\widetilde{\varphi_N}:\Big[\O(A)\odot\TMat\Big]^{\otimes n}\longrightarrow\O(A)\odot\TMat$ as the following composite
\begin{equation}
    \begin{tikzcd}
        \Big[\O(A)\odot\TMat\Big]^{\otimes n}\arrow[r,two heads, "\theta^{\otimes n}_N"] \arrow[rrr, bend right=10, "\widetilde{\varphi_N}"']& \orep^{\otimes n} \arrow[r,"\varphi_N"] & \orep \arrow[r,"\imath_1"] & \O(A)\odot\TMat
    \end{tikzcd}
\end{equation}

In other words, $\widetilde{\varphi_N}=\imath_1\circ\varphi_N\circ\theta^{\otimes n}_N$. What is most important to us now is that 
\begin{itemize}
    \item $\widetilde{\varphi_N}$ is $\GL$-equivariant
    \vspace{5pt}

    \item $\operatorname{ker}(\theta^{\otimes n}_N)\subset \operatorname{ker}(\widetilde{\varphi_N})$.
\end{itemize}

We will describe all such maps $\widetilde{\varphi_N}$ and then get back to $\varphi_N$ by using $\varphi_N=\theta_N\circ\widetilde{\varphi_N}\circ\imath_n$. Here we used identities $\theta^{\otimes n}_N\circ\imath_n=\operatorname{id}$ and $\theta_N\circ\imath_1=\operatorname{id}$. It might seem that we need to know the exact form of $\imath_n$, but that's unnecessary due to \eqref{f112} and the fact that $\widetilde{\varphi_N}$ vanishes on $\operatorname{ker}(\theta^{\otimes n}_N)$.

Let us pick a \textit{linear basis} in $\O(A)$ consisting of graded elements, say $\a_I$, for all $I$ belonging to an indexing set $\mathcal{I}$. This set is naturally graded by the degree of $\a_I$'s, which we denote simply by $|I|$. Due to \eqref{f123} we can view $\widetilde{\varphi_N}$ as an $\mathcal{I}\times\mathcal{I}^{n}$ matrix whose $(I,J)$-th entry for $I\in\mathcal{I},J\in\mathcal{I}^n$ belongs to 
\begin{align}
    \operatorname{Hom}_{\Bbbk}\left(\Mat^{\otimes |J|},\Mat^{\otimes |I|}\right)^{\GL},   
\end{align}
where $|J|$ stands for the total degree, namely for $J=(J_1,\ldots,J_n)\in\mathcal{I}^n$ the degree is the sum of degrees of the components $|J|:=|J_1|+\ldots+|J_n|$. 

Let us denote this matrix entry by $\eta_N^{I,J}$. Below we treat $\widetilde{\varphi_N}$ as a map from $\O(A)^{\otimes n}\odot\TMat$ to $\O(A)\odot\TMat$, so we can write
\begin{align}\label{f113}
    \widetilde{\varphi_N}(\a_J\otimes X)=\sum\limits_{I\in\mathcal{I}}\a_I\otimes \eta_N^{I,J}(X)
\end{align}
for any $J=(J_1,\ldots,J_n)\in\mathcal{I}^n$, $X\in\Mat^{\otimes |J|}$, and $\a_J=\a_{J_1}\otimes\ldots\otimes \a_{J_n}\in\O(A)^{\otimes n}$.

Now we need to describe the vector space of $\GL$-equivariant linear maps from $\Mat^{\otimes n}$ to $\Mat^{\otimes m}$ for any $n,m$. Let us describe a spanning set for it. Take any natural numbers $r$ and $s$ subject to $n+s-r=m$, and any permutations $\tau,w\in S(n+s)$. Let us define a linear maps $\nu^{(n)}_{w,\tau,r,s}$ as follows
    \begin{align}\label{f98}
        \nu^{(n)}_{w,\tau,r,s}:&\Mat^{\otimes n}\longrightarrow\Mat^{\otimes m},\\[5pt]
        &X\mapsto \varepsilon_1^r \operatorname{vect}(\tau)\operatorname{covect}(w)(\tr^{\otimes s}\otimes X).
    \end{align}
    
    Obviously, any such map is $\GL$-equivariant, because $\varepsilon$, $\operatorname{vect}(\tau)$, $\operatorname{covect}(w)$ are equivariant and $\tr\in\Mat$ is invariant. Let us check that any $\GL$-equivariant linear map is a linear combination of those. We have an obvious vector space isomorphism
    \begin{align}
    \operatorname{Hom}_{\Bbbk}\left(\operatorname{Mat}^*_N(\Bbbk)^{\otimes n},\operatorname{Mat}^*_N(\Bbbk)^{\otimes m}\right)^{\GL}\simeq\left(\operatorname{Mat}_N(\Bbbk)^{\otimes n+m}\right)^{\GL}.
\end{align}
By the first fundamental theorem of invariant theory the right-hand side equals the image of $\Bbbk[S(n+m)]$ in $\operatorname{Mat}_N(\Bbbk)^{\otimes n+m}\simeq \operatorname{End}_{\Bbbk}\Big((\Bbbk^N)^{\otimes n+m}\Big)$, where the symmetric group acts on the vector space $(\Bbbk^N)^{\otimes n+m}$ by permuting the tensor factors. Let us check that the image of a permutation $\sigma$ in $\operatorname{Hom}_{\Bbbk}\left(\operatorname{Mat}^*_N(\Bbbk)^{\otimes n},\operatorname{Mat}^*_N(\Bbbk)^{\otimes m}\right)^{\GL}$ is of the form $\nu^{(n)}_{w,\tau,r,s}$. This is a straightforward computation. 

Let $\sigma\in S(n+m)$, then its action on $\left(\Bbbk^N\right)^{\otimes n+m}$ is given by
\begin{equation}
    \sigma (v_1\otimes \ldots\otimes v_{n+m})=v_{\sigma^{-1}(1)}\otimes\ldots\otimes v_{\sigma^{-1}(n+m)}
\end{equation}
or in terms of matrix units
\begin{equation}
    \sigma=\sum\limits_{j_1,\ldots,j_{n+m}=1}^N E_{j_{\sigma^{-1}(1)} j_1}\otimes\ldots\otimes E_{j_{\sigma^{-1}(n+m)} j_{n+m}}\in \operatorname{Mat}_N(\Bbbk)^{\otimes n+m}.
\end{equation}
Then as an element of $\operatorname{Hom}_{\Bbbk}\left(\operatorname{Mat}^*_N(\Bbbk)^{\otimes n},\operatorname{Mat}^*_N(\Bbbk)^{\otimes m}\right)$ it is given by
\begin{align}
    \sigma=\sum\limits_{j_1,\ldots,j_{n+m}=1}^N E_{j_{\sigma^{-1}(1)} j_1}&\otimes\ldots\otimes E_{ j_{\sigma^{-1}(n)}j_{n}}\otimes E^*_{j_{n+1} j_{\sigma^{-1}(n+1)}}\otimes\ldots\otimes E^*_{j_{n+m}j_{\sigma^{-1}(n+m)}},
\end{align}
where $E_{ji}^*=\tr\left(E_{ij}\cdot-\right)\in\operatorname{Mat}^*_N(\Bbbk)$. 

Then the value of $\sigma$ on the element $E^*_{k_1l_1}\otimes\ldots\otimes E^*_{k_n l_n}$ of $\operatorname{Mat}^*_N(\Bbbk)^{\otimes n}$ is
\begin{align}\label{f85}
    \hspace{40pt}&\hspace{-40pt}\sigma\left(E^*_{k_1l_1}\otimes \ldots\otimes E^*_{k_n l_n}\right)\\
    &=\sum\limits_{j_1,\ldots,j_{n+m}=1}^N \delta\left(j_1=l_1,\ldots,j_n=l_n\right)\cdot\delta\left(j_{\sigma^{-1}(1)}=k_1,\ldots, j_{\sigma^{-1}(n)}=k_n\right)\\
    &\pushright{E^*_{j_{n+1} j_{\sigma^{-1}(n+1)}}\otimes\ldots\otimes E^*_{j_{n+m}j_{\sigma^{-1}(n+m)}}}\\
    &=\varepsilon^{2n}_1\Bigg(\sum\limits_{j_1,\ldots,j_{n+m}=1}^N E^*_{j_1l_1}\otimes\ldots\otimes E^*_{j_nl_n}\otimes E^*_{k_1j_{\sigma^{-1}(1)}}\otimes\ldots\otimes E^*_{k_nj_{\sigma^{-1}(n)}}\\
    &\pushright{E^*_{j_{n+1} j_{\sigma^{-1}(n+1)}}\otimes\ldots\otimes E^*_{j_{n+m}j_{\sigma^{-1}(n+m)}}}\Bigg).
\end{align}

Now it's clear that the right-hand side is of the form $\nu^{(n)}_{w,\tau,r,s}\left(E^*_{k_1l_1}\otimes \ldots\otimes E^*_{k_n l_n}\right)$ .

Let's come back to \eqref{f113} and consider a basis of $\operatorname{Hom}_{\Bbbk}\left(\Mat^{\otimes n},\Mat^{\otimes m}\right)^{\GL}$ consisting of linear maps of the form $\nu^{(n)}_{w,\tau,r,s}$. Let us denote this basis (as a set) by $\Omega_{n,m}$. Note that it depends on $N$ and from the Schur-Weyl duality it follows that for $N\geq n+m$ permutations from $S(n+m)$ provide such a basis while for $n+m>N$ the permutations become linear dependent. Let us also set $\Omega_{n,-}:=\bigsqcup\limits_{m\geq 0}\Omega_{n,m}$, which is a basis in $\operatorname{Hom}_{\Bbbk}\left(\Mat^{\otimes n},\TMat\right)^{\GL}$.

    So, we can write
    \begin{align}
        \eta_N^{I,J}=\sum\limits_{\nu\in \Omega_{|J|,|I|}} \xi^{I,J}(\nu)\nu
    \end{align}
    for some $\xi^{I,J}(\nu)\in\Bbbk$, which might depend on $N$ as well.
    
Now we substitute this expression in \eqref{f113}
    \begin{align}
        \widetilde{\varphi_N}(\a_J\otimes X)&=\sum\limits_{I\in\mathcal{I}} \ \sum\limits_{\nu\in \Omega_{|J|,|I|}}\xi^{I,J}(\nu)\cdot\a_I\otimes \nu(X)\\[5pt]
        &=\sum\limits_{m\geq 0}\ \ \sum\limits_{I\in\mathcal{I}:|I|=m} \ \sum\limits_{\nu\in \Omega_{|J|,m}}\xi^{I,J}(\nu)\cdot\a_I\otimes \nu(X)\tag{split the $1$-st sum}\\[5pt]
        &=\sum\limits_{m\geq 0} \ \sum\limits_{\nu\in \Omega_{|J|,m}}\ \ \sum\limits_{I\in\mathcal{I}:|I|=m}\xi^{I,J}(\nu)\cdot\a_I\otimes \nu(X)\tag{swap the $2$-nd and $3$-rd sums}%\\[5pt]
        %&=\sum\limits_{\nu\in \Omega_{|J|,-}}\ \ \sum\limits_{I\in\mathcal{I}:|I|=m}\xi^{I,J}(\nu)\cdot\a_I\otimes \nu(X_J)\tag{unite the $1$-st and $2$-nd sums}.
    \end{align}

    For any $\nu\in\Omega_{k,m}$ we define a linear map $\Phi_{\nu}:\O(A)^{\otimes n}\longrightarrow \O(A)$ by its values on the basis $\{a_{J}\}_{J\in\mathcal{I}}$
    \begin{align}
        \Phi_{\nu}(\a_J):=\begin{cases}
            \sum\limits_{I\in\mathcal{I}:|I|=m}\xi^{I,J}(\nu) \cdot\a_I,\ &\text{if}\ |J|=k,\\[5pt]
            0,\ &\text{otherwise}.
        \end{cases}
    \end{align}
    Note that for any particular $\a\in\O(A)^{\otimes n}$ of degree $k$ the value $\Phi_{\nu}(\a)$ is zero for any $\nu\in\Omega_{k,m}$ with $m$ large enough, because for any $J$ there exist only finitely many $I$ such that $\eta^{I,J}_N\neq 0$.

    Next, for $\widetilde{\varphi_N}$ and any homogeneous $\a\in\O(A)^{\otimes n}$, $X\in\TMat$ both of degree $k$ we have
    \begin{align}
        \widetilde{\varphi_N}(\a\otimes X)=\sum\limits_{\nu\in\Omega_{k,-}}\Phi_{\nu}(\a)\otimes \nu(X).
    \end{align}
Note that on the right-hand side we have an infinite sum, but for any fixed $\a$ this sum is actually finite, because in this case $\Phi_{\nu}(\a)$ is non-zero only for a finite number of $\nu$.

Thus, we can write
\begin{align}
        \varphi_N(\a|X)=\sum\limits_{\nu\in\Omega_{k,-}}\Big(\Phi_{\nu}(\a)\Big|\nu(X)\Big).
    \end{align}

    By \ref{cond14}, \ref{cond15}, \ref{cond10} from Lemma \ref{lemma5} we have
    \begin{align}
    \varphi_N(\a|X)=\Bigl(\Phi(\a)\, \Big|\, X\Bigr),
    \end{align}
where
    \begin{align}
        \Phi(\a)=\sum\limits_{\nu=\nu^{(n)}_{w,\tau,r,s}\in\Omega_{k,-}} \pi^{s} \left(w^{-1}\cdot \widehat{1}^{\, r}\Big(\Phi_{\nu}(\a)\Big)\cdot \tau\right)
    \end{align}
    and the proof is complete.
\end{proof}

Equivariant linear maps will not actually suffice. We will also have to deal with the followings linear maps
\begin{align}
    &\Hom_{\Bbbk}\Big(\orep^{\otimes n},\orep\Big)^{(p)}:=\operatorname{span}_{\Bbbk}\Bigg((\mathds{1}\otimes f)\circ\varphi\ \Bigg|\ f\in\Mat^{\otimes p},\ \\
    &\hspace{180pt}\varphi\in\Hom_{\Bbbk}\Big(\orep^{\otimes n},\orep\otimes \Matt^{\otimes p}\Big)^{\GL}\Bigg)
\end{align}
for various $p$. 

A typical element of this vector space is a linear map $\varphi_Y\in \Hom_{\Bbbk}\Big(\orep^{\otimes n},\orep\Big)^{\GL}$ indexed by $Y\in \Matt^{\otimes p}$ and such that $g.\varphi_Y=\varphi_{g.Y}$ for any $g\in\GL$.

\begin{theorem}\label{th2}
    A homogeneous linear map $\Phi:\O(A)^{\otimes n}\longrightarrow \O(A)[p]$ gives rise to a well-defined $\GL$-equivariant linear map $\varphi_N:\orep^{\otimes n}\longrightarrow \orep\otimes \Matt^{\otimes p}$ by the rule
    \begin{align}
        \varphi_N\big((\a_1|X_1),\ldots,(\a_n|X_n)\big)=\sum\limits_{\mathbf{i},\mathbf{j}\in[1,N]^p}\Big(\Phi(\a_1,\ldots,\a_n)\Big|E^*_{\mathbf{i}\mathbf{j}}\otimes X_1\otimes\ldots\otimes X_n\Big)\otimes E_{\mathbf{i}\mathbf{j}},
    \end{align}
    where $\a_1,\ldots,\a_n\in\O(A)$ and $X_1,\ldots,X_n\in\TMat$ are homogeneous elements such that $|\a_i|=|X_i|$, if and only if for any homogeneous elements $\a_1,\ldots,\a_n\in\O(A)$, any permutations $w_1\in S(|\a_1|),\ldots,w_n\in S(|\a_n|)$, and any $k$ the map $\Phi$ satisfies 
    \begin{align}
        &\hspace{-20pt}\Phi\big(\a_1,\ldots,\a_{k-1},\pi(\a_k),\a_{k+1},\ldots,\a_n\big)\equiv\pi\operatorname{Ad}((12)^{p+|\a_1|+\ldots+|\a_{k-1}|,1})\Phi(\a_1,\ldots,\a_n)\ 
        \operatorname{mod}\ \R_N(A),\label{f325}\\[10pt]
        &\hspace{-20pt}\Phi(\a_1,\ldots,\a_{k-1},\widehat{1}(\a_k),\a_{k+1},\ldots,\a_n)\equiv\operatorname{Ad}((12)^{1,p+|\a_1|+\ldots+|\a_{k-1}|})\widehat{1}\Phi(\a_1,\ldots,\a_n)\ \operatorname{mod}\ \R_N(A),\label{f339}\\[10pt]
        &\hspace{-20pt}\Phi\big(w_1\cdot \a_1, \ldots,w_n\cdot \a_n\big)\equiv(\id_p\times w_1\times\ldots\times w_n)\cdot\Phi(\a_1,\ldots,\a_n)\ 
        \operatorname{mod}\ \R_N(A),\label{f500}\\[5pt]
        &\hspace{-20pt}\Phi(\a_1\cdot w_1,\ldots,\a_n\cdot w_n)\equiv\Phi(\a_1\ldots,\a_n)\cdot \big(\id_p\times w_1\times\ldots\times w_n\big)\ 
        \operatorname{mod}\ \R_N(A).\label{f501}
    \end{align}

Moreover, any $\GL$-equivariant linear map $\varphi_N:\orep^{\otimes n}\longrightarrow \orep\otimes \Matt^{\otimes p}$ is of this form.
\end{theorem}
\begin{proof}%See also Image36
    The proof is similar to that of Theorem \ref{th3}, but we need to use the following description of $\GL$-equivariant linear maps from $\Mat^{\otimes n}$ to $\Mat^{\otimes m}\otimes \Matt^{\otimes p}$. They are spanned by the maps of the form 
    \begin{align}
        X\mapsto \sum\limits_{\mathbf{i},\mathbf{j}\in[1,N]^p}\varepsilon_1^r\vect(\tau)\covect(w)\left(\tr^{\otimes s}\otimes E^*_{\mathbf{i}\mathbf{j}}\otimes X\right)\otimes E_{\mathbf{i}\mathbf{j}}.
    \end{align}
    
    Let's check this. By the first fundamental theorem of invariant theory the equivariant maps we are interested in are spanned by the permutations $\sigma\in S(n+m+p)$. Similarly to \eqref{f85} we identify permutation $\sigma$ with its image in $\Hom_{\Bbbk}\Big(\Mat^{\otimes n},\Mat^{\otimes m}\otimes \Matt^{\otimes p}\Big)^{\GL}$ and write
    \begin{align}
    \hspace{20pt}&\hspace{-20pt}\sigma\left(E^*_{k_1l_1}\otimes \ldots\otimes E^*_{k_n l_n}\right)\\
    &=\varepsilon^{2n}_1\sum\limits_{j_1,\ldots,j_{n+m+p}=1}^N E^*_{j_1l_1}\otimes\ldots\otimes E^*_{j_nl_n}\otimes E^*_{k_1j_{\sigma^{-1}(1)}}\otimes\ldots\otimes E^*_{k_nj_{\sigma^{-1}(n)}}\\
    &\hspace{40pt}\otimes E^*_{j_{n+1} j_{\sigma^{-1}(n+1)}}\otimes\ldots\otimes E^*_{j_{n+m}j_{\sigma^{-1}(n+m)}}\otimes E_{j_{\sigma^{-1}(n+m+1)}j_{n+m+1}}\otimes\ldots\otimes E_{j_{\sigma^{-1}(n+m+p)}j_{n+m+p}}\\[5pt]
    &=\sum\limits_{\mathbf{i}\in[1,N]^p}\varepsilon_1^{2n}\sum\limits_{j_1,\ldots,j_{n+m+p}=1}^N\delta(i_1=j_{\sigma^{-1}(n+m+1)},\ldots,i_p=j_{\sigma^{-1}(n+m+p)})\\
    &\pushright{E^*_{j_1l_1}\otimes\ldots\otimes E^*_{j_nl_n}\otimes E^*_{k_1j_{\sigma^{-1}(1)}}\otimes\ldots\otimes E^*_{k_nj_{\sigma^{-1}(n)}}\otimes E^*_{j_{n+1} j_{\sigma^{-1}(n+1)}}\otimes\ldots\otimes E^*_{j_{n+m}j_{\sigma^{-1}(n+m)}}}\\
    &\pushright{\otimes E_{i_1j_{n+m+1}}\otimes\ldots\otimes E_{i_p j_{n+m+p}}}\\[5pt]
    &=\sum\limits_{\mathbf{i}\in[1,N]^p}\varepsilon_1^{2n+p}\sum\limits_{j_1,\ldots,j_{n+m+p}=1}^N E^*_{i_1j_{\sigma^{-1}(n+m+1)}}\otimes\ldots\otimes E^*_{i_pj_{\sigma^{-1}(n+m+p)}}\\
    &\pushright{\otimes E^*_{j_1l_1}\otimes\ldots\otimes E^*_{j_nl_n}\otimes E^*_{k_1j_{\sigma^{-1}(1)}}\otimes\ldots\otimes E^*_{k_nj_{\sigma^{-1}(n)}}\otimes E^*_{j_{n+1} j_{\sigma^{-1}(n+1)}}\otimes\ldots\otimes E^*_{j_{n+m}j_{\sigma^{-1}(n+m)}}}\\
    &\pushright{\otimes E_{i_1j_{n+m+1}}\otimes\ldots\otimes E_{i_p j_{n+m+p}}}
\end{align}
and the right-hand side is of the desired form.
\end{proof}

\subsection{Admissible maps}
    Corollary \ref{cor3} shows that we should work with $\ored(A)=\faktor{\O(A)}{\R(A)}$ instead of $\O(A)$, and we will do so.

    \begin{remark}\label{rem5}
        If a homogeneous linear map $\Phi:\O(A)^{\otimes n}\longrightarrow\O(A)[p]$ is as in Theorem \ref{th2}, i.e. $\Phi$ satisfies \eqref{f325},\eqref{f339},\eqref{f500},\eqref{f501}. Then $\Phi(\a_1,\ldots,\a_n)\in\R_N(A)$, when at least one of $\a$'s belongs to $\R_N(A)$. Indeed, the value of the induced linear map $\varphi_N:\orep^{\otimes n}\longrightarrow \orep\otimes \Matt^{\otimes p}$ is zero when at least one of the arguments equals $(\a|X)=0$ for $\a\in\R_N(A)$. Hence any such $\Phi$ descends to a diagonal $\S$-bimodule homomorphism $\O_N(A)^{\stimes n}\longrightarrow \O_N(A)[p]$ for which \eqref{f325},\eqref{f339} hold precisely, where $\O_N(A):=\faktor{\O(A)}{\R_N(A)}$.
    \end{remark}

    Let us temporarily denote by $\O$ any of the three algebras $\O(A)$, $\ored(A)$, or $\O_N(A)$.

\begin{definition}\label{def6}
     We say that a morphism of diagonal $\S$-bimodules $\Phi:\O^{\stimes n}\longrightarrow \O[p]$ is \textit{$p$-admissible}, if it satisfies \eqref{f325} and \eqref{f339}, i.e. for any homogeneous elements $\a_1,\ldots,\a_n\in\O$ and any $i$ one has 
    \begin{align}
        &\Phi\big(\a_1,\ldots,\a_{i-1},\pi(\a_i),\a_{i+1},\ldots,\a_n\big)=\pi\operatorname{Ad}((12)^{p+|\a_1|+\ldots+|\a_{i-1}|,1})\Phi(\a_1,\ldots,\a_n)\\
        &\Phi(\a_1,\ldots,\a_{i-1},\widehat{1}(\a_i),\a_{i+1},\ldots,\a_n)=\operatorname{Ad}((12)^{1,p+|\a_1|+\ldots+|\a_{i-1}|})\1\Phi(\a_1,\ldots,\a_n).
    \end{align}
    The set of $p$-admissible diagonal $\S$-bimodule homomorphisms $\O^{\stimes n}\longrightarrow \O[p]$ will be denoted by $\Hompadm{p}\big(\O^{\stimes n},\O\big)$. When $p=0$ we will simply say \textit{admissible} instead of $0$-admissible and write $\Homadm$ instead of $\Hompadm{0}$.
\end{definition}

For instance, in the case $n=1$, a diagonal $\S$-bimodule homomorphism $\Phi:\O\longrightarrow \O$ is admissible iff it commutes with $\pi$ and $\widehat{1}$.

And in the case $n=2$, a diagonal $\S$-bimodule homomorphism $\Phi:\O^{\stimes 2}\longrightarrow \O$ is admissible iff for any homogeneous $\alpha,\beta\in\O$ one has
\begin{align}
        \Phi(\pi(\a),\b)&=\pi\Phi(\a,\b),& \Phi(\a,\pi(\b))&=\pi\operatorname{Ad}((12)^{|\a|,1})\Phi(\a,\b),\\[5pt]
        \Phi(\widehat{1}(\a),\b)&=\widehat{1}\Phi(\a,\b),& \Phi(\a,\widehat{1}(\b))&= \operatorname{Ad}((12)^{1,|\a|})\widehat{1}\Phi(\a,\b).
    \end{align}
    
For instance, multiplications in $\O(A)$, $\ored(A)$, and $\O_N(A)$ are admissible.

\begin{lemma}\label{lemma19}
    Let $\O$ denote either $\ored(A)$ or $\O_N(A)$. 
    \begin{enumerate}[label=\theenumi),leftmargin=3ex]
        \item The natural action of $S(n)$ on diagonal $\S$-bimodule homomorphisms $\Hom_{\S}\big(\O^{\stimes n},\O[p]\big)$ by $(\sigma,\Phi)\mapsto \sigma.\Phi$, where 
        \begin{align}
        \sigma.\Phi(\a_1,\ldots,\a_n):=\Phi\big(\sigma^{-1}(\a_1\otimes\ldots\otimes\a_n)\big)=\operatorname{Ad}(\id_p\times \sigma^{|\a_{\sigma(1)}|,\ldots,|\a_{\sigma(n)}|})\Phi(\a_{\sigma(1)},\ldots,\a_{\sigma(n)}).
        \end{align}
        preserves $p$-admissible maps.

        \item Let $\Phi\in \Hompadm{p}\big(\O^{\stimes n},\O\big)$ and $\Psi\in \Hompadm{q}\big(\O^{\stimes m},\O\big)$. Then their $k$-th composition $\bullet_k$ defined by 
        \begin{align}
            &\Phi\bullet_k \Psi(\a_1,\ldots,\a_{n+m})\\
            &\hspace{60pt}:=\operatorname{Ad}(\id_{p}\times (12)^{|\a_1|+\ldots+|\a_{k-1}|,q})\Phi(\a_1,\ldots,\a_{k-1},\Psi(\a_k,\ldots,\a_{k+m-1}),\a_{k+m+1},\ldots,\a_{n+m})
        \end{align}
    for $k=1,\ldots,n$, and also by
        \begin{align}
            &\Phi\bullet_{n+1} \Psi(\a_1,\ldots,\a_{n+m}):=\operatorname{Ad}(\id_{p}\times (12)^{|\a_1|+\ldots+|\a_{n}|,q})\Phi(\a_1,\ldots,\a_{n})\Psi(\a_{n+1},\ldots,\a_{n+m})
        \end{align}
    for $k=n+1$, is $p+q$-admissible.

    \item For any homogeneous $\a\in\O$ the operator $L_{\a}:\O\longrightarrow\O$ of left multiplication by $\a$ is $|\a|$-admissible.

    \item For any $\Phi\in \Hompadm{p}\big(\O^{\stimes n},\O\big)$ and $w,\tau\in S(p)$ the map $w\cdot \Phi\cdot \tau:\O^{\stimes n}\longrightarrow \O(A)[p]$ defined by 
    \begin{align}
      w\cdot \Phi\cdot \tau(\a_1,\ldots,\a_n):=w\cdot \Phi(\a_1,\ldots,\a_n)\cdot \tau  
    \end{align}
    is $p$-admissible as well.
    \end{enumerate}
\end{lemma}
\begin{proof}
\phantom{a}
\begin{enumerate}[label=\theenumi),leftmargin=1ex]
    \item Let us denote by $\Phi_N$ the map induced by $\Phi$, i.e. $\Phi_N:\orep^{\otimes n}\longrightarrow\orep\otimes \Matt^{\otimes p}$. Then for any $\sigma\in S(n)$ one has $\sigma.\big(\Phi_N\big)=(\sigma.\Phi)_N$, hence by Theorem \ref{th2} (strictly speaking, by its obvious analog for $\O_N(A)$ or $\ored(A)$ instead of $\O(A)$) the map $\sigma.\Phi$ is $p$-admissible.

    \item One can consider the $k$-th composition of $\Hom_{\Bbbk}\Big(\orep^{\otimes n},\orep\otimes\Matt^{\otimes p}\Big)$ and $\Hom_{\Bbbk}\Big(\orep^{\otimes n},\orep\otimes\Matt^{\otimes q}\Big)$, which lands in $\Hom_{\Bbbk}\Big(\orep^{\otimes n},\orep\otimes\Matt^{\otimes p+q}\Big)$. Let's denote it by the same symbol $\bullet_k$. Then one clearly has $\Phi_N\bullet_k\Psi_N=(\Phi\bullet_k\Psi)_N$ and and we apply Theorem \ref{th2}.

    \item Consider the operator $L:\orep\longrightarrow\orep\otimes \Matt^{\otimes |\a|}$ given by 
    \begin{align}
        L((\b|X)):=\sum\limits_{\mathbf{i},\mathbf{j}\in[1,N]^{|\a|}}(\a|E^*_{\mathbf{i}\mathbf{j}})(\b|X)\otimes E_{\mathbf{i}\mathbf{j}}=\sum\limits_{\mathbf{i},\mathbf{j}\in[1,N]^{|\a|}}(\a\b|E^*_{\mathbf{i}\mathbf{j}}\otimes X)\otimes E_{\mathbf{i}\mathbf{j}}.
    \end{align}
    It is of the same form as in Theorem \ref{th2} for $\Phi(\b)=\a\b$, hence this $\Phi$ is $|\a|$-admissible.

    \item Let $\Phi_N:\orep^{\otimes n}\longrightarrow\orep\otimes \Matt^{\otimes p}$ the induced map. Then we consider the map $\Big(\mathds{1}_{\orep}\otimes \vect(w)\covect(\tau^{-1})\Big)\Phi_N:\orep^{\otimes n}\longrightarrow\orep\otimes \Matt^{\otimes p}$. Note that for $\Mat$ and $\Matt$ symbols $\vect$ and $\covect$ are switched. One has
    \begin{align}
        &\Big(\mathds{1}_{\orep}\otimes \vect(w)\covect(\tau^{-1})\Big)\Phi_N\big((\a_1|X_1),\ldots,(\a_n|X_n)\big)\\
        &\hspace{100pt}=\sum\limits_{\mathbf{i},\mathbf{j}\in[1,N]^p}\Big(\Phi(\a_1,\ldots,\a_n)\Big|E^*_{\mathbf{i}\mathbf{j}}\otimes X_1\otimes\ldots\otimes X_n\Big)\otimes E_{\tau^{-1}(\mathbf{i})w(\mathbf{j})}\\
        &\hspace{100pt}=\sum\limits_{\mathbf{i},\mathbf{j}\in[1,N]^p}\Big(\Phi(\a_1,\ldots,\a_n)\Big|\vect(\tau)\covect(w^{-1})(E^*_{\mathbf{i}\mathbf{j}})\otimes X_1\otimes\ldots\otimes X_n\Big)\otimes  E_{\mathbf{i}\mathbf{j}}\\
        &\hspace{100pt}=\sum\limits_{\mathbf{i},\mathbf{j}\in[1,N]^p}\Big(w\cdot\Phi(\a_1,\ldots,\a_n)\cdot \tau\Big|E^*_{\mathbf{i}\mathbf{j}}\otimes X_1\otimes\ldots\otimes X_n\Big)\otimes  E_{\mathbf{i}\mathbf{j}}.
    \end{align}
    We apply Theorem \ref{th2} once again to complete the proof.
\end{enumerate}
\end{proof}

\subsection{Double algebras over an operad}

Consider the usual endomorphism operad $\mathcal{E} nd_{\orep}$
\begin{align}
    \mathcal{E}nd_{\orep}=\Bigg(\operatorname{Hom}_{\Bbbk}\Big(\orep^{\otimes n},\orep\Big)\Bigg)_{n\geq 0}
\end{align}
and its suboperad of $\GL$-invariants
\begin{align}
    \mathcal{E}nd_{\orep}^{\operatorname{GL}_N}=\Bigg(\operatorname{Hom}_{\Bbbk}\Big(\orep^{\otimes n},\orep\Big)^{\GL}\Bigg)_{n\geq 0}
\end{align}

Below we denote by $\Com$ the operad of commutative unital algebras. The natural multiplication in $\orep$ gives a morphism of operads $\mu_N:\Com\longrightarrow \mathcal{E}nd_{\orep}^{\operatorname{GL}_N}$.

As we promised in the beginning of the present section, we focus on the case when the structure $\mathcal{S}$, double analog of which we wish to define, is given by a $\mathcal{P}$-algebra structure for an operad $\mathcal{P}$. Note that the multiplication in $\orep$ is always fixed and is never considered as a part of the additional structure $\mathcal{S}$. So, we should restrict ourselves to operads $\mathcal{P}$ containing $\operatorname{Com}$ and such that the morphism $p_N:\mathcal{P}\longrightarrow \mathcal{E}nd_{\orep}^{\operatorname{GL}_N}$, which is a part of the additional structure $\mathcal{S}$, restricts to $\mu_N$.

Now we wish to lift this ''commutative'' structure $\mathcal{S}=\big(\mathcal{P},\{p_N\}_{N\geq 1}\big)$ to the noncommutative world. Kontsevich-Rosenberg principle states that this lift should be a certain structure defined on an object related to $A$ and such that the representation functor maps this structure to $\mathcal{S}$. We will use the equivariant version of the Kontsevich-Rosenberg principle stated in the beginning of the present section, and assume that we are looking for an operad $\mathcalbb{E}_A$, which depends on the algebra $A$ but not on $N$, together with operad homomorphisms $\mathbbnew{e}_N:\mathcalbb{E}_A\longrightarrow \mathcal{E}nd_{\orep}^{\operatorname{GL}_N}$, $N\geq 1$, such that for any operad $\mathcal{P}$ and any sequence of operad homomorphisms $p_N:\mathcal{P}\longrightarrow \mathcal{E}nd_{\orep}^{\operatorname{GL}_N}$, $N\geq 1$ as above there is a unique operad homomorphism $\mathbbnew{p}$ such that the following diagram commutes for any $N$

\begin{equation}
    \begin{tikzcd}
        &  \mathcalbb{E}_A\arrow[dr,"\mathbbnew{e}_N"]& \\
        \mathcal{P} \arrow[rr,"p_N"']\arrow[ur,dashed,"\exists!\ \mathbbp"]& & \mathcal{E}nd_{\orep}^{\operatorname{GL}_N}.
    \end{tikzcd}
\end{equation}

Then we could treat $\mathcalbb{E}_A$ as the \textit{double endomorphism operad} and consider such homomorphisms $\mathbbnew{p}:\mathcal{P}\longrightarrow\mathcalbb{E}_A$ as \textit{double $\mathcal{P}$-algebra structures on $A$}. Such an operad $\mathcalbb{E}_A$ exists and is unique up to a unique isomorphism making all the diagrams commute. It is simply the direct product %the super standard proof is in Image35
\begin{align}
    \mathcalbb{E}_A=\prod\limits_{N\geq 1}\mathcal{E}nd_{\orep}^{\operatorname{GL}_N}.
\end{align}
The morphisms $\mathbbnew{e}_N:\mathcalbb{E}_A\longrightarrow \mathcal{E}nd_{\orep}^{\operatorname{GL}_N}$ are the natural projections, and for any sequence $p_N:\mathcal{P}\longrightarrow \mathcal{E}nd_{\orep}^{\operatorname{GL}_N}$, $N\geq 1$, the morphism $\mathbbnew{p}:\mathcal{P}\longrightarrow \mathcalbb{E}_A$ is the product $\mathbbnew{p}=(p_1,p_2,\ldots)$.

The main disadvantage of $\mathcalbb{E}_A$ is that it's very hard to work with and also the map $\mathbbnew{e}_N$, which is supposed to be induced by the representation functor according to the Kontsevich-Rosenberg principle, has nothing to do with the representation functor -- $\mathbbnew{e}_N$ is merely the projection to a direct factor. We intend to fix this by considering a suboperad of $\mathcalbb{E}_A$ for which the composition of the inclusion to $\mathcalbb{E}_A$ and $\mathbbnew{e}_N$ is naturally induced by the representation functor.

\begin{definition}
    By the \textit{double endomorphism operad} $\End_A$ we mean the following suboperad in the usual endomorphism operad $\mathcal{E}nd_{\ored(A)}$
    \begin{align}
        \End_A^{\operatorname{red}}:=\Big(\Homadm\left(\ored(A)^{\stimes n},\ored(A)\right)\Big)_{n\geq 0}.
    \end{align}
\end{definition}

We can also consider the following truncated version of it 
\begin{align}
    \End_A^N:=\Big(\Homadm\left(\O_N(A)^{\stimes n},\O_N(A)\right)\Big)_{n\geq 0}.
\end{align}

\begin{proposition}\label{prop19}
    \begin{enumerate}[label=\theenumi),leftmargin=3ex]
        \item $\End_A^{\operatorname{red}}$ and $\End_A^{N}$ are indeed operads.

        \item There are homomorphisms of operads
        \begin{equation}
            \begin{tikzcd}
                \End_A^{\operatorname{red}}\arrow[r] & \End_A^{N} \arrow[r,"\sim"] & \mathcal{E}nd_{\orep}^{\operatorname{GL}_N}
            \end{tikzcd}
        \end{equation}
        where the first arrow is induced by the canonical homomorphism $\ored(A)\longtwoheadrightarrow\O_N(A)$, and the second one is induced by the representation functor and is an isomorphism.
    \end{enumerate}
\end{proposition}
\begin{proof}
    The first item follows from Lemma \ref{lemma19}. The second item follows from Remark \ref{rem5} and Theorem \ref{th3}. More concretely, Theorem \ref{th3} provides an operad homomorphism $\mathcal{E}nd_{\orep}^{\operatorname{GL}_N}\longrightarrow \End_A^{N}$, which is clearly the inverse of the map induced by the representation functor.
\end{proof}

\begin{remark}\label{rem4}
    The second item of Proposition \ref{prop19} implies that
    \begin{align}
        \mathcalbb{E}_A\simeq\prod\limits_{N\geq 1} \End_A^{N}
    \end{align}
and we obtain an injective map $\End_A^{\operatorname{red}}\longhookrightarrow \mathcalbb{E}_A$ induced by the maps $\End_A^{\operatorname{red}}\longrightarrow \End_A^{N}$ from the second item of the proposition above.
\end{remark}

Note that the multiplication in $\orep$ lifts to $\mu:\Com\longrightarrow\Endred_A$ making the following diagram commute
\begin{equation}\label{f84}
    \begin{tikzcd}
        & \Endred_A \arrow[dr]& \\
        \Com \arrow[ur,dashed, "\mu"]\arrow[rr,"\mu_N"']& & \mathcal{E}nd_{\orep}^{\operatorname{GL}_N},
    \end{tikzcd}
\end{equation}
where the right diagonal arrow is given by the second item in Proposition \ref{prop19}.

We were lucky to lift the multiplication in $\orep$ to $\Endred_A$. For an arbitrary operad $\mathcal{P}$ and operad homomorphisms $p_N$ as in the beginning of the current subsection the similar result is slightly weaker. Using the isomorphism of operads from the second item of Proposition \ref{prop19}, we can uniquely lift $p_N$ to $\mathbbp_N:\mathcal{P}\longrightarrow\End_A^{N}$ in such a way that 
\begin{itemize}
    \item $\eval{\mathbbp_N}_{\Com}$ is the composition of $\mu$ defined in \eqref{f84} and the morphism $\Endred_A\longrightarrow\End_A^N$ from the second item of Proposition \ref{prop19};

    \item the following diagram commutes 
\begin{equation}
    \begin{tikzcd}
        & \End_A^N \arrow[dr]& \\
        \mathcal{P} \arrow[ur,dashed,"\exists!\ \mathbbp_N"]\arrow[rr,"p_N"]& & \mathcal{E}nd_{\orep}^{\operatorname{GL}_N}.
    \end{tikzcd}
\end{equation}
\end{itemize} 

\begin{definition}[Double algebra over an operad]\label{def8}
Let $\mathcal{P}$ be an operad containing $\Com$. By a \textit{double algebra over $\mathcal{P}$}, or a \textit{double $\mathcal{P}$-algebra}, we mean an associative algebra $A$ equipped with an operad homomorphism $\mathbbnew{p}:\mathcal{P}\longrightarrow \Endred_A$ such that the restriction of $\mathbbnew{p}$ to $\Com$ equals $\mu$ defined in \eqref{f84}.
\end{definition}

Note that double $\mathcal{P}$-algebras clearly satisfy the equivariant Kontsevich-Rosenberg principle. If $A$ is a $\mathcal{P}$-algebra, than the homomorphism $p_N:\mathcal{P}\longrightarrow \mathcal{E}nd_{\orep}^{\operatorname{GL}_N}$ is given by the composition of $\mathbbnew{p}$ and the homomorphism $\Endred_A\longrightarrow \mathcal{E}nd_{\orep}^{\operatorname{GL}_N}$ from Proposition \ref{prop19}, and also the restriction of $p_N$ to $\Com$ equals $\mu_N$ as desired.

\begin{remark}\label{rem7}
    Note that Definition \ref{def8} is maximal in the sense of Theorem \ref{th3} and Corollary \ref{cor3}. Namely, if one wishes to have a noncommutative structure on $A$ that induces a "universal" (= "stable" in $N$) $\GL$-equivariant $\mathcal{P}$-algebra structure on all $\rep$, then one has to deal with $\Endred_A$. So, one can regard $\Endred_A$ as a "universal" or "stable" part of $\mathcalbb{E}_A=\prod\limits_{N\geq 1}\mathcal{E}nd_{\orep}^{\operatorname{GL}_N}$.
\end{remark}

\begin{remark}\label{rem6}
    Definition \ref{def8} provides an answer to Question 36 from \cite{damas2024transposed}. 
\end{remark}

In the next section we consider the example $\mathcal{P}=\operatorname{Pois}$. 

\section{Di-twisted Poisson brackets and their deformations}\label{section_equivariant_linear_maps}

\begin{proposition}
A homogeneous linear map $\{-,-\}:\O(A)\otimes \O(A)\longrightarrow \O(A)$ of degree zero gives rise to a well-defined $\GL$-equivariant Poisson bracket $\{-,-\}_N:\orep\otimes \orep\longrightarrow \orep$ by the rule
    \begin{align}
        \big\{(\a|X),(\b|Y)\big\}_N=\Big(\{\a,\b\}\Big|X\otimes Y\Big),
    \end{align}
    where $\a,\b\in\O(A)$ and $X,Y\in\TMat$ are homogeneous elements such that $|\a|=|X|$ and $|\b|=|Y|$, if and only if for any graded elements $\a,\b\in\O(A)$ and any permutations $w_1\in S(|\a|),w_2\in S(|\b|)$ the map $\{-,-\}$ satisfies 
    \begin{align}
        \{\pi(\a),\b\}&\equiv\pi\{\a,\b\}\ 
        \operatorname{mod}\ \R_N(A),\\[5pt]
        \{\a,\pi(\b)\}&\equiv\pi\operatorname{Ad}((12)^{|\a|,1})\{\a,\b\}\ 
        \operatorname{mod}\ \R_N(A),\\[5pt]
        \{\widehat{1}(\a),\b\}&\equiv\widehat{1}\{\a,\b\}\ 
        \operatorname{mod}\ \R_N(A),\\[5pt]
        \{\a,\widehat{1}(\b)\}&\equiv \operatorname{Ad}((12)^{1,|\a|})\widehat{1}\{\a,\b\}\ 
        \operatorname{mod}\ \R_N(A),\\[5pt]
        \{w_1\cdot \a, w_2\cdot \b\}&\equiv(w_1\times w_2)\cdot\{\a,\b\}\ \operatorname{mod}\ \R_N(A),\\[5pt]
        \{\a\cdot w_1, \b\cdot w_2\}&\equiv\{\a,\b\}\cdot (w_1\times w_2)\ \operatorname{mod}\ \R_N(A),
    \end{align}
    and
    \begin{align}
    \{\b,\a\}&\equiv-\operatorname{Ad}((12)^{|\a|,|\b|})\{\a,\b\}\ \operatorname{mod}\ \R_N(A),\label{f48}\\[5pt]
    \{\a,\b\cc\}&\equiv\{\a,\b\}\cc+\operatorname{Ad}((12)^{|\b|,|\a|,|\cc|})\Bigl(\b\{\a,\cc\}\Bigr)\ \operatorname{mod}\ \R_N(A),\label{f154}\\[5pt]
    \{\a\b,\cc\}&\equiv\operatorname{Ad}((23)^{|\a|,|\cc|,|\b|})\Bigl(\{\a,\cc\}\b\Bigr)+\a\{\b,\cc\}\ \operatorname{mod}\ \R_N(A),\label{f155}\\[5pt]
    \hspace{-30pt}\bigl\{ \a,\{\b,\cc\}\bigr\} &+\operatorname{Ad}((123)^{|\b|,|\cc|,|\a|}) \bigl\{ \b,\{ \cc,\a\}\bigr\} +\operatorname{Ad}((132)^{|\cc|,|\a|,|\b|}) \bigl\{ \cc,\{ \a,\b\}\bigr\}\equiv0\ \operatorname{mod}\ \R_N(A).\label{f156}
    \end{align}

Moreover, any $\GL$-equivariant Poisson bracket on $\orep$ is of this form.
\end{proposition}
\begin{proof}
    This proposition is a corollary of Theorem \ref{th3} from which all the claims but \eqref{f48},\eqref{f154},\eqref{f155},\eqref{f156} follow. Their proof is similar to that of \eqref{f25}-\eqref{f58} from Theorem \ref{th3}.
\end{proof}

The next definition differs from its cousin Definition 2.6.1 from \cite{ginzburg2010differentialoperatorsbvstructures} only by the change $\ored(A)\leftrightarrow\O(A)$.

\begin{definition}\label{def12}
    By a di-twisted Poisson bracket on $\ored(A)$ we mean an admissible $\S$-bimodule homomorphism $\{-,-\}:\ored(A)\stimes \ored(A)\longrightarrow\ored(A)$ such that for any homogeneous $\a,\b,\cc\in\ored(A)$
    \begin{itemize}
        \item \textit{Skew-symmetry:\label{cond5}} 
    \begin{equation}\label{f49}
    \{\b,\a\}=-\operatorname{Ad}((12)^{|\a|,|\b|})\{\a,\b\}.    
    \end{equation}

    \item \textit{Leibniz rule:}\label{cond7}   
    \begin{align}
        \{\a,\b\cc\}&=\{\a,\b\}\cc+\operatorname{Ad}((12)^{|\b|,|\a|,|\cc|})\Bigl(\b\{\a,\cc\}\Bigr),\label{f50}\\[7pt]
        \{\a\b,\cc\}&=\operatorname{Ad}((23)^{|\a|,|\cc|,|\b|})\Bigl(\{\a,\cc\}\b\Bigr)+\a\{\b,\cc\}
    \end{align}

    \item \textit{Jacobi identity:}\label{cond8} 
    \begin{equation}\label{f52}
        \bigl\{ \a,\{\b,\cc\}\bigr\}+\operatorname{Ad}((123)^{|\b|,|\cc|,|\a|}) \bigl\{ \b,\{ \cc,\a\}\bigr\} +\operatorname{Ad}((132)^{|\cc|,|\a|,|\b|}) \bigl\{ \cc,\{ \a,\b\}\bigr\}=0.
    \end{equation}
    \end{itemize}
\end{definition}

\begin{remark}\label{rem2}
    Any double Poisson bracket on $A$ can be canonically extended to a di-twisted Poisson bracket on $\O(A)$, and hence, by Remark \ref{rem5}, it descends to a di-twisted Poisson bracket on $\ored(A)$. The part related to $\O(A)$ was mentioned without proof in \cite{ginzburg2010differentialoperatorsbvstructures}, see Remark 3.5.19, and proved in \cite{fernández2025symplecticwheelgebrasnoncommutativegeometry}, see Proposition 6.8. 
\end{remark}

    \subsection{Di-twisted deformations and star-products}

    \subsubsection{Deformations}

    \begin{definition}
        Let $\mathcal{A}$ be an associative algebra. By an associative $\Bbbk[[\hbar]]$-deformation of $\mathcal{A}$ we mean an associative $\Bbbk[[\hbar]]$-bilinear map $\star:\mathcal{A}[[\hbar]]\otimes \mathcal{A}[[\hbar]]\longrightarrow\mathcal{A}[[\hbar]]$ such that for any $a,b\in\mathcal{A}$ one has
    \begin{align}\label{f66}
        &a\star b=ab+O(\hbar).
    \end{align}
    We say that two $\Bbbk[[\hbar]]$-deformations $\star$, $\star'$ are \textit{gauge equivalent} if there is a $\kh$-linear map $f:\A[[\hbar]]\longrightarrow \A[[\hbar]]$ such that $f=\id+O(\hbar)$ and for any $a,b\in\A$ one has
    \begin{align}
        a\star' b=f^{-1}\big(f(a)\star f(b)\big).
    \end{align}
    \end{definition}

    Below we will consider deformations of $\orep$, and we will assume that they are $\GL$-equivariant along with the gauge transformations.

    Let $A$ be an associative algebra. By $\O(A)[[\hbar]]$ we denote the $\S$-bimodule whose $n$-th graded component is $\O(A)^{(n)}[[\hbar]]$ with the natural $S(n)$-bimodule action. To any homogeneous elements $\alpha(\hbar)=\sum\limits_{n\geq 0}\a_n\hbar^n\in\O(A)[[\hbar]]$ and $X\in\operatorname{T}(\Mat)$ such that $|\a(\hbar)|=|X|$ we can associate an element $(\a(\hbar)|X)\in\orep[[\hbar]]$ by the usual rule
    \begin{align}
        (\a(\hbar)|X)=\sum\limits_{n\geq 0}(\a_n|X)\hbar^n.
    \end{align}
    Note that this symbol is $\Bbbk[[\hbar]]$-linear, i.e. $\big(f(\hbar)\a(\hbar)|X\big)=f(\hbar)\big(\a(\hbar)|X\big)$ for $f(\hbar)\in\Bbbk[[\hbar]]$. 
    
    We can also extend operators $\pi,\widehat{1}:\O(A)\rightarrow\O(A)$ to $\O(A)[[\hbar]]$ by $\Bbbk[[\hbar]]$-linearity.

    %Let $\{-,-\}$ be a di-twisted Poisson bracket on $\O(A)$ and $\{-,-\}_N$ be the corresponding Poisson bracket on $\orep$. 

    \begin{proposition}
    A homogeneous $\Bbbk[[\hbar]]$-bilinear map $\dstar:\O(A)[[\hbar]]\otimes \O(A)[[\hbar]]\longrightarrow \O(A)[[\hbar]]$ of degree zero gives rise to a well-defined $\GL$-equivariant deformation $\star_N:\orep[[\hbar]]\otimes \orep[[\hbar]]\longrightarrow \orep[[\hbar]]$ by the rule
    \begin{align}
        \big(\a|X\big)\star_N\big(\b|Y\big)=\Big(\a\dstar\b\Big|X\otimes Y\Big),
    \end{align}
    where $\a,\b\in\O(A)$ and $X,Y\in\TMat$ are homogeneous elements such that $|\a|=|X|$ and $|\b|=|Y|$, if and only if for any graded elements $\a,\b,\cc\in\O(A)$ and any permutations $w_1\in S(|\a|),w_2\in S(|\b|)$ the map $\dstar$ satisfies 
    \begin{align}
        \pi\big(\a\big)\dstar\b&\equiv\pi\big(\a\dstar\b\big)\ 
        \operatorname{mod}\ \R_N(A)[[\hbar]],\\[5pt]
        \a\dstar\pi(\b)&\equiv\pi\operatorname{Ad}((12)^{|\a|,1})(\a\dstar\b)\ 
        \operatorname{mod}\ \R_N(A)[[\hbar]],\\[5pt]
        \widehat{1}(\a)\dstar\b&\equiv\widehat{1}(\a\dstar\b)\ 
        \operatorname{mod}\ \R_N(A)[[\hbar]],\\[5pt]
        \a\dstar\widehat{1}(\b)&\equiv \operatorname{Ad}((12)^{1,|\a|})\widehat{1}(\a\dstar\b)\ 
        \operatorname{mod}\ \R_N(A)[[\hbar]],\\[5pt]
        (w_1\cdot \a)\dstar (w_2\cdot \b)&\equiv(w_1\times w_2)\cdot(\a\dstar\b)\ \operatorname{mod}\ \R_N(A)[[\hbar]],\\[5pt]
        (\a\cdot w_1)\dstar (\b\cdot w_2)&\equiv(\a\dstar\b)\cdot (w_1\times w_2)\ \operatorname{mod}\ \R_N(A)[[\hbar]],
    \end{align}
    and
    \begin{align}
        \a\dstar \b&\equiv\a\b\ \operatorname{mod}\ \R_N(A)\oplus\hbar\O(A)[[\hbar]],\label{f44}\\
        %\a\dstar \b-\operatorname{Ad}((12)^{|\b|,|\a|})\big(\b\dstar \a\big)&\equiv\{\a,\b\}\hbar\ \operatorname{mod}\ \R_N(A)\oplus\hbar\R_N(A)\oplus\hbar^2\O(A)[[\hbar]],\label{f45}\\
        (\a\dstar\b)\dstar\cc&\equiv \a\dstar(\b\dstar\cc)\ \operatorname{mod}\ \R_N(A)[[\hbar]].\label{f46}
    \end{align}
    Moreover, any $\GL$-equivariant deformation of $\orep$ is of this form.
    \end{proposition}
    \begin{proof}
    All the claims but \eqref{f44},\eqref{f46} follow from Theorem \ref{th3} applied to Taylor's coefficients of $\star_N$ and $\dstar$ as series in $\hbar$. These identities follow from \eqref{f66} and associativity of $\star_N$.
    \end{proof}

    We call a $\Bbbk[[\hbar]]$-multilinear $\S$-bimodule homomorphism $\Phi:\big(\ored(A)[[\hbar]]\big)^{\stimes n}\longrightarrow\ored(A)[[\hbar]]$ \textit{admissible}, if all its Taylor's coefficients are admissible in the sense of Definition \ref{def6}.

    \begin{definition}\label{def7}
        By a \textit{di-twisted $\kh$-deformation} $\dstar$ of $\ored(A)$ we mean an associative admissible $\Bbbk[[\hbar]]$-bilinear $\S$-bimodule homomorphism $\dstar:\ored(A)[[\hbar]]\stimes \ored(A)[[\hbar]]\longrightarrow\ored(A)[[\hbar]]$ such that for any $\a,\b\in\ored(A)$ one has
        \begin{align}
        \a\dstar \b&=\a\b+O(\hbar).
    \end{align}
    We say that two di-twisted $\Bbbk[[\hbar]]$-deformations $\dstar$, $\dstar'$ are \textit{gauge equivalent} if there is an admissible $\kh$-linear $\S$-bimodule homomorphism $f:\ored(A)[[\hbar]]\longrightarrow \ored(A)[[\hbar]]$ such that $f=\id+O(\hbar)$ and for any $\a,\b\in\ored(A)$ one has
    \begin{align}
        \a\dstar' \b=f^{-1}\big(f(\a)\dstar f(\b)\big).
    \end{align}
    \end{definition}

    \begin{remark}%see Image37
        If $\dstar$ is a di-twisted $\kh$-deformation of $\ored(A)$, then the map $\{-,-\}:\ored(A)\stimes \ored(A)\longrightarrow\ored(A)$ defined by 
        \begin{align}
            \{\a,\b\}:=\cfrac{\a\dstar\b-\operatorname{Ad}((12)^{|\b|,|\a|})\b\dstar\a}{\hbar} \Bigg\vert_{\hbar=0}
        \end{align}
        is a di-twisted Poisson bracket.   
    \end{remark}

\subsubsection{Differential operators}

For commutative algebras we will use A.Grothendieck's definition of differential operators.

\begin{definition}
    Let $\mathcal{A}$ be a commutative algebra over $\Bbbk$. The only differential operator of order $-1$ is the zero. Differential operators of order $\leq m$ are those linear maps $D\in\Hom_{\Bbbk}(\mathcal{A},\mathcal{A})$ for which the map $[D,a]:b\mapsto D(ab)-aD(b)$ is a differential operator of order $\leq m-1$ for any $a\in\mathcal{A}$.
\end{definition}

We will apply this definition only to $\orep$, and we will actually replace arbitrary linear maps $\Hom_{\Bbbk}\Big(\orep,\orep\Big)$ with $\bigoplus\limits_{p\geq 0}\Hom_{\Bbbk}\Big(\orep,\orep\Big)^{(p)}$. Namely, we will say that a differential operator $D\in \Hom_{\Bbbk}\Big(\orep,\orep\Big)^{(p)}$ is of \textit{weight $p$}. Then it's clear that for any homogeneous $\a\in\O(A)$ and $X\in\Mat^{\otimes |\a|}$ the differential operator $\big[D,(\a|X)\big]$ is of weight $p+|\a|$, i.e. $\big[D,(\a|X)\big]\in \Hom_{\Bbbk}\Big(\orep,\orep\Big)^{(p+|\a|)}$.

This motivates us to define differential operators on $\ored(A)$ as follows.

\begin{definition}\label{def9}
    The only differential operator of order $-1$ is zero. We say that $\D\in\Hompadm{p}\Big(\ored(A),\ored(A)\Big)$ is a differential operator of order $\leq m$ and weight $p$ if for any homogeneous $\a\in\O(A)$ the operator $[\D,L_{\a}]^{(p)}:=\D\bullet L_{\a}-\operatorname{Ad}((12)^{|\a|,p})L_{\a}\bullet \D$ is a differential operator of order $\leq m-1$ and weight $p+|\a|$\footnote{It is clear from  Lemma \ref{lemma19} that for any $\D\in\Hompadm{p}\Big(\ored(A),\ored(A)\Big)$ and any homogeneous $\a\in\ored(A)$ one has $[D,L_{\a}]^{(p)}\in \Hompadm{p+|\a|}\Big(\ored(A),\ored(A)\Big)$.}, where $L_{\a}$ is the operator of left multiplication by $\a$ and the composition $\bullet$ was defined in Lemma \ref{lemma19}.
\end{definition}

This definition repeats the definition of differential operators studied \cite{ginzburg2010differentialoperatorsbvstructures}, see Section 2.4, with the only difference that we take into account $\1$, but the authors do not.

For any $\Phi\in\Hompadm{p}\Big(\ored(A)^{\stimes n},\ored(A)\Big)$ and $Y\in\Mat^{\otimes p}$ we denote by $(\Phi|Y)$ the induced map $(\Phi|Y)\in\Hom_{\Bbbk}\Big(\orep^{\otimes n},\orep\Big)^{(p)}$ defined by
\begin{align}
    (\a_1|X_1)\otimes\ldots\otimes(\a_n|X_n)\mapsto\Big(\Phi(\a_1,\ldots,\a_n)\, \Big|\, Y\otimes X_1\otimes\ldots\otimes X_n\Big).
\end{align}

\begin{proposition}\label{prop21}
    $\D\in\Hompadm{p}\Big(\ored(A),\ored(A)\Big)$ is a differential operator if and only if for any $N$ and any $Y\in\Mat^{\otimes p}$ the induced map $(\D|Y)$ is a differential operator of the same order as $\D$.
\end{proposition}
\begin{proof}
    The proof is by induction on the order of $\D$, which we will denote by $m$. The base $m=-1$ is clear from the definition of $\ored(A)$. Assume the claim holds for $m-1$. For any $N$ take any homogeneous $\a\in\ored(A)$, $X\in\Mat^{\otimes|\a|}$, and $Y\in\Mat^{\otimes p}$. Then for the induced operators we have
    \begin{align}
        &(\D\bullet L_{\a}|Y\otimes X)=(\D|Y)\circ L_{(a|X)},\\
        &\big(\operatorname{Ad}((12)^{|\a|,p})L_{\a}\bullet \D\big|Y\otimes X\big)=\big(L_{\a}\bullet \D\big|X\otimes Y\big)=L_{(\a|X)}\circ (\D|Y),
    \end{align}
    hence $\Big([\D,L_{\a}]^{(p)}\Big|Y\otimes X\Big)=\big[(\D|Y),L_{(a|X)}\big]$ and we are done.
\end{proof}

Let us extend the previous definition and proposition to poly-differential operators. 

\begin{definition}
    By an $n$-poly-differential operators on a commutative algebra $\A$ we mean $\Bbbk$-linear differential operators between $\A^{\otimes n}$ and $\A$ both viewed as $R=\A^{\otimes n}$-modules, where $\A$ is an $R$-module via the multiplication map $R\longrightarrow \A$. Explicitly, this means that the only $n$-poly-differential operator of order $-1$ is the zero operator, and an $n$-poly-differential operator of order $m\geq 0$ is a linear map $D:\A^{\otimes n}\longrightarrow \A$ such that for any $k=1,\ldots,n$ and any $a\in\A$ the operator $[D,L_{a}]_k:=D\circ_k L_a-L_a\circ D$ is an $n$-poly-differential operator of order $m-1$, where $\circ_k$ plugs $L_a$ into the $k$-th argument of $D$.
\end{definition}

\begin{definition}\label{def10}
    The only $n$-poly-differential operator of order $-1$ is zero. We say that $\D\in\Hompadm{p}\Big(\ored(A)^{\stimes n},\ored(A)\Big)$ is an $n$-poly-differential operator of order $\leq m$ and weight $p$ if for any $k=1,\ldots,n$ and any homogeneous $\a\in\O(A)$ the operator $[\D,L_{\a}]_k^{(p)}:=\D\bullet_k L_{\a}-\operatorname{Ad}((12)^{|\a|,p})L_{\a}\bullet_1 \D$ is a differential operator of order $\leq m-1$ and weight $p+|\a|$, where $\bullet_k$ is the composition defined in Lemma \ref{lemma19}.
\end{definition}

\begin{proposition}\label{prop20}
     $\D\in\Hompadm{p}\Big(\ored(A)^{\stimes n},\ored(A)\Big)$ is an $n$-poly-differential operator if and only if for any $N$ and any $Y\in\Mat^{\otimes p}$ the induced map $(\D|Y)$ is an $n$-poly-differential operator of the same order as $\D$.
\end{proposition}
\begin{proof}
    The proof is the same as of the proposition above and is based on the following identities
    \begin{align}
        &(\D\bullet_k L_{\a}|Y\otimes X)=(\D|Y)\circ_k L_{(a|X)},\\
        &\big(\operatorname{Ad}((12)^{|\a|,p})L_{\a}\bullet_1 \D\big|Y\otimes X\big)=\big(L_{\a}\bullet_1 \D\big|X\otimes Y\big)=L_{(\a|X)}\circ (\D|Y),
    \end{align}
    which yield $\Big([\D,L_{\a}]_k^{(p)}\Big|Y\otimes X\Big)=\big[(\D|Y),L_{(a|X)}\big]_k$.
\end{proof}

Particularly, the multiplication $\mu:\ored(A)\stimes\ored(A)\longrightarrow\ored(A)$ is a bi-differential operator of order zero.

\begin{corollary}\label{cor5}
    Let $\D$ be an $n$-poly-differential operator on $\ored(A)$ of order $\ord(\D)$ and weight $p$, and $\D'$ be an $m$-poly-differential operator on $\ored(A)$ of order $\ord(\D')$ and weight $q$. Then $\D\bullet_k\D'$ is an $n+m$-poly-differential operator of order $\ord(D)+\ord(D')$ and weight $p+q$ for any $k=1,\ldots,n+1$.
\end{corollary}

\subsubsection{Star-products}

\begin{definition}\label{def11}
        By a \textit{di-twisted star-product} $\dstar$ on $\ored(A)$ we mean a di-twisted deformation realized by bi-differential operators of weight zero, i.e. an associative admissible $\Bbbk[[\hbar]]$-bilinear $\S$-bimodule homomorphism $\dstar:\ored(A)[[\hbar]]\stimes \ored(A)[[\hbar]]\longrightarrow\ored(A)[[\hbar]]$ such that for any $\a,\b\in\ored(A)$ one has
        \begin{align}
        \a\dstar \b&=\a\b+O(\hbar),
    \end{align}
    and each Taylor coefficient of $\dstar$ is a bi-differential operator of weight zero $\ored(A)\stimes\ored(A)\longrightarrow\ored(A)$.
    We say that two di-twisted star-products $\dstar$, $\dstar'$ are \textit{gauge equivalent} if there is an admissible $\kh$-linear $\S$-bimodule homomorphism $f:\ored(A)[[\hbar]]\longrightarrow \ored(A)[[\hbar]]$ such that each its Taylor coefficient is a differential operator of weight zero, $f=\id+O(\hbar)$, and for any $\a,\b\in\ored(A)$ one has
    \begin{align}
        \a\dstar' \b=f^{-1}\big(f(\a)\dstar f(\b)\big).
    \end{align}
    \end{definition}

In Section \ref{section_example}, we present an example of a di-twisted star-product on a noncommutative affine space equipped with a double Poisson bracket, i.e. a star-product on $\O(\Bbbk\langle x_1,\ldots,x_d\rangle)$ with the di-twisted Poisson bracket coming from an arbitrary double Poisson bracket on $\Bbbk\langle x_1,\ldots,x_d\rangle$.

\subsection{Di-twisted Hochschild complex}

Let us recall the definition of the shifted Hochschild complex of an associative algebra $\A$ and the Gerstenhaber bracket on it along with the differential, turning it into a dg Lie algebra, see for example Section 3 in \cite{cattaneo2005deformation}. 

Let $L_{As}(\A)$ denote the shifted Hochschild complex of $\A$, whose $n$-th graded component is given by $\Hom_{\Bbbk}(\A^{\otimes n+1},\A)$. The Gerstenhaber bracket of two cochains $f\in L^n_{As}(\A)$ and $g\in L^m_{As}(\A)$ is given by 
\begin{align}
    [f,g]=f\bullet g-(-1)^{nm}g\bullet f,
\end{align}
where $\bullet$ is the following (nonassociative) product
\begin{align}\label{f51}
    (f\bullet g)(a_0,\ldots,a_{n+m}):=\sum_{i=0}^n(-1)^{mi}f(a_0,\ldots,a_{i-1},g(a_i,\ldots,a_{i+m}),a_{i+m+1},\ldots,a_{n+m}).
\end{align}

Then the differential $d$, which is of degree $+1$, is defined on homogeneous elements $f$ as $df=(-1)^{|f|}[\mu,f]$, where $\mu$ is the multiplication in $\A$. Explicitly,
\begin{align}
    \hspace{30pt}&\hspace{-30pt}(df)(a_0,\ldots,a_{n+1})\\
    &=a_0f(a_1,\ldots,a_{n+1})-\sum\limits_{i=0}^n(-1)^if(a_0,\ldots,a_{i-1},a_ia_{i+1},a_{i+2},\ldots,a_{n+1})+(-1)^nf(a_0,\ldots,a_n)a_{n+1}.
\end{align}

Recall that this dg Lie algebra $L_{As}(\A)$ controls the deformations of the multiplication in $\A$ modulo the gauge equivalence.

Consider now $L_{As}\big(\orep\big)$ for an associative unital algebra $A$ and a dg Lie subalgebra $L_{As}\Big(\orep\Big)^{\GL}$ consisting of $\GL$-equivariant cochains, i.e. 
\begin{align}
    L^n_{As}\Big(\orep\Big)^{\GL}:=\Hom_{\GL}\Big(\orep^{\otimes n+1},\orep\Big).
\end{align}

It's indeed a dg Lie subalgebra, as the multiplication in $\orep$ and the Gerstenhaber product $\bullet$ are $\GL$-equivariant, and it controls the $\GL$-equivariant deformations of $\orep$ modulo $\GL$-equivariant gauge equivalence. 

One can also consider the dg Lie subalgebra $L_{star}\Big(\orep\Big)^{\GL}$ of $L_{As}\Big(\orep\Big)^{\GL}$ consisting of poly-differential operators. This subalgebra controls associative deformations of $\orep$ realized by bi-differential operators, i.e. star-products, modulo gauge equivalence realized by differential operators.

\begin{definition}\label{def13}
    By  a \textit{di-twisted Hochschild complex} of $\ored(A)$, denoted by $L_{di\mh As}(\ored(A))$, we call the dg Lie subalgebra of $L_{As}(\ored(A))$ consisting of admissible, in the sense of Definition \ref{def6}, $\S$-bimodule homomorphisms, i.e.
    \begin{align}
        \Ldi^n(\ored(A))=\Hom_{\S}^{\operatorname{adm}}\big(\ored(A)^{\stimes n+1},\ored(A)\big).
    \end{align}
    We will denote by $\Ldistar(\ored(A))$ the dg Lie subalgebra of $\Ldi(\ored(A))$ consisting of poly-differential operators of weight zero in the sense of Definition \ref{def10}.
\end{definition}

Note that $\Ldi(\ored(A))$ and $\Ldistar(\ored(A))$ are indeed dg Lie algebras due to Lemma \ref{lemma19} and Corollary \ref{cor5}.

\begin{remark}\label{rem3}
\begin{enumerate}[label=\theenumi),leftmargin=3ex]
    \item The dg Lie algebra $\Ldi(\ored(A))$ controls the di-twisted deformations of the multiplication in $\ored(A)$ modulo the gauge equivalence realized by admissible $\S$-bimodule homomorphisms, and $\Ldistar(\ored(A))$ controls di-twisted star-products on $\ored(A)$ modulo the gauge equivalence realized by differential operators of weight zero.

    \item There is a commutative diagram in the category of dg Lie algebras
    \begin{equation}
        \begin{tikzcd}
            \Ldistar(\ored(A))\arrow[r,hook]\arrow[d] & \Ldi(\ored(A)) \arrow[d]  \\[10pt]            
            L_{star}\Big(\orep\Big)^{\GL}\arrow[r,hook] &  L_{As}\Big(\orep\Big)^{\GL},
        \end{tikzcd}
    \end{equation}
    where the horizontal arrows are natural inclusions and the vertical arrows are induced by the representation functor. 
\end{enumerate}    
\end{remark}

\subsubsection{Normalized cochains}\label{section_normalized_cochains}

By normalized cochains we will mean linear maps vanishing when at least one of their arguments equals $1$. So, we have complexes of normalized cochains $L_{1,star}\Big(\orep\Big)^{\GL}$, $L_{1,As}\Big(\orep\Big)^{\GL}$, $\Ldistarnorm(\ored(A))$, and $\Ldinorm(\ored(A))$, which control the deformations preserving the units. It would be convenient to organize them in the following commutative diagram, where all the vertical arrows are induced by the representation functor.

\begin{equation}\label{f87}
    \begin{tikzcd}
        &[25pt] &[-230pt] \Ldistar(\ored(A))\arrow[r,hook]\arrow[dd] &[25pt]  \Ldi(\ored(A))\arrow[dd]\\
        \Ldistarnorm(\ored(A))\arrow[r,hook,crossing over]\arrow[urr,hook,"\ii_{\star}"]\arrow[dd] & \Ldinorm(\ored(A))\arrow[urr,hook,"\ii"] & &
        \\[20pt]
        & & L_{star}\Big(\orep\Big)^{\GL}\arrow[r,hook] &  L_{As}\Big(\orep\Big)^{\GL} \\
        L_{1,star}\Big(\orep\Big)^{\GL}\arrow[r,hook]\arrow[urr,hook] & L_{1,As}\Big(\orep\Big)^{\GL}\arrow[urr,hook]\arrow[from=uu,crossing over] & &
    \end{tikzcd}
\end{equation}

The main goal of this subsection is to prove that $\ii_{\star}$ and $\ii$ on the top level of $\eqref{f87}$ are quasi-isomorphisms, which, particularly, implies that any di-twisted star-product or any di-twisted deformation on $\ored(A)$ admits a unit. We will explain that this quasi-isomorphism is merely a restriction of the well-known quasi-isomorphism between (non-twisted) cobar and normalized cobar complexes.  

\begin{proposition}\label{prop22}
    The natural inclusions $\ii$ and $\ii_{\star}$ from $\eqref{f87}$ are quasi-isomorphisms.
\end{proposition}
\begin{proof}
Below we briefly recall the construction from the proof of Theorem 6.1 in Section 6 of Chapter VIII from \cite{mac1963homology} which leads to a homotopy retract
\begin{equation}\label{f90}
    \begin{tikzcd}
        L_{As}(\A) 
            \arrow[loop left, looseness=6, "h"] 
            \arrow[r,two heads, shift left=0.8ex,"p"] 
        & L_{1,As}(\A) \arrow[l,hook', shift left=0.8ex,"i"]
    \end{tikzcd}
\end{equation}
\begin{align}
    pi=\id,\ \ \ \id-ip=dh+hd.
\end{align}
for any associative algebra $\A$. See also Corollary 2.2 in Chapter X of \cite{mac1963homology}. Note that here we work with the usual (non-twisted) complexes.

Retract \eqref{f90} is \textit{universal} in the sense that only the unity of $\A$ and the multiplication are used -- the maps $p$ and $h$ are given as certain sums with universal integer coefficients of the operations involving only composition of elements from $L_{As}(\A)$ with the multiplication and unity. This means that applied to $\A=\ored(A)$, this retract preserves admissible $\S$-bimodule homomorphisms and poly-differential operators due to Lemma \ref{lemma19} and the fact that for any admissible map/di-twisted poly-differential operator $\Phi$, the map $\Phi(\ldots,\mathds{1},\ldots)$ is admissible/di-twisted poly-differential operator as well, which is clear from the definition because $\mathds{1}\in\ored(A)$ is of degree zero. This proves the claim.

The construction in the proof of Theorem 6.1 in Section 6 of Chapter VIII from \cite{mac1963homology} is given for the bar complex of $\A$ instead of $L_{As}(\A)$. So, we will need to apply $\Hom_{\A-bimod}(-,\A)$ to it, which is easy. Let us set 
\begin{align}
    \b_n=\A\otimes \A^{\otimes n}\otimes \A,\ \ B_n=\A\otimes \overline{\A}^{\otimes n}\otimes \A,
\end{align}
where $\overline{\A}=\faktor{\A}{\Bbbk}$. Let us also denote by $d_i$ for $i=0,\ldots,n$ the linear map $ d_i:\b_n\longrightarrow\b_{n-1}$ that multiplies the $i$-th and $i+1$-st components, and by $s_i:\b_n\longrightarrow\b_{n+1}$ the linear map that inserts the unit in the $i+1$ position moving the rest to the right, i.e.
\begin{align}
    d_i(a_0\otimes a_1\otimes\ldots\otimes a_{n+1})&=a_0\otimes\ldots\otimes a_ia_{i+1}\otimes a_{n+1},\\
    s_i(a_0\otimes a_1\otimes\ldots\otimes a_{n+1})&=a_0\otimes \ldots a_i\otimes 1\otimes a_{i+1}\otimes \ldots\otimes a_{n+1}.
\end{align}
Then the differential in the complex $\b$ is given by
\begin{align}
    d=\sum\limits_{i=0}^n(-1)^id_i:\b_n\longrightarrow \b_{n-1},
\end{align}
which is also a well-defined differential for $B$. For any homogeneous element $a\in\b$ we set
\begin{align}
    h_k(a)=\begin{cases}
        (-1)^ks_k(a),\ \ &\text{if}\ |a|\geq k,\\
        0,\ \ &\text{otherwise}.
    \end{cases}
\end{align}
Then $t_k=1-dh_k-h_kd:\b\longrightarrow\b$\footnote{Warning: this notation diverges with the one used in \cite{mac1963homology}.} is a morphism of complexes homotopical to identity, and $h_k$ is the homotopy. Next, set $t=t_0t_1t_2\ldots=\overset{\longrightarrow}{\prod\limits_{k\geq 0}}t_k$. This is a well-defined morphism of complexes $\b\longrightarrow\b$ homotopical to identity due to $t_k(a)=a$ for $k$ large enough and the fact that each $t_k$ is homotopical to identity. Let's denote the homotopy between $t$ and the identity morphism by $h$. Let us denote the canonical projection $\b\longrightarrow B$ by $\rho$. Then one can check that the map $g:B\longhookrightarrow\b$ given by $g(a)=t(a')$, where $a'\in\b$ is any such that $\rho(a')=a$, is well-defined. This gives a homotopy retract
\begin{equation}
    \begin{tikzcd}
        \b
            \arrow[loop left, "h"] 
            \arrow[r,two heads, shift left=0.6ex,"\rho"] 
        & B \arrow[l,hook', shift left=0.6ex,"g"]
    \end{tikzcd}
\end{equation}
\begin{align}
    \rho g=\id,\ \ \ \id-g\rho=dh+hd,
\end{align}
which gives \eqref{f90} by dualizing.
\end{proof}

\subsection{Di-twisted Chevalley-Eilenberg complex}

Let us recall the definition of $L_{Lie}(\g)$, the shifted Chevalley-Eilenberg complex of a Lie algebra $\g$. The $n$-th cochains are given by 
\begin{align}
    L^n_{Lie}(\g)=\Hom_{\Bbbk}(\wedge^{n+1}\g,\g).
\end{align}

The Richardson-Nijenhuis bracket $[-,-]_{RN}$, turning $L_{Lie}(\g)$ into a graded Lie algebra, is given by 
\begin{align}
    [f,g]_{RN}=f\bulletrn g-(-1)^{nm}g\bulletrn f
\end{align}
for $f\in L^n_{Lie}(\g)$, and $g\in L^m_{Lie}(\g)$, where
\begin{align}
    f\bulletrn g(x_0,\ldots,x_{n+m})=\sum\limits_{\sigma\in\Sh(m+1,n)}(-1)^{\sigma}f\Big(g(x_{\sigma(0)},\ldots,x_{\sigma(m)}),x_{\sigma(m+1)},\ldots,x_{\sigma(n+m)}\Big),
\end{align}
and $\Sh(m+1,n)$ is the set of permutations \mbox{$\sigma\in S(n+m+1)$} such that $\sigma(0)<\ldots<\sigma(m)$ and $\sigma(m+1)<\ldots<\sigma(n+m)$.

The differential $d$ turning $L_{Lie}(\g)$ into a dg Lie algebra is given by $df=(-1)^{|f|}[\b,f]_{RN}$, where $\b$ stands for the Lie bracket in $\g$. This dg Lie algebra controls the deformations of the Lie bracket in $\g$.

Let $\Big(\A,\{-,-\}\Big)$ be a Poisson algebra. Consider a dg Lie subalgebra $L_{Pois}(\A)$ of $L_{Lie}(\A)$, where $\A$ equipped with $\{-,-\}$ is viewed as a Lie algebra, consisting of cochains that are derivations in each argument. This dg Lie algebra $L_{Pois}(\A)$ controls the deformations of the Poisson bracket in $\A$.

Suppose now that $\orep$ is equipped with a $\GL$-invariant Poisson bracket and we wish to study its $\GL$-equivariant deformations. Then one has to deal with the dg Lie subalgebra of $L_{Pois}(\orep)$ consisting of $\GL$-equivariant cochains. We will denote this dg Lie algebra by $L_{Pois}\Big(\orep\Big)^{\GL}$.

\begin{proposition}\label{prop15}
    A homogeneous linear map $\Phi:\O(A)^{\otimes n}\longrightarrow \O(A)$ of degree zero gives rise to a well-defined $\GL$-equivariant $n$-vector field $\varphi_N$ on $\rep$, i.e. a $\GL$-equivariant linear map $\varphi_N:\bigwedge^n\orep\longrightarrow \orep$ which is a derivation in each argument, by the rule
    \begin{align}
        \varphi_N\big((\a_1|X_1)\wedge\ldots\wedge(\a_n|X_n)\big)=\Big(\Phi(\a_1,\ldots,\a_n)\Big|X_1\otimes\ldots\otimes X_n\Big),
    \end{align}
    where $\a_1,\ldots,\a_n\in\O(A)$ and $X_1,\ldots,X_n\in\TMat$ are homogeneous elements such that $|\a_i|=|X_i|$, if and only if the map $\Phi$ satisfies \eqref{f25}-\eqref{f58} from Theorem \ref{th3}, and, for any $\sigma\in S(n)$ and any homogeneous $\a_1,\ldots,\a_n\in\O(A)$, and $\b_1\in\O(A)$, one has
    \begin{align}
        &\Phi(\a_{\sigma^{-1}(1)},\ldots,\a_{\sigma^{-1}(n)})\equiv (-1)^{\sigma}\operatorname{Ad}(\sigma^{|\a_1|,\ldots,|\a_n|})\Phi(\a_1,\ldots,\a_n)\ 
        \operatorname{mod}\ \R_N(A)\\[10pt]
        &\Phi(\a_1,\ldots,\a_{i-1},\a_i\b_i,\a_{i+1},\ldots,\a_n)\equiv\operatorname{Ad}((12)^{|\a_i|,|\a_1|+\ldots+|\a_{i-1}|})\Big(\a_i\Phi(\a_1,\ldots,\a_{i-1},\b_i,\a_{i+1},\ldots,\a_n)\Big)\\
        &\hspace{180pt}+\operatorname{Ad}((12)^{|\b_i|,|\a_1|+\ldots+|\a_i|})\Big(\b_i\Phi(\a_1,\ldots,\a_n)\Big)\ 
        \operatorname{mod}\ \R_N(A).
    \end{align}

Moreover, any $\GL$-equivariant $n$-vector field $\varphi_N$ on $\rep$ is of this form.
\end{proposition}%see Image9 for the proof.

The partial case $n=1$ of this proposition suggests the following definition.

\begin{definition}
    By a \textit{di-twisted derivation} of $\ored(A)$ we mean an admissible $\S$-bimodule homomorphism $\partial:\ored(A)\longrightarrow\ored(A)$ which is a derivation, i.e. such that for any $\a,\b\in\ored(A)$ one has $\partial(\a\b)=\partial(\a)\b+\a\partial(\b)$. 
\end{definition}

\begin{definition}\label{def14}
    Let us denote by $L^n_{di\mh Pois}(\ored(A))\subset \Homadm\Big(\bigwedge\nolimits_{\S}^{n+1}\ored(A),\ored(A)\Big)$ the subspace of maps $\Phi$ satisfying 
\begin{align}
    &\Phi(\a_0,\ldots,\a_{i-1},\a_i\b_i,\a_{i+1},\ldots,\a_n)=\operatorname{Ad}((12)^{|\a_i|,|\a_0|+\ldots+|\a_{i-1}|})\Big(\a_i\Phi(\a_0,\ldots,\a_{i-1},\b_i,\a_{i+1},\ldots,\a_n)\Big)\\
        &\pushright{+\operatorname{Ad}((12)^{|\b_i|,|\a_0|+\ldots+|\a_i|})\Big(\b_i\Phi(\a_0,\ldots,\a_n)\Big)\hspace{40pt}}\label{f47}
\end{align}
for any homogeneous $\a_0,\ldots,\a_n,\b_i\in\ored(A)$ and any $i$.
\end{definition}
Note that it's enough to claim \eqref{f47} only for $i=n$: the rest follows from skew-symmetry. %see Image10

For $i=n$ one has
\begin{align}
    &\Phi(\a_0,\ldots,\a_{n-1},\a_n\b_n)=\Phi(\a_0,\ldots,\a_n)\b_n+\operatorname{Ad}((12)^{|\a_n|,|\a_0|+\ldots+|\a_{n-1}|})\Big(\a_n\Phi(\a_0,\ldots,\a_{n-1},\b_n)\Big).
\end{align}

Assume that $\{-,-\}$ is a fixed di-twisted Poisson bracket on $\ored(A)$. We are going to define a bracket $[-,-]_{di\mh RN}$ and a differential $d$ on $L_{di\mh Pois}(\ored(A))$ which turn it into a dg Lie algebra that controls the deformations of $\{-,-\}$.

The bracket $[-,-]_{di\mh RN}$ is given by 
\begin{align}
    [\Phi,\Psi]_{di\mh RN}=\Phi\bullet_{di\mh\mathsmaller{RN}} \Psi-(-1)^{nm}\Psi\bullet_{di\mh\mathsmaller{RN}} \Phi
\end{align}
for $\Phi\in L^n_{di\mh Pois}(\ored(A))$, and $\Psi\in L^m_{di\mh Pois}(\ored(A))$, where
\begin{align}
    &\Phi\bullet_{di\mh\mathsmaller{RN}} \Psi(\a_0,\ldots,\a_{n+m})\\
    &\hspace{50pt}=\sum\limits_{\sigma\in\Sh(m+1,n)}(-1)^{\sigma}\operatorname{Ad}(\sigma^{|\a_{\sigma(0)}|,\ldots,|\a_{\sigma(n+m)}|})\Big(\Phi\big(\Psi(\a_{\sigma(0)},\ldots,\a_{\sigma(m)}),\a_{\sigma(m+1)},\ldots,\a_{\sigma(n+m)}\big)\Big).
\end{align}

The differential $d$ is given by $d\Phi=(-1)^{\Phi}[\mathcal{B},\Phi]_{di\mh RN}$, where $\mathcal{B}=\{-,-\}\in L^1_{di\mh Pois}(\ored(A))$.

\begin{proposition}\label{prop18}
    \begin{enumerate}[label=\theenumi), leftmargin=3ex]
        \item The bracket $[-,-]_{di\mh RN}$ and differential $d$ defined above turn $L_{di\mh Pois}(\ored(A))$ into a dg Lie algebra.

        \item The dg Lie algebra $L_{di\mh Pois}(\ored(A))$ controls the deformations of the di-twisted Poisson bracket $\{-,-\}$ on $\ored(A)$ modulo the gauge equivalence infinitesimally given by di-twisted derivations of $\ored(A)$.

        \item The representation functor induces a dg Lie algebra homomorphism
        \begin{align}
            L_{di\mh Pois}(\ored(A))\longrightarrow L_{Pois}\Big(\orep\Big)^{\GL}.
        \end{align}
    \end{enumerate}
\end{proposition}
\begin{proof}
    The third item holds by the very definition of the bracket and differential. The first item follows then from Lemma \ref{lemma19} and Theorem \ref{th3}. The second item is straightforward.
\end{proof}

\section{Di-twisted star-product \texorpdfstring{$\dstar$}{} on \texorpdfstring{$N\amsmathbb{A}^d$}{NAd}}\label{section_example}

Throughout this section we assume that $A=\Bbbk\langle x_1,\ldots,x_d\rangle$ and $\{\!\!\{-,-\}\!\!\}$ is an arbitrary double Poisson on it.

\subsection{Combinatorial data: double admissible graphs and principal splittings}

Recall the definition of admissible graphs, see \cite{kontsevich2003deformation}.

\begin{definition}
    The set $G_n$ consists of labeled oriented finite graphs $\Gamma$ with simple edges (no multiple edges are allowed) without loops satisfying the following conditions
    \begin{enumerate}[label=\arabic*)]
        \item $\Gamma$ has $n+2$ vertices and $2n$ edges;

        \item The vertices are labeled by $\{1,\ldots,n\}\sqcup\{\mathbf{L},\mathbf{R}\}$;

        \item The edges are labeled by $e_1^1,e_1^2,\ldots,e_n^1,e_n^2$;

        \item For each $k=1,\ldots,n$ edges $e_k^1$ and $e_k^2$ start at the vertex $k$ (and necessarily end at different vertices).
    \end{enumerate}
    So, $G_n$ is a finite set with $(n(n+1))^n$ elements. Its elements are called \textit{admissible graphs}.
\end{definition}

For any admissible graph $\Gamma$ we denote by $E_{\Gamma}(k)$ the set of all edges terminating at vertex $k$. Similarly, we write $E_{\Gamma}(\mathbf{L})$ and $E_{\Gamma}(\mathbf{R})$ for the edges terminating at $\mathbf{L}$ and $\mathbf{R}$ respectively.

Now we are going to introduce two combinatorial structures which serve as a proper analog of admissible graphs: double admissible graphs and their principal splittings. 

\begin{definition}\label{def3}
    For any pairs of non-negative integers $l=(l_1,l_2)$ and $r=(r_1,r_2)$ we let $\amsmathbb{G}_n^{(l,r)}$ be the set of labelled oriented graphs $\mathbb{\Gamma}$ with simple edges (no multiple edges regardless of their direction are allowed) and satisfying the following conditions
    \begin{enumerate}[label=\arabic*)]
        \item $\mathbb{\Gamma}$ has $2n+l_1+l_2+r_1+r_2$ vertices and $4n+l_2+r_2$ edges;

        \item The vertices of $\mathbb{\Gamma}$ are of the following kinds 
        \begin{enumerate}[label=\alph*)]
            \item $2n$ ''numbered vertices'': $1.1$, $1.2$, $2.1$, $2.2,\ldots$, $n.1$, $n.2$;

            \item $l_1$ ''left vertices'': $\mathbf{L_1},\ldots,\mathbf{L_{l_1}}$;

            \item $r_1$ ''right vertices'': $\mathbf{R_1},\ldots,\mathbf{R_{r_1}}$;

            \item $l_2$ ''left loop vertices'': $\mathbf{O^{l}_{1}},\ldots,\mathbf{O^{l}_{l_2}}$;

            \item $r_2$ ''right loop vertices'': $\mathbf{O^{r}_{1}},\ldots,\mathbf{O^{r}_{r_2}}$.
        \end{enumerate}
        Below we consider pairs of numbered vertices of the form $b_k=\{k.1,k.2\}$, which we will call \textit{blocks}.

        \item The edges of $\mathbb{\Gamma}$ are of the following four kinds:
        \begin{enumerate}[label=\alph*)]
            \item $2n$ ''proper edges'': $e_1^1,e_1^2,\ldots,e_n^1,e_n^2$;

            \item $2n$ ''reflected edges'': $^*\!e_1^1,\phantom{}^*\!e_1^2,\ldots,\phantom{}^*\!e_n^1,\phantom{}^*\!e_n^2$;

            \item $l_2$ ''left loops'': $o_1^l,\ldots,o_{l_2}^l$ which are loops starting and ending at the corresponding left loop vertices;

            \item $r_2$ ''right loops'': $o_1^r,\ldots,o_{r_2}^r$ which are loops starting and ending at the corresponding right loop vertices.
        \end{enumerate}

        So, the loops are always fixed and are in bijection with the loop vertices.

        \item For each $k=1,\ldots,n$ 
        \begin{itemize}
            \item reflected edges are completely determined by the proper edges: $\phantom{}^*\!e_k^1$ starts at the terminal point of $e_k^1$ and ends at $k.2$, while $\phantom{}^*\!e_k^2$ starts at the terminal point of $e_k^2$ and ends at $k.1$;

            \item edges $e_k^1$ and $e_k^2$ start at $k.1$ and $k.2$ respectively, and 
            \begin{itemize}
                \item end at numbered vertices belonging to different blocks both distinct from $b_k$ or

                \item one of the edges ends at a left, right, left loop, or right loop vertex and the other ends at a numbered vertex from a block distinct from $b_k$ or

                \item one of the edges ends at a left or left loop vertex and the other ends at a right or right loop vertex.
            \end{itemize}
        \end{itemize}
    \end{enumerate}
    We call elements of $\amsmathbb{G}_n^{(l,r)}$ \textit{double admissible graphs}. 
\end{definition}

To fix a double admissible graph $\mathbb{\Gamma}$ it's enough to specify the proper edges. So, $\amsmathbb{G}_n^{(l,r)}$ is a finite set with $(2v)^{n}$ elements, where $v:=|r|\cdot|l|+2(n-1)(|r|+|l|+n-2)$ and $|l|=l_1+l_2$, $|r|=r_1+r_2$\footnote{For each block there is $2v=|l|(|r|+2(n-1))+|r|(|l|+2(n-1))+2(n-1)(|l|+|r|+2(n-2))$ options for outgoing edges.} for $n\geq 1$, and $|\amsmathbb{G}_0^{(l,r)}|=1$. For instance, $|\amsmathbb{G}_2^{(l,r)}|=100$ and $|\amsmathbb{G}_3^{(l,r)}|=17576$ for $|l|=|r|=1$, $|\amsmathbb{G}_3^{(l,r)}|=46656$ for $|l|=1$, $|r|=2$, and $|\amsmathbb{G}_3^{(l,r)}|=110592$ for $|l|=|r|=2$.

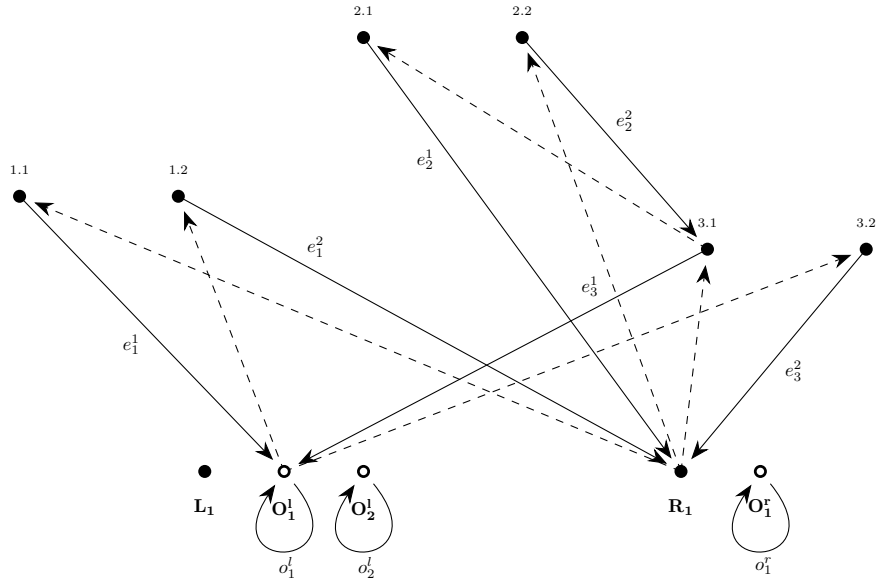
\begin{figure}[ht]
    \centering 
    \begin{tikzpicture}[scale=0.7, every node/.style={transform shape},
dot/.style = {circle, fill, minimum size=#1, inner sep=0pt, outer sep=5pt},
emptydot/.style = {circle, fill=white, minimum size=#1, inner sep=0pt, outer sep=5pt}]
%%%%%%%%%%%%%%%%%%%%%%%%%%%%%%%% preliminaries %%%%%%%%%%%%%%%%%%%%%%%%%%%%%%%%
\def\radd{7pt}% radius of vertices
\def\raddd{4pt}% radius of empty vertices
%%%%%%%%%%%%%%%%%%%%%%%%%%%%%%%% blocks and vertices %%%%%%%%%%%%%%%%%%%%%%%%%%%%%%%%

%block 1
\coordinate (center1) at (-3,0);

\path (center1) ++(-1.5,0.2) coordinate (1_1);
\node[dot=\radd ,label=$\mathsmaller{1.1}$] (p1_1) at (1_1) {};

\path (center1) ++(1.5,0.2) coordinate (1_2);
\node[dot=\radd, label=$\mathsmaller{1.2}$] (p1_2) at (1_2) {};

%block 2
\coordinate (center2) at (3.5,3);

\path (center2) ++(-1.5,0.2) coordinate (2_1);
\node[dot=\radd, label=$\mathsmaller{2.1}$] (p2_1) at (2_1) {};

\path (center2) ++(1.5,0.2) coordinate (2_2);
\node[dot=\radd, label=$\mathsmaller{2.2}$] (p2_2) at (2_2) {};

%block 3
\coordinate (center3) at (10,-1);

\path (center3) ++(-1.5,0.2) coordinate (3_1);
\node[dot=\radd,label=$\mathsmaller{3.1}$] (p3_1) at (3_1) {};

\path (center3) ++(1.5,0.2) coordinate (3_2);
\node[dot=\radd,label=$\mathsmaller{3.2}$] (p3_2) at (3_2) {};

%left
\coordinate (left) at (-1,-5);
\node[dot=\radd, label={[label distance=3pt]-90 :$\mathbf{L_1}$}] (l) at (left) {};

\path (left) ++(1.5,0) coordinate (left_loop_1);
\node[dot=\radd, label={[label distance=3pt]-90 :$\mathbf{O^l_1}$}] (ll1) at (left_loop_1) {};

\path (left_loop_1) ++(1.5,0) coordinate (left_loop_2);
\node[dot=\radd, label={[label distance=3pt]-90:$\mathbf{O^l_2}$}] (ll2) at (left_loop_2) {};

%right
\coordinate (right) at (8,-5);
\node[dot=\radd, label={[label distance=3pt]-90 :$\mathbf{R_1}$}] (r) at (right) {};

\path (right) ++(1.5,0) coordinate (right_loop_1);
\node[dot=\radd, label={[label distance=3pt]-90:$\mathbf{O^r_1}$}] (rl1) at (right_loop_1) {};
%%%%%%%%%%%%%%%%%%%%%%%%%%%%%%%% edges %%%%%%%%%%%%%%%%%%%%%%%%%%%%%%%%
%edges from block 1
\draw[-{Stealth[scale=1.5]}, black] (p1_1.center) -- (ll1) node[midway, below left] {$e_1^1$};
\draw[-{Stealth[scale=1.5]}, black] (p1_2.center) -- (r) node[near start, above right] {$e_1^2$};
\draw[-{Stealth[scale=1.5]}, black, dashed] (ll1.center) -- (p1_2);
\draw[-{Stealth[scale=1.5]}, black, dashed] (r.center) -- (p1_1) ;

%edges from block 2
\draw[-{Stealth[scale=1.5]}, black] (p2_1.center) -- (r) node[near start, below left] {$e_2^1$};
\draw[-{Stealth[scale=1.5]}, black] (p2_2.center) -- (p3_1) node[midway, above right] {$e_2^2$};
\draw[-{Stealth[scale=1.5]}, black, dashed] (r.center) -- (p2_2);
\draw[-{Stealth[scale=1.5]}, black, dashed] (p3_1.center) -- (p2_1);

%edges from block 3
\draw[-{Stealth[scale=1.5]}, black] (p3_1.center) -- (ll1) node[near start, above left] {$e_3^1$};
\draw[-{Stealth[scale=1.5]}, black] (p3_2.center) -- (r) node[midway, below right] {$e_3^2$};
\draw[-{Stealth[scale=1.5]}, black, dashed] (ll1.center) -- (p3_2);
\draw[-{Stealth[scale=1.5]}, black, dashed] (r.center) -- (p3_1);

%left loops
\path[-{Stealth[scale=1.5]}, black, every loop/.style={looseness=15}] (ll1) edge [in=-124,out=-51, loop] node[below] {$o_1^l$} (); 

\path[-{Stealth[scale=1.5]}, black, every loop/.style={looseness=15}] (ll2) edge  [in=-124,out=-51, loop] node[below] {$o_2^l$} ();

%right loops
\path[-{Stealth[scale=1.5]}, black, every loop/.style={looseness=15}] (rl1) edge  [in=-124,out=-51, loop] node[below] {$o_1^r$} ();

%to draw the empty dots above all edges
\node[emptydot=\raddd] at (left_loop_1) {};
\node[emptydot=\raddd] at (left_loop_2) {};
\node[emptydot=\raddd] at (right_loop_1) {};

\end{tikzpicture}
    \caption{A double admissible graph from $\amsmathbb{G}^{(l,r)}_3$ for $l=(1,2)$, $r=(1,1)$. Reflected edges are represented by dashed arrows. Labels of reflected edges are omitted for simplicity. The loop vertices are depicted as empty dots, all other vertices are represented by black dots.}
    \label{fig8}
\end{figure}

\begin{notation}
    For a given double admissible graph $\mathbb{\Gamma}$ we denote the set of
\begin{itemize}
    \item all edges of $\mathbb{\Gamma}$ by $\amsmathbb{E}_{\mathbb{\Gamma}}$;

    \item proper edges of $\mathbb{\Gamma}$ by $E_{\mathbb{\Gamma}}$;

    \item proper edges of $\mathbb{\Gamma}$ terminating at vertex $v$ by $E_{\mathbb{\Gamma}}(v)$;

    \item dashed edges of $\mathbb{\Gamma}$ starting at vertex $v$ by $\phantom{}^*\!E_{\mathbb{\Gamma}}(v)$.
\end{itemize}
\end{notation}

Note that the set $E_{\mathbb{\Gamma}}(v)$ is naturally linearly ordered by the first number of the starting point of every edge, which is a number from $1$ to $n$.

Clearly, there is a map $\amsmathbb{G}_n^{(l,r)}\rightarrow G_n$ which 
\begin{enumerate}[label=\arabic*)]
    \item deletes all dashed edges;

    \item merges vertices $k.1$ and $k.2$ for each $k=1,\ldots,n$;

    \item collapses the vertices $\mathbf{L_1},\ldots,\mathbf{L_{l_1}}$ and $\mathbf{O^l_{1}},\ldots,\mathbf{O^l_{l_2}}$ to $\mathbf{L}$, and $\mathbf{R_1},\ldots,\mathbf{R_{r_1}}$ and $\mathbf{O^r_{1}},\ldots,\mathbf{O^r_{r_2}}$ to $\mathbf{R}$.
\end{enumerate}
We will denote this map by $m_n^{(l,r)}:\amsmathbb{G}_n^{(l,r)}\rightarrow G_n$. 

\begin{definition}\label{def2}
    By an \textit{oriented path} in a directed graph we mean a sequence of vertices such that any two consecutive vertices are connected with an edge and all these edges are equally oriented, i.e. there is no vertex in the path which is the endpoint of two consecutive edges. Equivalently, an oriented path is a sequence of edges such that any two consecutive edges have a common vertex which is the endpoint of the first edge and the starting point of the second edge. By a \textit{cycle} in a directed graph we mean a cyclic oriented path considered up to a cyclic permutation of the vertices it consists of.
\end{definition}

\begin{definition}\label{def5}
    By a \textit{splitting} of a double admissible graph $\mathbb{\Gamma}\in\amsmathbb{G}^{(l,r)}_n$ we mean a collection of oriented paths $\tau_1,\ldots,\tau_{l_1}$, $\overline{\tau}_1,\ldots,\overline{\tau}_{r_1}$, and cycles $\gamma_1,\gamma_2,\ldots$ in $\mathbb{\Gamma}$, denoted by $\mathcal{S}$. The number of $\gamma$'s, which may be zero, and their order are not fixed. All these paths are subject to the following constraints
    \begin{enumerate}[label=\arabic*)]
        \item for any $t=1,\ldots,l_1$ path $\tau_t$ starts at the vertex $\mathbf{L_t}$ and ends at a left or right vertex, the length of $\tau_t$ is zero if and only if there are no proper edges terminating at $\mathbf{L_t}$;

        \item for any $s=1,\ldots,r_1$ path $\overline{\tau}_s$ starts at the vertex $\mathbf{R_s}$ and ends at a left or right vertex, the length of $\overline{\tau}_s$ is zero if and only if there are no proper edges terminating at $\mathbf{R_s}$;

        \item\label{cond4} endpoints of paths $\tau_1,\ldots,\tau_{l_1}$ and $\overline{\tau}_1,\ldots,\overline{\tau}_{r_1}$ are pairwise distinct;

        \item\label{cond3} each edge of $\mathbb{\Gamma}$ occurs in $\mathcal{S}$ exactly once;

        \item\label{cond2} patterns like $e\phantom{}^*\!e$, where $e$ is a proper edge of $\mathbb{\Gamma}$, are not allowed in paths or cycles from $\mathcal{S}$;

        \item\label{cond1} if a path or a
        cycle from $\mathcal{S}$ goes through a numbered vertex $v$ admitting at least one outgoing reflected edge, i.e. $|\phantom{}^*\!E_{\mathbb{\Gamma}}(v)|>0$, and comes to this vertex $v$ along a reflected edge, then it goes out of $v$ along a reflected edge too.
    \end{enumerate}
    Note that we do not forbid paths $\tau_1,\ldots,\tau_{l_1}$, $\overline{\tau}_1,\ldots,\overline{\tau}_{r_1}$ to be cyclic, however we do not consider them as cycles in the sense of Definition \ref{def2}.
\end{definition}

For the double admissible graph represented in Figure \ref{fig8} condition \ref{cond1} means that any path and any cycle passing through the vertex $3.1$, which comes to it along the reflected edge $\phantom{}^*\!e_3^2$, goes further along the reflected edge $\phantom{}^*\!e_2^2$ but not $e_3^1$. 

For any splitting $\mathcal{S}$ of $\mathbb{\Gamma}$ we denote by $P(\mathcal{S})$ the paths from $\mathcal{S}$ starting at a left or right vertex and ending at a left or right vertex, and by $C(\mathcal{S})$ all the cycles from $\mathcal{S}$ (''P'' from ''path'' and ''C'' from ''cycle''). Note that $|P(\mathcal{S})|$ is fixed and equals $l_1+r_1$.

To any splitting $\mathcal{S}$ of $\mathbb{\Gamma}$ we associate a collection of integer numbers $\{d_{\mathcal{S}}(v)\}_{v\in\mathbb{\Gamma}}$ indexed by the vertices $v$ of $\mathbb{\Gamma}$. We construct it by reordering the incoming proper edges of a vertex $v$ with the help of the algorithm below. The output of this algorithm is a collection of words $\{w_v\}_{v\in\mathbb{\Gamma}}$ indexed by the vertices of $\mathbb{\Gamma}$. The word $w_v$ is a word in the alphabet $E_{\mathbb{\Gamma}}(v)$. Then $d_{\mathcal{S}}(v)$ is the number of letters (or edges) in this word $w_v$. The algorithm works as follows. First we determine the initial edges by the rules 

\begin{itemize}
    \item If $v$ is a left or right vertex, then there is a unique path in $P(\mathcal{S})$ starting at $v$. By definition $w_v(1)$ is the proper edge corresponding to the reflected edge that path starts with.

    \item If $v$ is a left loop or right loop vertex, then there is a unique path or cycle that contains the loop corresponding to this vertex. The edge in this path or cycle which is next to the loop edge is necessarily reflected. By definition $w_v(1)$ is the proper edge corresponding to this reflected edge.

    \item Let $v$ be a numbered vertex. If there are no proper edges terminating at $v$, then $w_v$ permutes the empty set and is trivial. If there is at least one proper edge terminating at $v$, we consider the path or cycle from $\mathcal{S}$ which contains the reflected edge, say $\phantom{}^*\!e$, terminating at $v$. Then due to condition \ref{cond1} of Definition \ref{def1} the edge of the path which is next to $\phantom{}^*\!e$ is necessarily reflected. We set $w_v(1)$ to be the proper edge corresponding to this reflected edge.
\end{itemize}

Now, we know the first edges and start applying to them the algorithm described below. It either breaks or produces another proper edge. If the latter, we apply the algorithm to this new proper edge again and continue doing this until the algorithm breaks. 
\vspace{\baselineskip}

\textbf{Algorithm.}\label{algorithm}
Suppose a vertex of $\mathbb{\Gamma}$ and a proper edge terminating at this vertex are given. This edge belongs to a unique path or cycle from the splitting $\mathcal{S}$.
\begin{itemize}
    \item We break the algorithm, if \begin{itemize}
        \item[---]  in the case of a right or left vertex the proper edge belongs to a path which terminates at the chosen vertex;

        \item[---] in the case of a numbered vertex the path goes out of the vertex along the proper edge corresponding to this vertex; 

        \item[---] in the case of a loop vertex the edge in the path that is next to the given one is the loop.
    \end{itemize}

    \item Otherwise the path goes out of the chosen vertex along a reflected edge. We say that the proper edge corresponding to this reflected edge is the immediate follower of the proper edge we started with.
\end{itemize}

Note that conditions \ref{cond3} and \ref{cond2} from Definition \ref{def1} ensure that the words we obtain after application of the algorithm contain at most one copy of each edge. So, these words are almost permutational (each letter is used exactly once) except the thing that some edges may be missing. 

\begin{definition}
    We say that a splitting $\mathcal{S}$ of $\mathbb{\Gamma}$ is \textit{principal}, if the words $\{w_v\}_{v\in\mathbb{\Gamma}}$ obtained with the help of the algorithm above are permutational or equivalently the collection of numbers $d_{\mathcal{S}}$ is maximal in the sense that $d_{\mathcal{S}}(v)=|E_{\mathbb{\Gamma}}(v)|$ for any $v\in\mathbb{\Gamma}$. 
\end{definition}

The set of all principal splittings of $\mathbb{\Gamma}$ will be denoted by $\operatorname{split}(\mathbb{\Gamma})$. 

\begin{example}
    Consider two principal splittings given in Figure \ref{fig23}. For the reader's convenience we give two presentations for paths and cycles below: as collections of edges and as sequences of vertices.
    
    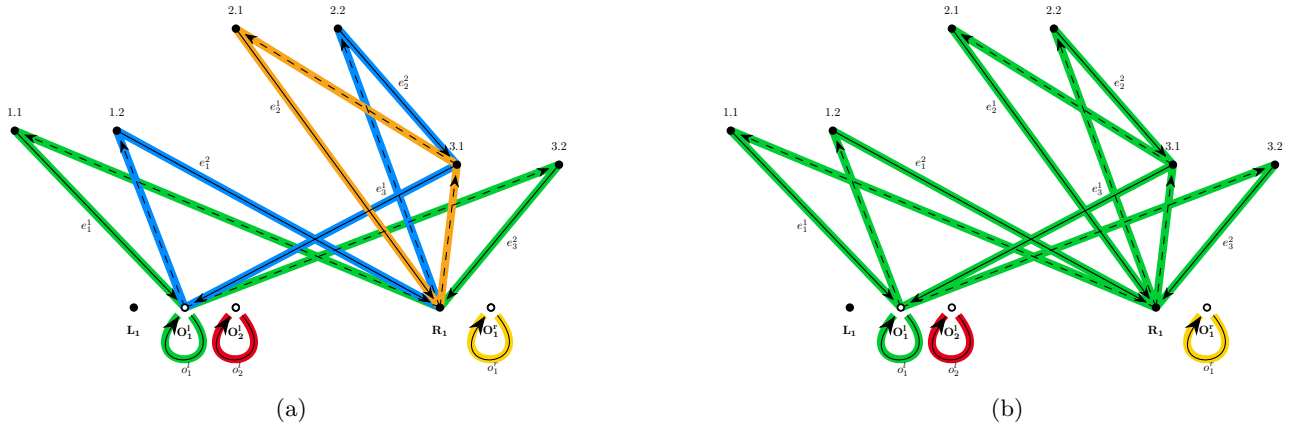
\begin{figure}[ht]
    \centering
    \begin{subfigure}[b]{0.45\textwidth}
         \centering 
         \begin{tikzpicture}[scale=0.45, every node/.style={transform shape},
dot/.style = {circle, fill, minimum size=#1, inner sep=0pt, outer sep=5pt},
emptydot/.style = {circle, fill=white, minimum size=#1, inner sep=0pt, outer sep=5pt}]
%%%%%%%%%%%%%%%%%%%%%%%%%%%%%%%% preliminaries %%%%%%%%%%%%%%%%%%%%%%%%%%%%%%%%
\def\radd{7pt}% radius of vertices
\def\raddd{4pt}% radius of empty vertices
%initialize vertices, so that we can draw colored arrows on the background
%block 1
\coordinate (center1) at (-3,0);
\path (center1) ++(-1.5,0.2) coordinate (1_1);
\path (center1) ++(1.5,0.2) coordinate (1_2);

%block 2
\coordinate (center2) at (3.5,3);
\path (center2) ++(-1.5,0.2) coordinate (2_1);
\path (center2) ++(1.5,0.2) coordinate (2_2);

%block 3
\coordinate (center3) at (10,-1);
\path (center3) ++(-1.5,0.2) coordinate (3_1);
\path (center3) ++(1.5,0.2) coordinate (3_2);

%left
\coordinate (left) at (-1,-5);
\path (left) ++(1.5,0) coordinate (left_loop_1);
\path (left_loop_1) ++(1.5,0) coordinate (left_loop_2);

\node[dot=\radd] (ll1) at (left_loop_1) {};
\node[dot=\radd] (ll2) at (left_loop_2) {};

%right
\coordinate (right) at (8,-5);
\path (right) ++(1.5,0) coordinate (right_loop_1);

\node[dot=\radd] (rl1) at (right_loop_1) {};
%%%%%%%%%%%%%%%%%%%%%%% colored arrows corresponding to splitting of the graph %%%%%%%%%%%%%%%%%%%%%%%
\draw[blue!20!green, line width=3.0pt] (left_loop_1) -- (3_2);
\draw[blue!20!green, line width=3.0pt] (3_2) -- (right);
\draw[blue!20!green, line width=3.0pt] (right) -- (1_1);
\draw[blue!20!green, line width=3.0pt] (1_1) -- (left_loop_1);
\path[blue!20!green, line width=3.0pt, every loop/.style={looseness=15}] (ll1) edge  [in=-131,out=-52, loop] (); 

\draw[blue!45!cyan, line width=3.0pt, opacity=1] (left_loop_1) -- (1_2);
\draw[blue!45!cyan, line width=3.0pt, opacity=1] (1_2) -- (right);
\draw[blue!45!cyan, line width=3.0pt, opacity=1] (right) -- (2_2);
\draw[blue!45!cyan, line width=3.0pt, opacity=1] (2_2) -- (3_1);
\draw[blue!45!cyan, line width=3.0pt, opacity=1] (3_1) -- (left_loop_1);

\draw[yellow!20!orange, line width=3.0pt] (right) -- (3_1);
\draw[yellow!20!orange, line width=3.0pt] (3_1) -- (2_1);
\draw[yellow!20!orange, line width=3.0pt] (2_1) -- (right);

\path[blue!10!red, line width=3.0pt, every loop/.style={looseness=15}] (ll2) edge  [in=-131,out=-52, loop] (); 

\path[yellow!70!orange, line width=3.0pt, every loop/.style={looseness=15}] (rl1) edge  [in=-131,out=-52, loop] (); 
%%%%%%%%%%%%%%%%%%%%%%%%%%%%%%%% blocks and vertices %%%%%%%%%%%%%%%%%%%%%%%%%%%%%%%%
%block 1
\node[dot=\radd ,label=$1.1$] (p1_1) at (1_1) {};
\node[dot=\radd, label=$1.2$] (p1_2) at (1_2) {};

%block 2
\node[dot=\radd, label=$2.1$] (p2_1) at (2_1) {};
\node[dot=\radd, label=$2.2$] (p2_2) at (2_2) {};

%block 3
\node[dot=\radd,label=$3.1$] (p3_1) at (3_1) {};
\node[dot=\radd,label=$3.2$] (p3_2) at (3_2) {};

%left
\node[dot=\radd, label={[label distance=3pt]-90 :$\mathbf{L_1}$}] (l) at (left) {};

\node[dot=\radd, label={[label distance=3pt]-90 :$\mathbf{O_1^l}$}] (ll1) at (left_loop_1) {};

\node[dot=\radd, label={[label distance=3pt]-90 :$\mathbf{O_2^l}$}] (ll2) at (left_loop_2) {};

%right
\node[dot=\radd,label={[label distance=3pt]-90 :$\mathbf{R_1}$}] (r) at (right) {};

\node[dot=\radd, label={[label distance=3pt]-90 :$\mathbf{O_1^r}$}] (rl1) at (right_loop_1) {};

%%%%%%%%%%%%%%%%%%%%%%%%%%%%%%%% edges %%%%%%%%%%%%%%%%%%%%%%%%%%%%%%%%
%edges from block 1
\draw[-{Stealth[scale=1]}, black] (p1_1.center) -- (ll1) node[midway, below left] {$e_1^1$};
\draw[-{Stealth[scale=1]}, black] (p1_2.center) -- (r) node[near start, above right] {$e_1^2$};
\draw[-{Stealth[scale=1]}, black, dashed] (ll1.center) -- (p1_2);
\draw[-{Stealth[scale=1]}, black, dashed] (r.center) -- (p1_1) ;

%edges from block 2
\draw[-{Stealth[scale=1]}, black] (p2_1.center) -- (r) node[near start, below left] {$e_2^1$};
\draw[-{Stealth[scale=1]}, black] (p2_2.center) -- (p3_1) node[midway, above right] {$e_2^2$};
\draw[-{Stealth[scale=1]}, black, dashed] (r.center) -- (p2_2);
\draw[-{Stealth[scale=1]}, black, dashed] (p3_1.center) -- (p2_1);

%edges from block 3
\draw[-{Stealth[scale=1]}, black] (p3_1.center) -- (ll1) node[near start, above left] {$e_3^1$};
\draw[-{Stealth[scale=1]}, black] (p3_2.center) -- (r) node[midway, below right] {$e_3^2$};
\draw[-{Stealth[scale=1]}, black, dashed] (ll1.center) -- (p3_2);
\draw[-{Stealth[scale=1]}, black, dashed] (r.center) -- (p3_1);

%left loops
\path[-{Stealth[scale=1.5]}, black, every loop/.style={looseness=15}] (ll1) edge  [in=-124,out=-51, loop] node[below] {$o_1^l$} (); 

\path[-{Stealth[scale=1.5]}, black, every loop/.style={looseness=15}] (ll2) edge  [in=-124,out=-51, loop] node[below] {$o_2^l$} (); 

%right loops

\path[-{Stealth[scale=1.5]}, black, every loop/.style={looseness=15}] (rl1) edge  [in=-124,out=-51, loop] node[below] {$o_1^r$} ();

%to draw the empty dots above all edges
\node[emptydot=\raddd] at (left_loop_1) {};
\node[emptydot=\raddd] at (left_loop_2) {};
\node[emptydot=\raddd] at (right_loop_1) {};

\end{tikzpicture}
         \caption{}
         \label{fig9}
     \end{subfigure}
     \hfill
     \begin{subfigure}[b]{0.45\textwidth}
         \centering
         \begin{tikzpicture}[scale=0.45, every node/.style={transform shape},
dot/.style = {circle, fill, minimum size=#1, inner sep=0pt, outer sep=5pt},
emptydot/.style = {circle, fill=white, minimum size=#1, inner sep=0pt, outer sep=5pt}]
%%%%%%%%%%%%%%%%%%%%%%%%%%%%%%%% preliminaries %%%%%%%%%%%%%%%%%%%%%%%%%%%%%%%%
\def\radd{7pt}% radius of vertices
\def\raddd{4pt}% radius of empty vertices
%initialize vertices, so that we can draw colored arrows on the background
%block 1
\coordinate (center1) at (-3,0);
\path (center1) ++(-1.5,0.2) coordinate (1_1);
\path (center1) ++(1.5,0.2) coordinate (1_2);

%block 2
\coordinate (center2) at (3.5,3);
\path (center2) ++(-1.5,0.2) coordinate (2_1);
\path (center2) ++(1.5,0.2) coordinate (2_2);

%block 3
\coordinate (center3) at (10,-1);
\path (center3) ++(-1.5,0.2) coordinate (3_1);
\path (center3) ++(1.5,0.2) coordinate (3_2);

%left
\coordinate (left) at (-1,-5);
\path (left) ++(1.5,0) coordinate (left_loop_1);
\path (left_loop_1) ++(1.5,0) coordinate (left_loop_2);

%right
\coordinate (right) at (8,-5);
\path (right) ++(1.5,0) coordinate (right_loop_1);
%%%%%%%%%%%%%%%%%%%%%%% colored arrows corresponding to splitting of the graph %%%%%%%%%%%%%%%%%%%%%%%
\draw[blue!20!green, line width=3.0pt] (right) -- (3_2);
\draw[blue!20!green, line width=3.0pt] (right) -- (1_2);
\draw[blue!20!green, line width=3.0pt] (right) -- (2_1);
\draw[blue!20!green, line width=3.0pt] (right) -- (2_2);
\draw[blue!20!green, line width=3.0pt] (right) -- (3_1);
\draw[blue!20!green, line width=3.0pt] (left_loop_1) -- (1_2);
\draw[blue!20!green, line width=3.0pt] (3_1) -- (2_1);
\draw[blue!20!green, line width=3.0pt] (3_1) -- (2_2);
\draw[blue!20!green, line width=3.0pt] (left_loop_1) -- (3_1);
\draw[blue!20!green, line width=3.0pt] (left_loop_1) -- (3_2);
\draw[blue!20!green, line width=3.0pt] (right) -- (1_1);
\draw[blue!20!green, line width=3.0pt] (1_1) -- (left_loop_1);
\path[blue!20!green, line width=3.0pt, every loop/.style={looseness=15}] (ll1) edge  [in=-131,out=-52, loop] ();

\path[blue!10!red, line width=3.0pt, every loop/.style={looseness=15}] (ll2) edge  [in=-131,out=-52, loop] (); 

\path[yellow!70!orange, line width=3.0pt, every loop/.style={looseness=15}] (rl1) edge  [in=-131,out=-52, loop] (); 
%%%%%%%%%%%%%%%%%%%%%%%%%%%%%%%% blocks and vertices %%%%%%%%%%%%%%%%%%%%%%%%%%%%%%%%
%block 1
\node[dot=\radd ,label=$1.1$] (p1_1) at (1_1) {};
\node[dot=\radd, label=$1.2$] (p1_2) at (1_2) {};

%block 2
\node[dot=\radd, label=$2.1$] (p2_1) at (2_1) {};
\node[dot=\radd, label=$2.2$] (p2_2) at (2_2) {};

%block 3
\node[dot=\radd,label=$3.1$] (p3_1) at (3_1) {};
\node[dot=\radd,label=$3.2$] (p3_2) at (3_2) {};

%left
\node[dot=\radd, label={[label distance=3pt]-90 :$\mathbf{L_1}$}] (l) at (left) {};

\node[dot=\radd, label={[label distance=3pt]-90 :$\mathbf{O_1^l}$}] (ll1) at (left_loop_1) {};

\node[dot=\radd, label={[label distance=3pt]-90 :$\mathbf{O_2^l}$}] (ll2) at (left_loop_2) {};

%right
\node[dot=\radd, label={[label distance=3pt]-90 :$\mathbf{R_1}$}] (r) at (right) {};

\node[dot=\radd, label={[label distance=3pt]-90 :$\mathbf{O_1^r}$}] (rl1) at (right_loop_1) {};
%%%%%%%%%%%%%%%%%%%%%%%%%%%%%%%% edges %%%%%%%%%%%%%%%%%%%%%%%%%%%%%%%%
%edges from block 1
\draw[-{Stealth[scale=1]}, black] (p1_1.center) -- (ll1) node[midway, below left] {$e_1^1$};
\draw[-{Stealth[scale=1]}, black] (p1_2.center) -- (r) node[near start, above right] {$e_1^2$};
\draw[-{Stealth[scale=1]}, black, 
dashed] (ll1.center) -- (p1_2);
\draw[-{Stealth[scale=1]}, black, dashed] (r.center) -- (p1_1) ;

%edges from block 2
\draw[-{Stealth[scale=1]}, black] (p2_1.center) -- (r) node[near start, below left] {$e_2^1$};
\draw[-{Stealth[scale=1]}, black] (p2_2.center) -- (p3_1) node[midway, above right] {$e_2^2$};
\draw[-{Stealth[scale=1]}, black, dashed] (r.center) -- (p2_2);
\draw[-{Stealth[scale=1]}, black, dashed] (p3_1.center) -- (p2_1);

%edges from block 3
\draw[-{Stealth[scale=1]}, black] (p3_1.center) -- (ll1) node[near start, above left] {$e_3^1$};
\draw[-{Stealth[scale=1]}, black] (p3_2.center) -- (r) node[midway, below right] {$e_3^2$};
\draw[-{Stealth[scale=1]}, black, dashed] (ll1.center) -- (p3_2);
\draw[-{Stealth[scale=1]}, black, dashed] (r.center) -- (p3_1);

%left loops
\path[-{Stealth[scale=1.5]}, black, every loop/.style={looseness=15}] (ll1) edge  [in=-124,out=-51, loop] node[below] {$o_1^l$} (); 

\path[-{Stealth[scale=1.5]}, black, every loop/.style={looseness=15}] (ll2) edge  [in=-124,out=-51, loop] node[below] {$o_2^l$} (); 

%right loops
\path[-{Stealth[scale=1.5]}, black, every loop/.style={looseness=15}] (rl1) edge  [in=-124,out=-51, loop] node[below] {$o_1^r$} ();

%to draw the empty dots above all edges
\node[emptydot=\raddd] at (left_loop_1) {};
\node[emptydot=\raddd] at (left_loop_2) {};
\node[emptydot=\raddd] at (right_loop_1) {};

\end{tikzpicture}
         \caption{}
     \end{subfigure}
    \caption{Two principal splittings of the graph $\mathbb{\Gamma}$ represented in Figure \ref{fig8}. The only path from $P(\mathcal{S})$ is highlighted in green. All other colors represent the cycles from $C(\mathcal{S})$.}
    \label{fig23}
\end{figure}

    \begin{align}
        (a): \mathcal{S}&=(\tau,\overline{\tau}, \gamma_1,\gamma_2,\gamma_3,\gamma_4), &(b): \mathcal{S}&=(\tau,\overline{\tau}, \gamma_1,\gamma_2),\\
        \tau&=\diameter=(\mathbf{L_1}), &\phantom{aaaa}\tau&=\diameter=(\mathbf{L_1}),\\
        \phantom{aaaa}\overline{\tau}&=(\phantom{}^*\!e_1^2, e_1^1,o_1^l,\phantom{}^*\!e_3^1,e_3^2) &\phantom{aaaa}\overline{\tau}&=(\phantom{}^*\!e_1^2, e_1^1,o_1^l,\phantom{}^*\!e_3^1, e_3^2, \phantom{}^*\!e_2^1, e_2^2, e_3^1, \phantom{}^*\!e_1^1,e_1^2,\phantom{}^*\!e_3^2,\phantom{}^*\!e_2^2,e_2^1),\\
        \phantom{aaaa\overline{\tau}}&=(\mathbf{R_1},1.1,\mathbf{O_1^l},\mathbf{O_1^l},3.2,\mathbf{R_1}), &\phantom{aaaa}&=(\mathbf{R_1},1.1,\mathbf{O_1^l},\mathbf{O_1^l},3.2,\mathbf{R_1},2.2,3.1,\mathbf{O_1^l},\\
        & &\phantom{aaaa}&\phantom{aaaaaaaaaaaaaaaaaaaaaaa}1.2,\mathbf{R_1},3.1,2.1,\mathbf{R_1}),\\
        \phantom{aaaa}\gamma_1&=(\phantom{}^*\!e_3^2,\phantom{}^*\!e_2^2,e_2^1)&\phantom{aaa_2}\gamma_1&=(o_2^l)=(\mathbf{O_2^l},\mathbf{O_2^l}),\\
        \phantom{aaaac_1}&=(\mathbf{R_1},3.1,2.1,\mathbf{R_1}), \\
        \phantom{aaaa}\gamma_2&=(\phantom{}^*\!e_2^1,e_2^2, e_3^1, \phantom{}^*\!e_1^1,e_1^2) &\phantom{aaa_2}\gamma_2&=(o_1^r)=(\mathbf{O_1^r},\mathbf{O_1^r}).\\
        \phantom{aaaac_2}&=(\mathbf{R_1},2.2,3.1,\mathbf{O_1^l},1.2,\mathbf{R_1}),\\
        \phantom{aaaa}\gamma_3&=(o_2^l)=(\mathbf{O_2^l},\mathbf{O_2^l}),\\
        \phantom{aaaa}\gamma_4&=(o_1^r)=(\mathbf{O_1^r},\mathbf{O_1^r}).
    \end{align}

The only nontrivial permutations corresponding to these splittings are given below. All other permutations are trivial due to the fact that they permute the empty set.

    \begin{align}
        (a): w_{\mathbf{O_1^l}}&=(e_3^1,e_1^1), &(b): w_{\mathbf{O_1^l}}&=(e_3^1,e_1^1)\\
        w_{\mathbf{R_1}}&=(e_1^2,e_2^1,e_3^2), & w_{\mathbf{R_1}}&=(e_1^2,e_3^2,e_2^1)\\
        w_{3.1}&=(e_2^2).& w_{3.1}&=(e_2^2). 
    \end{align}

Note that $w_{3.1}$ is identical in both cases because $|E_{\mathbb{\Gamma}}(3.1)|=1$.
\end{example}

\begin{example}
Here is an example of a non-principal splitting of a double admissible graph, see Figure \ref{fig11}.
\begin{figure}[ht]
    \centering 
    \begin{tikzpicture}[scale=0.45, every node/.style={transform shape},
dot/.style = {circle, fill, minimum size=#1, inner sep=0pt, outer sep=5pt}]
%%%%%%%%%%%%%%%%%%%%%%%%%%%%%%%% preliminaries %%%%%%%%%%%%%%%%%%%%%%%%%%%%%%%%
\def\radd{7pt}% radius of vertices
%initialize vertices, so that we can draw colored arrows on the background
%block 1
\coordinate (center1) at (-3,0);
\path (center1) ++(-1.5,0.2) coordinate (1_1);
\path (center1) ++(1.5,0.2) coordinate (1_2);

%block 2
\coordinate (center2) at (3.5,3);
\path (center2) ++(-1.5,0.2) coordinate (2_1);
\path (center2) ++(1.5,0.2) coordinate (2_2);

%block 3
\coordinate (center3) at (10,-1);
\path (center3) ++(-1.5,0.2) coordinate (3_1);
\path (center3) ++(1.5,0.2) coordinate (3_2);

%left
\coordinate (left) at (-1,-5);

%right
\coordinate (right) at (8,-5);
%%%%%%%%%%%%%%%%%%%%%%% colored arrows corresponding to splitting of the graph %%%%%%%%%%%%%%%%%%%%%%%
\draw[blue!20!green, line width=3.0pt] (left) -- (3_2);
\draw[blue!20!green, line width=3.0pt] (3_2) -- (right);
\draw[blue!20!green, line width=3.0pt] (right) -- (2_1);
\draw[blue!20!green, line width=3.0pt] (2_1) -- (left);
\draw[blue!20!green, line width=3.0pt] (left) -- (1_2);
\draw[blue!20!green, line width=3.0pt] (1_2) -- (right);

\draw[blue!10!red, line width=3.0pt, opacity=1] (right) -- (3_1);
\draw[blue!10!red, line width=3.0pt, opacity=1] (3_1) -- (left);

\draw[blue!45!cyan, line width=3.0pt] (right) -- (1_1);
\draw[blue!45!cyan, line width=3.0pt] (1_1) -- (left);
\draw[blue!45!cyan, line width=3.0pt] (left) -- (2_2);
\draw[blue!45!cyan, line width=3.0pt] (2_2)-- (right);
%%%%%%%%%%%%%%%%%%%%%%%%%%%%%%%% blocks and vertices %%%%%%%%%%%%%%%%%%%%%%%%%%%%%%%%
%block 1
\node[dot=\radd ,label=$1.1$] (p1_1) at (1_1) {};
\node[dot=\radd, label=$1.2$] (p1_2) at (1_2) {};

%block 2
\node[dot=\radd, label=$2.1$] (p2_1) at (2_1) {};
\node[dot=\radd, label=$2.2$] (p2_2) at (2_2) {};

%block 3
\node[dot=\radd,label=$3.1$] (p3_1) at (3_1) {};
\node[dot=\radd,label=$3.2$] (p3_2) at (3_2) {};

%left
\node[dot=\radd, label=below :$\mathbf{L_1}$] (l) at (left) {};

%right
\node[dot=\radd,label=below :$\mathbf{R_1}$] (r) at (right) {};

%%%%%%%%%%%%%%%%%%%%%%%%%%%%%%%% edges %%%%%%%%%%%%%%%%%%%%%%%%%%%%%%%%
%edges from block 1
\draw[-{Stealth[scale=1]}, black] (p1_1.center) -- (l) node[midway, below left] {$e_1^1$};
\draw[-{Stealth[scale=1]}, black] (p1_2.center) -- (r) node[near start, above right] {$e_1^2$};
\draw[-{Stealth[scale=1]}, black, dashed] (l.center) -- (p1_2);
\draw[-{Stealth[scale=1]}, black, dashed] (r.center) -- (p1_1) ;

%edges from block 2
\draw[-{Stealth[scale=1]}, black] (p2_1.center) -- (l) node[near start, below right] {$e_2^1$};
\draw[-{Stealth[scale=1]}, black] (p2_2.center) -- (r) node[midway, above right] {$e_2^2$};
\draw[-{Stealth[scale=1]}, black, dashed] (r.center) -- (p2_1);
\draw[-{Stealth[scale=1]}, black, dashed] (l.center) -- (p2_2);

%edges from block 3
\draw[-{Stealth[scale=1]}, black] (p3_1.center) -- (l) node[near start, above left] {$e_3^1$};
\draw[-{Stealth[scale=1]}, black] (p3_2.center) -- (r) node[midway, below right] {$e_3^2$};
\draw[-{Stealth[scale=1]}, black, dashed] (l.center) -- (p3_2);
\draw[-{Stealth[scale=1]}, black, dashed] (r.center) -- (p3_1);

\end{tikzpicture}
    \caption{An example of a non-principal splitting $\mathcal{S}$ of a double admissible graph from $\amsmathbb{G}_3^{(l,r)}$ for $l=(1,0)$ and $r=(1,0)$. Two paths from $P(\mathcal{S})$ are highlighted in green and red. The only cycle is in blue.}
    \label{fig11}
\end{figure}

\begin{align}
    &\mathcal{S}=(\tau,\overline{\tau},\gamma), &&\tau=(\phantom{}^*\!e_1^1,e_1^2,\phantom{}^*\!e_2^2, e_2^1,\phantom{}^*\!e_3^1,e_3^2)=(\mathbf{L_1},1.2,\mathbf{R_1},2.1,\mathbf{L_1},3.2,\mathbf{R_1}),\\
    & &&\overline{\tau}=(\phantom{}^*\!e_3^2,e_3^1)=(\mathbf{R_1},3.1,\mathbf{L_1}),\\
    & &&\gamma=(\phantom{}^*\!e_2^1,e_2^2,\phantom{}^*\!e_1^2,e_1^1)=(\mathbf{L_1},2.2,\mathbf{R_1},1.1,\mathbf{L_1}).
\end{align}
Then $w_{\mathbf{R_1}}=(e_3^2)$, hence $d_{\mathcal{S}}(\mathbf{R_1})=1$, but $|E_{\mathbb{\Gamma}}(\mathbf{R_1})|=3$.
\end{example}

\begin{proposition}\label{prop2}
    The algorithm on page \pageref{algorithm} defines a bijection between $\operatorname{split}(\mathbb{\Gamma})$ and the product of permutation groups $\prod\limits_{v\in \mathbb{\Gamma}} S\left(E_{\mathbb{\Gamma}}(v)\right)$. Particularly, $\operatorname{split}(\Gamma)$ has $\prod\limits_{v\in \mathbb{\Gamma}} |E_{\mathbb{\Gamma}}(v)|!$ elements.
\end{proposition}

Recall that the set $E_{\mathbb{\Gamma}}(v)$ is naturally linearly ordered by the first number of the starting point of every edge, which is a number from $1$ to $n$, so we can identify the group of all permutations $S\left(E_{\mathbb{\Gamma}}(v)\right)$ with the usual symmetric group on $|E_{\mathbb{\Gamma}}(v)|$ letters.

\begin{proof}[Proof of Proposition \ref{prop2}]
By the very definition of principal splittings we have a map from $\operatorname{split}(\mathbb{\Gamma})$ to the product of symmetric groups. Let's denote it by $\varphi$. The inverse map, denoted by $\psi$, is obtained with the help of the inverse algorithm described below. This algorithm may seem a bit involved -- it contains many cases and subcases depending on the vertex considered, but the main idea is simple: for each vertex we invert the step of the algorithm given on page \pageref{algorithm}. We explain how to check that $\psi$ is indeed the inverse of $\varphi$ below.

Assume we have a collection of permutations $\{w_v\}_{v\in E_{\mathbb{\Gamma}}}$ and want to construct a principal splitting $\mathcal{S}$ of $\mathbb{\Gamma}$. To construct a particular path from $\mathcal{S}$ we start applying the inverse algorithm to the initial vertex of the path. This initial vertex is given in advance in the case of paths from $P(\mathcal{S})$ and can be arbitrarily chosen in the case of cycles. We construct the paths from $P(\mathcal{S})$ first and after that we pick any vertex and any unoccupied edge terminating at that vertex to start constructing the first cycle. When it is constructed, we pick another vertex with an unoccupied edged to construct the second cycle and so on.

\textbf{Inverse algorithm.}
Suppose that we already have a path $\tau$, that terminates at a vertex $v$ of the double admissible graph $\mathbb{\Gamma}$. The following rules tells us how to prolong it according to the collection of permutations $\{w_v\}_{v\in \mathbb{\Gamma}}$. When we construct a path from $P(\mathcal{S})$ we consider the case when the given path $\tau$ is of length zero too. But when we construct a cycle, after all $l_1+r_1$ paths from $P(\mathcal{S})$ are constructed, we assume that the length of the given path $\tau$ is strictly greater than zero due to the discussion at the end of the previous paragraph. We have to consider several cases depending on the type of the vertex $v$, which may be numbered, left, or right. Below we identify the set $E_{\mathbb{\Gamma}}(v)$ with $\{1,\ldots,|E_{\mathbb{\Gamma}}(v)|\}$ using the comment after Proposition \ref{prop2}.
\vspace{\baselineskip}

\textit{Left or right vertex $v$}
\begin{enumerate}[label=\arabic*.]
    \item If the length of $\tau$ is zero, we assume that $\tau$ is going to be a path from $P(\mathcal{S})$. Then by definition $\tau$ starts with the reflected edge $\phantom{}^*\!w_v(1)$. 
    
    \item If the length of $\tau$ is greater than zero, then the last edge of $\tau$ is proper and we find such $p$ that this edge equals $w_v(p)$. 
    \begin{enumerate}[label*=\arabic*.]
        \item If $p=|E_{\mathbb{\Gamma}}(v)|$, then the path $\tau$ terminates and we do not prolong it anymore. 
        
        \item If $p<|E_{\mathbb{\Gamma}}(v)|$, then we add the reflected edge $\phantom{}^*\!w_v(p+1)$ to $\tau$.
    \end{enumerate}
\end{enumerate}
\vspace{\baselineskip}

\textit{Loop vertex $v$}
\begin{enumerate}[label=\arabic*.]
    \item If the last edge of $\tau$ is the loop, then we add the reflected edge $\phantom{}^*\!w_v(1)$ to $\tau$. 
    
    \item If the last edge of $\tau$ is proper, then we find such $p$ that this edge equals $w_v(p)$.
    \begin{enumerate}[label*=\arabic*.]
        \item If $p=|E_{\mathbb{\Gamma}}(v)|$, then we add the loop to $\tau$.
        
        \item If $p<|E_{\mathbb{\Gamma}}(v)|$, then we add the reflected edge $\phantom{}^*\!w_v(p+1)$ to $\tau$.
    \end{enumerate}
\end{enumerate}

\vspace{\baselineskip}

\textit{Numbered vertex $v$.}
\begin{enumerate}[label=\arabic*.]
    \item If there are no proper edges terminating at $v$, then there are no reflected edges going out of $v$ and we add the only proper edge going out of $v$ to our path $\tau$.

    \item Suppose there is at least one proper edge terminating at $v$. 
    \begin{enumerate}[label*=\arabic*.]
        \item If the last edge of $\tau$ is reflected, then we add the reflected edge $\phantom{}^*\!w_v(1)$ to $\tau$.

        \item If the last edge of $\tau$ is proper, then we find such $h$ that this last edge equals $w_v(h)$.
        \begin{enumerate}[label*=\arabic*.]
            \item If $h=|E_{\mathbb{\Gamma}}(v)|$, then we add to $\tau$ the only proper edge going out of $v$.

            \item If $h<|E_{\mathbb{\Gamma}}(v)|$, then we add the reflected edge $\phantom{}^*\!w_v(h+1)$ to $\tau$.
        \end{enumerate}
    \end{enumerate}
\end{enumerate} 
We should add two comments to this algorithm.

First, from the given rules its clear under what conditions we stop prolonging paths from $P(\mathcal{S})$. For the cycles we use the following intuitive rule: if the cycle becomes complete at some point of our manipulations, we do not prolong it anymore and proceed with another cycle or terminate the whole procedure.

Second, the result of application of this inverse algorithm does not depend on the order meaning that we first construct paths from $P(\mathcal{S})$ in any order we want and then switch to the cycles. The resulting collection of paths and cycles does not depend on the chosen order. This is explained by the fact that in the inverse algorithm, each edge of the graph $\mathbb{\Gamma}$ which is not an initial or terminal edge of a path from $P(\mathcal{S})$ is addressed exactly twice during construction of a path or a cycle. Once as an incoming edge of a vertex and once as an outgoing edge of a different vertex, but during construction of all other paths and cycles different from that one, this particular edge is not considered at all. Initial or terminal edges of paths from $P(\mathcal{S})$ are addressed only once in the inverse algorithm. From this we also see that all the paths constructed accordingly to the inverse algorithm after the paths from $P(\mathcal{S})$ have been constructed, have to be cyclic, because all the vertices, where a path could terminate, are already occupied by the paths from $P(\mathcal{S})$.

Let us check that $\psi$ is the inverse of $\varphi$. For this we must prove that $\varphi\psi$ and $\psi\varphi$ are the identity maps. The former is a map from the product of symmetric groups indexed by the vertices of the graph $\mathbb{\Gamma}$ to itself. The verification that $\varphi\psi=\mathds{1}$ can be done ''locally'', i.e. for each vertex independently. And this in turn can be easily checked because for each vertex the step of the inverse algorithm is literally the inversion of the step of the algorithm given on page \pageref{algorithm}. The verification of the identity $\psi\varphi=\mathds{1}$, which is an identity for maps from $\operatorname{split}(\mathbb{\Gamma})$ to itself, can be organised as follows. We take $\mathcal{S}\in\operatorname{split}(\mathbb{\Gamma})$ and check that any path $\tau\in P(\mathcal{S})$ coincides with $\psi\varphi(\tau)$. Then we do the same with the cycles. Let us illustrate this for paths. To check this property for any path $\tau\in\mathcal{S}$ we go along that path and compare it to the result of application of the map $\psi\varphi$ to $\tau$. During this we compute $\psi\varphi(\tau)$ step by step according to the inverse algorithm by adding edges to the path we obtained on the previous step. And after we add an edge, we compare the result with the same part of the path $\tau$. And this can be easily done due to the fact that each step of the inverse algorithm is the inversion of the similar step of the algorithm from page \pageref{algorithm}.
\end{proof}

\subsection{Double bidifferential operators \texorpdfstring{$\amsmathbb{B}_{\Gamma}$}{BG}}
\phantom{}

Recall that $A=\Bbbk\langle x_1,\ldots,x_d\rangle$ and $\{\!\!\{-,-\}\!\!\}$ is an arbitrary double Poisson on it.

\subsubsection{Building blocks: \texorpdfstring{$X_{I,\mathcal{S}}(v)$}{XISv}, \texorpdfstring{$X_{I,\mathcal{S}}(\tau)$}{XISt}, \texorpdfstring{$X_{I,\mathcal{S}}(\gamma)$}{XISy}}

Suppose we have 
\begin{itemize}
    \item two pairs of non-negative integers $l=(l_1,l_2)$, $r=(r_1,r_2)$;

    \item $l_1+l_2+r_1+r_2$ elements of $A$, say $a_1,\ldots,a_{l_1}$, $f_1,\ldots,f_{l_2}$, $b_1,\ldots,b_{r_1}$, and $g_1,\ldots,g_{r_2}$;

    \item a double admissible graph $\mathbb{\Gamma}\in\amsmathbb{G}_n^{(l,r)}$;

    \item and a principal splitting $\mathcal{S}\in\operatorname{split}(\mathbb{\Gamma})$.
\end{itemize}

  We will assume that the elements of $A$ are attached to the left, right, and loop vertices of $\mathbb{\Gamma}$, namely, $a$'s are attached to the left vertices, $f$'s are attached to the left loop vertices, $b$'s are attached to the right vertices, and $g$'s are attached to the right loop vertices.
  
  We are going to associate an expression $X_{I,\mathcal{S}}(v)$, which depends on elements of $A$ given above, the principal splitting $\mathcal{S}$ and an arbitrary map $I:E_{\mathbb{\Gamma}}\rightarrow [1,d]$, to each vertex $v$ of the graph $\mathbb{\Gamma}$. We do this in two steps, first we consider preliminary expression $\overline{X}_{I,\mathcal{S}}(v)$ and in terms of it we define $X_{I,\mathcal{S}}(v)$. Recall that we denote the permutation corresponding to the vertex $v$ and obtained from $\mathcal{S}$ with the help of Proposition \ref{prop2} by $w_v$.

  Let us denote by $\partial_k\in\amsmathbb{D}er(A)$ the double derivation of $A$ that acts as $x_i\mapsto\delta_{i,k}(1\otimes 1)$. Then for any $k_1,\ldots,k_m$ we set $\partial_{k_1,\ldots,k_m}:A\longrightarrow A^{\otimes m+1}$, $\partial_{k_1,\ldots,k_m}=\partial_{k_1}^{(1)}\ldots\partial_{k_m}^{(1)}$, where each $\partial_k^{(1)}$ means that we apply this double derivation to the first tensor factor. One can extend this definition in an obvious way replacing the set $k_1,\ldots,k_m$ with any linearly ordered finite set of indices ranging from $1$ to $d$.

  With all this in mind we set
\begin{equation}\label{f13}
    \overline{X}_{I,\mathcal{S}}(v):=\partial_{I(w_v)}\begin{cases}
        a_t,\ \text{if}\ v=\mathbf{L_t}\ \text{for some}\ t=1,\ldots,l_1,\\[4pt]
        f_p,\ \text{if}\ v=\mathbf{O_p^l}\ \text{for some}\ p=1,\ldots,l_2,\\[4pt]
        \left\{\!\!\!\left\{x_{I(e_k^1)},x_{I(e_k^2)}\right\}\!\!\!\right\}',\ \text{if}\ v=k.1\ \text{for some}\ k=1,\ldots,n,\\[5pt]
        \left\{\!\!\!\left\{x_{I(e_k^1)},x_{I(e_k^2)}\right\}\!\!\!\right\}'',\ \text{if}\ v=k.2\ \text{for some}\ k=1,\ldots,n,\\[5pt]
        b_s,\ \text{if}\ v=\mathbf{R_s}\ \text{for some}\ s=1,\ldots,r_1,\\[4pt]
        g_q,\ \text{if}\ v=\mathbf{O_q^r}\ \text{for some}\ q=1,\ldots,r_2.\\
    \end{cases}
\end{equation}

Here we apply operator $\partial_{I(w_v)}$ to expressions, which do not depend on the splitting $\mathcal{S}$, and are associated to the vertices of graph $\mathbb{\Gamma}$. Recall that the permutation $w_v$ may be viewed as a permutational word if we forget the group structure, which in turn is a linearly ordered set. The operator $\partial_{I(w_v)}$ is defined according to it.

In the second and third lines of \eqref{f13} we referred to the first and second Sweedler's components of the double Poisson bracket without referring to the other one. It is an illegal operation -- Sweedler's components are not well-defined separately. However, we will not use the second line of this definition without the third one, and the final expression will make sense. Equivalently, we can say that we fix somehow Sweedler's components and define the expression $\overline{X}_{I,\mathcal{S}}(v)$ with respect to this choice. But it will be clear that the final expression will not depend on this choice.

One easily sees that $\overline{X}_{I,\mathcal{S}}(v)\in A^{\otimes \left(|E_{\mathbb{\Gamma}}(v)|+1\right)}$, so it would be convenient to index Sweedler's components of $\overline{X}_{I,\mathcal{S}}(v)$ by the elements of $E_{\mathbb{\Gamma}}(v)$. This can be done due to that fact that the set $E_{\mathbb{\Gamma}}(v)$ is naturally linearly ordered by starting points of the edges. So, we can write 
\begin{equation}
    \overline{X}_{I,\mathcal{S}}(v)=\overline{X}_{I,\mathcal{S}}^{(0)}(v)\otimes \bigotimes\limits_{e\in E_{\mathbb{\Gamma}}(v)} \overline{X}_{I,\mathcal{S}}^{(e)}(v),
\end{equation}
where the tensor product is taken with respect to the linear order on $E_{\mathbb{\Gamma}}(v)$. Of course, summation is implicitly assumed. One may think of the zero Sweedler's component $\overline{X}_{I,\mathcal{S}}^{(0)}(v)$ as of the vertex $v$ itself, if $v$ is a left or right vertex, as of the loop edge, if the vertex $v$ is a loop vertex, and as of the only reflected edge terminating at the vertex $v$, if $v$ is a numbered vertex. To simplify the future notation we would like to twist the non-zero part of $\overline{X}_{I,\mathcal{S}}(v)$ by the permutation $w_v$. So, we set
\begin{equation}
    X_{I,\mathcal{S}}(v):=\left(\operatorname{id}\otimes w_v\right) \overline{X}_{I,\mathcal{S}}(v)
\end{equation}
and index the Sweedler's components of this new expression by the elements of $E_{\mathbb{\Gamma}}(v)$. 

To any path or cycle in a principal splitting $\mathcal{S}$ of $\mathbb{\Gamma}$ we associate an expression according to the rules listed below. This expression is an element of $A$ in the case of a path, and it is an element of $A_{\natural}$ in the case of a cycle. The rules we use to produce these elements are organised as follows. We go along a path or a cycle and multiply on the right an initial expression by some expressions assigned to the vertices belonging to the path or cycle according to the multiplication rules below.
\vspace{\baselineskip}

\textit{Initial expressions.}
Take any path or any cycle $\tau$ from $\mathcal{S}$. Then its initial expression is $X_{I,\mathcal{S}}^{(0)}(v)$, if $\tau\in P(\mathcal{S})$ and $v$ is the starting point of $\tau$. Suppose now $\tau$ is a cycle. Take any vertex $v$ belonging to $\tau$ and pick any edge belonging to $\tau$ that terminates at $v$. It is important that we fix these vertex and edge and will not change them during the whole procedure -- multiplication rules depend on them. Then the initial expression is 
\begin{itemize}
    \item  $X_{I,\mathcal{S}}^{(0)}(v)$, if \begin{itemize}
        \item $v$ is a numbered vertex and the chosen edge is the only reflected edge terminating at $v$ or

        \item $v$ is a loop vertex and the chosen edge is the loop
    \end{itemize} 

    \item $X_{I,\mathcal{S}}^{(e)}(v)$, if the chosen edge is a proper edge $e$.
\end{itemize}

We have defined the initial expressions for all paths and cycles from $\mathcal{S}$. Now we start multiplying them by certain elements of $A$. 
\vspace{\baselineskip}

\textit{Multiplication rules.}
Take a path or a cycle from $\mathcal{S}$. We go along the fist edge of this path or cycle and multiply the initial expression on the right by a certain element of $A$ related to the endpoint of this edge. Then we go along the second edge and multiply the resulting expression by a certain expression related to the endpoint of this second edge and so on. In the case of a cycle we start this walking at the vertex chosen during the definition of the initial expression of the cycle. Suppose at some point of this procedure we came to a vertex $v$. Then we multiply the expression we obtained after the previous step by 
 \begin{itemize}
    \item $X_{I,\mathcal{S}}^{(0)}(v)$, if \begin{itemize}
        \item $v$ is a numbered vertex and we came to $v$ along the only reflected edge terminating at $v$ or

        \item $v$ is a loop vertex and we came to $v$ along the loop
    \end{itemize} 

    \item $X_{I,\mathcal{S}}^{(e)}(v)$, if we came to $v$ along a proper edge $e$.
\end{itemize}
 
 After we walked along the whole path or cycle we obtain an element of $A$, which is given as a product of several terms, the number of which equals the number of vertices in the path or cycle. In the case of a cycle, we also apply the natural projection $A\rightarrow A_{\natural}$, $a\mapsto\overline{a}$ to the resulting expression we obtained after we walked along the whole cycle. For any $\tau\in\mathcal{S}$ we denote the result of these manipulations by $X_{I,\mathcal{S}}(\tau)$. So, $X_{I,\mathcal{S}}(\tau)\in A$, if $\tau$ starts at a left or right vertex, and $X_{I,\mathcal{S}}(\tau)\in A_{\natural}$, if $\tau$ is a cycle. It is clear that $X_{I,\mathcal{S}}(\tau)$ does not depend on the vertex and edge chosen during the definition of initial expressions in the case when $\tau$ is a cycle due to projection $A\rightarrow A_{\natural}$. Note that, if there are no proper edges terminating at a left, right, or loop vertex, then the corresponding path or cycle $\tau$ is trivial, but the expression $X_{I,\mathcal{S}}(\tau)$ is not trivial, it equals the initial expression (projected to $A_{\natural}$ in the case of a loop vertex), which is some $a_t\in A$, $b_s\in A$, $\overline{f_p}\in A_{\natural}$, or $\overline{g_q}\in A_{\natural}$.

\begin{example}
Let us work out in details an example of the principal splitting $\mathcal{S}$ represented in Figure \ref{fig9}, see also Figure \ref{fig7} below.
\begin{figure}[ht]
    \centering
    \begin{tikzpicture}[scale=0.45, every node/.style={transform shape},
dot/.style = {circle, fill, minimum size=#1, inner sep=0pt, outer sep=5pt},
emptydot/.style = {circle, fill=white, minimum size=#1, inner sep=0pt, outer sep=5pt}]
%%%%%%%%%%%%%%%%%%%%%%%%%%%%%%%% preliminaries %%%%%%%%%%%%%%%%%%%%%%%%%%%%%%%%
\def\radd{7pt}% radius of vertices
\def\raddd{4pt}% radius of empty vertices
%initialize vertices, so that we can draw colored arrows on the background
%block 1
\coordinate (center1) at (-3,0);
\path (center1) ++(-1.5,0.2) coordinate (1_1);
\path (center1) ++(1.5,0.2) coordinate (1_2);

%block 2
\coordinate (center2) at (3.5,3);
\path (center2) ++(-1.5,0.2) coordinate (2_1);
\path (center2) ++(1.5,0.2) coordinate (2_2);

%block 3
\coordinate (center3) at (10,-1);
\path (center3) ++(-1.5,0.2) coordinate (3_1);
\path (center3) ++(1.5,0.2) coordinate (3_2);

%left
\coordinate (left) at (-1,-5);
\path (left) ++(1.5,0) coordinate (left_loop_1);
\path (left_loop_1) ++(1.5,0) coordinate (left_loop_2);

\node[dot=\radd] (ll1) at (left_loop_1) {};
\node[dot=\radd] (ll2) at (left_loop_2) {};

%right
\coordinate (right) at (8,-5);
\path (right) ++(1.5,0) coordinate (right_loop_1);

\node[dot=\radd] (rl1) at (right_loop_1) {};
%%%%%%%%%%%%%%%%%%%%%%% colored arrows corresponding to splitting of the graph %%%%%%%%%%%%%%%%%%%%%%%
\draw[blue!20!green, line width=3.0pt] (left_loop_1) -- (3_2);
\draw[blue!20!green, line width=3.0pt] (3_2) -- (right);
\draw[blue!20!green, line width=3.0pt] (right) -- (1_1);
\draw[blue!20!green, line width=3.0pt] (1_1) -- (left_loop_1);
\path[blue!20!green, line width=3.0pt, every loop/.style={looseness=15}] (ll1) edge  [in=-131,out=-52, loop] (); 

\draw[blue!45!cyan, line width=3.0pt, opacity=1] (left_loop_1) -- (1_2);
\draw[blue!45!cyan, line width=3.0pt, opacity=1] (1_2) -- (right);
\draw[blue!45!cyan, line width=3.0pt, opacity=1] (right) -- (2_2);
\draw[blue!45!cyan, line width=3.0pt, opacity=1] (2_2) -- (3_1);
\draw[blue!45!cyan, line width=3.0pt, opacity=1] (3_1) -- (left_loop_1);

\draw[yellow!20!orange, line width=3.0pt] (right) -- (3_1);
\draw[yellow!20!orange, line width=3.0pt] (3_1) -- (2_1);
\draw[yellow!20!orange, line width=3.0pt] (2_1) -- (right);

\path[blue!10!red, line width=3.0pt, every loop/.style={looseness=15}] (ll2) edge  [in=-131,out=-52, loop] (); 

\path[yellow!70!orange, line width=3.0pt, every loop/.style={looseness=15}] (rl1) edge  [in=-131,out=-52, loop] (); 
%%%%%%%%%%%%%%%%%%%%%%%%%%%%%%%% blocks and vertices %%%%%%%%%%%%%%%%%%%%%%%%%%%%%%%%
%block 1
\node[dot=\radd ,label=$1.1$] (p1_1) at (1_1) {};
\node[dot=\radd, label=$1.2$] (p1_2) at (1_2) {};

%block 2
\node[dot=\radd, label=$2.1$] (p2_1) at (2_1) {};
\node[dot=\radd, label=$2.2$] (p2_2) at (2_2) {};

%block 3
\node[dot=\radd,label=$3.1$] (p3_1) at (3_1) {};
\node[dot=\radd,label=$3.2$] (p3_2) at (3_2) {};

%left
\node[dot=\radd, label={[label distance=3pt]-90 :$\mathbf{L_1}$}] (l) at (left) {};

\node[dot=\radd, label={[label distance=3pt]-90 :$\mathbf{O_1^l}$}] (ll1) at (left_loop_1) {};

\node[dot=\radd, label={[label distance=3pt]-90 :$\mathbf{O_2^l}$}] (ll2) at (left_loop_2) {};

%right
\node[dot=\radd,label={[label distance=3pt]-90 :$\mathbf{R_1}$}] (r) at (right) {};

\node[dot=\radd, label={[label distance=3pt]-90 :$\mathbf{O_1^r}$}] (rl1) at (right_loop_1) {};

%%%%%%%%%%%%%%%%%%%%%%%%%%%%%%%% edges %%%%%%%%%%%%%%%%%%%%%%%%%%%%%%%%
%edges from block 1
\draw[-{Stealth[scale=1]}, black] (p1_1.center) -- (ll1) node[midway, below left] {$e_1^1$};
\draw[-{Stealth[scale=1]}, black] (p1_2.center) -- (r) node[near start, above right] {$e_1^2$};
\draw[-{Stealth[scale=1]}, black, dashed] (ll1.center) -- (p1_2);
\draw[-{Stealth[scale=1]}, black, dashed] (r.center) -- (p1_1) ;

%edges from block 2
\draw[-{Stealth[scale=1]}, black] (p2_1.center) -- (r) node[near start, below left] {$e_2^1$};
\draw[-{Stealth[scale=1]}, black] (p2_2.center) -- (p3_1) node[midway, above right] {$e_2^2$};
\draw[-{Stealth[scale=1]}, black, dashed] (r.center) -- (p2_2);
\draw[-{Stealth[scale=1]}, black, dashed] (p3_1.center) -- (p2_1);

%edges from block 3
\draw[-{Stealth[scale=1]}, black] (p3_1.center) -- (ll1) node[near start, above left] {$e_3^1$};
\draw[-{Stealth[scale=1]}, black] (p3_2.center) -- (r) node[midway, below right] {$e_3^2$};
\draw[-{Stealth[scale=1]}, black, dashed] (ll1.center) -- (p3_2);
\draw[-{Stealth[scale=1]}, black, dashed] (r.center) -- (p3_1);

%left loops
\path[-{Stealth[scale=1.5]}, black, every loop/.style={looseness=15}] (ll1) edge  [in=-124,out=-51, loop] node[below] {$o_1^l$} (); 

\path[-{Stealth[scale=1.5]}, black, every loop/.style={looseness=15}] (ll2) edge  [in=-124,out=-51, loop] node[below] {$o_2^l$} (); 

%right loops

\path[-{Stealth[scale=1.5]}, black, every loop/.style={looseness=15}] (rl1) edge  [in=-124,out=-51, loop] node[below] {$o_1^r$} ();

%to draw the empty dots above all edges
\node[emptydot=\raddd] at (left_loop_1) {};
\node[emptydot=\raddd] at (left_loop_2) {};
\node[emptydot=\raddd] at (right_loop_1) {};

\end{tikzpicture} 
    \caption{The principal splitting $\mathcal{S}=(\tau,\overline{\tau}, \gamma_1,\gamma_2,\gamma_3,\gamma_4)$ represented in Figure \ref{fig9}. The only non-trivial path $\overline{\tau}$ is highlighted in green and cycles $\gamma_1$ and $\gamma_2$ are highlighted in orange and blue respectively.}
    \label{fig7}
\end{figure}
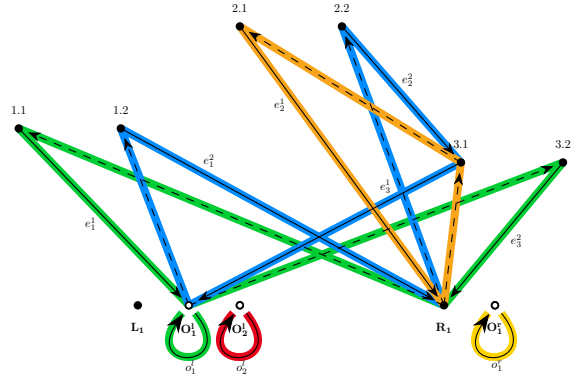

Recall that $\mathcal{S}=(\tau,\overline{\tau}, \gamma_1,\gamma_2,\gamma_3,\gamma_4)$, where
\begin{align}
        &&\tau&=\diameter=(\mathbf{L_1}), &\overline{\tau}&=(\phantom{}^*\!e_1^2, e_1^1,o_1^l,\phantom{}^*\!e_3^1,e_3^2)=(\mathbf{R_1},1.1,\mathbf{O_1^l},\mathbf{O_1^l},3.2,\mathbf{R_1}),\\
        &&\gamma_1&=(\phantom{}^*\!e_3^2,\phantom{}^*\!e_2^2,e_2^1)=(\mathbf{R_1},3.1,2.1,\mathbf{R_1}), &\gamma_2&=(\phantom{}^*\!e_2^1,e_2^2, e_3^1, \phantom{}^*\!e_1^1,e_1^2)=(\mathbf{R_1},2.2,3.1,\mathbf{O_1^l},1.2,\mathbf{R_1}),\\
        &&\gamma_3&=(o_2^l)=(\mathbf{O_2^l},\mathbf{O_2^l}), &\gamma_4&=(o_1^r)=(\mathbf{O_1^r},\mathbf{O_1^r}).
    \end{align}

The corresponding non-trivial permutations are $w_{\mathbf{O_1^l}}=(e_3^1,e_1^1)$, $w_{\mathbf{R_1}}=(e_1^2,e_2^1,e_3^2)$, $w_{3.1}=(e_2^2)$. Note that the last two permutation are identical and the first one is the transposition $(21)$.

Let us choose vertex $3.1$ as a starting point for cycles $\gamma_1$ and $\gamma_2$. This vertex admits a unique incoming edge $\phantom{}^*\!e_3^2$ belonging to $\gamma_1$ and a unique incoming edge $e_2^2$ belonging to $\gamma_2$, so we have no choice in picking edges to define the initial expressions. According to the definition above, we have
\begin{align}
    X_{I,\mathcal{S}}(\tau)&=a_1\in A,\\[5pt]
    X_{I,\mathcal{S}}(\overline{\tau})&=X_{I,\mathcal{S}}^{(0)}(\mathbf{R_1})X_{I,\mathcal{S}}^{(0)}(1.1)X_{I,\mathcal{S}}^{(e_1^1)}(\mathbf{O_1^l})X_{I,\mathcal{S}}^{(0)}(\mathbf{O_1^l})X_{I,\mathcal{S}}^{(0)}(3.2)X_{I,\mathcal{S}}^{(e_3^2)}(\mathbf{R_1})\in A,\\[7pt]
    X_{I,\mathcal{S}}(\gamma_1)&=\overline{X_{I,\mathcal{S}}^{(0)}(3.1)X_{I,\mathcal{S}}^{(0)}(2.1)X_{I,\mathcal{S}}^{(e_2^1)}(\mathbf{R_1})}\in A_{\natural},\\[7pt]
    X_{I,\mathcal{S}}(\gamma_2)&=\overline{X_{I,\mathcal{S}}^{(e_2^2)}(3.1)X_{I,\mathcal{S}}^{(e_3^1)}(\mathbf{O_1^l})X_{I,\mathcal{S}}^{(0)}(1.2)X_{I,\mathcal{S}}^{(e_1^2)}(\mathbf{R_1})X_{I,\mathcal{S}}^{(0)}(2.2)}\in A_{\natural},\\[5pt]
    X_{I,\mathcal{S}}(\gamma_3)&=\overline{f_2}\in A_{\natural},\\[5pt]
    X_{I,\mathcal{S}}(\gamma_4)&=\overline{g_1}\in A_{\natural}.
\end{align}

Let us expand these expressions. By the very definition we have 
\begin{multline}
    X_{I,\mathcal{S}}(\mathbf{O_1^l})=(\operatorname{id}\otimes w_{\mathbf{O_1^l}})\partial_{I(w_{\mathbf{O_1^l}})}(f_1)=(\operatorname{id}\otimes w_{\mathbf{O_1^l}})\partial_{I(e_3^1)}^{(1)}\partial_{I(e_1^1)}(f_1)\\
    =(\operatorname{id}\otimes w_{\mathbf{O_1^l}})\left[\partial_{I(e_3^1)}\left(\partial_{I(e_1^1)}(f_1)'\right)'\otimes \partial_{I(e_3^1)}\left(\partial_{I(e_1^1)}(f_1)'\right)''\otimes \partial_{I(e_1^1)}(f_1)''\right]\\
    =\partial_{I(e_3^1)}\left(\partial_{I(e_1^1)}(f_1)'\right)'\otimes \partial_{I(e_1^1)}(f_1)''\otimes \partial_{I(e_3^1)}\left(\partial_{I(e_1^1)}(f_1)'\right)''.
\end{multline}
So, we set 
\begin{equation}
    \begin{split}
        &X_{I,\mathcal{S}}^{(0)}(\mathbf{O_1^l})=\partial_{I(e_3^1)}\left(\partial_{I(e_1^1)}(f_1)'\right)',\ \ \ \ \  X_{I,\mathcal{S}}^{(e_1^1)}(\mathbf{O_1^l})=\partial_{I(e_1^1)}(f_1)'',\ \ \ \ \ X_{I,\mathcal{S}}^{(e_3^1)}(\mathbf{O_1^l})=\partial_{I(e_3^1)}\left(\partial_{I(e_1^1)}(f_1)'\right)''.
    \end{split}
\end{equation}
Similarly,
\begin{multline}
     X_{I,\mathcal{S}}(\mathbf{R_1})=\partial_{I(w_{\mathbf{R_1}})}(b_1)=\partial_{I(e_1^2)}^{(1)}\partial_{I(e_2^1)}^{(1)}\partial_{I(e_3^2)}(b_1)=\partial_{I(e_1^2)}^{(1)}\left[\partial_{I(e_2^1)}\left(\partial_{I(e_3^2)}(b_1)'\right)\otimes \partial_{I(e_3^2)}(b_1)''\right]\\
     =\partial_{I(e_1^2)}\left(\partial_{I(e_2^1)}\left(\partial_{I(e_3^2)}(b_1)'\right)'\right)\otimes \partial_{I(e_2^1)}\left(\partial_{I(e_3^2)}(b_1)'\right)''\otimes \partial_{I(e_3^2)}(b_1)''
\end{multline}
and
\begin{align}
    &X_{I,\mathcal{S}}^{(0)}(\mathbf{R_1})=\partial_{I(e_1^2)}\left(\partial_{I(e_2^1)}\left(\partial_{I(e_3^2)}(b_1)'\right)'\right)',\ \ \ \ \ \ \ &&X_{I,\mathcal{S}}^{(e_1^2)}(\mathbf{R_1})=\partial_{I(e_1^2)}\left(\partial_{I(e_2^1)}\left(\partial_{I(e_3^2)}(b_1)'\right)'\right)'',\\
    &X_{I,\mathcal{S}}^{(e_2^1)}(\mathbf{R_1})=\partial_{I(e_2^1)}\left(\partial_{I(e_3^2)}(b_1)'\right)'', &&X_{I,\mathcal{S}}^{(e_3^2)}(\mathbf{R_1})=\partial_{I(e_3^2)}(b_1)''.
\end{align}

Next, $X_{I,\mathcal{S}}(3.1)=\partial_{I(e_2^2)}\left(\left\{\!\!\!\left\{x_{I(e_3^1)},x_{I(e_3^2)}\right\}\!\!\!\right\}'\right)$, so
    \begin{align}
        &X_{I,\mathcal{S}}^{(0)}(3.1)=\partial_{I(e_2^2)}\left(\left\{\!\!\!\left\{x_{I(e_3^1)},x_{I(e_3^2)}\right\}\!\!\!\right\}'\right)', &X_{I,\mathcal{S}}^{(e_2^2)}(3.1)=\partial_{I(e_2^2)}\left(\left\{\!\!\!\left\{x_{I(e_3^1)},x_{I(e_3^2)}\right\}\!\!\!\right\}'\right)''.
    \end{align}

For all other vertices the expressions are even more simple: 
\begin{equation}
    X_{I,\mathcal{S}}(k.i)=X_{I,\mathcal{S}}^{(0)}(k.i)=\left\{\!\!\!\left\{x_{I(e_k^1)},x_{I(e_k^2)}\right\}\!\!\!\right\}^{(i)},
\end{equation}
where $k=1,2,3$, $i=1,2$, but $k.i\neq 3.1$ and the superscript $(i)$ on the right-hand side stands for the $i$-th Sweedler's component.

Thus, we arrive at
\begin{multline}
        X_{I,\mathcal{S}}(\overline{\tau})=\partial_{I(e_1^2)}\left(\partial_{I(e_2^1)}\left(\partial_{I(e_3^2)}(b_1)'\right)'\right)'\times \left\{\!\!\!\left\{x_{I(e_1^1)},x_{I(e_1^2)}\right\}\!\!\!\right\}'\times \partial_{I(e_1^1)}(f_1)''\times\partial_{I(e_3^1)}\left(\partial_{I(e_1^1)}(f_1)'\right)'\\
        \shoveright{\times\left\{\!\!\!\left\{x_{I(e_3^1)},x_{I(e_3^2)}\right\}\!\!\!\right\}''\times\partial_{I(e_3^2)}(b_1)''\in A,}\\
        \shoveleft{X_{I,\mathcal{S}}(\gamma_1)=\overline{\partial_{I(e_2^2)}\left(\left\{\!\!\!\left\{x_{I(e_3^1)},x_{I(e_3^2)}\right\}\!\!\!\right\}'\right)'\times\left\{\!\!\!\left\{x_{I(e_2^1)},x_{I(e_2^2)}\right\}\!\!\!\right\}'\times \partial_{I(e_2^1)}\left(\partial_{I(e_3^2)}(b_1)'\right)''}\in A_{\natural},}\\
        \shoveleft{X_{I,\mathcal{S}}(\gamma_2)=\overline{\partial_{I(e_2^2)}\left(\left\{\!\!\!\left\{x_{I(e_3^1)},x_{I(e_3^2)}\right\}\!\!\!\right\}'\right)''\times \partial_{I(e_3^1)}\left(\partial_{I(e_1^1)}(f_1)'\right)''\times \left\{\!\!\!\left\{x_{I(e_1^1)},x_{I(e_1^2)}\right\}\!\!\!\right\}''}}\\
        \overline{\times\partial_{I(e_1^2)}\left(\partial_{I(e_2^1)}\left(\partial_{I(e_3^2)}(b_1)'\right)'\right)''\times \left\{\!\!\!\left\{x_{I(e_2^1)},x_{I(e_2^2)}\right\}\!\!\!\right\}''}\in A_{\natural}.
\end{multline}
\end{example}

\subsubsection{Definition of the double bidifferential operators \texorpdfstring{$\amsmathbb{B}_{\Gamma}$}{BG}}

Now we are going to define preliminary expressions $\amsmathbb{B}_{\mathbb{\Gamma},\mathcal{S}}$, a certain sum of which will be later denoted by $\amsmathbb{B}_{\Gamma}$.

Let $l=(l_1,l_2)$ and $r=(r_1,r_2)$ be pairs of non-negative integers as before. To any principal splitting $\mathcal{S}\in\operatorname{split}(\mathbb{\Gamma})$ of a double admissible graph $\mathbb{\Gamma}\in\amsmathbb{G}_n^{(l,r)}$ we associate a permutation $\sigma\in S(P(\mathcal{S}))$ defined by endpoints of the paths. Namely, any path $\tau\in P(\mathcal{S})$ is mapped by $\sigma$ to the path starting at the endpoint of $\tau$. We will denote this permutation by $\sigma(\mathcal{S})$. Note that $P(\mathcal{S})$ is linearly ordered by starting points, which are in a natural bijection with $\{1,\ldots,l_1+r_1\}$. This bijection identifies left vertices with $\{1,\ldots,l_1\}$ and right vertices with $\{l_1+1,\ldots,l_1+r_1\}$, so $\sigma(\mathcal{S})$ may be viewed as an element of $S(l_1+r_1)$. In this interpretation $\sigma$ maps any number $k$ from $\{1,\ldots,l_1+r_1\}$ to the number corresponding to the endpoint of the path starting at the vertex corresponding to $k$.

 Let us set
 \begin{align}
    &\amsmathbb{B}_{\mathbb{\Gamma},\mathcal{S}}:\Bigl(A^{\otimes l_1}\otimes A^{\otimes l_2}\otimes \Bbbk\left[S(l_1)\right]\Bigr)\otimes \Bigl(A^{\otimes r_1}\otimes A^{\otimes r_2}\otimes \Bbbk\left[S(r_1)\right]\Bigr)\longrightarrow \O(A)^{(l_1+r_1)},\\[5pt] \label{f23}
    &\amsmathbb{B}_{\mathbb{\Gamma},\mathcal{S}}(a\otimes f\otimes u,b\otimes g\otimes v)\\
    &\hspace{70pt}=\sum\limits_{I:E_{\mathbb{\Gamma}}\rightarrow [1,d]}\sigma(\mathcal{S})\cdot\left[\left(\bigotimes\limits_{\tau\in P(\mathcal{S})} X_{I,\mathcal{S}}(\tau)\right)\otimes \left(\prod\limits_{\gamma\in C(\mathcal{S})}X_{I,\mathcal{S}}(\gamma)\right)\otimes (u\times v)\right],
 \end{align}
 where 
 \begin{itemize}    
     \item $a=a_1\otimes \ldots\otimes a_{l_1}\in A^{\otimes l_1}$, $b=b_1\otimes\ldots\otimes b_{r_1}\in A^{\otimes r_1}$;

     \item $u\in S(l_1)$ and $v\in S(r_1)$;

     \item $f=f_1\otimes \ldots\otimes f_{l_2}\in A^{\otimes l_2}$ and $g=g_1\otimes\ldots\otimes g_{r_2}\in A^{\otimes r_2}$;

     \item $\sigma(\mathcal{S})\cdot$ stands for the left action of $S(l_1+r_1)$ on $\O(A)^{(l_1+r_1)}$;

     \item elements $X_{I,\mathcal{S}}$'s are defined according to these $a$'s, $b$'s, $f$'s, and $g$'s;

     \item the tensor product $\bigotimes\limits_{\tau\in P(\mathcal{S})}$ is taken over the totally ordered set $P(\mathcal{S})$;

     \item the symbol $\prod\limits_{\gamma\in C(\mathcal{S})}$ stands for the product in the symmetric algebra $\operatorname{S}^{\bullet}(A_{\natural})$.
 \end{itemize}

 Note that the product $\prod\limits_{\gamma\in C(\mathcal{S})}$ may be empty, in this case it is equal to the unit in the symmetric algebra by definition. 

Now we are ready to define the double bidifferential operators $\amsmathbb{B}_{\Gamma}$. They are supposed to be maps from $\O(A)\otimes\O(A)$ to $\O(A)$. Here we define them simply as

\begin{align}
&\amsmathbb{B}_{\Gamma}:\Bigl(\bigoplus\limits_{l_1\geq 0}A^{\otimes l_1}\otimes \operatorname{T}(A)\otimes \Bbbk\left[S(l_1)\right]\Bigr)\otimes \Bigl(\bigoplus\limits_{r_1\geq 0} A^{\otimes r_1}\otimes \operatorname{T}(A)\otimes \Bbbk\left[S(r_1)\right]\Bigr)\longrightarrow \O(A),\\ \label{f53}
&\phantom{aaaaaaaaaaaaaaaaaaaa}\amsmathbb{B}_{\Gamma}=\sum\limits_{l,r\in \amsmathbb{Z}^2_{\geq 0}}\ \sum\limits_{\substack{\mathbb{\Gamma}\in \amsmathbb{G}_n^{(l,r)}:\\ m_n^{(l,r)}(\mathbb{\Gamma})=\Gamma}}\ \ \sum\limits_{\mathcal{S}\in\operatorname{split}(\mathbb{\Gamma})} \amsmathbb{B}_{\mathbb{\Gamma},\mathcal{S}},
\end{align}
where $\operatorname{T}(A)$ stands for the tensor algebra of the vector space $A$. There is an infinite sum on the right-hand side of \eqref{f53}, which is essentially finite when applied to any element of the domain of definition of $\amsmathbb{B}_{\Gamma}$.

If we want $\amsmathbb{B}_{\Gamma}$ to be a map $\O(A)\otimes\O(A)\rightarrow \O(A)$, then the domain of definition of $\amsmathbb{B}_{\Gamma}$ is almost correct except the factors $\operatorname{T}(A)$, which should be $\operatorname{S}(A_{\natural})$. We are going to prove that $\amsmathbb{B}_{\Gamma}$ actually factors through $\O(A)\otimes\O(A)$. 

\begin{proposition}\label{prop17}
    For any admissible graph $\Gamma\in G_n$ the linear map $\amsmathbb{B}_{\Gamma}$ defined by \eqref{f53} gives rise to a well-defined map $\O(A)\otimes \O(A)\rightarrow \O(A)$. Moreover, for any homogeneous elements $\a,\b\in\O(A)$, a natural number $N$, and tuples of integers $\mathbf{i_1},\mathbf{j_1}, \mathbf{i_2}, \mathbf{j_2}$ ranging from $1$ to $N$ of appropriate lengths the following equality holds
    \begin{equation}\label{f55}
        B_{\Gamma}(\a_{\mathbf{i_1}\mathbf{j_1}},\b_{\mathbf{i_2}\, \mathbf{j_2}})=\amsmathbb{B}_{\Gamma}(\a,\b)_{\mathbf{i_1}\sqcup\mathbf{i_2},\,\mathbf{j_1}\sqcup \mathbf{j_2}},
    \end{equation}
    where $B_{\Gamma}$ is the bidifferential operator from Kontsevich's universal formula, see Section 2 in \cite{kontsevich2003deformation}.
\end{proposition}
\begin{proof}
    Note that $\O(A)$ is a quotient of $\bigoplus\limits_{n\geq 0}A^{\otimes n}\otimes \operatorname{T}(A)\otimes \Bbbk\left[S(n)\right]$ by the two-sided ideal generated by $\mathds{1}\otimes (a\otimes b-b\otimes a)\otimes \mathds{1}$ and $\mathds{1}\otimes (ab-ba)\otimes \mathds{1}$ for $a,b\in A$. Then we must prove that if at least one of the arguments of $\amsmathbb{B}_{\Gamma}$ belongs to this ideal, then the result is zero. Let's do that for the first argument. Note that the ideal is graded. Then take any its homogeneous element, say $\widetilde{a}$, and any element $\widetilde{b}\in A^{\otimes n}\otimes \operatorname{T}(A)\otimes \Bbbk\left[S(n)\right]$ for some $n$. Then by the very definition of $\amsmathbb{B}_{\Gamma}$ the expression $\amsmathbb{B}_{\Gamma}(\widetilde{a},\widetilde{b})_{\mathbf{i_1}\sqcup\mathbf{i_2},\,\mathbf{j_1}\sqcup \mathbf{j_2}}$ equals the right-hand side of the expression from Corollary \ref{cor2} below, which we will prove independently of course. Then by Corollary \ref{cor2} we have 
    \begin{equation}\label{f54}
        \amsmathbb{B}_{\Gamma}(\widetilde{a},\widetilde{b})_{\mathbf{i_1}\sqcup\mathbf{i_2},\,\mathbf{j_1}\sqcup \mathbf{j_2}}=B_{\Gamma}(\widetilde{a}_{\mathbf{i_1}\mathbf{j_1}},\widetilde{b}_{\mathbf{i_2} \, \mathbf{j_2}}),
    \end{equation}
where on the right-hand side we apply the indices to elements of $\bigoplus\limits_{n\geq 0}A^{\otimes n}\otimes \operatorname{T}(A)\otimes \Bbbk\left[S(n)\right]$ exactly in the same way as in the case of $\O(A)$, that is we apply $\tr$ to the tensor algebra component. But the right-hand side of \eqref{f54} is zero, since $\widetilde{a}_{\mathbf{i_1}\mathbf{j_1}}=0$ due to the fact that $\widetilde{a}$ belongs to the ideal mentioned above. Hence $\amsmathbb{B}_{\Gamma}(\widetilde{a},\widetilde{b})_{\mathbf{i_1}\sqcup\mathbf{i_2},\,\mathbf{j_1}\sqcup \mathbf{j_2}}=0$ for any tuples of integers $\mathbf{i_1},\mathbf{i_2},\mathbf{j_1}, \mathbf{j_2}$. Then by Proposition \ref{prop5} we have $\amsmathbb{B}_{\Gamma}(\widetilde{a},\widetilde{b})=0\in\O(A)$. Thus, $\amsmathbb{B}_{\Gamma}:\O(A)\otimes\O(A)\rightarrow\O(A)$ is a well-defined map. Then \eqref{f55} is just a trivial consequence of Corollary \ref{cor2}.    
\end{proof}

\begin{remark}
    From Theorem \ref{th3} and Proposition \ref{prop5} it follows that $\amsmathbb{B}_{\Gamma}:\O(A)\otimes\O(A)\rightarrow\O(A)$ is an admissible morphism of $\S$-modules. 
\end{remark}

\subsection{Definition of the di-twisted star-product \texorpdfstring{$\dstar$}{} on \texorpdfstring{$N\amsmathbb{A}^d$}{NAd}}
\phantom{a}

Let us define $\dstar:\O(A)[[\hbar]]\stimes \O(A)[[\hbar]]\longrightarrow \O(A)[[\hbar]]$ by
\begin{equation}\label{f56}
    \a\dstar \b=\sum\limits_{n=0}^{\infty}\hbar^n \sum\limits_{\Gamma\in G_n} w_{\Gamma}\amsmathbb{B}_{\Gamma}(\a,\b), \ \ \ \ \a,\b\in \O(A),
\end{equation}
where $w_{\Gamma}$ are coefficients from Kontsevich's universal formula, see \cite{kontsevich2003deformation}, and extend this definition to the whole $\O(A)[[\hbar]]\otimes \O(A)[[\hbar]]$ by $\Bbbk[[\hbar]]$-bilinearity.

Let us denote by $\{-,-\}$ the canonical extension of the double Poisson bracket $\{\!\!\{-,-\}\!\!\}$ on $A$ to a di-twisted Poisson bracket on $\O(A)$ provided by Remark \ref{rem2}. Let us denote the induced Poisson bracket on $\orep$ by $\{-,-\}_N$.

\begin{theorem}
    The map $\dstar:\O(A)[[\hbar]]\stimes \O(A)[[\hbar]]\longrightarrow \O(A)[[\hbar]]$ defined by \eqref{f56} is a di-twisted deformation on $\big(\O(A),\{-,-\}\big)$ in the sense of Definition \ref{def7}. Moreover, for any $N$, any homogeneous elements $\a,\b\in\O(A)$ and $X,Y\in\operatorname{T}(\Mat)$ such that $|\a|=|X|$ and $|\b|=|Y|$, one has
        \begin{equation}\label{f59}
            (\a|X)\star_N(\b|Y)=\big(\a\dstar\b\, \big|\, X\otimes Y\big),
        \end{equation}
        where $\star_N$ stand for the Kontsevich's universal star-product from \cite{kontsevich2003deformation} computed for $\{-,-\}_N$.
\end{theorem}
\begin{proof}
    Note that \eqref{f59} immediately follows from Proposition \ref{prop17}. Associativity of $\dstar$ follows then from \eqref{f59}, Proposition \ref{prop5}, and Theorem 2.3 from \cite{kontsevich2003deformation}. Identities from Definition \ref{def7} follow from similar identites for $\star_N$, \eqref{f59}, and Proposition \ref{prop5}.
\end{proof}

\subsection{Proof of Corollary \ref{cor2}}
Corollary \ref{cor2} is formulated and proved below. We shall start with some preparation. 

We keep assuming that $A=\Bbbk\langle x_1,\ldots,x_d\rangle$ and $\{\!\!\{-,-\}\!\!\}$ is a double Poisson bracket on $A$.

By $\amsmathbb{D}er(A)$ we denote the double derivations of $A$ with respect to the outer bimodule structure; it is an $A$-bimodule with respect to the inner bimodule structure and this bimodule is freely generated by $\partial_k$ for $k=1,\ldots,d$, which acts as $x_i\mapsto \delta_{i,k}(1\otimes 1)$.

Recall that any double derivation $\Delta\in\amsmathbb{D}er(A)$ determines a matrix of derivations $\left(\Delta_{ij}\right)_{ij=\overline{1,N}}$ of $\mathcal{O}(\operatorname{Rep}(A,N))$ by the rule 
\begin{equation}
    \Delta_{ij}(a_{kl}):=(\Delta(a))^{kj}_{il}.
\end{equation}

We will identify $\mathcal{O}(\operatorname{Rep}(A,N))=\Bbbk[x_{kij}]_{\substack{k=1,\ldots,d\\ i,j=1,\ldots,N}}$, where $x_{kij}=(x_k)_{ij}$. Denote by $\partial_{(k,p,q)}$ the derivation corresponding to $x_{kpq}$, i.e. $\partial_{(k,p,q)}=\frac{\partial}{\partial x_{kpq}}$. One easily sees that 
    \begin{equation}\label{f1}
        \partial_{(k,p,q)}=(\partial_k)_{qp},
    \end{equation}
    where $\partial_k$ on the right-hand side stands for the double derivation. Indeed, both sides of equality \eqref{f1} define derivations of $\mathcal{O}(\operatorname{Rep}(A,N))$ and they agree on the generators. 

    Now we would like to obtain an analog of \eqref{f1} for the product of $m>1$ derivations of the form $\partial_{(k,p,q)}$.

    We will often view elements of $S(m)$ as permutational words, which are words in $1,\ldots, m$ such that each letter is used exactly once. For any $w\in S(m)$ we set $i(w)=(i(w)_1,\ldots,i(w)_{m})$, where $i(w)_s$ equals the position of $s$ in the permutational word obtained from $w$ by removing all letters less then $s$. For instance,
    $w=(1,2,3,4,5)$, then $i(w)=(1,1,1,1,1)$ or $w=(1,2,5,4,3)$ and $i(w)=(1,1,3,2,1)$. Note that $w\mapsto i(w)$ is a bijective map between $S(m)$ and $m$-tuples of positive integers, where the first number is taken from $\{1,\ldots,m\}$, the second --- from $\{1,\ldots, m-1\}$ and so on. For instance, the last number is always equal to $1$, but the penultimate number is taken from the set $\{1,2\}$.

    We will also need permutational words with letters from an arbitrary linearly ordered finite set, say $X$. We denote them by $S(X)$ and the tuple of integers $i(w)$ is defined similarly.

    For any permutational word $w\in S(m)$ and any m-tuple $(k_1,\ldots, k_m)$ of integer numbers from $1$ to $d$ we define a linear map $\partial^{(w)}_{k_1,\ldots,k_m}:A\rightarrow A^{\otimes m+1}$ by setting
    \begin{equation}
        \partial^{(w)}_{k_1,\ldots,k_m}(a):=\partial_{k_1}^{(i(w)_1)}\ldots \partial_{k_{m}}^{(i(w)_{m})}(a),
    \end{equation}
    where on the right-hand side by $\partial_{k}^{(j)}$ we mean the operator between appropriate tensor powers of $A$ which acts on the $j$-th component of the domain by $\partial_k$ and does not change other tensor factors. So, $\partial_{k}^{(j)}$ increases the number of tensor factors by $1$ and is given by
    \begin{equation}
        \begin{gathered}
            \partial_k^{(j)}:A^{\otimes l}\rightarrow A^{\otimes l+1}\\
            a_1\otimes\ldots\otimes a_l\mapsto a_1\otimes \ldots a_{j-1}\otimes \partial_k(a_j)\otimes a_{j+1}\otimes\ldots a_l.
        \end{gathered}
    \end{equation}

    For example, let $m=3$ and $w=321$. Then one has $i(w)=(3,2,1)$ and
    \begin{equation}
         \partial^{(w)}_{k_1,k_2,k_3}(a)=\partial_{k_1}^{(3)}\partial_{k_2}^{(2)}\partial_{k_3}^{(1)}(a)=\partial_{k_3}(a)'\otimes \partial_{k_2}\left(\partial_{k_3}(a)''\right)'\otimes \partial_{k_1}\left(\partial_{k_2}\left(\partial_{k_3}(a)''\right)''\right)\in A^{\otimes 4},\ \ \ a\in A.
    \end{equation}

Let us set $\partial_{k_1,\ldots,k_2}:=\partial_{k_1}^{(1)}\ldots \partial_{k_m}^{(1)}$. 
    \begin{proposition}\label{prop1}
        One has
        \begin{equation}\label{f3}
        \partial^{(w)}_{k_1,\ldots,k_m}=\partial_{k_{w(1)},\ldots,k_{w(2)}}
        \end{equation}
    \end{proposition}
    \begin{proof}
        It is straightforward to check that
        \begin{equation}\label{f82}
            \partial_{i}^{(k)}\partial_{j}^{(1)}=\partial_{j}^{(1)}\partial_{i}^{(k-1)}
        \end{equation}
        as operators between some tensor powers of $A$, say $A^{\otimes n}\rightarrow A^{\otimes n+2}$, where $n\geq k-1$. Indeed, it's trivial if $k\geq3$, since in this case $\partial_{i}^{(k)}$ and $\partial_{j}^{(1)}$ do not interact at all. If $k=2$, it is sufficient to consider only the case $n=1$, for which \eqref{f82} follows from a simple computation.
        
        Having \eqref{f82} proven, we can proceed with \eqref{f3}. We will prove it by induction on $m$. The base of induction is $m=1$, which is trivial. Let $m>1$. Then 
        \begin{equation}
            \partial_{k_1,\ldots,k_m}^{(w)}=\partial_{k_1}^{(i(w)_1)}\partial_{k_2,\ldots,k_m}^{(\widehat{w})},
        \end{equation}
        where $\widehat{w}\in S([2,m])$ is obtained from $w$ by removing the symbol $1$. By induction hypothesis applied to $\partial_{k_2,\ldots,k_m}^{(\widehat{w})}$ we have
        \begin{equation}
            \partial_{k_1,\ldots,k_m}^{(w)}=\partial_{k_1}^{(i(w)_1)}\partial^{(1)}_{k_{\widehat{w}(2)}}\ldots \partial^{(1)}_{k_{\widehat{w}(m)}}.
        \end{equation}
        Next, we use \eqref{f82} $i(w)_1$ times and obtain
        \begin{equation}
            \partial_{k_1,\ldots,k_m}^{(w)}=\partial^{(1)}_{k_{\widehat{w}(2)}}\ldots \partial_{k_1}^{(1)}\ldots \partial^{(1)}_{k_{\widehat{w}(m)}},
        \end{equation}
        where $\partial_{k_1}^{(1)}$ is placed on the $i(w)_1$-st position. So, the resulting permutational word equals $\widehat{w}(2)\ldots 1 \ldots \widehat{w}(m)$, where $1$ is placed on the $i(w)_1$-st position. This word equals $w$ itself, which completes the proof.
    \end{proof}

    Sometimes it would be convenient for us to write the pair of indices of an element $a_{ij}$ of $\mathcal{O}(\operatorname{Rep}(A,N))$ in parentheses, like $a(ij)$. Moreover, we will use tuples of pairs of indices, and in such cases, it would be crucial for us to write these indices in parentheses but not as a subscript. We will separate the pairs of indices by the symbol $\otimes$ in order to indicate that they correspond to different matrix components. For instance, we write $(i_1j_1)\otimes\ldots \otimes (i_mj_m)$ for $E^*_{i_1j_1}\otimes\ldots\otimes E^*_{i_mj_m}$.

    \begin{lemma}\label{lemma1}
    For any $a\in A$, $k_1,\ldots,k_m\in[1,d]$, and $i,j,p_1,q_1,\ldots, p_m,q_m\in [1,N]$ we have 
        \begin{equation}
            \partial_{(k_1,p_1,q_1)}\ldots \partial_{(k_m,p_m,q_m)}a_{ij}=\sum\limits_{w\in S(m)}\left(\partial_{k_1,\ldots,k_m}^{(w)}(a)\mathlarger{\mathlarger{|}} ip(w)q(w)j\right),
        \end{equation}
        where $ip(w)q(w)j:=(ip_{w(1)})\otimes (q_{w(1)}p_{w(2)})\otimes \ldots \otimes (q_{w(m-1)}p_{w(m)})\otimes (q_{w(m)}j)$, see above for notation.
    \end{lemma}
    \begin{proof}
        Induction on $m$. The base of induction, which is $m=1$, follows from \eqref{f1}. Let us perform the induction step. By the induction hypothesis we have
        \begin{equation}
            \partial_{(k_2,p_2,q_2)}\ldots \partial_{(k_m,p_m,q_m)}a_{ij}=\sum\limits_{w\in S([2,m])}\left(\partial_{k_2,\ldots,k_m}^{(w)}(a)\mathlarger{\mathlarger{|}} ip(w)q(w)j\right),
        \end{equation}
        where $ip(w)q(w)j:=(ip_{w(2)})\otimes (q_{w(2)}p_{w(3)})\otimes \ldots \otimes (q_{w(m-1)}p_{w(m)})\otimes (q_{w(m)}j)$.
        Then we apply the derivation $\partial_{(k_1,p_1,q_1)}$ to each summand of this sum. We need Sweedler's notation $\partial_{k_2,\ldots,k_m}^{(w)}(a)=X^{(1)}\otimes \ldots \otimes X^{(m)}\in A^{\otimes m}$. Then by the Leibniz rule we have
        \begin{align}
            \partial_{(k_1,p_1,q_1)}&\left(\partial_{k_2,\ldots,k_m}^{(w)}(a)\mathlarger{\mathlarger{|}} ip(w)q(w)j\right)\\
            &=\partial_{(k_1,p_1,q_1)}\left[\left(X^{(1)}\mathlarger{|}ip_{w(2)}\right)\left(X^{(2)}\mathlarger{|}q_{w(2)}p_{w(3)}\right)\ldots\left(X^{(m-1)}\mathlarger{|}q_{w(m-1)}p_{w(m)}\right)\left(X^{(m)}\mathlarger{|}q_{w(m)}j\right)\right]\\
            &=\sum\limits_{s=1}^m\left(X^{(1)}\mathlarger{|}ip_{w(2)}\right)\ldots\partial_{(k_1,p_1,q_1)}\left[\left(X^{(s)}\mathlarger{|}q_{w(s)}p_{w(s+1)}\right)\right]\ldots\left(X^{(m)}\mathlarger{|}q_{w(m)}j\right).
        \end{align}
        By \eqref{f1} it equals
        \begin{equation}
            \sum\limits_{s=1}^m\left(X^{(1)}\mathlarger{|}ip_{w(2)}\right)\ldots\\
            \left(\partial_{k_1}(X^{(s)})\mathlarger{|}q_{w(s)}p_1\otimes q_1p_{w(s+1)}\right)\ldots\left(X^{(m)}\mathlarger{|}q_{w(m)}j\right)
        \end{equation}
        which in turn equals
        \begin{equation}
            \sum \left(\partial_{k_1,\ldots,k_m}^{(\widehat{w})}(a)\mathlarger{\mathlarger{|}} ip(\widehat{w})q(\widehat{w})j\right)
        \end{equation}
        where the sum runs over all permutational words $\widehat{w}$ in $m$ letters $\widehat{w}\in S(m)$ such that after removing the letter $1$ we obtain the word $w$. This completes the proof.
    \end{proof}

    Now we would like to obtain an analog of Lemma \ref{lemma1} for the product of $>1$ terms of the form $a_{ij}$. First we rewrite the Leibniz rule in a suitable form, see Lemma \ref{lemma2} below. For this we need a few notation. 
    
    Let $\mathtt{S}$ be a set. For any natural number $m$ we define the set of \textit{m-splittings} of $\mathtt{S}$ as follows
    \begin{equation}
        \operatorname{m-split}(\mathtt{S})=\{(s_1,\ldots,s_m)\mid s_i\subset \mathtt{S},\ s_i\cap s_j=\diameter\ \text{for}\ i\neq j,\ s_1\cup\ldots\cup s_m=\mathtt{S}\}.
    \end{equation}

    \begin{lemma}[Leibniz rule]\label{lemma2}
        Let $\partial_1,\ldots, \partial_l$ be commuting derivations of a commutative algebra and let $\alpha_1,\ldots,\alpha_m$ be elements of this algebra. One has
        \begin{equation}
            \partial_1\ldots \partial_l(\alpha_1\ldots \alpha_m)=\sum\limits_{s\in\operatorname{m-split}([1,l])}\prod_{t=1}^m\left[\left(\prod_{i\in s_t}\partial_i\right)\alpha_t\right],
        \end{equation}
        where on the right-hand side (in the second product) we assume that $s=(s_1,\ldots,s_m)\in\operatorname{m-split}([1,l])$.
    \end{lemma}

    Recall that $\partial_{k_1,\ldots,k_2}:=\partial_{k_1}^{(1)}\ldots \partial_{k_m}^{(1)}$. The next result follows immediately from Lemma \ref{lemma1}, Lemma \ref{lemma2}, and Proposition \ref{prop1}.

    \begin{corollary}\label{cor1}
        For any $a_1,\ldots,a_m\in A$, $i_1,j_1,\ldots,i_m,j_m\in [1,N]$, $k_1,\ldots,k_l\in[1,d]$, and $p_1,q_1,\ldots, p_l,q_l\in [1,N]$ we have
        \begin{equation}
            \partial_{(k_1,p_1,q_1)}\ldots \partial_{(k_l,p_l,q_l)}\left[(a_1)_{i_1j_1}\ldots (a_m)_{i_mj_m}\right]=\!\!\sum\limits_{\substack{s\in\operatorname{m-split}([1,l]) \\ w_1,\ldots,w_m:\, w_t\in S(s_t)}}\!\!\prod_{t=1}^m\left(\partial_{\{k_{w_t(i)}\}_{i\in s_t}}(a_t)\mathlarger{\mathlarger{|}} i_t p(w_t) q(w_t) j_t\right).
        \end{equation}
    \end{corollary}
    Since the set $s_t\subset [1,l]$ is naturally totally ordered, so is $\{k_{w_t(i)}\}_{i\in s_t}$. The expression $\partial_{\{k_{w_t(i)}\}_{i\in s_t}}$ is defined according to this total order, i.e. $\partial_{\{k_{w_t(i)}\}_{i\in s_t}}=\partial_{k_{w_t(1)},\ldots,k_{w_t(|s_t|)}}$.

We have just finished with the preparatory part, and we can now switch to the main object discussed in Corollary \ref{cor2} -- the bidifferential operator $B_{\Gamma}$ from Kontsevich's universal formula. In our case of the representation space $\rep$ 
equipped with a Poisson bracket $\{-,-\}$ coming from a double Poisson bracket $\{\!\!\{-,-\}\!\!\}$ on $A$ the operator $B_{\Gamma}$ is given by
\begin{multline}\label{f12}
    B_{\Gamma}(F,G)=\sum\limits_{\mathcal{I}:E_{\Gamma}\rightarrow [1,d]\times [1,N]^2}\left[\prod_{k=1}^{n}\left(\prod_{e\in E_{\Gamma}(k)}\partial_{\mathcal{I}(e)}\right)\left\{x_{\mathcal{I}(e_k^1)},x_{\mathcal{I}(e_k^2)}\right\}\right]\\
    \times\left(\prod_{e\in E_{\Gamma}(\mathbf{L})}\partial_{\mathcal{I}(e)}\right)F\times \left(\prod_{e\in E_{\Gamma}(\mathbf{R})}\partial_{\mathcal{I}(e)}\right)G
\end{multline}
for any $F,G\in\orep$.

 We will usually write $\mathcal{I}$ as a triple $\mathcal{I}=(I,J,K)$, where $I:E_{\Gamma}\rightarrow [1,d]$ and $J,K:E_{\Gamma}\rightarrow [1,N]$.

By the very definition of the Poisson bracket on $\mathcal{O}(\operatorname{Rep}(A,N))$ we have
\begin{equation}
    \left\{x_{\mathcal{I}(e_k^1)},x_{\mathcal{I}(e_k^2)}\right\}=\left\{\!\!\!\left\{x_{I(e_k^1)},x_{I(e_k^2)}\right\}\!\!\!\right\}_{J(e_k^2),K(e_k^1)\sqcup J(e_k^1),K(e_k^2)}.
\end{equation}

Let $l$ and $r$ be non-negative integers. They correspond to the pairs $(l,0)$ and $(r,0)$ in the sense of our notation related to double admissible graphs. Consider the following elements $F=(a_1)_{i_1j_1}\ldots (a_l)_{i_lj_l}$ and $G=(b_1)_{\overline{i_1}\,\overline{j_1}}\ldots (b_r)_{\overline{i_r}\,\overline{j_r}}$ of $\mathcal{O}(\operatorname{Rep}(A,N))$ with arbitrary $a_1,\ldots,a_l\in A$, $b_1,\ldots,b_r\in A$, and $i_1,\ldots,i_r,j_1,\ldots,j_l\in[1,N]$, $\overline{i_1},\ldots, \overline{i_l},\overline{j_1},\ldots,\overline{j_r}\in[1,N]$. 

Let us temporarily denote $\amsmathbb{G}_n^{((l,0),(r,0))}$ by $\widetilde{\amsmathbb{G}_n^{(l,r)}}$. For any double admissible graph $\mathbb{\Gamma}\in\widetilde{\amsmathbb{G}_n^{(l,r)}}$ and any path $\tau$ belonging to $P(\mathcal{S})$ for a principal splitting $\mathcal{S}$ of $\mathbb{\Gamma}$ we set
\begin{align}\label{f20}
    &i(\tau)=\begin{cases}
        i_t,\ \text{if}\ \tau\ \text{starts at}\ \mathbf{L_t}\ \text{for}\ t=1,\ldots,l,\\
        \overline{i_s},\ \text{if}\ \tau\ \text{starts at}\ \mathbf{R_s}\ \text{for}\ s=1,\ldots,r;
    \end{cases}& j(\tau)=\begin{cases}
        j_t,\ \text{if}\ \tau\ \text{ends at}\ \mathbf{L_t}\ \text{for}\ t=1,\ldots,l,\\
        \overline{j_s},\ \text{if}\ \tau\ \text{ends at}\ \mathbf{R_s}\ \text{for}\ s=1,\ldots,r.
    \end{cases}
\end{align}

\begin{lemma}\label{lemma3}
For $f,g\in\mathcal{O}(\operatorname{Rep}(A,N))$ as above we have
    \begin{equation}\label{f14}
        B_{\Gamma}(F,G)=\sum\limits_{\substack{\mathbb{\Gamma}\in \widetilde{\amsmathbb{G}_n^{(l,r)}}:\\ \widetilde{m_n^{(l,r)}}(\mathbb{\Gamma})=\Gamma}}\ \ \sum\limits_{\mathcal{S}\in\operatorname{split}(\mathbb{\Gamma})}\ \ \sum\limits_{I:E_{\mathbb{\Gamma}}\rightarrow [1,d]}\ \ \prod\limits_{\tau\in P(\mathcal{S})} \left(X_{I,\mathcal{S}}(\tau)\right)_{i(\tau)j(\tau)}\prod\limits_{\gamma\in C(\mathcal{S})}\operatorname{tr}(X_{I,\mathcal{S}}(\gamma)),
\end{equation}
where all $X_{I,\mathcal{S}}(\tau)$'s and $X_{I,\mathcal{S}}(\gamma)$'s are defined according to $a$'s and $b$'s, and by $\widetilde{m_n^{(l,r)}}$ we mean the map $m_n^{((l,0),(r,0))}$ introduced below Definition \ref{def3}.
\end{lemma}

\begin{proof}
    The proof is straightforward. By Corollary \ref{cor1} we have
\begin{equation}
    \left(\prod_{e\in E_{\Gamma}(\mathbf{L})}\partial_{\mathcal{I}(e)}\right)F=\sum\limits_{\substack{P\in\operatorname{l-split}(E_{\Gamma}(\mathbf{L})) \\ w_1,\ldots,w_l:\, w_t\in S(P_t)}}\ \ \ \prod_{t=1}^l\left(\partial_{I(w_t)}(a_t)\mathlarger{\mathlarger{|}} i_t J(w_t) K(w_t) j_t\right),
\end{equation}
where
\begin{multline}
    i_t J(w_t) K(w_t) j_t:=\Bigl(i_t J(w_t(1))\Bigr)\otimes \Bigl(K(w_t(1)) J(w_t(2))\Bigr)\otimes\\
    \ldots \otimes \Bigl(K(w_t(|P_t|-1)) J(w_t(|P_t|))\Bigr)\otimes \Bigl(K(w_t(|P_t|)) j_t\Bigr).
\end{multline}
Here we identified $P_t$ and $\{1,\ldots,|P_t|\}$ as usually.

Similarly,
\begin{equation}
    \left(\prod_{e\in E_{\Gamma}(\mathbf{R})}\partial_{\mathcal{I}(e)}\right)G=\sum\limits_{\substack{Q\in\operatorname{r-split}(E_{\Gamma}(\mathbf{R})) \\ u_1,\ldots,u_r:\, u_r\in S(Q_r)}}\ \ \ \prod_{s=1}^r\left(\partial_{I(u_s)}(b_s)\mathlarger{\mathlarger{|}} \overline{i_s} J(u_s) K(u_s) \overline{j_s}\right).
\end{equation}

Applying Corollary \ref{cor1} once again, we obtain
\begin{multline}
        \left(\prod_{e\in E_{\Gamma}(k)}\partial_{\mathcal{I}(e)}\right)\left\{x_{\mathcal{I}(e_k^1)},x_{\mathcal{I}(e_k^2)}\right\}\\
        \shoveright{=\sum\limits_{\substack{R^k\in\operatorname{2-split}(E_{\Gamma}(k)) \\ v_1^k\in S(R_1^k)\\ v_2^k\in S(R_2^k)}}\ \ \ \left(\partial_{I(v_1^k)}\left(\left\{\!\!\!\left\{x_{I(e_k^1)},x_{I(e_k^2)}\right\}\!\!\!\right\}'\right)\mathlarger{\mathlarger{\mathlarger{|}}} J(e_k^2) J(v_1^k) K(v_1^k) K(e_k^1)\right)}\\
        \times \left(\partial_{I(v_2^k)}\left(\left\{\!\!\!\left\{x_{I(e_k^1)},x_{I(e_k^2)}\right\}\!\!\!\right\}''\right)\mathlarger{\mathlarger{\mathlarger{|}}} J(e_k^1) J(v_2^k) K(v_2^k) K(e_k^2)\right).
\end{multline}

Putting together all the pieces we obtain

\begin{align}\label{f9}
        B_{\Gamma}(F,G)=&\sum\limits_{\substack{I:E_{\Gamma}\rightarrow [1,d] \\ J:E_{\Gamma}\rightarrow[1,N]\\ K:E_{\Gamma}\rightarrow [1,N]}}\ \ \ \sum\limits_{\substack{P\in\operatorname{l-split}(E_{\Gamma}(\mathbf{L})) \\ Q\in\operatorname{r-split}(E_{\Gamma}(\mathbf{R})) \\ w_1,\ldots,w_l:\, w_t\in S(P_t)\\ u_1,\ldots,u_r:\, u_s\in S(Q_s)}}\ \ \ \sum\limits_{\substack{R^k\in\operatorname{2-split}(E_{\Gamma}(k)) \\ v_1^k\in S(R_1^k)\\ v_2^k\in S(R_2^k)\\ \forall k=1,\ldots,n}}\notag\\
        &\phantom{aa}\prod_{k=1}^{n}\left[\left(\partial_{I(v_1^k)}\left(\left\{\!\!\!\left\{x_{I(e_k^1)},x_{I(e_k^2)}\right\}\!\!\!\right\}'\right)\mathlarger{\mathlarger{\mathlarger{|}}} J(e_k^2) J(v_1^k) K(v_1^k) K(e_k^1)\right)\right.\\
        &\phantom{aaaaaaaaaaaaaaaa}\left.\times\left(\partial_{I(v_2^k)}\left(\left\{\!\!\!\left\{x_{I(e_k^1)},x_{I(e_k^2)}\right\}\!\!\!\right\}''\right)\mathlarger{\mathlarger{\mathlarger{|}}} J(e_k^1) J(v_2^k) K(v_2^k) K(e_k^2)\right)\right]\\
        &\phantom{aa}\prod_{t=1}^l\left(\partial_{I(w_t)}(a_t)\mathlarger{\mathlarger{|}} i_t J(w_t) K(w_t) j_t\right)\times \prod_{s=1}^r\left(\partial_{I(u_s)}(b_s)\mathlarger{\mathlarger{|}} \overline{i_s} J(u_s) K(u_s) \overline{j_s}\right).
\end{align}

Now we are going to replace the sum over $P$, $Q$, and all $R^k$'s in \eqref{f9} with a sum over some double admissible graphs. It is possible to do due to a bijection between the following sets
\begin{equation}
     \left\{\mathbb{\Gamma}\subset\widetilde{\amsmathbb{G}_n^{(l,r)}}\ \middle|\  \widetilde{m_n^{(l,r)}}(\mathbb{\Gamma})=\Gamma\right\}
\end{equation}
and
\begin{equation}
    \left\{\left(P,Q,R^1,\ldots,R^n\right)\middle|\  \substack{P\in\operatorname{l-split}(E_{\Gamma}(\mathbf{L})),\\ Q\in\operatorname{r-split}(E_{\Gamma}(\mathbf{R})),}\ \substack{R^k\in\operatorname{2-split}(E_{\Gamma}(k)) \\ \forall k=1,\ldots,n} \right\},
\end{equation}
which is defined by $\mathbb{\Gamma}\mapsto \left(P,Q,R^1,\ldots,R^n\right)$, where
\begin{itemize}
    \item $P=(P_1,\ldots,P_l)$, $P_t$ is the set of all edges of $\mathbb{\Gamma}$ terminating at $\mathbf{L_t}$;

    \item $Q=(Q_1,\ldots,Q_r)$, $Q_s$ is the set of all edges of $\mathbb{\Gamma}$ terminating at $\mathbf{R_s}$;

    \item $R^k=(R^k_1,R^k_2)$, $R^k_1$ is the set of all \textit{proper} edges terminating at $k.1$ and $R^k_2$ is the set of all \textit{proper} edges terminating at $k.2$, for each $k=1,\ldots,n$.
\end{itemize}

Now we shall rewrite formula \eqref{f9}. Below for any permutational word $w\in  S(E_{\mathbb{\Gamma}}(v))$ we denote by $\phantom{}^*\!w\in S(\phantom{}^*\!E_{\mathbb{\Gamma}}(v))$ the induced permutational word which is obtained from $w$ by replacing each proper edge by the corresponding reflected one.

With all these notation in mind we arrive at

\begin{align}\label{f15}
        B_{\Gamma}(F,G)=&\sum\limits_{\substack{\mathbb{\Gamma}\in \widetilde{\amsmathbb{G}_n^{(l,r)}}:\\ \widetilde{m_n^{(l,r)}}(\mathbb{\Gamma})=\Gamma}}\ \ \ \sum\limits_{\substack{w_1,\ldots,w_l:\\ w_t\in S(E_{\mathbb{\Gamma}}(\mathbf{L_t}))\\ u_1,\ldots,u_r:\\ u_s\in S(E_{\mathbb{\Gamma}}(\mathbf{R_s}))}}\ \ \ \sum\limits_{\substack{v_1^k\in S(E_{\mathbb{\Gamma}}(k.1))\\ v_2^k\in S(E_{\mathbb{\Gamma}}(k.2))\\ \forall k=1,\ldots,n}}\ \ \ \sum\limits_{I:E_{\mathbb{\Gamma}}\rightarrow [1,d]}\ \ \ \sum\limits_{J:\amsmathbb{E}_{\mathbb{\Gamma}}\rightarrow[1,N]}\notag\\
        &\phantom{aa}\prod_{k=1}^{n}\left[\left(\partial_{I(v_1^k)}\left(\left\{\!\!\!\left\{x_{I(e_k^1)},x_{I(e_k^2)}\right\}\!\!\!\right\}'\right)\mathlarger{\mathlarger{\mathlarger{|}}} J(\phantom{}^*\!e_k^2) J(\phantom{}^*\!v_1^k) J(v_1^k) J(e_k^1)\right)\right.\\
        &\phantom{aaaaaaaaaaaaaaaa}\times\left.\left(\partial_{I(v_2^k)}\left(\left\{\!\!\!\left\{x_{I(e_k^1)},x_{I(e_k^2)}\right\}\!\!\!\right\}''\right)\mathlarger{\mathlarger{\mathlarger{|}}} J(\phantom{}^*\!e_k^1) J(\phantom{}^*\!v_2^k) J(v_2^k) J(e_k^2)\right)\right]\\
        &\phantom{aa}\prod_{t=1}^l\left(\partial_{I(w_t)}(a_t)\mathlarger{\mathlarger{|}} i_t J(\phantom{}^*\!w_t) J(w_t) j_t\right)\times \prod_{s=1}^r\left(\partial_{I(u_s)}(b_s)\mathlarger{\mathlarger{|}} \overline{i_s} J(\phantom{}^*\!u_s) J(u_s) \overline{j_s}\right).
\end{align}

Note that we not only rewrote formula \eqref{f9} in new notation, but also changed summation indices. Now we do not have K-indices, but we have twice more J-indices. We came to this as follows. First, we swapped J-indices and K-indices. After that we replaced all the K-indices with new J-indices. So, in new notation J-indices attached to proper edges correspond to K-indices in the old notation, and J-indices attached to reflected edges correspond to J-indices in the old notation.

Now we rewrite \eqref{f15} with the help of principal splittings of a double admissible graph $\mathbb{\Gamma}$. Namely, we replace the sum over permutations $\{w_t\}_{t=1,\ldots,l}$, $\{u_s\}_{s=1,\ldots,r}$, $\{v_1^k\}_{k=1,\ldots,n}$, $\{v_2^k\}_{k=1,\ldots,n}$ with the sum over the principal splittings of $\mathbb{\Gamma}$ with the help of Proposition \ref{prop2}. Then we switch to notation \eqref{f13} in order to simplify the subsequent explanations. 

\begin{align}\label{f16}
        B_{\Gamma}(F,G)=&\sum\limits_{\substack{\mathbb{\Gamma}\in \widetilde{\amsmathbb{G}_n^{(l,r)}}:\\ \widetilde{m_n^{(l,r)}}(\mathbb{\Gamma})=\Gamma}}\ \sum\limits_{\mathcal{S}\in\operatorname{split}(\mathbb{\Gamma})}\ \sum\limits_{I:E_{\mathbb{\Gamma}}\rightarrow [1,d]}\ \sum\limits_{J:\amsmathbb{E}_{\mathbb{\Gamma}}\rightarrow[1,N]}\notag\\
        &\phantom{aa}\prod_{k=1}^{n}\left[\left(\overline{X}_{I,\mathcal{S}}(k.1)\mathlarger{\mathlarger{\mathlarger{|}}} J(\phantom{}^*\!e_k^2) J(\phantom{}^*\!w_{k.1}) J(w_{k.1}) J(e_k^1)\right)\right.\\
        &\pushright{\times\left.\left(\overline{X}_{I,\mathcal{S}}(k.2)\mathlarger{\mathlarger{\mathlarger{|}}} J(\phantom{}^*\!e_k^1) J(\phantom{}^*\!w_{k.2}) J(w_{k.2}) J(e_k^2)\right)\right]}\\
        &\phantom{aa}\prod_{t=1}^l\left(\overline{X}_{I,\mathcal{S}}(\mathbf{L_t})\mathlarger{\mathlarger{|}} i_t J(\phantom{}^*\!w_{\mathbf{L_t}}) J(w_{\mathbf{L_t}}) j_t\right)\times \prod_{s=1}^r\left(\overline{X}_{I,\mathcal{S}}(\mathbf{R_s})\mathlarger{\mathlarger{|}} \overline{i_s} J(\phantom{}^*\!w_{\mathbf{R_s}}) J(w_{\mathbf{R_s}}) \overline{j_s}\right).
\end{align}

The only thing that remains to be performed is the summation over $J:\amsmathbb{E}_{\mathbb{\Gamma}}\rightarrow[1,N]$. For this we shall expand each of the factors and convince ourselves that the summation produces exactly the paths from $\mathcal{S}$. In \eqref{f16} we have three products. They can be united in a bigger product which runs over all the vertices of the graph $\mathbb{\Gamma}$. We expand each of the factors and consider two cases depending on the type of the vertex of $\mathbb{\Gamma}$ the factor corresponds to.
\vspace{\baselineskip}

\textit{Left or right vertex.} Let $v=\mathbf{L_t}$ or $v=\mathbf{R_s}$. Set $i=i_t$, $j=j_t$ or $i=\overline{i_s}$, $j=\overline{j_s}$. Then we have

\begin{align}\label{f18}
    &\left(\overline{X}_{I,\mathcal{S}}(v)\mathlarger{\mathlarger{|}} i J(\phantom{}^*\!w_v) J(w_v) j\right)=\left(\left(\operatorname{id}\otimes w_v^{-1}\right)X_{I,\mathcal{S}}(v)\mathlarger{\mathlarger{|}} i J(\phantom{}^*\!w_v) J(w_v) j\right)\notag\\
    &\phantom{aaaaa}=\left(X_{I,\mathcal{S}}^{(0)}(v)\mathlarger{\mathlarger{|}} i J(\phantom{}^*\!w_v(1))\right)\\
    &\phantom{aaaaaaaaaaaa}\times \left(X_{I,\mathcal{S}}^{(w_v(1))}(v)\mathlarger{\mathlarger{|}} J(w_v(1)) J(\phantom{}^*\!w_v(2)) \right)\ldots \left(X_{I,\mathcal{S}}^{(w_v(p-1))}(v)\mathlarger{\mathlarger{|}} J(w_v(p-1)) J(\phantom{}^*\!w_v(p)) \right)\\
    &\pushright{\times\left(X_{I,\mathcal{S}}^{(w_v(p))}(v)\mathlarger{\mathlarger{|}} J(w_v(p)) j\right),}
\end{align}
where $p=|E_{\mathbb{\Gamma}}(v)|$ and we identified $E_{\mathbb{\Gamma}}(v)$ with $\{1,\ldots,p\}$ as linearly ordered sets. 

\vspace{\baselineskip}

\textit{Numbered vertex.} Let $v$ be a numbered vertex $k.1$ or $k.2$. Let us denote by $\phantom{}^*\overline{e}_v$ the only reflected edge terminating at $v$ and by $e_v$ the only proper edge starting at $v$. Then we have

\begin{align}\label{f19}
    &\left(\overline{X}_{I,\mathcal{S}}(v)\mathlarger{\mathlarger{|}} J(\phantom{}^*\overline{e}_v) J(\phantom{}^*\!w_v) J(w_v) J(e_v)\right)=\left(\left(\operatorname{id}\otimes w_v^{-1}\right)X_{I,\mathcal{S}}(v)\mathlarger{\mathlarger{|}} J(\phantom{}^*\overline{e}_v) J(\phantom{}^*\!w_v) J(w_v) J(e_v)\right)\notag\\
    &\phantom{aaaaaaa}=\left(X_{I,\mathcal{S}}^{(0)}(v)\mathlarger{\mathlarger{|}} J(\phantom{}^*\overline{e}_v) J(\phantom{}^*\!w_v(1))\right)\\
    &\phantom{aaaaaaaaaaaa}\times \left(X_{I,\mathcal{S}}^{(w_v(1))}(v)\mathlarger{\mathlarger{|}} J(w_v(1)) J(\phantom{}^*\!w_v(2)) \right)\ldots \left(X_{I,\mathcal{S}}^{(w_v(q-1))}(v)\mathlarger{\mathlarger{|}} J(w_v(q-1)) J(\phantom{}^*\!w_v(q)) \right)\\
    &\pushright{\times\left(X_{I,\mathcal{S}}^{(w_v(q))}(v)\mathlarger{\mathlarger{|}} J(w_v(q)) J(e_v)\right),}
\end{align}
where $q=|E_{\mathbb{\Gamma}}(v)|$.

Let us analyze these expressions. If we multiply all the factors of the form of \eqref{f18} and \eqref{f19}, we obtain a product of $4n+l+r$ terms. Let us denote it by $\Pi$. It depends on many things, like $\mathbb{\Gamma}$, $\mathcal{S}$, $I$ etc. but it is not important for us now. By a trivial calculation of ''J-indices'' corresponding to proper edges in the factors of $\Pi$ we see that there are 
\begin{equation}
    \sum\limits_{\text{left or right}\, v\in\mathbb{\Gamma}}|E_{\mathbb{\Gamma}}(v)|+\sum\limits_{\text{numbered}\, v\in\mathbb{\Gamma}}\left(|E_{\mathbb{\Gamma}}(v)|+1\right)=4n
\end{equation}
possibly coinciding ''J-indices'' corresponding to proper edges and the same number of possibly coinciding ''J-indices'' corresponding to reflected edges. It's readily seen that due to the structure of the factors of $\Pi$ each edge occurs in $\Pi$ at least four times: 
\begin{itemize}
    \item once in a ''J-index'' corresponding to a proper edge standing in the first position e.g., the last term in \eqref{f18};
    
    \item once in a ''J-index'' corresponding to a proper edge standing in the second position e.g., the last term in \eqref{f19};

    \item once in a ''J-index'' corresponding to a reflected edge standing in the first position e.g., the first term in \eqref{f19};

    \item once in a ''J-index'' corresponding to a reflected edge standing in the second position e.g., the first term in \eqref{f18}.
\end{itemize}

So, each edge occurs in the product exactly four times, because the number of edges in $\Gamma$ is $2n$ and the total number of ''J-indices'' is $8n$. Thus, after summation over all ''J-indices'' we obtain with the help of relation \eqref{f17} a product of several terms $l+r$ of which are of the form $(\alpha|i_{t_1}j_{t_2})$ or $(\alpha|\overline{i_{s_1}}\,\overline{j_{s_2}})$, or $(\alpha|\overline{i_{s}}j_t)$, or $(\alpha|i_t\, \overline{j_s})$, while all other terms are of the form $\tr(\beta)=\sum\limits_{i=1}^N\beta_{ii}$; for some $\alpha,\beta\in A$ depending on all $X_{I,\mathcal{S}}$'s. Then one can readily check that the rules according to which we defined elements $X_{I,\mathcal{S}}(\tau)$ mimic the procedure of summation over the ''J-indices'' and application of relation \eqref{f17} -- just look at \eqref{f18} and \eqref{f19}. This completes the proof.
\end{proof}

Now we are going to generalize Lemma \ref{lemma3} to the case $F=\a_{\mathbf{i}\mathbf{j}}$ and $G=\b_{\overline{\mathbf{i}}\, \overline{\mathbf{j}}}$ for arbitrary $\a,\widehat{b}\in\O(A)$.

Let $l=(l_1,l_2)$ and $r=(r_1,r_2)$ be pairs of non-negative integers and $\a\in\O(A)_{l_1}$ and $\b\in\O(A)_{r_1}$ be given by
\begin{align}
    \a&=a_1\otimes\ldots\otimes a_{l_1}\otimes u\otimes f_1\cdot\ldots\cdot f_{l_2}\in\O(A)_{l_1}, &a_1,\ldots,a_{l_1}\in A,\ \ u\in S(l_1),\ \ f_1,\ldots, f_{l_2}\in A_{\natural},\\
    \b&=b_1\otimes\ldots\otimes b_{r_1}\otimes v\otimes g_1\cdot\ldots\cdot g_{r_2}\in \O(A)_{r_1}, &b_1,\ldots,b_{r_1}\in A,\ \ v\in S(r_1),\ \ g_1,\ldots, g_{r_2}\in A_{\natural}.
\end{align}

Let $\mathbf{i}$, $\mathbf{j}$ and $\overline{\mathbf{i}}$, $\overline{\mathbf{j}}$ be $l_1$- and $r_1$-tuples of integers ranging from $1$ to $N$.

For any double admissible graph $\mathbb{\Gamma}\in\amsmathbb{G}_n^{(l,r)}$ and any path $\tau$ belonging to $P(\mathcal{S})$ for a principal splitting $\mathcal{S}$ of $\mathbb{\Gamma}$ we use a notation similar to \eqref{f20}:

\begin{flalign}
    &i(\tau)=\begin{cases}
        i_{u^{-1}(t)},\ \text{if}\ \tau\ \text{starts at}\ \mathbf{L_t}\ \text{for}\ t=1,\ldots,l_1,\\
        \overline{i}_{v^{-1}(s)},\ \text{if}\ \tau\ \text{starts at}\ \mathbf{R_s}\ \text{for}\ s=1,\ldots,r_1;
    \end{cases}& j(\tau)=\begin{cases}
        j_t,\ \text{if}\ \tau\ \text{ends at}\ \mathbf{L_t}\ \text{for}\ t=1,\ldots,l_1,\\
        \overline{j}_s,\ \text{if}\ \tau\ \text{ends at}\ \mathbf{R_s}\ \text{for}\ s=1,\ldots,r_1.
    \end{cases}
\end{flalign}

\begin{corollary}\label{cor2}
With the notation above we have
    \begin{equation}
        B_{\Gamma}(\a_{\mathbf{i}\mathbf{j}},\b_{\overline{\mathbf{i}}\, \overline{\mathbf{j}}})=\sum\limits_{\substack{\mathbb{\Gamma}\in \amsmathbb{G}_n^{(l,r)}:\\ m_n^{(l,r)}(\mathbb{\Gamma})=\Gamma}}\ \ \sum\limits_{\mathcal{S}\in\operatorname{split}(\mathbb{\Gamma})}\ \ \sum\limits_{I:E_{\mathbb{\Gamma}}\rightarrow [1,d]}\ \ \prod\limits_{\tau\in P(\mathcal{S})} \left(X_{I,\mathcal{S}}(\tau)\right)_{i(\tau)j(\tau)}\prod\limits_{\gamma\in C(\mathcal{S})}\operatorname{tr}(X_{I,\mathcal{S}}(\gamma)),
\end{equation}
where all $X_{I,\mathcal{S}}(\tau)$'s and $X_{I,\mathcal{S}}(\gamma)$'s are defined according to $a$'s, $b$'s, $f$'s, and $g$'s.
\end{corollary}
\begin{proof}
    By the very definition we have
    \begin{equation}\label{f27}
        B_{\Gamma}(\a_{\mathbf{i}\mathbf{j}},\b_{\overline{\mathbf{i}}\, \overline{\mathbf{j}}})\\
        =\sum\limits_{n_1,\ldots,n_{l_2}=1}^N\ \ \sum\limits_{m_1,\ldots,m_{r_2}=1}^N \ \ B_{\Gamma}\left(F_{n_1,\ldots,n_{l_2}},G_{m_1,\ldots,m_{r_2}}\right),
    \end{equation}
    where
    \begin{align}
        F_{n_1,\ldots,n_{l_2}}&=(a_1)_{i_{u^{-1}(1)}j_1}\ldots (a_{l_1})_{i_{u^{-1}(l_1)}j_{l_1}}(f_1)_{n_1n_1}\ldots (f_{l_2})_{n_{l_2}n_{l_2}},\\
        G_{m_1,\ldots,m_{r_2}}&=(b_1)_{\overline{i}_{v^{-1}(1)}\ \overline{j}_1}\ldots (b_{r_1})_{\overline{i}_{v^{-1}(r_1)}\overline{j}_{r_1}}(g_1)_{m_1m_1}\ldots (g_{r_2})_{m_{r_2}m_{r_2}}.
    \end{align}

    Now we can apply Lemma \ref{lemma3} to the right-hand side of \eqref{f27}. Then one can easily see that the extra summations in \eqref{f27} correspond exactly to the loops in double admissible graphs. Note that there is a bijection between $\operatorname{split}(\mathbb{\Gamma})$ and $\operatorname{split}(\widetilde{\mathbb{\Gamma}})$, where $\mathbb{\Gamma}$ is as in Corollary \ref{cor2} and $\widetilde{\mathbb{\Gamma}}$ is the same double admissible graph as $\mathbb{\Gamma}$, but without loop vertices, which are replaced by left and right vertices, and all the loop edges are removed. This bijection is given by gluing paths or cycles from $\operatorname{split}(\widetilde{\mathbb{\Gamma}})$ which end and start at the extra left and right vertices, which are going to be the loop vertices in $\mathbb{\Gamma}$. The inverse map exists because we always know where we glued two paths or cycles -- the corresponding loop edge is responsible for that. The bijection yields the desired result.
\end{proof}

\section{Double formality theorem for \texorpdfstring{$N\amsmathbb{A}^d$}{NAd}}\label{section_double_formality}

We set $A=\Bbbk\langle x_1,\ldots,x_d\rangle$ for the rest of the section. We do not distinguish $\ored(A)$ and $\O(A)$ below due to Proposition \ref{prop5}.

\subsection{Differential operators on \texorpdfstring{$N\amsmathbb{A}^d$}{NAd}}

We should rewrite the formula from Lemma \ref{lemma1} in our notation related to $\O(A)$.

Let $k_1,\ldots,k_m=1,\ldots,d$. We set
\begin{align}
    \tildeD_{k_1,\ldots,k_m}&:A\longrightarrow A^{\otimes m+1}\otimes \s \otimes\Bbbk[S(m+1)]\\
    \tildeD_{k_1,\ldots,k_m}(a)&:=\sum\limits_{w\in S(m)}\operatorname{Ad}\big(w^{-1}\times\id_1\big)\Big(\partial_{k_{w^{-1}(1)}}^{(1)}\ldots\partial_{k_{w^{-1}(m)}}^{(1)}(a)\otimes\mathds{1}\otimes (12)^{m,1}\Big).
\end{align}
Recall that  
    \begin{equation}
        \begin{gathered}
            \partial_k^{(1)}:A^{\otimes l}\rightarrow A^{\otimes l+1}\\
            a_1\otimes\ldots\otimes a_l\mapsto \partial_k(a_1)\otimes a_{2}\otimes\ldots a_l.
        \end{gathered}
    \end{equation}

Next, we extend $\tildeD_{k_1,\ldots,k_m}$ to a map $\D_{k_1,\ldots,k_m}:\O(A)\longrightarrow \O(A)[m]$, where $m$ stands for shift of diagonal $\S$-bimodules to the left, as follows
\begin{align}
    \phantom{20pt}&\hspace{-20pt}\D_{k_1,\ldots,k_m}(\a)\\
    &:=\sum\limits_{s\in\operatorname{l+r-split}([1,m])}\left[\pi^r\operatorname{Ad}((12)^{m+l,r}(\sigma_s^{-1}\times \id_{l+r}))\left(\prod\limits_{t=1}^l\tildeD_{\{k_i\}_{i\in s_t}}(a_t)\prod\limits_{t=1}^r\tildeD_{\{k_i\}_{i\in s_{t+l}}}(f_t^*) \right)\right]\cdot (\id_m\times u)
\end{align}
for $\a=(a_1\otimes \ldots\otimes a_l)\otimes  (f_1\ldots f_r)\otimes u\in\O(A)^{(l)}$, where 
\begin{itemize}
    \item $\sigma_s\in S(m)$ is the permutation that maps $E^*_{q_1p_1}\otimes\ldots\otimes E^*_{q_mp_m}$ to $\bigotimes\limits_{t=1}^{l}E^*_{(q_ip_i)_{i\in s_t}}$;

    \item $f_1,\ldots,f_r\in A_{\natural}$ and $f_i^*\in A$ are such that $\overline{f_i^*}=f_i$.
\end{itemize}

This definition does not depend on $f_i^*$ due to the main argument with representation spaces\footnote{Even for different $f^*_t$ we obtain the same expression after passing to representation spaces, so the result does not depend on the choice made as $\R(\Bbbk\langle x_1,\ldots,x_d\rangle)=0$, see Corollary \ref{cor3}} and Lemma \ref{lemma12} below.

Note that if $\a=a\otimes\mathds{1}\otimes\id_1\in\O(A)^{(1)}$, one has $\D_{k_1,\ldots,k_m}(\a)=\tildeD_{k_1,\ldots,k_m}(a)$.

Then next lemma together with Proposition \ref{prop21} show that $\D_{k_1,\ldots,k_m}$ is a differential operator of order $m$ and weight $m$ in the sense of Definition \ref{def9}.

\begin{lemma}\label{lemma12}
Let $k_1,\ldots,k_m=1,\ldots,d$; $\a\in\O(A)^{(n)}$, and $X\in\Mat^{\otimes n}$. Then for any $p_1,\ldots,p_m,q_1,\ldots,q_m\in[1,N]$ one has
    \begin{align}
    \partial_{(k_1,p_1,q_1)}\ldots \partial_{(k_m,p_m,q_m)}(\a|X)=\Big(\D_{k_1,\ldots,k_m}(\a)\,\Big|\,E^*_{\mathbf{q}\mathbf{p}}\otimes X\Big),
\end{align}
where we denoted $E^*_{q_1p_1}\otimes \ldots\otimes E^*_{q_mp_m}$ by $E^*_{\mathbf{q}\mathbf{p}}$.
\end{lemma}
\begin{proof}
Let us start with the case $n=1$. Then by Lemma \ref{lemma1} one has
    \begin{align}
    \hspace{0pt}&\hspace{-0pt}\partial_{(k_1,p_1,q_1)}\ldots \partial_{(k_m,p_m,q_m)}(a)_{ij}\\
    &=\sum\limits_{w\in S(m)}\left(\partial_{k_1,\ldots,k_m}^{(w)}(a)\, \middle|\, ip(w)q(w)j\right)\\
    &=\sum\limits_{w\in S(m)}\left(\partial_{k_{w(1)},\ldots,k_{w(m)}}(a)\, \middle|\, ip(w)q(w)j\right)\tag{Proposition \ref{prop1}}\\
    &=\sum\limits_{w\in S(m)}\left(\partial_{k_{w(1)},\ldots,k_{w(m)}}(a) \, \middle|\, \vect((12)^{m,1})\big(q(w)p(w)\otimes (ij)\big)\right)\\
    &=\sum\limits_{w\in S(m)}\left(\partial_{k_{w(1)},\ldots,k_{w(m)}}(a)\otimes\mathds{1}\otimes (12)^{m,1} \, \middle|\, q(w)p(w)\otimes (ij)\right)\tag{new notation}\\
    &=\sum\limits_{w\in S(m)}\left(\partial_{k_{w^{-1}(1)},\ldots,k_{w^{-1}(m)}}(a)\otimes\mathds{1}\otimes (12)^{m,1} \, \middle|\, \vect(w)\covect(w)E^*_{\mathbf{q}\mathbf{p}}\otimes (ij)\right)\\
    &=\sum\limits_{w\in S(m)}\left((w^{-1}\times\id_1)\cdot\big(\partial_{k_{w^{-1}(1)},\ldots,k_{w^{-1}(m)}}(a)\otimes\mathds{1}\otimes (12)^{m,1}\big)\cdot (w\times \id_1) \, \middle|\, E^*_{\mathbf{q}\mathbf{p}}\otimes (ij)\right)
\end{align}

Recall that $\partial_{k_1,\ldots,k_m}=\partial_{k_1}^{(1)}\ldots\partial_{k_m}^{(1)}$, which finishes the proof of the partial case $n=1$. To pass from it to the general case $n\geq 1$, we use the Leibniz rule from Lemma \ref{lemma2}.

\begin{align}
    \hspace{20pt}&\hspace{-20pt}\partial_{(k_1,p_1,q_1)}\ldots \partial_{(k_m,p_m,q_m)}\big((a_1)_{i_1j_1}\ldots (a_l)_{i_lj_l}\big)\\
    &=\sum\limits_{s\in\operatorname{l-split}([1,m])}\prod_{t=1}^l\left[\left(\prod_{i\in s_t}\partial_{(k_i,p_i,q_i)}\right)(a_t)_{i_tj_t}\right]\\
    &=\sum\limits_{s\in\operatorname{l-split}([1,m])}\prod_{t=1}^l\Big(\tildeD_{\{k_i\}_{i\in s_t}}(a_t)\big)\,\Big|\,(q_ip_i)_{i\in s_t}\otimes (i_tj_t)\Big)\tag{by the partial case $n=1$}\\
    &=\sum\limits_{s\in\operatorname{l-split}([1,m])}\left(\prod\limits_{t=1}^l\tildeD_{\{k_i\}_{i\in s_t}}(a_t)\,\middle|\,\bigotimes\limits_{t=1}^{l} (q_ip_i)_{i\in s_t}\otimes (i_tj_t)\right)\tag{$\prod\limits_{t=1}^l$ is in $\O(A)$}
\end{align}

Using the definition of $\sigma_{s}$, one obtains 
\begin{align}
    \hspace{20pt}&\hspace{-20pt}\partial_{(k_1,p_1,q_1)}\ldots \partial_{(k_m,p_m,q_m)}\big((a_1)_{i_1j_1}\ldots (a_l)_{i_lj_l}\big)\\
    &=\sum\limits_{s\in\operatorname{l-split}([1,m])}\left(\prod\limits_{t=1}^l\tildeD_{\{k_i\}_{i\in s_t}}(a_t)\,\middle|\,\vect(\sigma_{s})\covect(\sigma_{s})(q_1p_1)\otimes\ldots\otimes (q_mp_m)\otimes (i_1j_1)\otimes\ldots \otimes (i_lj_l)\right)\\
    &=\sum\limits_{s\in\operatorname{l-split}([1,m])}\left(\operatorname{Ad}(\sigma_s^{-1}\times \id_l)\left(\prod\limits_{t=1}^l\tildeD_{\{k_i\}_{i\in s_t}}(a_t)\right)\,\middle|\,(q_1p_1)\otimes\ldots\otimes (q_mp_m)\otimes (i_1j_1)\otimes\ldots \otimes (i_lj_l)\right)\\
    &=\left(\sum\limits_{s\in\operatorname{l-split}([1,m])}\operatorname{Ad}(\sigma_s^{-1}\times \id_l)\left(\prod\limits_{t=1}^l\tildeD_{\{k_i\}_{i\in s_t}}(a_t)\right)\,\middle|\,(q_1p_1)\otimes\ldots\otimes (q_mp_m)\otimes (i_1j_1)\otimes\ldots \otimes (i_lj_l)\right)
\end{align}
Next, let $f_1,\ldots,f_r\in A_{\natural}$ and $f_i^*\in A$ are such that $\overline{f_i^*}=f_i$. 

 \begin{align}
    \hspace{20pt}&\hspace{-20pt}\partial_{(k_1,p_1,q_1)}\ldots \partial_{(k_m,p_m,q_m)}\big((a_1)_{i_1j_1}\ldots (a_l)_{i_lj_l}\tr(f_1)\ldots\tr(f_r)\big)\\
    &=\Bigg(\sum\limits_{s\in\operatorname{l+r-split}([1,m])}\operatorname{Ad}(\sigma_s^{-1}\times \id_{l+r})\left(\prod\limits_{t=1}^l\tildeD_{\{k_i\}_{i\in s_t}}(a_t)\prod\limits_{t=1}^r\tildeD_{\{k_i\}_{i\in s_{t+l}}}(f_t^*) \right)\,\Bigg|\,(q_1p_1)\otimes\ldots\otimes (q_mp_m)\otimes\\
    &\pushright{(i_1j_1)\otimes\ldots \otimes (i_lj_l)\otimes \tr^{\otimes r}\Bigg)}\\
    &=\Bigg(\sum\limits_{s\in\operatorname{l+r-split}([1,m])}\operatorname{Ad}(\sigma_s^{-1}\times \id_{l+r})\left(\prod\limits_{t=1}^l\tildeD_{\{k_i\}_{i\in s_t}}(a_t)\prod\limits_{t=1}^r\tildeD_{\{k_i\}_{i\in s_{t+l}}}(f_t^*) \right)\\
    &\hspace{60pt}\,\Bigg|\,\vect((12)^{r,m+l})\covect((12)^{r,m+l})\tr^{\otimes r} \otimes (q_1p_1)\otimes\ldots\otimes (q_mp_m)\otimes(i_1j_1)\otimes\ldots \otimes (i_lj_l)\Bigg)\\
    &=\Bigg(\sum\limits_{s\in\operatorname{l+r-split}([1,m])}\operatorname{Ad}((12)^{m+l,r}(\sigma_s^{-1}\times \id_{l+r}))\left(\prod\limits_{t=1}^l\tildeD_{\{k_i\}_{i\in s_t}}(a_t)\prod\limits_{t=1}^r\tildeD_{\{k_i\}_{i\in s_{t+l}}}(f_t^*) \right)\\
    &\pushright{\,\Bigg|\,\tr^{\otimes r} \otimes (q_1p_1)\otimes\ldots\otimes (q_mp_m)\otimes(i_1j_1)\otimes\ldots \otimes (i_lj_l)\Bigg)}\\
    &=\Bigg(\sum\limits_{s\in\operatorname{l+r-split}([1,m])}\pi^r\operatorname{Ad}((12)^{m+l,r}(\sigma_s^{-1}\times \id_{l+r}))\left(\prod\limits_{t=1}^l\tildeD_{\{k_i\}_{i\in s_t}}(a_t)\prod\limits_{t=1}^r\tildeD_{\{k_i\}_{i\in s_{t+l}}}(f_t^*) \right)\\
    &\pushright{\,\Bigg|\, (q_1p_1)\otimes\ldots\otimes (q_mp_m)\otimes(i_1j_1)\otimes\ldots \otimes (i_lj_l)\Bigg)}
\end{align}

Finally,
\begin{align}
    \hspace{20pt}&\hspace{-20pt}\partial_{(k_1,p_1,q_1)}\ldots \partial_{(k_m,p_m,q_m)}\big((a_1)_{i_{u^{-1}(1)}j_1}\ldots (a_l)_{i_{u^{-1}(l)}j_l}\tr(f_1)\ldots\tr(f_r)\big)\\
    &=\Bigg(\sum\limits_{s\in\operatorname{l+r-split}([1,m])}\pi^r\operatorname{Ad}((12)^{m+l,r}(\sigma_s^{-1}\times \id_{l+r}))\left(\prod\limits_{t=1}^l\tildeD_{\{k_i\}_{i\in s_t}}(a_t)\prod\limits_{t=1}^r\tildeD_{\{k_i\}_{i\in s_{t+l}}}(f_t^*) \right)\\
    &\hspace{180pt}\,\Bigg|\, \vect(\id_m\times u)(q_1p_1)\otimes\ldots\otimes (q_mp_m)\otimes(i_1j_1)\otimes\ldots \otimes (i_lj_l)\Bigg)\\
    &=\Bigg(\sum\limits_{s\in\operatorname{l+r-split}([1,m])}\!\left[\pi^r\operatorname{Ad}((12)^{m+l,r}(\sigma_s^{-1}\times \id_{l+r}))\left(\prod\limits_{t=1}^l\tildeD_{\{k_i\}_{i\in s_t}}(a_t)\prod\limits_{t=1}^r\tildeD_{\{k_i\}_{i\in s_{t+l}}}(f_t^*) \right)\!\right]\!\cdot\! (\id_m\times u)\\
    &\pushright{\,\Bigg|\, (q_1p_1)\otimes\ldots\otimes (q_mp_m)\otimes(i_1j_1)\otimes\ldots \otimes (i_lj_l)\Bigg)}
\end{align}
And the proof is complete.
\end{proof}

Let us extend Lemma \ref{lemma12} to poly-differential operators.

Let $\mathbf{k_1},\ldots,\mathbf{k_n}$ be tuples of integers from $[1,d]$. We will denote the length of $\mathbf{k_i}$ by $\#\mathbf{k_i}$, i.e. $\mathbf{k_i}\in[1,d]^{\#\mathbf{k_i}}$. To each of then we can associate the differential operator $\D_{\mathbf{k_i}}:\O(A)\longrightarrow\O(A)[\#k_i]$ as defined above Lemma \ref{lemma12}.

For any tuples $\mathbf{p_1},\mathbf{q_1},\ldots,\mathbf{p_n},\mathbf{q_n}$ of integers from $1$ to $N$ such that $\#\mathbf{p_i}=\#\mathbf{q_i}=\#\mathbf{k_i}$ we can consider the differential operator $\partial_{(\mathbf{k_i},\mathbf{p_i},\mathbf{q_i})}:\orep\longrightarrow\orep$ of order $\#\mathbf{k_i}$. Then we can consider the $n$-poly-differential operator $\partial_{(\mathbf{k_1},\mathbf{p_1},\mathbf{q_1)}}\otimes\ldots \otimes \partial_{(\mathbf{k_n},\mathbf{p_n},\mathbf{q_n})}:\orep^{\otimes n}\longrightarrow\orep$.

\begin{corollary}
Let $\a_1,\ldots,\a_n\in\O(A)$ be homogeneous elements. With the notation above one has
    \begin{align}
    \hspace{50pt}&\hspace{-50pt}\partial_{(\mathbf{k_1},\mathbf{p_1},\mathbf{q_1)}}\otimes\ldots \otimes \partial_{(\mathbf{k_n},\mathbf{p_n},\mathbf{q_n})}\Big((\a_1|X_1)\otimes\ldots\otimes (\a_n|X_n)\Big)\\
    &=\Big(\operatorname{Ad}(\sigma_{\a,\#\mathbf{k}})\D_{\mathbf{k_1}}(\a_1)\ldots\D_{\mathbf{k_n}}(\a_n)\, \Big|\,E^*_{\mathbf{q_1}\mathbf{p_1}}\otimes\ldots\otimes E^*_{\mathbf{q_n}\mathbf{p_n}}\otimes X_1\otimes\ldots\otimes \otimes X_n \Big),
\end{align}
where $\sigma_{\a,\#\mathbf{k}}\in S(\#\mathbf{k_1}+\ldots+\#\mathbf{k_n}+|\a_1|+\ldots+|\a_n|)$ is the permutation that permutes the blocks of sizes $\#\mathbf{k_1},\ldots,\#\mathbf{k_n}$ and $|\a_1|,\ldots,|\a_n|$ and maps the block decomposition $\#\mathbf{k_1},|\a_1|,\ldots,\#\mathbf{k_n},|\a_n|$ to $\#\mathbf{k_1},\ldots,\#\mathbf{k_1},|\a_1|,\ldots,|\a_n|$.
\end{corollary}
\begin{proof}
    By Lemma \ref{lemma12} the left-hand side equals
    \begin{align}
        \hspace{30pt}&\hspace{-30pt}\Big(\D_{\mathbf{k_1}}(\a_1)\,\Big|\,E^*_{\mathbf{q_1}\mathbf{p_1}}\otimes X_1\Big)\ldots \Big(\D_{\mathbf{k_n}}(\a_n)\,\Big|\,E^*_{\mathbf{q_n}\mathbf{p_n}}\otimes X_n\Big)\\
        &=\Big(\D_{\mathbf{k_1}}(\a_1)\ldots\D_{\mathbf{k_n}}(\a_n)\, \Big|\,E^*_{\mathbf{q_1}\mathbf{p_1}}\otimes X_1\otimes\ldots\otimes E^*_{\mathbf{q_n}\mathbf{p_n}}\otimes X_n \Big)\\
        &=\Big(\operatorname{Ad}(\sigma_{\a,\#\mathbf{k}})\D_{\mathbf{k_1}}(\a_1)\ldots\D_{\mathbf{k_n}}(\a_n)\, \Big|\,E^*_{\mathbf{q_1}\mathbf{p_1}}\otimes\ldots\otimes E^*_{\mathbf{q_n}\mathbf{p_n}}\otimes X_1\otimes\ldots\otimes \otimes X_n \Big).
    \end{align}
\end{proof}

Let us denote the $n$-poly-differential operator $\a_1\otimes\ldots\otimes\a_n\mapsto \operatorname{Ad}(\sigma_{\a,\#\mathbf{k}})\D_{\mathbf{k_1}}(\a_1)\ldots\D_{\mathbf{k_n}}(\a_n)$ by $\D_{\mathbf{k_1},\ldots,\mathbf{k_n}}:\O(A)^{\stimes n}\longrightarrow\O(A)[\#\mathbf{k_1}+\ldots+\#\mathbf{k_n}]$.

\subsection{The \texorpdfstring{$L_{\infty}$}{L}-quasi-isomorphism}
Let us recall the definition of admissible graphs, see Section 6.1 in \cite{kontsevich2003deformation}.

\begin{definition}
    By $G_{n,m}$ we will denote the set of \textit{admissible graphs}, which are by definition labeled graphs with $n+m$ vertices and $2n+m-2\geq 0$ edges without loops or multiple edges, and such that 
    \begin{enumerate}[label=\theenumi),leftmargin=3ex]
        \item the set of vertices is $\{1,\ldots,n\}\sqcup\{\overline{1},\ldots,\overline{m}\}$, vertices from $\{1,\ldots,n\}$ are called \textit{vertices of the first type} and vertices from $\{\overline{1},\ldots,\overline{m}\}$ are called vertices of the second type;

        \item the set of edges is denoted by $E_{\Gamma}$, and every edge starts at a vertex of first type;

        \item for every vertex of first type $k\in\{1,\ldots,n\}$, the set of edges starting at $k$:
        \begin{align}
            star(k):=\{e\in E_{\Gamma}\mid e=(k,*)\}
        \end{align}
        is labeled by symbols $e_k^1,\ldots,e_{k}^{\#star(k)}$.
    \end{enumerate}
     By $E_{\Gamma}(k)$ we will denote the set of edges terminating at vertex $k$, which may also be of second type.
\end{definition}

\begin{definition}
    Let $\Gamma\in G_{n,m}$ be any admissible graph, and $k_i=\#star(i)-1$ for $i=1,\ldots,n$. For any $\Phi_1\in L^{k_1}_{di\mh Pois}(\O(A)),\ldots,\Phi_n\in L^{k_n}_{di\mh Pois}(\O(A))$, and homogeneous $\a_1,\ldots,\a_m\in\O(A)$ we set
    \begin{align}
    &\U_{\Gamma}(\Phi_1,\ldots,\Phi_n)(\a_1,\ldots,\a_m)\\
    &\hspace{8pt}:=\sum\limits_{I:E_{\Gamma}\longrightarrow[1,d]}\pi^{2|E_{\Gamma}|}\left[(\sigma_{\Gamma}^2)^{-1}\cdot\left(\prod\limits_{i=1}^{n}\D_{\{I(e)\}_{e\in E_{\Gamma}(i)}}\Phi_i\Big(x_{I(e_i^1)},\ldots,x_{I(e_i^{k_i+1})}\Big)\prod\limits_{j=1}^{m}\D_{\{I(e)\}_{e\in E_{\Gamma}(\overline{j})}}(\a_j)\right)\cdot \sigma_{\Gamma}^1\right],
\end{align}
where permutations $\sigma_{\Gamma}^1,\sigma_{\Gamma}^2\in S(2|E_{\Gamma}|+|\a_1|+\ldots+|\a_m|)$ are such that for any $J,K:E_{\Gamma}\longrightarrow [1,N]$
\begin{align}
    \hspace{120pt}&\hspace{-120pt}\vect(\sigma_{\Gamma}^1)\covect(\sigma_{\Gamma}^2)\left[\left(\bigotimes\limits_{e\in E_{\Gamma}}J(e)J(e)\otimes K(e)K(e)\right)\otimes X_1\otimes\ldots\otimes X_m\right]\\
    &=\bigotimes\limits_{i=1}^{n}\left[\bigotimes\limits_{e\in E_{\Gamma}(i)}K(e)J(e)\otimes J(e_i)K(e_i)\right]\otimes \bigotimes\limits_{j=1}^{m}\left[\bigotimes\limits_{e\in E_{\Gamma}(\overline{j})}K(e)J(e)\otimes X_j\right],
\end{align}
where $J(e_i)K(e_i):=J(e_i^1)K(e_i^1)\otimes\ldots\otimes J(e_i^{k_i+1})K(e_i^{k_i+1})$ and the total order used for $\bigotimes\limits_{e\in E_{\Gamma}(i)}$ and $\bigotimes\limits_{e\in E_{\Gamma}(\overline{j})}$ is the natural one -- the edges are ordered by their starting points. For $\bigotimes\limits_{e\in E_{\Gamma}}$ we first take all edges going out of vertex $1$, then the edges going out of vertex $2$, and so on. Recall that edges going out of a particular vertex are totally ordered by definition.
\end{definition}
So, $\U_{\Gamma}(\Phi_1,\ldots,\Phi_n)(\a_1,\ldots,\a_m)\in\O(A)^{|\a_1|+\ldots+|\a_m|}$.

\begin{lemma}\label{lemma13}
    For any $\Phi_i\in L^{k_i}_{di\mh Pois}(\O(A))$, one has
    \begin{align}
        \mathcal{U}_{\Gamma}\Big((\Phi_1)_N,\ldots,(\Phi_n)_N\Big)=\Big(\U_{\Gamma}(\Phi_1,\ldots,\Phi_n)\Big)_N,
    \end{align}
    where $(-)_N$ stands for the maps induced by the representation functor, and $\mathcal{U}_{\Gamma}$ is defined in Section 6.3 from \cite{kontsevich2003deformation}.
\end{lemma}
\begin{proof}
    By the very definition of $\mathcal{U}_{\Gamma}$ one has 
\begin{align}
    \hspace{0pt}&\hspace{-0pt}\mathcal{U}_{\Gamma}\Big((\Phi_1)_N,\ldots,(\Phi_n)_N\Big)\big((\a_1|X_1),\ldots,(\a_m|X_m)\big)\\
    &=\sum\limits_{\mathcal{I}:E_{\Gamma}\longrightarrow[1,d]\times [1,N]^2}\prod\limits_{i=1}^{n}\left[\left(\prod_{e\in E_{\Gamma}(i)}\partial_{\mathcal{I}(e)}\right)(\Phi_i)_N\Big(x_{\mathcal{I}(e_i^1)},\ldots,x_{\mathcal{I}(e_i^{k_i+1})}\Big)\right]\prod\limits_{j=1}^m\left[\left(\prod_{e\in E_{\Gamma}(\overline{j})}\partial_{\mathcal{I}(e)}\right)(\a_j|X_j)\right]
\end{align}

We will write $\mathcal{I}$ as a triple $\mathcal{I}=(I,J,K)$, where $I:E_{\Gamma}\rightarrow [1,d]$ and $J,K:E_{\Gamma}\rightarrow [1,N]$.

Then by the very definition of $(\Phi_i)_N$ one has
\begin{align}
    (\Phi_i)_N\Big(x_{\mathcal{I}(e_i^1)},\ldots,x_{\mathcal{I}(e_i^{k_i+1})}\Big)=\Bigg(\Phi_i\Big(x_{I(e_i^1)},\ldots,x_{I(e_i^{k_i+1})}\Big)\,\Bigg|\,J(e_i^1)K(e_i^1)\otimes\ldots\otimes J(e_i^{k_i+1})K(e_i^{k_i+1})\Bigg),
\end{align}
where $E^*_{J(e_i^t)K(e_i^1)}$ is denoted simply by $J(e_i^t)K(e_i^t)$. Then by Lemma \ref{lemma12} one has
\begin{align}
    &\left(\prod_{e\in E_{\Gamma}(i)}\partial_{\mathcal{I}(e)}\right)(\Phi_i)_N\Big(x_{\mathcal{I}(e_i^1)},\ldots,x_{\mathcal{I}(e_i^{k_i+1})}\Big)\\
    &\hspace{30pt}=\Bigg(\D_{\{I(e)\}_{e\in E_{\Gamma}(i)}}\Phi_i\Big(x_{I(e_i^1)},\ldots,x_{I(e_i^{k_i+1})}\Big)\,\Bigg|\,\bigotimes\limits_{e\in E_{\Gamma}(i)}K(e)J(e)\otimes J(e_i)K(e_i)\Bigg),
\end{align}
where we set $K(e_i)J(e_i):=K(e_i^1)J(e_i^1)\otimes\ldots\otimes K(e_i^{k_i+1})J(e_i^{k_i+1})$ for brevity. Recall that the set $E_{\Gamma}(i)$ is naturally totally ordered by the starting points of edges. We use the same total order in $\{I(e)\}_{e\in E_{\Gamma}(i)}$ and $\bigotimes\limits_{e\in E_{\Gamma}(i)}$.

Next,

\begin{align}
    &\prod\limits_{i=1}^{n}\left[\left(\prod_{e\in E_{\Gamma}(i)}\partial_{\mathcal{I}(e)}\right)(\Phi_i)_N\Big(x_{\mathcal{I}(e_i^1)},\ldots,x_{\mathcal{I}(e_i^{k_i+1})}\Big)\right]\\
    &\hspace{30pt}=\left(\prod\limits_{i=1}^{n}\D_{\{I(e)\}_{e\in E_{\Gamma}(i)}}\Phi_i\Big(x_{I(e_i^1)},\ldots,x_{I(e_i^{k_i+1})}\Big)\,\Bigg|\,\bigotimes\limits_{i=1}^{n}\left[\bigotimes\limits_{e\in E_{\Gamma}(i)}K(e)J(e)\otimes J(e_i)K(e_i)\right]\right),
\end{align}
and similarly,
\begin{align}
    &\prod\limits_{j=1}^m\left[\left(\prod_{e\in E_{\Gamma}(\overline{j})}\partial_{\mathcal{I}(e)}\right)(\a_j|X_j)\right]=\left(\prod\limits_{j=1}^{m}\D_{\{I(e)\}_{e\in E_{\Gamma}(\overline{j})}}(\a_j)\,\Bigg|\,\bigotimes\limits_{j=1}^{m}\left[\bigotimes\limits_{e\in E_{\Gamma}(\overline{j})}K(e)J(e)\otimes X_j\right]\right).
\end{align}
Then
\begin{align}
    &\prod\limits_{i=1}^{n}\left[\left(\prod_{e\in E_{\Gamma}(i)}\partial_{\mathcal{I}(e)}\right)(\Phi_i)_N\Big(x_{\mathcal{I}(e_i^1)},\ldots,x_{\mathcal{I}(e_i^{k_i+1})}\Big)\right]\prod\limits_{j=1}^m\left[\left(\prod_{e\in E_{\Gamma}(\overline{j})}\partial_{\mathcal{I}(e)}\right)(\a_j|X_j)\right]\\
    &\hspace{30pt}=\Bigg(\prod\limits_{i=1}^{n}\D_{\{I(e)\}_{e\in E_{\Gamma}(i)}}\Phi_i\Big(x_{I(e_i^1)},\ldots,x_{I(e_i^{k_i+1})}\Big)\prod\limits_{j=1}^{m}\D_{\{I(e)\}_{e\in E_{\Gamma}(\overline{j})}}(\a_j)\,\Bigg|\,\\
    &\hspace{140pt}\bigotimes\limits_{i=1}^{n}\left[\bigotimes\limits_{e\in E_{\Gamma}(i)}K(e)J(e)\otimes J(e_i)K(e_i)\right]\otimes \bigotimes\limits_{j=1}^{m}\left[\bigotimes\limits_{e\in E_{\Gamma}(\overline{j})}K(e)J(e)\otimes X_j\right]\Bigg).
\end{align}
Hence by the definition of $\sigma_{\Gamma}^1$ and $\sigma_{\Gamma}^2$ one has
\begin{align}
    &\prod\limits_{i=1}^{n}\left[\left(\prod_{e\in E_{\Gamma}(i)}\partial_{\mathcal{I}(e)}\right)(\Phi_i)_N\Big(x_{\mathcal{I}(e_i^1)},\ldots,x_{\mathcal{I}(e_i^{k_i+1})}\Big)\right]\prod\limits_{j=1}^m\left[\left(\prod_{e\in E_{\Gamma}(\overline{j})}\partial_{\mathcal{I}(e)}\right)(\a_j|X_j)\right]\\
    &\hspace{30pt}=\Bigg((\sigma_{\Gamma}^2)^{-1}\cdot\left[\prod\limits_{i=1}^{n}\D_{\{I(e)\}_{e\in E_{\Gamma}(i)}}\Phi_i\Big(x_{I(e_i^1)},\ldots,x_{I(e_i^{k_i+1})}\Big)\prod\limits_{j=1}^{m}\D_{\{I(e)\}_{e\in E_{\Gamma}(\overline{j})}}(\a_j)\right]\cdot \sigma_{\Gamma}^1\,\Bigg|\,\\
    &\hspace{240pt}\left(\bigotimes\limits_{e\in E_{\Gamma}}J(e)J(e)\otimes K(e)K(e)\right)\otimes X_1\otimes\ldots\otimes X_m\Bigg).
\end{align}
Finally, summing over $I$, $J$, and $K$ indices we obtain the desired result.

\end{proof}

\begin{corollary}
    For any $\Phi_i\in L^{k_i}_{di\mh Pois}(\O(A))$, where $k_i=\#star(i)-1$, one has $\U_{\Gamma}(\Phi_1,\ldots,\Phi_n)\in \Ldistar^{k_1+\ldots+k_n+1-n}(\O(A))$.
\end{corollary}
\begin{proof}
    Combine Lemma \ref{lemma13} and Proposition \ref{prop20}.
\end{proof}

Let us extend $\U_{\Gamma}: L^{k_1}_{di\mh Pois}(\O(A))\otimes\ldots\otimes L^{k_n}_{di\mh Pois}(\O(A))\longrightarrow \Ldistar^{k_1+\ldots+k_n+1-n}(\O(A))$ for $k_i=\#star(i)-1$ to a map
\begin{align}
    \U_{\Gamma}:\Lpoisdi(\O(A))^{\otimes n}\longrightarrow\Ldistar(\O(A))[1-n],
\end{align}
assuming that it equals zero on all graded components but $(k_1,\ldots,k_n)$.

Next, we define $\U_n:\Lpoisdi(\O(A))^{\otimes n}\longrightarrow\Ldistar(\O(A))[1-n]$ by
\begin{align}
    \U_n:=\sum\limits_{m\geq 0}\sum\limits_{\Gamma\in G_{n,m}}W_{\Gamma}\,\U_{\Gamma},
\end{align}
where $W_{\Gamma}\in\amsmathbb{R}$\footnote{Starting from here we assume that the base field $\Bbbk$ is the field of complex numbers.} are defined in Section 6.2 of \cite{kontsevich2003deformation}.

\begin{theorem}[Double Formality Theorem for $N\amsmathbb{A}^d$]\label{th4}
Let $A=\Bbbk\langle x_1,\ldots,x_d\rangle$. Then
\begin{equation}
    \U=\big(\U_n\big)_{n\geq 1}:\begin{tikzcd}
        \Lpoisdi(\O(A))\arrow[r,rightsquigarrow]& \Ldistar(\O(A))
    \end{tikzcd}
\end{equation}
is an $L_{\infty}$-quasi-isomorphism.
\end{theorem}
\begin{proof}
    Let's prove the $L_{\infty}$-morphism part. We need to check that $\U_n$ is graded anti-symmetric
    \begin{align}\label{f83}
        \U_n(\Phi_1,\ldots,\Phi_i,\Phi_{i+1},\ldots,\Phi_n)+(-1)^{|\Phi_i|\cdot|\Phi_{i+1}|}\U_n(\Phi_1,\ldots,\Phi_{i+1},\Phi_{i},\ldots,\Phi_n)=0
    \end{align}
    and satisfies a quadratic relation of the form
    \begin{align}\label{f81}
        &(-1)^n\sum\limits_{i<j}(-1)^s\U_{n-1}\Big([\Phi_i,\Phi_j]_{di\mh RN},\Phi_1,\ldots,\widehat{\Phi_i},\ldots,\widehat{\Phi_j},\ldots,\Phi_n\Big)\\
        &-\cfrac{1}{2}\sum\limits_{p+q=n}\sum\limits_{\sigma\in\Sh(p,q)}(-1)^{pn+t}\Big[\U_p\big(\Phi_{\sigma(1)},\ldots,\Phi_{\sigma(p)}\big),\U_q\big(\Phi_{\sigma(p+1)},\ldots,\Phi_{\sigma(n)}\big)\Big]-d\Big(\U_n\big(\Phi_1,\ldots,\Phi_n\big)\Big)=0
    \end{align}
    for any homogeneous $\Phi_1,\ldots,\Phi_n\in\Lpoisdi(\O(A))$, where $s$ and $t$ are explicitly defined integers which do not depend on $\Phi_i$'s, see Section 3.6 in \cite{cattaneo2005deformation} for details. Note that the left-hand sides are admissible maps for any homogeneous $\Phi_1,\ldots,\Phi_n\in\Lpoisdi(\O(A))$.

    For any $\Psi\in\Ldi^n(\O(A))$ we will denote by $(\Psi)_N$ the induced linear map
    \begin{align}
        (\Psi)_N:\orep^{\otimes n+1}\longrightarrow\orep.
    \end{align}
    Due to Proposition \ref{prop5} it's enough to prove that $(L_{\U}(\Phi_1,\ldots,\Phi_n))_N=0$ for any $N$, where $L_{\U}(\Phi_1,\ldots,\Phi_n)$ stands for the left-hand sides of \eqref{f83} and \eqref{f81}. This follows from Lemma \ref{lemma13}, Remark \ref{rem3}, Proposition \ref{prop18}, and Theorem 6.3 from \cite{kontsevich2003deformation}. Indeed, by Lemma \ref{lemma13} one has
    \begin{align}
        \Big(\U_n(\Phi_1,\ldots,\Phi_n)\Big)_N=\mathcal{U}_n\Big((\Phi_1)_N,\ldots,(\Phi_n)_N\Big),
    \end{align}
    where $\mathcal{U}_n=\sum\limits_{m\geq 0}\sum\limits_{\Gamma\in G_{n,m}}W_{\Gamma}\,\mathcal{U}_{\Gamma}$. Hence, by last items in Remark \ref{rem3} and Proposition \ref{prop18} one has
    \begin{align}
        \Big(L_{\U}(\Phi_1,\ldots,\Phi_n)\Big)_N=L_{\mathcal{U}}\Big((\Phi_1)_N,\ldots,(\Phi_n)_N\Big),
    \end{align}
    where $L_{\mathcal{U}}$ on the right-hand side stands for a similar expression as $L_{\U}$ but computed for $\mathcal{U}$ instead of $\U$. Finally, by Theorem 6.3 from \cite{kontsevich2003deformation} the right-hand side equals zero and the proof is complete.

    Let's prove the quasi-isomorphism part. We need to check that
    \begin{align}
        \U_1:\Lpoisdi(\O(A))\longrightarrow\Ldistar(\O(A))    
    \end{align}
    is a quasi-isomorphism. One easily sees that this map is simply the skew-symmetrization, i.e. it's induced by 
    \begin{align}
        &\O(A)^{\stimes n}\xtwoheadrightarrow{\hspace{20pt}} \bigwedge\nolimits_{\S}^{n}\O(A),\\
        &\a\mapsto \sum\limits_{\sigma\in S(n)}(-1)^{\sigma}\sigma(\a).
    \end{align}
    So, $\U_1$ actually lands in $\Ldistarnorm(\O(A))$. Then by Proposition \ref{prop22}, it's enough to prove that 
    \begin{align}
        \U_1:\Lpoisdi(\O(A))\longrightarrow\Ldistarnorm(\O(A))    
    \end{align}
    is a quasi-isomorphism. It's enough to have a "universal (=stable in $N$)" homotopy retract
    \begin{equation}\label{f91}
    \begin{tikzcd}
        L_{1,star}\Big(\orep\Big)^{\GL} 
            \arrow[loop left, out=-172, in=172, distance=4em, "h_N"] 
            \arrow[r,two heads, shift left=0.8ex,"p_N"] 
        &L_{Pois}\Big(\orep\Big)^{\GL}, \arrow[l,hook', shift left=0.8ex,"i_N"]
    \end{tikzcd}
\end{equation}
which can be lifted to $\ored(A)=\O(A)$ due to the universal nature of maps $h_N$ and $p_N$. Such a retract does exist, as follows from Section 4 in \cite{de1995homotopy}. More concretely, in \cite{de1995homotopy} the authors constructed an explicit homotopy retract
\begin{equation}
    \begin{tikzcd}
        L_{1,star}\Big(\Bbbk[x_1,\ldots,x_m]\Big) 
            \arrow[loop left, out=-172, in=172, distance=4em, "H_m"] 
            \arrow[r,two heads, shift left=0.8ex,"P_m"] 
        &L_{Pois}\Big(\Bbbk[x_1,\ldots,x_m]\Big), \arrow[l,hook', shift left=0.8ex,"I_m"]
    \end{tikzcd}
\end{equation}
where $I_m$ is the skew-symmetrization, and $H_m$ and $P_m$ are sums with universal (rational) coefficients of certain expressions involving 
\begin{itemize}
    \item summation over $\{x_1,\ldots,x_m\}$,
    
    \item the bracket $[-,-]$ in $L_{star}\Big(\Bbbk[x_1,\ldots,x_m]\Big)$,

    \item evaluation of elements of $L_{star}\Big(\Bbbk[x_1,\ldots,x_m]\Big)$ on $x_1,\ldots,x_m$,

    \item compositions of elements of $L_{star}\Big(\Bbbk[x_1,\ldots,x_m]\Big)$ with derivatives $\partial_i:=\cfrac{\partial}{\partial x_i}$.
\end{itemize}

Although the summation over the set $\{x_1, \ldots, x_m\}$ is involved, the maps $H_m$, $P_m$, and, clearly, $I_m$ are $\operatorname{GL}_m$-invariant, since each $x_i$ appears together with the corresponding derivative $\partial_i$.

When we apply this construction to our case $m=dN^2$ with $\orep=\Bbbk[x_{(k,i,j)}]$, where $k=1,\ldots,d$ and $i,j=1,\ldots,N$, the only obstacle to lifting $I_m$, $H_m$, and $P_m$ to $\O(A)$ is the summation over all the indices $k=1,\ldots,d$, $i,j=1,\ldots,N$. But that's not a problem as long as the formulas are $\GL\subset\operatorname{GL}_{dN^2}$-invariant. We do the same thing we have already done numerous times in this paper -- we forget about the indices corresponding to the noncommutative variables $x_1,\ldots,x_d$ and sum over all remaining indices ranging from $1$ to $N$. By careful looking at formulas in Section 4 of \cite{de1995homotopy} one realizes that such a summation can be performed because these remaining indices enter all the formulas only as a certain power of $\tr$ twisted by appropriate permutations, i.e. $\vect(\tau)\covect(w)\tr^{\otimes r}$. In this way we obtain "universal" retract \eqref{f91} and its lift to $\O(A)$. All the necessary properties of this lift follow from Theorem \ref{th3} and Proposition \ref{prop5}. Now the proof is complete.
\end{proof}

    \addtocontents{toc}{\protect\setcounter{tocdepth}{0}}
    \printbibliography
    \addtocontents{toc}{\protect\setcounter{tocdepth}{3}}
    \addcontentsline{toc}{section}{References}

\end{document}